%% file: grain20160331.tex
\documentclass[a4paper,11pt]{amsart}
\usepackage{latexsym,amsmath,amssymb,amsthm}
\newtheorem{thm}{Theorem}[section]
\newtheorem{lemma}[thm]{Lemma}
\newtheorem{cor}[thm]{Corollary}
\newtheorem{rem}[thm]{Remark} 

\newtheorem{define}[thm]{Definition}
\newtheorem{prop}[thm]{Proposition}

\numberwithin{equation}{section}

\usepackage{constants}
\usepackage{color}
\input{amssym.def}

\usepackage{graphicx}
\usepackage{epstopdf}
\usepackage{epsfig} 

\newconstantfamily{eps}{symbol=\epsilon}
\newconstantfamily{c}{symbol=c}
\newcommand{\e}{\varepsilon}

\title{On the mean curvature flow of grain boundaries}
\author{Lami Kim}
\author{Yoshihiro Tonegawa}
\thanks{Both authors are partially supported by JSPS Grant-in-aid for scientific research (A) $\#$25247008 and (S) $\#$26220702. They thank the anonymous referee for a number of valuable comments.}
\address{Department of Mathematics, Tokyo Institute of Technology, 
Ookayama, 152-8550,  JAPAN}
\email{lmkim@math.titech.ac.jp}
\email{tonegawa@math.titech.ac.jp}
\allowdisplaybreaks[1]
\begin{document}
\begin{abstract}
Suppose that $\Gamma_0\subset\mathbb R^{n+1}$ is a closed countably 
$n$-rectifiable set whose complement
$\mathbb R^{n+1}\setminus \Gamma_0$ consists of more than one connected 
component.
Assume that the
$n$-dimensional Hausdorff measure of $\Gamma_0$ is finite or 
grows at most exponentially near infinity.
Under these assumptions, we prove a global-in-time existence of
mean curvature flow in the sense of Brakke starting from 
$\Gamma_0$. There exists a finite family of open sets which move continuously with respect to the Lebesgue measure, and whose boundaries
coincide with the space-time support of the mean curvature flow. 
\end{abstract}

\maketitle
\section{Introduction}
A family of $n$-dimensional surfaces $\{\Gamma(t)\}_{t\geq 0}$ in $\mathbb R^{n+1}$ is called the mean curvature
flow (hereafter abbreviated by MCF) if the velocity is equal to its mean curvature at
each point and time. Since the 1970's, the MCF has been studied by numerous researchers
as it is one of the fundamental geometric evolution problems (see \cite{BellettiniT,Colding,Ecker,Giga1,MantegazzaT} for the overview and 
references related to the MCF) appearing in fields such as 
differential geometry, general relativity, image processing and materials science. Given a smooth surface 
$\Gamma_0$, one can find a smoothly moving MCF starting from $\Gamma_0$ 
until some singularities such as
vanishing or pinching occur. The theory of MCF inclusive of such occurrence of 
singularities started with the pioneering work of 
Brakke in his seminal work \cite{Brakke}. 
He formulated a notion of MCF in the
setting of geometric measure theory and discovered a number of striking
measure-theoretic properties in this general setting. It is often called the
Brakke flow and we call the flow by this name hereafter. It is a family of 
varifolds representing generalized surfaces which satisfy the motion law
of MCF in a distributional sense. His aim was to allow
a broad class of singular surfaces to move by the MCF which can
undergo topological changes. 
Quoting from \cite[p.1]{Brakke}: ``A
physical system exhibiting this behavior is the motion of grain boundaries in an 
annealing pure metal [...] It is experimentally observed that these grain boundaries 
move with a velocity proportional to their mean curvature.'' 
One of Brakke's major achievements
is his general existence theorem \cite[Chapter 4]{Brakke}. 
Given a general integral
varifold as an initial data, he proved a global-in-time existence of Brakke flow
with an ingenious approximation scheme and 
delicate compactness-type theorems on varifolds. One serious uncertainty on his existence theorem, however, is 
that there is no guarantee that the MCF he obtained
is nontrivial. That is, since the definition of Brakke flow is flexible enough to
allow sudden loss of measure at any time, 
whatever the initial $\Gamma_0$ is,
setting $\Gamma(t)=\emptyset$ for all $t>0$, we obtain a Brakke flow satisfying
the definition trivially. 
The proof of existence in \cite{Brakke} does not preclude the unpleasant
possibility of getting this trivial flow when one takes the limit of approximate
sequence.
The idea of such ``instantaneous vanishing'' may appear unlikely, but
the very presence of singularities of $\Gamma_0$ 
may potentially cause such catastrophe in his
approximation scheme. For this reason,
rigorous global-in-time existence
of MCF of grain boundaries has been considered completely open among 
the specialists. 

In this regard, we have two aims in this paper. The first aim is to reformulate 
and modify the
approximation scheme so that we always 
obtain a nontrivial MCF even if $\Gamma_0$ is singular. We 
prove for the first time a rigorous global-in-time existence theorem of the MCF of
grain boundaries which was not known even for the 1-dimensional case. 
The main existence theorem of the 
present paper may be stated roughly as follows. 
\begin{thm}
Let $n$ be a natural number and
suppose that $\Gamma_0\subset\mathbb R^{n+1}$ is a closed countably 
$n$-rectifiable set whose complement
$\mathbb R^{n+1}\setminus \Gamma_0$ is not connected. 
Assume that the $n$-dimensional Hausdorff measure of $\Gamma_0$ 
is finite or grows at most exponentially near infinity. 
Let $E_{0,1},\ldots,E_{0,N}\subset\mathbb R^{n+1}$ be mutually disjoint non-empty 
open sets with $N\geq 2$ such that $\mathbb R^{n+1}\setminus \Gamma_0=
\cup_{i=1}^N E_{0,i}$. 
Then, for each $i=1,\ldots,N$, there exists a family of open sets $\{E_i(t)\}_{t\geq 0}$
with $E_i(0)=E_{0,i}$ such that $E_1(t),\ldots E_N(t)$ are mutually disjoint for 
each $t\geq 0$ and $\Gamma(t):=\mathbb R^{n+1}\setminus \cup_{i=1}^N E_i(t)$ 
is a MCF with $\Gamma(0)=\Gamma_0$, in the sense that $\Gamma(t)$ 
coincides with the space-time support of a Brakke flow 
starting from $\Gamma_0$. Each $E_i(t)$ moves continuously in time with respect 
to the Lebesgue measure. 
\end{thm}
We may regard each $E_i(t)\subset\mathbb R^{n+1}$ 
as a region of ``$i$-th grain'' at time $t$, and 
$\Gamma(t)$ as the ``grain boundaries'' which move by their mean 
curvature. Some of $E_i(t)$ shrink and vanish, and some may grow
and may even occupy the whole $\mathbb R^{n+1}$ in finite time. We may 
also consider a periodic setting, and in that case, a typical phenomenon is a 
grain coarsening. As a framework,
loosely speaking, instead of working only with
varifolds as Brakke did, we perceive the varifolds as boundaries of a finite number of 
open sets $E_i(t)$ at each time. The open sets are designed to 
move continuously with respect to the Lebesgue measure, so that the boundaries do
not vanish instantaneously at $t=0$. Sudden loss of measure may still occur when some
``interior boundaries'' inside $E_i(t)$ appear, but otherwise, one cannot vanish certain portion of 
boundaries arbitrarily. The resulting MCF as boundaries of open sets 
is more or less in accordance with the 
MCF of physical grain boundaries originally envisioned by Brakke. If $\mathbb R^{n+1}\setminus
\Gamma_0$ consists of $N$ connected components, 
we naturally define them to
be $E_{0,1},\ldots,E_{0,N}$. 
If there are infinitely many connected components, 
we need to group them to be finitely many mutually disjoint open sets
$E_{0,1},\ldots,E_{0,N}$ for some arbitrary $N\geq 2$, hence there is 
already non-uniqueness of grouping at this point in our scheme. 
Even if they are finitely many, simple examples indicate that 
the flow is non-unique in general, even though it is not clear how
generic the non-uniqueness prevails. 

The second aim of the paper is to clarify the content of \cite{Brakke} with a
number of modifications and simplifications. Despite the potential importance 
of the claim, there have been no review on the existence theory of
\cite{Brakke} so far. 
Also, we need to provide different definitions and proofs working in the 
framework of sets of boundaries. 
Here, we present a mostly self-contained 
proof which should be accessible to interested researchers versed in the basics of
geometric measure theory. A good working knowledge on rectifiability
\cite{AFP,Evans-Gariepy, Federer} and basics on the theory of
varifolds in \cite{Allard,Simon} are assumed. 

Next, we briefly describe and compare the known results on the 
existence of Brakke flow to that of the present paper. Given a smooth
compact embedded hypersurface in general dimension, one has a smooth MCF
until the first time singularities occur. For $n=1$, it is well known that the curves stay
embedded until they become convex and shrink to a point by the results due to 
Gage-Hamilton \cite{Gage} and Grayson \cite{Grayson} (see 
also \cite{Andrew} for the elegant and short proof). 
For general dimensions,
one has the notion of viscosity solution \cite{CGG,ES1} which gives a family
of closed sets as a unique weak solution of the MCF even after the occurrence of 
singularities. It is possible that the closed set may develop nontrivial interior afterwords,
a phenomenon called fattening, 
and it is not clear if the set is Brakke flow after singularities appear 
in general. On the other hand, Evans and
Spruck proved that almost all level sets of viscosity solution are
unit density Brakke flows \cite{ES2}.
As a different track, there are other methods such as elliptic regularization \cite{Ilmanen1}
and phase field approximation via the Allen-Cahn equation \cite{Ilmanen2,Tonegawa1} to obtain rigorous global-in-time existence results of Brakke flow. 
All of the above results use the ansatz that the MCF is
represented as a boundary of a single time-parametrized set, so that
it is not possible to handle grain boundaries with more than two grains in general.
For more general cases such as triple junction figure on a plane and the 
higher dimensional analogues, all known results up to this point are based more or less 
on a certain parametrized framework and the existence results cannot be extended
past topological changes in general. For three regular curves meeting at a
triple junction of 120 degrees, Bronsard and Reitich \cite{Bronsard} proved
short-time existence and uniqueness using a theory of system of parabolic PDE \cite{Solo}. 
There are numerous results studying existence, uniqueness (or non-uniqueness) and
stability under various boundary conditions as well as studies on the 
self-similar shrinking/expanding solutions, and we mention 
\cite{Bellettini,Chen,Chen2,Freire,Freire2,Garcke,Ikota2005,Neves,
Kinderlehrer,Magni,Mantegazza,Saez,Sch}. 
Compared to the above known results, our existence theorem does not require any
parametrization and there is no restriction on the dimension or configuration. 
The regularity assumption put on the closed set $\Gamma_0$ is 
countable $n$-rectifiability, 
which allows wide variety of singularities, and the solution can undergo past topological
changes. In this sense, even the results for the 1-dimensional case are new in an
essential way. 

On the computational side of the MCF of grain boundaries, there are
enormous number of works on the simulations and algorithms, which are far
beyond the scope of this paper. Here we simply mention
for a point of reference that Brakke developed
an interactive software {\it Surface Evolver} \cite{Brakke2} which handles 
variety of geometric flow problems including the MCF. See video clips
of 1-dimensional MCF of grain boundaries of as many as N=10,000
in Brakke's home page \cite{Brakke3}.  

We end the introduction by describing  
the organization of the paper. In Section 2, we state our basic 
notation and present preliminary materials from geometric measure theory. 
In Section 3, we state the main existence results and give an overview of
the proof. Section 4 introduces notions of open partition and a certain class of 
admissible functions as well as some preliminary materials concerning varifold smoothing.
Section 5 contains a number of estimates on the approximation of smoothed
mean curvature vector essential to the construction of approximate solutions. 
Section 6 gives the actual construction of approximate solutions with 
good estimates derived in Section 5. Section 7 and 8 are mostly independent 
from the previous sections and prove
compactness-type theorems for rectifiability and integrality, respectively, of the
limit varifold. Gathering all the results up to this point,  
Section 9 proves that the family of limit measures is a Brakke flow.
Section 10 proves a certain continuity property of domains of ``grains''.
Section 11 gives additional comments on the property of the solution.

\section{Notation and preliminaries}
\subsection{Basic notation}
$\mathbb N$, $\mathbb Q$, $\mathbb R$ are the sets of natural numbers, rational
numbers, real numbers, respectively. We set $\mathbb R^+:=\{x\in \mathbb R : 
x\geq 0\}$. We reserve $n\in \mathbb N$ for the dimension of hypersurface and 
$\mathbb R^{n+1}$ is the $n+1$-dimensional Euclidean space. 
For $r\in (0,\infty)$ and $a\in \mathbb R^{n+1}$ define
\begin{equation*}
\begin{split}
& B_r(a):=\{x\in \mathbb R^{n+1} : |x-a|\leq r\}, \, B_r^n (a):=\{x\in \mathbb R^n : 
|x-a|\leq r\}, \\ &U_r(a):=\{x\in \mathbb R^{n+1} : |x-a|<r\},\, U_r^n(a):=\{x\in \mathbb R^n : 
|x-a|<r\}
\end{split}
\end{equation*}
and when $a=0$ define $B_r:=B_r(0)$, $B_r^n:=B_r^n(0)$, $U_r:=U_r(0)$ and $U_r^n
:=U_r^n(0)$. For a subset $A\subset \mathbb R^{n+1}$, ${\rm int}\,A$ is the set of interior
points of $A$, and ${\rm clos}\,A$ denotes the closure of $A$. ${\rm diam}\,A$ is the
diameter of $A$. For two subsets $A,B
\subset \mathbb R^{n+1}$, define $A\triangle B:=(A\setminus B)\cup (B\setminus A)$. 
For an open subset $U\subset\mathbb R^{n+1}$ let $C_c(U)$ be the set of all 
compactly supported continuous functions defined on $U$ and let $C_c(U ; 
\mathbb R^{n+1})$ be the set of all compactly supported continuous vector fields. 
Indices $l$ of $C_c^l(U)$ and $C_c^l(U;\mathbb R^{n+1})$ indicate continuous 
$l$-th order differentiability. For $g\in C^1(U;\mathbb R^{n+1})$, we regard
$\nabla g(x)$ as an element of ${\rm Hom}(\mathbb R^{n+1};\mathbb R^{n+1})$.
Similarly for $g\in C^2(U)$, we regard the Hessian matrix $\nabla^2 g(x)$ as an
element of ${\rm Hom}(\mathbb R^{n+1};\mathbb R^{n+1})$. For a Lipschitz
function $f$, ${\rm Lip}\,(f)$ is the Lipschitz constant. 

\subsection{Notation related to measures}
$\mathcal L^{n+1}$ denotes the Lebesgue measure on 
$\mathbb R^{n+1}$ and $\mathcal H^n$ denotes the $n$-dimensional Hausdorff measure
on $\mathbb R^{n+1}$. $\mathcal H^0$ denotes the counting measure.  
We use $\omega_n:=\mathcal H^n(U_1^n)$ and  
$\omega_{n+1}:=\mathcal L^{n+1}(U_1)$. The restriction of $\mathcal H^n$ to a set $A$ is denoted by
$\mathcal H^n\lfloor_A$, and when $f$ is a $\mathcal H^n$ measurable function 
defined on $\mathbb R^{n+1}$, 
$\mathcal H^n\lfloor_f$ is the weighted $\mathcal H^n$ by $f$. Let ${\bf B}_{n+1}$ 
be the constant appearing in Besicovitch's
covering theorem (see \cite[1.5.2]{Evans-Gariepy}) on $\mathbb R^{n+1}$.  

For a Radon measure $\mu$ on $\mathbb R^{n+1}$ and $\phi\in 
C_c(\mathbb R^{n+1})$, we often write $\mu(\phi)$ for $\int_{\mathbb R^{n+1}}
\phi\, d\mu$. Let ${\rm spt}\,\mu$ be the support of $\mu$, i.e., ${\rm spt}\,\mu
:=\{x\in \mathbb R^{n+1} : \mu(B_r(x))>0 \mbox{ for all $r>0$}\}$. By definition, 
it is a closed set. Let $\theta^{* n}(\mu,x)$ be defined by $\limsup_{r\rightarrow 0+}
\mu(B_r(x))/(\omega_n r^n)$ and let $\theta^n (\mu,x)$ be defined as 
$\lim_{r\rightarrow 0+} \mu(B_r(x))/(\omega_n r^n)$ when the limit exists. The set of
$\mu$ measurable and (locally) square integrable functions as well as 
vector fields is denoted by $L^2(\mu)$
($L^2_{loc}(\mu)$). For a set $A\subset \mathbb R^{n+1}$, $\chi_A$ is the characteristic function of $A$. If $A$ is a set of finite perimeter,
$\|\nabla\chi_A\|$ is the total variation measure of the distributional derivative
$\nabla\chi_A$. 

\subsection{The Grassmann manifold and varifold}
Let ${\bf G}(n+1,n)$ be the space of $n$-dimensional subspaces of ${\mathbb R}^{n+1}$.
For $S\in {\bf G}(n+1,n)$, we identify $S$ with the corresponding orthogonal 
projection of ${\mathbb R}^{n+1}$ onto $S$. Let $S^{\perp}\in {\bf G}(n+1,1)$ be the orthogonal
complement of $S$. For two  elements $A$
and $B$ of ${\rm Hom}\, ({\mathbb R}^{n+1};{\mathbb R}^{n+1})$, define a scalar product $A\cdot B:={\rm trace}\, (A^{\top}\circ B)$ where $A^{\top}$ is the transpose of $A$ and $\circ$ indicates the composition. The identity of ${\rm Hom}\, ({\mathbb R}^{n+1};{\mathbb R}^{n+1})$ is denoted by $I$. Let $a\otimes
b\in {\rm Hom}\, ({\mathbb R}^{n+1};{\mathbb R}^{n+1})$ be the tensor product of $a,\, b\in {\mathbb R}^{n+1}$, i.e., as an $(n+1)\times(n+1)$ matrix, the $(i,j)$-component
is given by $a_i b_j$ where $a=(a_1,\ldots,a_{n+1})$ and similarly for $b$. For $A\in {\rm Hom}\, ({\mathbb R}^{n+1};{\mathbb R}^{n+1})$ define
\begin{equation*}
|A|:=\sqrt{A\cdot A},\hspace{1cm}
\|A\|:=\sup\{|A(x)|\, :\, x\in {\mathbb R}^{n+1},\, |x|=1\}.
\end{equation*}
For $A\in {\rm Hom}\,({\mathbb R}^{n+1};{\mathbb R}^{n+1})$ and 
$S\in {\bf G}(n+1,n)$, let $|\Lambda_n(A\circ S)|$ be the absolute value of the
Jacobian of the map
$A\lfloor_S$. If $S$ is spanned by a set of orthonormal basis $v_1,\ldots,v_n$,
then $|\Lambda_n (A\circ S)|$ is the $n$-dimensional volume of the
parallelepiped formed by $A(v_1),\ldots,A(v_n)$. If we form a $(n+1)\times n$
matrix $B$ 
with these vectors as the columns, we may compute $|\Lambda_n (A\circ S)|$ as 
the square root of the sum of the squares of the determinants of the 
$n\times n$ submarices of $B$, or we may compute it as $\sqrt{{\rm det}\,(B^{\top}\circ B)}$.

We recall some notions related to varifolds and refer to \cite{Allard,Simon} for more
details. Define ${\bf G}_n(\mathbb R^{n+1}):=\mathbb R^{n+1}\times {\bf G}(n+1,n)$. For any subset $C\subset \mathbb R^{n+1}$, we similarly 
define ${\bf G}_n(C):=C\times{\bf G}(n+1,n)$. 
A general $n$-varifold in $\mathbb R^{n+1}$ is a Radon measure on 
${\bf G}_n(\mathbb R^{n+1})$. The set of all general $n$-varifolds in $\mathbb R^{n+1}$ is denoted by ${\bf V}_n(\mathbb R^{n+1})$. 
For $V\in {\bf V}_n(\mathbb R^{n+1})$, let $\|V\|$ be the weight measure of $V$, 
namely, for all $\phi\in C_c(\mathbb R^{n+1})$,
\begin{equation*}
\|V\|(\phi):=\int_{{\bf G}_n(\mathbb R^{n+1})}\phi(x)\, dV(x,S).
\end{equation*}
For a
proper map $f\in C^1({\mathbb R}^{n+1};{\mathbb R}^{n+1})$ define
$f_{\#}V$ as the push-forward of varifold $V\in {\bf V}_n(\mathbb R^{n+1})$ (see \cite[3.2]{Allard} for the definition).
Given any ${\mathcal H}^n$ measurable 
countably $n$-rectifiable set $\Gamma \subset \mathbb R^{n+1}$ with
locally finite ${\mathcal H}^n$ measure, there
is a natural $n$-varifold $|\Gamma |\in {\bf V}_n(\mathbb R^{n+1})$ defined by
\begin{equation*}
|\Gamma |(\phi):=\int_{\Gamma }\phi(x,{\rm Tan}^n(\Gamma,x))\, d{\mathcal H}^n(x)
\end{equation*}
for all $\phi\in C_c({\bf G}_n(\mathbb R^{n+1}))$.
Here, ${\rm Tan}^n( \Gamma,x)\in {\bf G}(n+1,n)$ is the approximate tangent space which exists
${\mathcal H}^n$ a.e$.$ on $\Gamma$ (see \cite[2.2.11]{AFP}). In this case, $\||\Gamma|\|={\mathcal H}^n \lfloor_\Gamma$.

We say $V\in {\bf V}_n(\mathbb R^{n+1})$ is {\it rectifiable} if for all $\phi \in C_c({\bf G}_n(\mathbb R^{n+1}))$,
\begin{equation*}
V(\phi)=\int_{\Gamma}\phi(x,{\rm Tan}^n( \Gamma,x))\theta(x)\, d{\mathcal H}^n(x)
\end{equation*}
for some ${\mathcal H}^n$ measurable 
countably $n$-rectifiable set $\Gamma \subset \mathbb R^{n+1}$ and locally ${\mathcal H}^n$ integrable 
non-negative function $\theta$ defined on $\Gamma$. The set of all rectifiable $n$-varifolds
is denoted by ${\bf RV}_n(\mathbb R^{n+1})$.
Note that for such varifold, $\theta^n(\|V\|,x)
=\theta(x)$, approximate tangent space as varifold 
exists and is equal to ${\rm Tan}^n(\Gamma,x)$, 
${\mathcal H}^n$ a.e$.$ on $\Gamma$. The approximate tangent space is denoted by ${\rm Tan}^n (\|V\|,x)$. 
In addition, if $\theta(x)\in \mathbb N$ for
$\mathcal H^n$ a.e$.$ on $\Gamma$, we say $V$ is {\it integral}. The set of all 
integral $n$-varifolds in $\mathbb R^{n+1}$ is denoted by ${\bf IV}_n(\mathbb R^{n+1})$.
We say $V$ is a {\it unit density}
$n$-varifold if $V$ is integral and $\theta=1$ $\mathcal H^n$ a.e$.$ on $\Gamma$, i.e., $V=|\Gamma|$. 

\subsection{First variation and generalized mean curvature}
For $V\in {\bf V}_n(\mathbb R^{n+1})$ let $\delta V$ be the first variation of $V$, namely,
\begin{equation}
\delta V(g):=\int_{{\bf G}_n(\mathbb R^{n+1})}\nabla g(x)\cdot S\, dV(x,S)
\label{fvf0}
\end{equation}
for $g\in C_c^1(\mathbb R^{n+1};{\mathbb R}^{n+1})$. 
Let $\|\delta V\|$ be the total variation measure when it exists. If 
$\|\delta V\|$ is absolutely continuous with respect to $\|V\|$, by the 
Radon-Nikodym theorem, we have
for some $\|V\|$ measurable vector field $h(\cdot,V)$
\begin{equation}
\delta V(g)=-\int_{\mathbb R^{n+1}}g(x)\cdot h(x,V)\, d\|V\|(x).
\label{fvf}
\end{equation}
The vector field $h(\cdot,V)$ is called the {\it generalized mean curvature} of $V$.
For any $V\in {\bf IV}_n(\mathbb R^{n+1})$ with bounded first variation (so in particular
when $h(x,V)$ exists), Brakke's 
perpendicularity theorem of generalized mean curvature
\cite[Chapter 5]{Brakke} says that we have for $V$ a.e$.$ $(x,S)\in {\bf G}_n(\mathbb
R^{n+1})$
\begin{equation}
S^{\perp}(h(x,V))=h(x,V).
\label{fvf2}
\end{equation}
One may also understand
this property in connection with $C^2$ rectifiability of varifold 
established in \cite{Menne}. 

\subsection{The right-hand side of the MCF equation}
For any $V\in {\bf V}_n(\mathbb R^{n+1})$, $\phi\in C^1_c(\mathbb R^{n+1};{\mathbb R}^+)$ and $g\in C^1_c(\mathbb R^{n+1};\mathbb R^{n+1})$, define
\begin{equation}
\delta(V,\phi)(g):=\int_{{\bf G}_n(\mathbb R^{n+1})} \phi(x)\nabla g(x)\cdot S\,\, dV(x,S)+
\int_{\mathbb R^{n+1}} g(x)\cdot \nabla \phi(x)\, d\|V\|(x).
\label{Bdef}
\end{equation}
As explained in \cite[2.10]{Brakke}, $\delta(V,\phi)(g)$ may be considered as a 
$\phi$-weighted first variation of $V$ in the direction of $g$. 
Using $\phi\nabla g=\nabla(\phi g)-g\otimes \nabla\phi$ and \eqref{fvf0}, we have
\begin{equation}
\begin{split}
\delta(V,\phi)(g)&=\delta V(\phi g)+\int_{{\bf G}_n(\mathbb R^{n+1})}
 g(x)\cdot (I-S)(\nabla\phi(x))\, dV(x,S) \\ &
 =\delta V(\phi g)+\int_{{\bf G}_n(\mathbb R^{n+1})}
 g(x)\cdot S^{\perp}(\nabla\phi(x))\, dV(x,S).
 \end{split}
 \label{Bdefs}
 \end{equation}
When $\|\delta V\|$ is locally finite and absolutely continuous
with respect to $\|V\|$,  \eqref{fvf} and \eqref{Bdefs} show
\begin{equation}
\label{Bdefs2}
\delta(V,\phi)(g)= -\int_{\mathbb R^{n+1}} \phi(x)g(x)\cdot h(x,V)\,d\|V\|(x)
+\int_{{\bf G}_n(\mathbb R^{n+1})}
 g(x)\cdot S^{\perp}(\nabla\phi(x))\, dV(x,S).
 \end{equation}
 Furthermore, if $V\in {\bf IV}_n(\mathbb R^{n+1})$ with $h(\cdot,V)\in L^2_{loc}(\|V\|)$,
 by approximating each component of $h(\cdot,V)$ by a sequence of smooth functions,
 we may naturally define
 \begin{equation}
 \label{Bdefs3}
 \delta(V,\phi)(h(\cdot,V)):=\int_{\mathbb R^{n+1}} -\phi(x)|h(x,V)|^2+ h(x,V)\cdot\nabla\phi(x)\, d\|V\|(x).
 \end{equation}
 Here, we also used \eqref{fvf2}. 
 It is convenient to define $\delta(V,\phi)(h(\cdot,V))$ when some of the 
 conditions above are not satisfied. Thus, we define (even if $h(\cdot, V)$ does not 
 exist)
 \begin{equation}
 \label{Bdefs4}
 \delta(V,\phi)(h(\cdot,V)):=-\infty
 \end{equation}
unless $\|\delta V\|$ is locally finite, absolutely continuous with respect to $\|V\|$,
$V\in {\bf IV}_n(\mathbb R^{n+1})$ and $h(\cdot,V)\in L^2_{loc}(\|V\|)$. 
Formally, if a family of smooth $n$-dimensional surfaces $\{\Gamma(t)\}_{t\in\mathbb R^+}$ moves
by the velocity equal to the mean curvature, then one
can check that $V_t=|\Gamma(t)|$ satisfies
\begin{equation}
\frac{d}{dt}\|V_t\|(\phi(\cdot,t))\leq \delta(V_t,\phi(\cdot,t))(h(\cdot,V_t))+\|V_t\|\big(\frac{\partial\phi}{\partial t}(\cdot,t)\big)
\label{formal}
\end{equation}
for all $\phi=\phi(x,t)\in C_c^1(\mathbb R^{n+1}\times\mathbb R^+;{\mathbb R}^+)$.
In fact, \eqref{formal}
holds with equality. Conversely, if \eqref{formal} is satisfied for all such $\phi$, 
then one can 
prove that the velocity of motion is equal to the mean curvature. The inequality
in \eqref{formal} allows the sudden loss of measure and it is the source of 
general non-uniqueness of Brakke's formulation. 
\section{Main results}
\subsection{Weight function $\Omega$}
\label{defomega}
To include unbounded sets which may have infinite measures 
in ${\mathbb R}^{n+1}$,
we choose a weight function $\Omega\in C^2({\mathbb R}^{n+1})$ satisfying
\begin{equation}
0<\Omega(x)\leq 1,\,\, |\nabla \Omega(x)|\leq \Cl[c]{c_1}\Omega(x),\,\,\|\nabla^2
\Omega(x)\|\leq \Cr{c_1}\Omega(x)\,\mbox{for all }x\in {\mathbb R}^{n+1}
\label{omegacon1}
\end{equation}
where $\Cr{c_1}\in \mathbb R^+$ 
is a constant depending on the choice of $\Omega$. 
If one is interested in
sets of finite $\mathcal H^n$ measure, one may choose
\begin{equation*}
\Omega(x)=1
\end{equation*}
and $\Cr{c_1}=0$ in this case. 
Another example is
\begin{equation*}
\Omega(x)
=e^{-\sqrt{1+|x|^2}}.
\end{equation*} 
Note that the second condition of \eqref{omegacon1} restricts the behavior of $\Omega$
at infinity in the sense that $e^{-\Cr{c_1}|x|}\Omega(0) \leq \Omega(x)$ with
$\Cr{c_1}$ as in \eqref{omegacon1}. Thus we cannot choose 
arbitrarily fast decaying $\Omega$. Depending on the choice of $\Omega$, 
we may have different
solutions in the end. Note that we are not so concerned with the uniqueness of 
the flow in this paper. 

We often use the following
\begin{lemma}
\label{omegapro}
Let $\Cr{c_1}$ be as in \eqref{omegacon1}. Then for $x,y\in {\mathbb R}^{n+1}$, we have
\begin{equation}
\Omega(x)\leq \Omega(y)\exp(\Cr{c_1}|x-y|).
\label{omegaproeq}
\end{equation}
\end{lemma}
\subsection{Main existence theorems}

The first theorem states that there exists a Brakke flow starting from $\Gamma_0$.
The nontriviality is described subsequently. 
\begin{thm}
\label{sthm}
Suppose that $\Gamma_0\subset \mathbb R^{n+1}$ is a closed countably
$n$-rectifiable set whose complement $\mathbb R^{n+1}\setminus \Gamma_0$
consists of more than one connected component and suppose
\begin{equation}
\label{sthm1}
\mathcal H^n\lfloor_{\Omega}(\Gamma_0)\Big(=\int_{\Gamma_0} \Omega(x)\, d\mathcal H^n(x)\Big)<\infty.
\end{equation}
For some $N\geq 2$, choose a finite collection of non-empty open sets $E_{0,1},\ldots,E_{0,N}$ 
such that
they are disjoint and $\cup_{i=1}^{N} E_{0,i}=\mathbb R^{n+1}\setminus \Gamma_0$. 
Then there exists a family $\{V_t\}_{t\in \mathbb R^+}\subset{\bf V}_n(\mathbb R^{n+1})$
with the following property. 
\begin{itemize}
\item[(1)] $V_0=|\Gamma_0|$.
\item[(2)] For $\mathcal L^1$ a.e$.$ $t\in \mathbb R^+$, $V_t\in {\bf IV}_n
(\mathbb R^{n+1})$ and $h(\cdot, V_t)\in L^2(\|V_t\|\lfloor_{\Omega})$.
\item[(3)] For all $t>0$, $\|V_t\|(\Omega)\leq \mathcal H^n\lfloor_{\Omega}(\Gamma_0)
\exp(\Cr{c_1}^2 t/2)$ and 
$\int_0^t \int_{\mathbb R^{n+1}} |h(\cdot,V_s)|^2 \Omega\, d\|V_s\|ds<\infty$.
\item[(4)] For any $0\leq t_1<t_2<\infty$ and $\phi\in C_c^1(\mathbb R^{n+1}
\times \mathbb R^+;\mathbb R^+)$, we have
\begin{equation}
\label{sthm2}
\|V_t\|(\phi(\cdot,t))\Big|_{t=t_1}^{t_2}\leq \int_{t_1}^{t_2} \delta(V_t,\phi(\cdot,t))
(h(\cdot,V_t))+\|V_t\|\big(\frac{\partial\phi}{\partial t}(\cdot,t)\big)\, dt.
\end{equation}
\end{itemize}
\end{thm}
The choice of $E_{0,1},\ldots,E_{0,N}$ appears irrelevant here but there are
more properties as explained in Theorem \ref{sthm3}. 
The assumption \eqref{sthm1} allows various possibilities for the choice of $\Gamma_0$. 
If $\mathcal H^n(\Gamma_0)
<\infty$, then, we may work with $\Omega=1$ and $\Cr{c_1}=0$ as stated before. 
If $\mathcal H^n(\Gamma_0\cap B_r)
\leq c e^r$ for some $c>0$ and for all $r>0$, 
we may choose $\Omega(x)=e^{-2\sqrt{1+|x|^2}}$
with a suitable $\Cr{c_1}>0$ and
we may satisfy \eqref{sthm1}. By (2), for a.e$.$ $t$, $\delta(V_t,\phi(\cdot,t))(h(\cdot,
V_t))$ in \eqref{sthm2} is expressed by \eqref{Bdefs3}. 
Note that \eqref{sthm2} is an integral version of \eqref{formal}. 

For above $\{V_t\}_{t
\in \mathbb R^+}$, we define the corresponding space-time Radon measure $\mu$:
\begin{define}
Define a Radon measure $\mu$ on $\mathbb R^{n+1}\times \mathbb R^+$ by
$d\mu=d\|V_t\|dt$, i.e.,
\begin{equation*}
\int_{\mathbb R^{n+1}\times\mathbb R^+}\phi(x,t)\, d\mu(x,t)=\int_{\mathbb R^+}
\int_{\mathbb R^{n+1}} \phi(x,t)\, d\|V_t\|(x)dt\,\,\mbox{ for }\phi\in C_c(\mathbb R^{n+1}\times\mathbb R^+).
\end{equation*}
\label{sptime}
\end{define}
We have the following relations between $\|V_t\|$ and 
$\mu$ as well as a finiteness 
of the support. 
\begin{prop}
\label{sptimeato}
For all $t>0$ and $r>0$, 
\begin{equation}
{\rm spt}\,\|V_t\|\subset \{x : (x,t)\in {\rm spt}\,\mu\} \mbox{ and } \mathcal H^n(B_r\cap \{x : (x,t)\in {\rm spt}\,\mu\})<\infty.
\label{thsup1}
\end{equation}
\end{prop}

We next state the
existence of open complements, which may be considered as moving grains and which prevent arbitrary loss of measure of $\|V_t\|$.
\begin{thm}
\label{sthm3}
Under the same assumptions of Theorem \ref{sthm}, there exists a family of
open sets $\{E_i(t)\}_{t\in \mathbb R^+}$ for each $i=1,\ldots,N$ with the
following property. Define $\Gamma(t):=\cup_{i=1}^N \partial E_i(t)$. 
\begin{itemize}
\item[(1)] $E_i(0)=E_{0,i}$ for $i=1,\ldots,N$ and $\Gamma_0=\Gamma(0)$. 
\item[(2)] $E_1(t),\ldots,E_N(t)$ are disjoint for each $t\in \mathbb R^+$.
\item[(3)] $\{x : (x,t)\in {\rm spt}\,\mu\}=\mathbb R^{n+1}\setminus \cup_{i=1}^N E_i(t)
=\Gamma(t)$ for each $t>0$.
\item[(4)] $\|V_t\|\geq \|\nabla \chi_{E_i(t)}\|$ for each $t\in \mathbb R^+$ and $i=1,\ldots,N$. 
\item[(5)] $S(i):=\{(x,t) : x\in E_i(t), t\in \mathbb R^+\}$ is open in $\mathbb R^{n+1}\times\mathbb R^+$ for each $i=1,\ldots, N$. 
\item[(6)] Fix $i=1,\ldots,N$, $t\in \mathbb R^+$, $x\in \mathbb R^{n+1}$ and $r>0$, and define 
\begin{equation*}
g(s):=\mathcal L^{n+1}
((E_i(t)\triangle E_i(s))\cap B_r(x))
\end{equation*}
for $s\in[0,\infty)$. Then $g\in C^{0,\frac12}((0,\infty))\cap C([0,\infty))$.
\end{itemize}
\end{thm}
Since the Lebesgue measure of $E_i(t)$ changes locally continuously by (6), and the 
boundary measure bounds $\|V_t\|$ from below by (4), one may conclude that
$\|V_t\|$ remains non-zero at least for some positive time. 
If $\Gamma_0$ is compact, $\|V_t\|$ will vanish
in finite time. If unbounded, it may stay non-zero for all time. 

We say that $\{V_t\}_{t\in \mathbb R^+}$ is a {\it unit density flow} if 
$V_t$ is a unit density varifold for a.e$.$ $t\in \mathbb R^+$. 
Under this unit density assumption, 
the results of partial regularity theory of \cite{Brakke,KT,Tonegawa2} (see also 
\cite{Lahiri}) apply to this flow.
\begin{thm}
\label{ktthm}
Let $\{V_t\}_{t\in \mathbb R^+}$ be as in Theorem \ref{sthm} and additionally 
assume that it is a unit density flow. Then, for a.e$.$
$t\in \mathbb R^+$, there exists a closed set $S_t\subset\mathbb R^{n+1}$ with 
the following property. We have $\mathcal H^n(S_t)=0$, and for any $x\in 
\mathbb R^{n+1}\setminus S_t$, there exists a space-time neighborhood $O_{(x,t)}$
of $(x,t)$ such that, either ${\rm spt}\,\mu\cap O_{(x,t)}=\emptyset$ or ${\rm spt}\, \mu$ is
a $C^{\infty}$ embedded $n$-dimensional MCF in $O_{(x,t)}$. 
\end{thm}
For further properties of $\{V_t\}_{t\in \mathbb R^+}$, see Section \ref{question}. In 
particular, under a mild measure-theoretic condition on $\Gamma_0$ 
(see Section \ref{astr}), Theorem
\ref{ktthm} is always applicable for an initial short time interval. Such 
general short-time
existence of partially regular flow is also new in all dimensions.  
\subsection{Heuristic description of the proof}
It is worthwhile to summarize the main steps to prove the existence of Brakke flow
at this point. The proof may be roughly divided into two phases, the first is the 
construction of 
sequence of time-discrete approximate flows, and the second is to prove that 
the limit satisfies the desired properties of Brakke flow. 
\subsubsection{Construction of approximate flows}
Starting from $\{E_{0,i}\}_{i=1}^N$
where $\Gamma_0=\cup_{i=1}^N \partial E_{0,i}$, time-discrete approximate
flows are constructed by alternating two steps. Let $\Delta t_j$ be a small time
grid size which goes to 0 as $j\rightarrow\infty$. The very first step is to map $\{E_{0,i}\}_{i=1}^N$ by a Lipschitz map so that the image under this map almost minimizes 
$n$-dimensional measure of boundaries
in a small length scale of order $j^{-2}$ but at the same time, keeping the structure of 
``$\Omega$- finite open partition'' (Definition~\ref{LFOP}). We introduce a certain
admissible class of Lipschitz functions called $\mathcal E$-admissible functions for 
this purpose (Definition~\ref{LFOP3}). This ``Lipschitz deformation step'' (1st step) has a
regularization effect in a small length scale, which is essential to prove 
the rectifiability and integrality of the limit flow. The map should also have an effect
of de-singularizing certain unstable singularities, even though we do not know how
to utilize it so far. After this first step, we next move the open partition 
by a smooth approximate
mean curvature which is computed by smoothing the varifold. The length scale of 
smoothing is much smaller than that of Lipschitz deformation, and the time step
$\Delta t_j$ is even much smaller than the smoothing length scale, so that the 
motion of this step remains very small and the map is a
diffeomorphism. We need to estimate how close the
approximations are for various quantities  and this takes up all of Section \ref{SMC}. 
We obtain a number of estimates which are expected to hold for the limit flow and
this is a general guideline to keep in mind. After this ``mean curvature motion step'' 
(2nd step),
we go back and do the 1st step, and then the 2nd step and we keep moving open partitions by
repeating these two steps alternatingly. We make sure that we have the right estimates by
an inductive argument (Proposition~\ref{exapp}). 
\subsubsection{Proof of properties of Brakke flow}
Once we have a sequence of approximate flows with proper estimates, such as the
time semi-decreasing property and approximate motion law, we see that there exists
a subsequence which converges as measures on $\mathbb R^{n+1}$ (not as 
varifolds at this point) for all $t\in \mathbb R^+$ (Proposition~\ref{promecone}). We then proceed to prove that the limit
measures are rectifiable first (Section~\ref{secrec}), and then integral next (Section~\ref{secint}), for a.e$.$ $t$. Because of the way
they are constructed, for a.e$.$ $t$, we know that the approximate mean curvatures
are $L^2$ bounded and they are almost minimizing in a small length scale.  
The latter gives a uniform lower density ratio bound for the limit measure
(Proposition~\ref{rect0}), and since the
$L^2$ norm of mean curvature is lower-semicontinuous under measure convergence, 
we are in a setting where Allard's rectifiability theorem applies. This gives 
rectifiability of the limit measure. Once this is done, we can focus on generic
points where the approximate tangent space exists. Since we only have a control 
of $L^2$ norms of approximate mean curvature, not the exact mean curvature, 
some extra information on a small length scale has to come in. This is provided by
small tilt excess and almost minimizing properties, which show that the hypersurfaces 
look like a finite number of layered hyperplanes in term of measure in a small 
length scale (Lemma~\ref{itglem}). 
This combined with some argument of Allard's compactness theorem
of integral varifold shows that the density of the limit flow is integer-valued wherever the
approximate tangent space exists. Since an approximate motion law is available,
we show the limit flow satisfies the exact motion law of Brakke flow (Section~\ref{secBra}).
We in addition
need to analyze the behavior of open partitions using Huisken's monotonicity formula
and the relative isoperimetric inequality in 
the end to make sure that the desired properties in Theorem \ref{sthm3} hold
(Section~\ref{secBas}). 
\section{Further preliminaries for construction of approximate flows}
\subsection{$\Omega$-finite open partition}
\begin{define}
\label{LFOP} 
A finite and ordered collection of sets
${\mathcal E}=\{E_i\}_{i=1}^{N}$ in ${\mathbb R}^{n+1}$ is called an $\Omega$-finite open partition of $N$ elements if
\begin{itemize}
\item[(a)] $E_1,\ldots,E_N$ are open and mutually disjoint,
\item[(b)] ${\mathcal H}^{n}\lfloor_{\Omega}(\mathbb R^{n+1}\setminus 
\cup_{i=1}^{N} E_i)<\infty$,
\item[(c)] $\cup_{i=1}^{N} \partial E_i$ is countably 
$n$-rectifiable.
\end{itemize}
The set of all $\Omega$-finite open partitions of $N$ elements is denoted by $\mathcal{OP}_{\Omega}^N$.
\end{define}
\begin{rem}
\label{LFOP1}
Since $\Omega(x)\geq e^{-\Cr{c_1}|x|}\Omega(0)$, (b) implies 
that, for any compact set $K\subset {\mathbb R}^{n+1}$, we have 
${\mathcal H}^{n}\lfloor_K
(\mathbb R^{n+1}\setminus \cup_{i=1}^{N} E_i)<\infty$. Also, this implies
\begin{equation}
\label{siqe}
\mathbb R^{n+1}\setminus\cup_{i=1}^N E_i=\cup_{i=1}^N\partial E_i.
\end{equation}
\end{rem}
Any open set $E\subset {\mathbb R}^{n+1}$ with ${\mathcal H}^{n}(
\partial E)<\infty$ has finite perimeter and $\|\nabla\chi_{E}\|
\leq {\mathcal H}^{n}\lfloor_{\partial E}$ (see \cite[Proposition 3.62]{AFP}). 
By De Giorgi's theorem, the reduced boundary of $E$ is 
countably $n$-rectifiable.
On the other hand, it may differ from the topological boundary 
$\partial E$ in general and the assumption
(c) is not redundant. 

Given ${\mathcal E}=\{E_i\}_{i=1}^{N}\in \mathcal{OP}_{\Omega}^N$, we define
\begin{equation}
\label{LFOP2}
\partial {\mathcal E}:=|\cup_{i=1}^{N} \partial E_i|\in {\bf IV}_{n}({\mathbb R}^{n+1})
\end{equation}
which is a unit density varifold induced naturally from the countably $n$-rectifiable 
set $\cup_{i=1}^N\partial E_i$. 
By (b), \eqref{siqe} and \eqref{LFOP2}, we have $\| \partial {\mathcal E} \|(\Omega)=\mathcal H^{n}
\lfloor_{\Omega}(\cup_{i=1}^N \partial E_i)<\infty$
for ${\mathcal E}\in \mathcal{OP}_{\Omega}^N$. 
\subsection{${\mathcal E}$-admissible function and its push-forward map $f_{\star}$}
\begin{define}
\label{LFOP3}
For ${\mathcal E}=\{E_i\}_{i=1}^{N} \in \mathcal{OP}_{\Omega}^N$, 
a function $f:{\mathbb R}^{n+1}\rightarrow
{\mathbb R}^{n+1}$ is called ${\mathcal E}$-admissible if it is Lipschitz and 
satisfies the following. 
Define ${\tilde E}_i:={\rm int}\,(f(E_i))$ for each $i$. Then
\begin{itemize}
\item[(a)] $\{{\tilde E}_i\}_{i=1}^{N}$ are mutually disjoint, 
\item[(b)] ${\mathbb R}^{n+1}\setminus \cup_{i=1}^{N} {\tilde E}_i\subset
f(\cup_{i=1}^{N} \partial E_i)$,
\item[(c)] $\sup_{x\in {\mathbb R}^{n+1}} |f(x)-x|<\infty$. 
\end{itemize}
\end{define}
\begin{lemma}
\label{LFOP4}
For ${\mathcal E}=\{E_i\}_{i=1}^{N}\in \mathcal{OP}_{\Omega}^N$, suppose
that $f$ is ${\mathcal E}$-admissible. 
Define $\tilde{\mathcal E}:=\{{\tilde E}_i\}_{i=1}^{N}$
with ${\tilde E}_i
:={\rm int}\, (f(E_i))$. Then we have $\tilde{\mathcal E}\in \mathcal{OP}_{\Omega}^N$.
\end{lemma}
{\it Proof}. We need to check that $\tilde{\mathcal E}$ satisfies Definition~\ref{LFOP} (a)-(c).
$\{\tilde E_{i}\}_{i=1}^{N}$ are open and mutually 
disjoint by Definition~\ref{LFOP3} (a). 
By Definition~\ref{LFOP3} (b), we have
\begin{equation}
\label{siba1}
\begin{split}
{\mathcal H}^{n}\lfloor_{\Omega}(\mathbb R^{n+1}\setminus \cup_{i=1}^{N} {\tilde E}_i) &
\leq {\mathcal H}^{n}\lfloor_{\Omega}(f(\cup_{i=1}^{N} \partial E_i))
 \leq ({\rm Lip}(f))^{n}\int_{\cup_{i=1}^{N} \partial E_i} \Omega(f(y))\, d
{\mathcal H}^{n}(y) \\
& \leq ({\rm Lip}(f))^{n} \exp(\Cr{c_1}\sup_{y\in {\mathbb R}^{n+1}}
|f(y)-y|) \, {\mathcal H}^{n}\lfloor_{\Omega}(\cup_{i=1}^{N} \partial E_i),
\end{split}
\end{equation}
where we used \eqref{omegaproeq}. The last quantity is finite due to 
Definition~\ref{LFOP} (b) and \eqref{siqe} for $\mathcal E$ and Definition~\ref{LFOP3} (c) for $f$. Since $\cup_{i=1}^{N}\partial E_i$ is
countably $n$-rectifiable, so is the Lipschitz image $f(\cup_{i=1}^{N}\partial E_i)$. Any subset of countably $n$-rectifiable set is again countably 
$n$-rectifiable, thus by Definition~\ref{LFOP3} (b) and \eqref{siqe} for 
$\tilde{\mathcal E}$, $\tilde{\mathcal E}$ satisfies
Definition~\ref{LFOP} (c) as well. This concludes the proof. 
\hfill{$\Box$}
\begin{define}
For $\mathcal E =\{E_i\}_{i=1}^{N}\in \mathcal{OP}_{\Omega}^N$ and 
$\mathcal E$-admissible function $f$, let $\mathcal{\tilde E}$ be defined as in
Lemma \ref{LFOP4}. We
define $\tilde{\mathcal E}$ to be the push-forward of 
$\mathcal E$ by $f$ and define
\begin{equation*}
f_{\star}{\mathcal E}:=\tilde{\mathcal E}\in \mathcal{OP}_{\Omega}^N.
\end{equation*}
\end{define}

Under the definition of $f_{\star}$, the unit density varifold
$\partial f_{\star}\mathcal E$ (cf.~\eqref{LFOP2}) is $|\cup_{i=1}^{N}
\partial\tilde E_i|$ and 
is in general different from the usual push-forward of varifold $f_{\sharp}
\partial\mathcal E=f_{\sharp}|\cup_{i=1}^{N}\partial E_i|$ in that
it does not count the multiplicity of image under the map. Moreover,
$\partial f_{\star}\mathcal E$ is not defined as the varifold induced from the
set $f(\cup_{i=1}^{N}\partial E_i)$ in general. For example,  if $({\rm int}\,f(E_i))\cap f(\partial E_i)$ is non-empty (whose possibility is not 
excluded by Definition~\ref{LFOP3}), 
it does not belong to 
$\cup_{i=1}^{N}\partial \tilde E_i$ and thus $f(\cup_{i=1}^{N}\partial E_i)
\neq \cup_{i=1}^{N}\partial \tilde E_i$ in this case. 
\subsection{Examples of $f_{\star}\mathcal E$}
\label{exfe}
It is worthwhile to see some simple examples of $\mathcal E$ and 
$\mathcal E$-admissible
functions to see what to expect. The choice of this particular admissible class 
characterizes general tendencies of what would happen to singularities. As 
we explain in Section \ref{arld}, we are interested in maps which reduce
$\mathcal H^n$ measure of boundaries. 
\subsubsection{Two lines crossing with four different open sets} \label{tlc4}
Consider the following Figure \ref{fig_1}, where two lines are intersecting, 
and $\mathcal E$ consists of four open sectors as shown. To reduce length of 
boundaries, one may consider
a Lipschitz map $f$ which vertically crashes triangle areas of $E_1$ and $E_3$
to a horizontal line segment, shrinks the neighboring areas next to them, 
and stretches some portion 
of $E_2$ and $E_4$ so that the map is Lipschitz. The map reduces the length of
boundaries, and also $\mathcal E$-admissible since $f(\cup_{i=1}^4 \partial E_i)
=\cup_{i=1}^4 \partial \tilde E_i$. This example indicates that, if we choose $f$
which locally reduces boundary measure, junctions of more than
three curves are likely to break up into triple junctions.
\begin{figure}[h]
\centering
\def\svgwidth{1\textwidth}
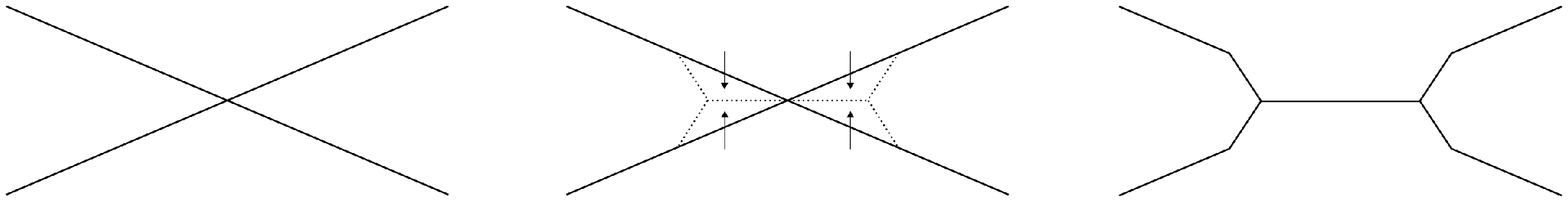
\caption{}
\label{fig_1}
\end{figure}
\subsubsection{Interior boundary} \label{ibry}
\begin{figure}[h]
\centering
\def\svgwidth{.9\textwidth}
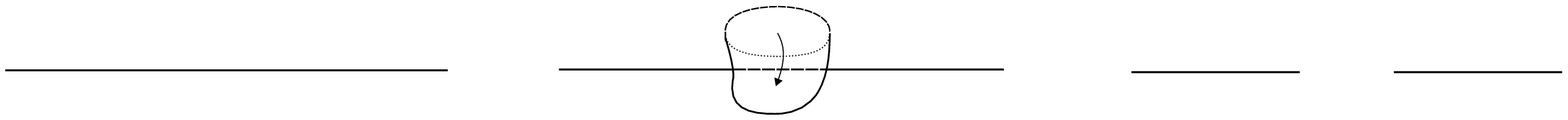
\caption{}
\label{fig_2}
\end{figure}

\begin{figure}[h]
\centering
\def\svgwidth{.8\textwidth}
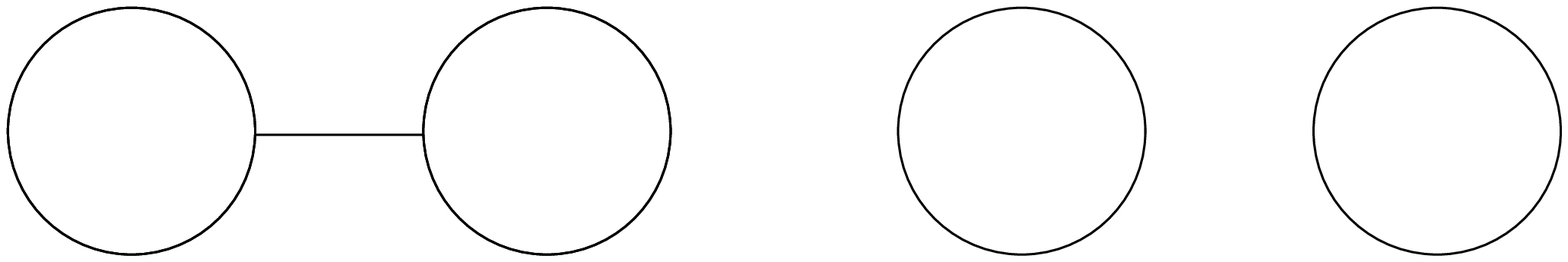
\caption{}
\label{fig_new}
\end{figure}

Suppose that we have only $E_1$ as shown in Figure \ref{fig_2}, which is the complement
of $x$-axis. For $f$, we may take a smooth map such that the 
dotted region of the second figure is stretched downwards to
hang over the lower half. Then the portion of $x$-axis covered by this map will be
interior points of the image of $E_1$ under $f$, thus we have $\tilde E_1$ as shown
in the third figure. By considering such ``stretching map'', we may even 
eliminate the whole $x$-axis with arbitrarily small deformation. This example indicates
that interior boundary is likely to be eliminated under measure reducing $f$. 
For this reason, as illustrated in Figure
\ref{fig_new}, if (a) is the initial data, the line segment connecting two circles is likely to vanish
instantly. 
\subsubsection{Two lines crossing with two different open sets}\label{tlc2}
\begin{figure}[h]
\centering
\def\svgwidth{.9\textwidth}
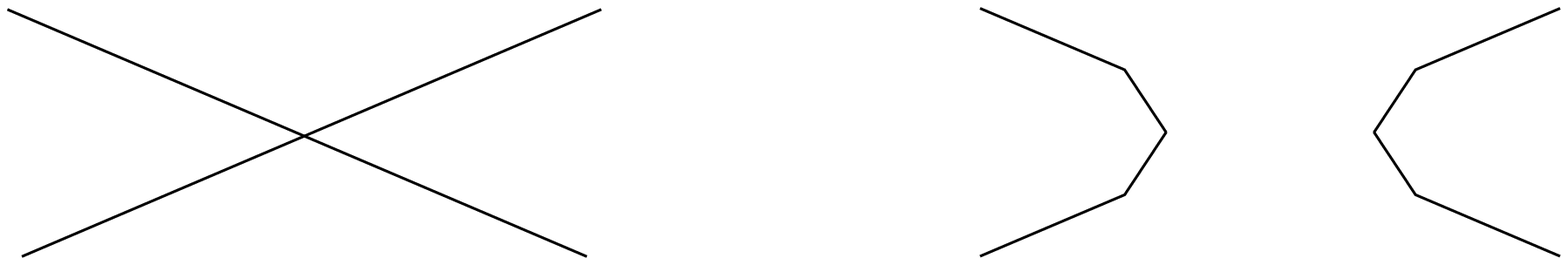
\caption{}
\label{fig_3}
\end{figure}

The next example is similar to \ref{tlc4}, but labeling is different as shown
in Figure \ref{fig_3}. By using a Lipschitz map of \ref{tlc4} and then composing a map 
of \ref{ibry} to eliminate the 
horizontal line segment appearing in Figure \ref{fig_1} (c), we can obtain Figure \ref{fig_3}
(b). Thus, depending on the combination of domains, we expect to have different
behaviors. 
\subsubsection{Radial projection}\label{rpn}
\begin{figure}[h]
\centering
\def\svgwidth{.9\textwidth}
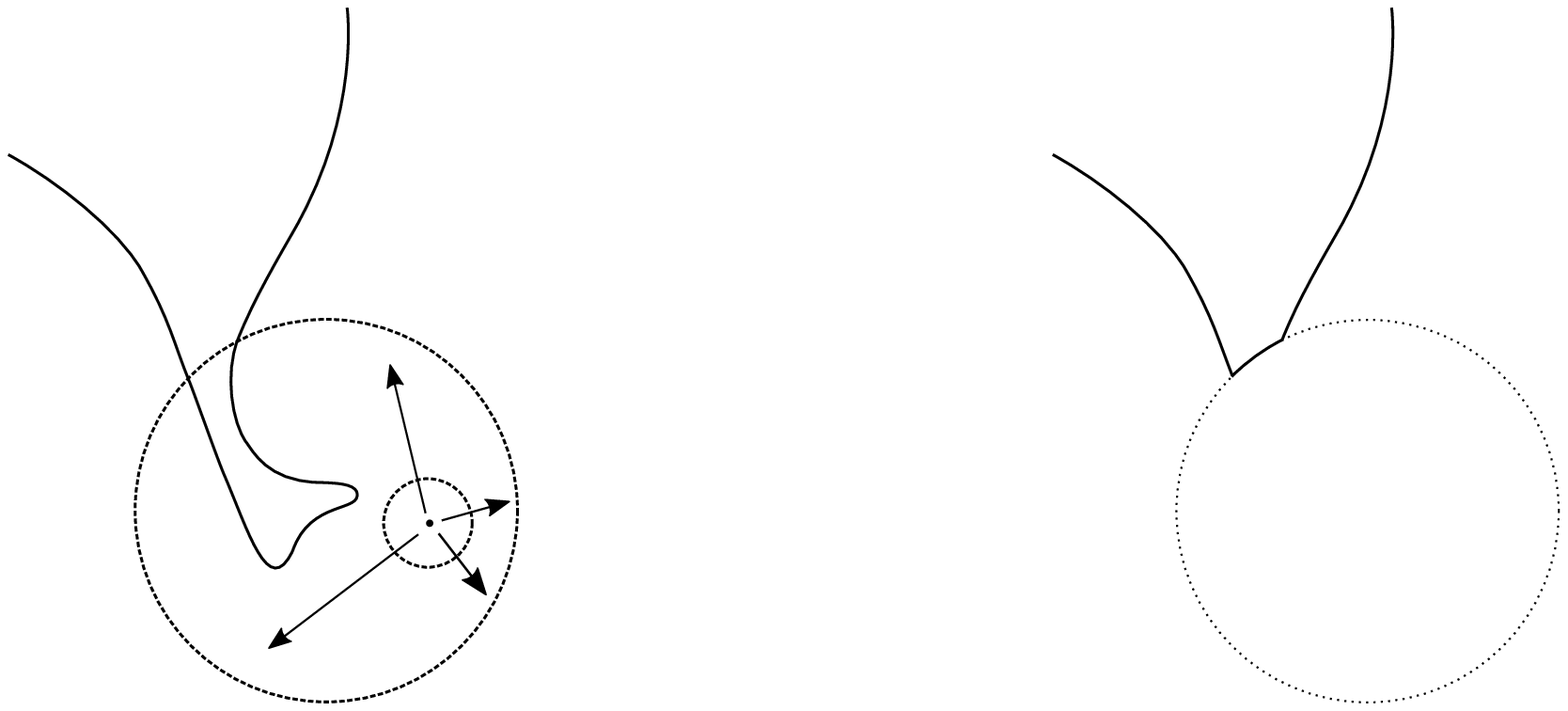
\caption{}
\label{fig_4}
\end{figure}

As in Figure \ref{fig_4}, consider a Lipschitz map which radially projects the annular region 
bounded by two dotted circles to the larger circle, with the trace of map being 
a radial line emanating from $x_0$. $f$ expands
the smaller disc to fill the larger disc one-to-one. Outside the larger disc, $f$ is
identity.  This map is $\mathcal E$-admissible since the new boundary is in 
$f(\partial E_1)$. Note that some portion of $f(\partial E_1)$ does not become
part of $\partial \tilde E_1$ because it is mapped to the interior of $\tilde E_2$. 
Depending on how much length there is inside the disc, the map
reduces the length. This type of projection map is used when we prove the rectifiability and integrality of the
limit flow. 
\subsection{Families ${\mathcal A}_j$ and ${\mathcal B}_j$ of test functions and vector fields}
We define sets of test functions ${\mathcal A}_j$ and vector fields ${\mathcal B}_j$
for $j\in {\mathbb N}$ as
\begin{equation}
\begin{split}
{\mathcal A}_j:=\{\phi\in& C^2({\mathbb R}^{n+1};{\mathbb R}^+)\,:\,
\phi(x)\leq \Omega(x),\,|\nabla\phi(x)|\leq j\,\phi(x),\, \\ &
\|\nabla^2\phi(x)\|\leq j\,\phi(x) \mbox{ for all }x\in {\mathbb R}^{n+1}\},
\end{split}
\label{adef}
\end{equation}
\begin{equation}
\begin{split}
{\mathcal B}_j:=&\{g\in C^2({\mathbb R}^{n+1};{\mathbb R}^{n+1})\,:\,
|g(x)|\leq j\Omega(x),\, \|\nabla g(x)\|\leq j \Omega(x),\\& \|\nabla^2 g(x)\|\leq j\Omega(x)
\mbox{ for all }x\in {\mathbb R}^{n+1}\mbox{ and }\|\Omega^{-1} g\|_{L^2 ({\mathbb R}^{n+1})}
\leq j\}.
\end{split}
\label{bdef}
\end{equation}
Note that $\Omega\in {\mathcal A}_{j}$ if $j\geq \max\{1,\Cr{c_1}\}$. 
Elements of ${\mathcal A}_j$
are strictly positive on ${\mathbb R}^{n+1}$ unless identically equal to $0$. For 
$V\in {\bf V}_n({\mathbb R}^{n+1})$ with $\|V\|(\Omega)<\infty$, we 
have $\|V\|(\phi)<\infty$ for $\phi\in {\mathcal A}_j$ since 
$\phi\leq \Omega$ from \eqref{adef}. 
For $g\in {\mathcal B}_j$, we naturally define
$\delta V(g)$ as 
\begin{equation*}
\delta V(g):=\int_{{\bf G}_n({\mathbb R}^{n+1})} S\cdot \nabla g(x)
\, dV(x,S)
\end{equation*} 
which is finite and well-defined due to $\|\nabla g\|\leq j\Omega$ of 
\eqref{bdef}.

Using \eqref{adef}, \eqref{omegaproeq} and \eqref{bdef}, the following can be 
seen easily. 
\begin{lemma}
For all $x,y\in {\mathbb R}^{n+1}$,
$j\in {\mathbb N}$ and $\phi\in {\mathcal A}_j$, we have
\begin{equation}
\phi(x)\leq \phi(y)\exp(j|x-y|),
\label{phiproa}
\end{equation}
\begin{equation}
|\phi(x)-\phi(y)|\leq j |x-y|\phi(x)\exp(j|x-y|),
\label{phipro}
\end{equation}
\begin{equation}
|\phi(x)-\phi(y)-\nabla\phi(y)\cdot(x-y)|\leq j|x-y|^2 \phi(y)\exp(j|x-y|).
\label{phipro2}
\end{equation}
\end{lemma}
\begin{lemma}
Let $\Cr{c_1}$ be as in \eqref{omegacon1}. Then for all $x,y\in {\mathbb R}^{n+1}$,
$j\in {\mathbb N}$ and $g\in {\mathcal B}_j$, we have
\begin{equation}
|g(x)-g(y)|\leq j |x-y|\Omega(x)\exp(\Cr{c_1}|x-y|).
\label{phiprog}
\end{equation}
\end{lemma}
As these inequalities indicate, within a small distance of order $1/j$, 
minimum and maximum values of $\phi$ are compatible up to some fixed constant
and this fact is used quite heavily in the following. 
\subsection{Area reducing Lipschitz deformation}
\label{arld}
\begin{define}
\label{mrld}
For ${\mathcal E}=\{E_i\}_{i=1}^N \in \mathcal{OP}_{\Omega}^N$
and $j\in{\mathbb N}$, define
${\bf E}(\mathcal E,j)$ to be the set of all $\mathcal E$-admissible
functions $f:{\mathbb R}^{n+1}\rightarrow{\mathbb R}^{n+1}$ such that
\begin{itemize}
\item[(a)] $|f(x)-x|\leq 1/j^2$ for all $x\in {\mathbb R}^{n+1}$,
\item[(b)] ${\mathcal L}^{n+1} (\tilde E_i\triangle E_i)\leq 1/j$ for all
$i=1,\ldots,N$ and where $\{\tilde E_i\}_{i=1}^N=f_{\star}
\mathcal E$,
\item[(c)] $\| \partial f_{\star} \mathcal E \|(\phi)\leq \| \partial\mathcal E \|(\phi)$
for all $\phi\in {\mathcal A}_j$.
\end{itemize}
\end{define}
${\bf E}(\mathcal E,j)$ includes the identity map $f(x)=x$, thus it is not
empty. We are interested in this class with large $j$, so that 
(a) and (b) restrict $f$ to be a very small deformation. 
Since $\Omega\in {\mathcal A}_j$ for all $j\geq\max\{1,\Cr{c_1}\}$, if $f\in 
{\bf E}(\mathcal E,j)$ with $j\geq \max\{1,\Cr{c_1}\}$, then we have
\begin{equation}
\label{mrld1}
\|\partial f_{\star}\mathcal E\|(\Omega)\leq \| \partial\mathcal E \|(\Omega).
\end{equation}
\begin{define}
\label{mrld2def}
For $\mathcal E\in \mathcal{OP}_{\Omega}^N$ and $j\in \mathbb N$,
we define
\begin{equation}
\label{mrld2}
\Delta_{j}\|\partial \mathcal E\| (\Omega)
:=\inf_{f\in {\bf E}(\mathcal E,j)} (\|\partial f_{\star}\mathcal E\|(\Omega)
-\|\partial\mathcal E\|(\Omega)).
\end{equation}
\end{define}
In addition, for localized deformations, we define
for a compact set $C\subset{\mathbb R}^{n+1}$
\begin{equation}
\label{mrld2.2}
{\bf E}(\mathcal E, C, j):=\{f\in {\bf E}(\mathcal E,j) : \{x : f(x)\neq x\}\cup
\{f(x) : f(x)\neq x\}\subset C\}
\end{equation} 
and 
\begin{equation}
\label{mrld3}
\Delta_{j}\|\partial \mathcal E\| (C)
:=\inf_{f\in {\bf E}(\mathcal E,C,j)} (\|\partial f_{\star}\mathcal E\|(C)
-\|\partial\mathcal E\|(C)).
\end{equation}
Since the identity map is in ${\bf E}(\mathcal E,j)$ and ${\bf E}(\mathcal E, C,j)$, 
$\Delta_j\|\partial\mathcal E\|(\Omega)$ and $\Delta_j\|\partial\mathcal E\|
(C)$ are 
always non-positive. They measure the extent to which $\|\partial \mathcal E\|$ 
can be reduced under the Lipschitz deformation in the $\mathcal E$-admissible
class. For $\mathcal E\in \mathcal{OP}_{\Omega}^N$ and $j\in \mathbb N$, 
we state their basic properties. 
\begin{lemma}
\label{escov}
For compact sets $C\subset \tilde C$, we have
\begin{equation}
\Delta_j\|\partial\mathcal E\|(\tilde C)\leq \Delta_j\|\partial\mathcal E\|(C)
\label{escov1}
\end{equation}
and 
\begin{equation}
\label{escov1.1}
\Delta_j\|\partial\mathcal E\|(\Omega)\leq (\max_C \Omega)
\{\Delta_j\|\partial\mathcal E\|(C)+(1-\exp(-\Cr{c_1}{\rm diam}\,C))
\|\partial\mathcal E\|(C)\}.
\end{equation}
\end{lemma}
{\it Proof}. By \eqref{mrld2.2}, ${\bf E}(\mathcal E,C,j)\subset {\bf E}(\mathcal E,
\tilde C, j)$. For any $f\in {\bf E}(\mathcal E,C,j)$, $\|\partial f_\star\mathcal E\|
(\tilde C)-\|\partial\mathcal E\|(\tilde C)=\|\partial f_\star\mathcal E\|
(C)-\|\partial\mathcal E\|(C)$ since $f\lfloor_{\tilde C\setminus C}$ is 
identity and $f(C)\subset C$. Then \eqref{escov1} follows from \eqref{mrld3}.
For \eqref{escov1.1}, take arbitrary $f\in {\bf E}(\mathcal E,C,j)$ and since
$f\in {\bf E}(\mathcal E,j)$, \eqref{mrld2} and \eqref{mrld2.2} give
\begin{equation}
\begin{split}
\Delta_j \|\partial\mathcal E\|(\Omega)&\leq \|\partial f_\star\mathcal E\|(\Omega)
-\|\partial\mathcal E\|(\Omega)=\|\partial f_\star\mathcal E\|\lfloor_C(\Omega)
-\|\partial\mathcal E\|\lfloor_C(\Omega) \\
&\leq (\max_C \Omega)\|\partial f_\star\mathcal E\|(C)-(\min_C \Omega)
\|\partial\mathcal E\|(C) \\
& \leq (\max_C \Omega)\{\|\partial f_\star\mathcal E\|(C)-\|\partial\mathcal E\|
(C)+(1-\exp(-\Cr{c_1}{\rm diam}\,C))\|\partial\mathcal E\|(C)\}
\end{split}
\end{equation}
where we used $(\min_C \Omega)/(\max_C\Omega)\geq \exp(-\Cr{c_1}
{\rm diam}\, C)$ which follows from \eqref{omegaproeq}. By taking $\inf$ over
${\bf E}(\mathcal E,C,j)$, we obtain \eqref{escov1.1}.
\hfill{$\Box$}
\begin{lemma}
\label{escov2}
Suppose that $\{C_i\}_{i=1}^{\infty}$ is a sequence of compact sets which are 
mutually disjoint and suppose that
$C$ is a compact set with $\cup_{i=1}^{\infty}C_i\subset C$ and $\mathcal L^{n+1}
(C)<1/j$. Then
\begin{equation}
\Delta_j\|\partial\mathcal E\|(C)\leq \sum_{i=1}^{\infty} \Delta_j\|\partial\mathcal 
E\|(C_i).
\label{escov3}
\end{equation}
\end{lemma}
{\it Proof}. By Lemma \ref{escov}, if $\Delta_j \|\partial\mathcal E\|(C)>-\infty$, 
then $\Delta_j\|\partial\mathcal E\|(C_i)>-\infty$ for all $i$. 
Let $m\in \mathbb N$ and $\e\in (0,1)$ be arbitrary. For all $i\leq m$, choose $f_i
\in {\bf E}(\mathcal E,C_i,j)$ such that $\Delta_j \|\partial\mathcal E\|(C_i)+\e
\geq \|\partial (f_i)_\star \mathcal E\|(C_i)-\|\partial\mathcal E\|(C_i)$. We define a
map $f:\mathbb R^{n+1}\rightarrow\mathbb R^{n+1}$ by setting $f\lfloor_{C_i}(x)=(f_i)\lfloor_{C_i}(x)$ and $f\lfloor_{\mathbb R^{n+1}\setminus \cup_{i=1}^m C_i}(x)=x$.
Since $\{C_i\}_{i=1}^m$ are disjoint, $f$ is well-defined, Lipschitz and 
$\mathcal E$-admissible. Using $\mathcal L^{n+1}(C)<1/j$, one checks that
$f\in {\bf E}(\mathcal E,C,j)$. Thus we have
\begin{equation}
\label{escov4}
\begin{split}
\Delta_j\|\partial\mathcal E\|(C)&\leq \|\partial f_\star\mathcal E\|(C)-\|\partial\mathcal
E\|(C)=\sum_{i=1}^{m}\|\partial (f_i)_\star\mathcal E\|(C_i)-\|\partial\mathcal E
\|(C_i) \\
&\leq m\e+\sum_{i=1}^m \Delta_j\|\partial\mathcal E\|(C_i).
\end{split}
\end{equation}
By letting $\e\rightarrow 0$ first and letting $m\rightarrow\infty$, we obtain \eqref{escov3}.
\hfill{$\Box$}

\begin{lemma}(\cite[4.10]{Brakke})
\label{eslip}
If $\mathcal E=\{E_i\}_{i=1}^N \in \mathcal{OP}_{\Omega}^N$, $j\in {\mathbb N}$,
$C$ is a compact set of ${\mathbb R}^{n+1}$, 
$f:{\mathbb R}^{n+1}\rightarrow{\mathbb R}^{n+1}$
is a $\mathcal E$-admissible function such that
\begin{itemize}
\item[(a)] $\{x : f(x)\neq x\} \cup \{ f(x) : f(x)\neq x\}\subset C$,
\item[(b)] $|f(x)-x|\leq 1/j^2$ for all $x\in {\mathbb R}^{n+1}$,
\item[(c)] ${\mathcal L}^{n+1} (\tilde E_i\triangle E_i)\leq 1/j$ for all
$i=1,\ldots,N$ and where $\{\tilde E_i\}_{i=1}^N=f_{\star}
\mathcal E$,
\item[(d)] $\| \partial f_{\star}\mathcal E\|(C)\leq \exp(-j\, {\rm diam}\,C) \, \| \partial 
\mathcal E \|(C)$,
\end{itemize}
then $f\in {\bf E}(\mathcal E,C,j)$.
\end{lemma}
{\it Proof}. We only need to check Definition~\ref{mrld} (c). By (a), $\|\partial f_{\star}
\mathcal E\|\lfloor_{{\mathbb R}^{n+1}\setminus C}=\|\partial{\mathcal E}\|\lfloor_{
{\mathbb R}^{n+1}\setminus C}$. Suppose $\phi\in {\mathcal 
A}_j$. Then by \eqref{phiproa}
\begin{equation*}
\begin{split}
\| \partial f_{\star}\mathcal E\|(\phi)-\| \partial \mathcal E\|(\phi)
&=\| \partial f_{\star}\mathcal E\|\lfloor_C(\phi)-\|  \partial \mathcal E\|\lfloor_C(\phi)\\
&\leq \max_C \phi\, \|\partial f_{\star}\mathcal E\|(C)-\min_C \phi\, \| \partial\mathcal
E\|(C) \\
&\leq \min_C \phi\, \big(\exp(j\,{\rm diam}\, C)\|\partial f_{\star}\mathcal E\|(C)- \|\partial\mathcal
E\|(C)) \leq 0
\end{split}
\end{equation*}
where (d) is used in the last line. 
\hfill{$\Box$}
\subsection{Smoothing function $\Phi_{\e}$}
Let $\psi\in C^{\infty}({\mathbb R}^{n+1})$ be a radially symmetric function such that
\begin{equation}
\begin{split}
&\psi(x)=1\mbox{ for }|x|\leq 1/2,\,\, \psi(x)=0\mbox{ for }|x|\geq 1, \\
&0\leq \psi(x)\leq 1,\,\,  |\nabla\psi(x)|\leq 3,\,\, \|\nabla^2\psi(x)\|\leq 9 \mbox{ for all }x\in
{\mathbb R}^{n+1}.
\end{split}
\label{defpsi}
\end{equation}
Define for each $\e\in(0,1)$
\begin{equation}
\hat\Phi_{\e}(x):=\frac{1}{(2\pi\e^2)^{\frac{n+1}{2}}}\exp\big(-\frac{|x|^2}{2\e^2}\big),\,\,
\Phi_{\e}(x):=c(\e)\psi(x)\hat\Phi_{\e}(x),
\label{defPsi}
\end{equation}
where the constant $c(\e)$ is chosen so that 
\begin{equation}
\int_{{\mathbb R}^{n+1}}\Phi_{\e}(x)\, dx=1.
\label{propPsi1}
\end{equation}
Since $\int_{{\mathbb R}^{n+1}}\hat\Phi_{\e}(x)\,dx=1$ for any $\e>0$
and $\hat\Phi_{\e}$ converges to the delta function as $\e\rightarrow 0+$,
there exists a constant $c(n)$ depending only on $n$ such that 
\begin{equation}
1< c(\e)\leq c(n)\mbox{ for }\e\in (0,1)\,\,\mbox{and } \lim_{\e\rightarrow 0+} c(\e)=1.
\label{propce}
\end{equation}
From the definitions of $\psi$ and $\Phi_{\e}$, we have the following estimates.
\begin{lemma}
There exists a constant $c$ depending only on $n$ such that, for $\e\in (0,1)$, we have
\begin{equation}
|\nabla\Phi_{\e}(x)|\leq \frac{|x|}{\e^2}\Phi_{\e}(x)+ c \chi_{B_1\setminus B_{1/2}}(x)
\exp(-\e^{-1}),
\label{Phies1}
\end{equation}
\begin{equation}
\|\nabla^2 \Phi_{\e}(x)\|\leq \frac{|x|^2}{\e^4} \Phi_{\e}(x)+\frac{c}{\e^2}\Phi_{\e}(x)
+c\chi_{B_1\setminus B_{1/2}} (x)\exp(-\e^{-1}).
\label{Phies2}
\end{equation}
\label{lemma1.2}
\end{lemma}
\begin{lemma}
With $c(\e)$ as in \eqref{defPsi}, we have
\begin{equation}
 x\,\Phi_{\e}(x)+\e^2 \nabla \Phi_{\e}(x)= \e^2 c(\e) \nabla\psi(x){\hat \Phi}_{\e}(x).
\label{Phies3}
\end{equation}
\label{lemma1.3}
\end{lemma}
The exponential smallness of the right-hand side of \eqref{Phies3}
will be of critical importance in Proposition~\ref{prop_mc2}.
\subsection{Smoothing of varifold \cite[4.3]{Brakke}}
In this subsection, we consider a smoothing of varifold and derive various
estimates. For general distribution $T$, there is a notion of smoothing of $T$
using a duality, i.e., $\Phi_{\e}\ast T(\phi)=T(\Phi_{\e}\ast \phi)$ for any 
$\phi\in C^{\infty}_c(\mathbb R^{n+1})$.
Here, given a varifold $V
\in {\bf V}_n(\mathbb R^{n+1})$, we smooth out with respect to only
the space variables and not
the Grassmannian part. 
\begin{define}
For $V\in {\bf V}_n({\mathbb R}^{n+1})$, we define $\Phi_{\e}\ast V\in {\bf V}_n({\mathbb R}^{n+1})$ through 
\begin{equation}
(\Phi_{\e}\ast V)(\phi):=V(\Phi_{\e}\ast\phi)=\int_{{\bf G}_n({\mathbb R}^{n+1})} 
\int_{{\mathbb R}^{n+1}} \phi(x-y,S)\Phi_{\e}(y)\, dy\,dV(x,S)
\label{defPhist}
\end{equation}
for $\phi\in C_c({\bf G}_n({\mathbb R}^{n+1}))$. 
\end{define}
If $\|V\|(\Omega)<\infty$, we have $\|\Phi_{\e}\ast V\|(\Omega)
<\infty$ since
\begin{equation}
\|\Phi_{\e}\ast V\|(\Omega)\leq \int_{{\bf G}_n({\mathbb R}^{n+1})}
\int_{{\mathbb R}^{n+1}} e^{\Cr{c_1}} \Omega(x)\Phi_{\e}(y)\,dy\,dV(x,S)
= e^{\Cr{c_1}} \|V\|(\Omega)
\label{toma3}
\end{equation}
by \eqref{omegaproeq} and \eqref{propPsi1}. Thus we 
have $\|\Phi_{\e}\ast V\|(\phi)<\infty$ for $\phi\in {\mathcal A}_j$ as
well. 
For a general Radon measure $\mu$ on ${\mathbb R}^{n+1}$, we similarly define 
a Radon measure $\Phi_{\e}\ast \mu$. $\Phi_{\e}\ast \mu$ may be identified with a smooth
function on ${\mathbb R}^{n+1}$ via the $L^2$ inner product, because, for $\phi\in C_c({\mathbb R}^{n+1})$, 
\begin{equation}
\begin{split}
(\Phi_{\e}\ast \mu)(\phi)& =\int_{{\mathbb R}^{n+1}}\int_{{\mathbb R}^{n+1}}
\phi(y) \Phi_{\e}(x-y)\, dy\,d\mu(x) \\ &=\int_{{\mathbb R}^{n+1}} \phi(y)\int_{{\mathbb R}^{n+1}}
\Phi_{\e}(x-y)\, d\mu (x) \,dy=<\Phi_{\e}\ast\mu,\phi>_{L^2(\mathbb R^{n+1})},
\end{split}
\label{musm0}
\end{equation}
and we may identify $\Phi_{\e}\ast\mu\in C^{\infty}({\mathbb R}^{n+1})$ with
\begin{equation}
(\Phi_{\e}\ast \mu)(x)=\int_{{\mathbb R}^{n+1}} \Phi_{\e}(y-x)\, d\mu(y).
\label{musm}
\end{equation}

In a similar way, for general $V\in {\bf V}_n({\mathbb R}^{n+1})$, we may
define $\Phi_{\e}\ast \delta V$ as a $C^{\infty}$ vector field as follows. 
Note that $V$ may not have a 
bounded first variation in general. 
 For $g\in C^1_c({\mathbb R}^{n+1};
{\mathbb R}^{n+1})$, $\Phi_{\e}\ast\delta V$ should be defined to satisfy
\begin{equation}
\begin{split}
\int_{{\mathbb R}^{n+1}} (\Phi_{\e}\ast \delta V)(x)\cdot g(x)\, dx&=\delta V(\Phi_{\e}\ast g) 
 = \int_{{\bf G}_n({\mathbb R}^{n+1})} S\cdot( (\nabla\Phi_{\e}\ast g)(x))\, dV(x,S) \\
& = \int_{{\mathbb R}^{n+1}} g(y)\cdot
\int_{{\bf G}_n({\mathbb R}^{n+1})} S(\nabla\Phi_{\e}(x-y))\, dV(x,S)dy.
\end{split}
\label{musm1}
\end{equation}
The equality \eqref{musm1} motivate the definition of $\Phi_{\e}\ast\delta V$ as 
a $C^{\infty}$ vector field
\begin{equation}
(\Phi_{\e}\ast \delta V)(x):=\int_{{\bf G}_n({\mathbb R}^{n+1})} S(\nabla\Phi_{\e}(y-x))\, dV(y,S).
\label{musm2}
\end{equation}
\begin{lemma}
For $V\in {\bf V}_n({\mathbb R}^{n+1})$, we have
\begin{equation}
\Phi_{\e}\ast\|V\|=\|\Phi_{\e}\ast V\|,
\label{musm3}
\end{equation}
\begin{equation}
\Phi_{\e}\ast \delta V=\delta (\Phi_{\e}\ast V).
\label{musm4}
\end{equation}
\end{lemma}
{\it Proof}. For $\phi\in C_c({\mathbb R}^{n+1})$, we have
\begin{equation}
\begin{split}
\int_{{\mathbb R}^{n+1}} \phi\, d\|\Phi_{\e}\ast V\| &=\int_{{\bf G}_n({\mathbb R}^{n+1})}
\phi(x)\, d(\Phi_{\e}\ast V)(x,S) \\
&=\int_{{\bf G}_n({\mathbb R}^{n+1})}(\Phi_{\e}\ast \phi )(x) \,dV(x,S)\,\,\,(\mbox{by \eqref{defPhist}}) \\
&=\int_{{\mathbb R}^{n+1}}\phi(y)\int_{{\mathbb R}^{n+1}} \Phi_{\e}(x-y)\, d\|V\|(x)\, dy  \\
&= \int_{{\mathbb R}^{n+1}} \phi\, d(\Phi_{\e}\ast\|V\|)\,\,\,(\mbox{by \eqref{musm0}}).
\end{split}
\label{musm5}
\end{equation}
Thus we proved \eqref{musm3}. For $g\in C_c^1({\mathbb R}^{n+1};{\mathbb R}^{n+1})$,
by \eqref{musm1}, 
\begin{equation}
(\Phi_{\e}\ast\delta V)(g)=\delta V(\Phi_{\e}\ast g)=\int_{{\bf G}_n({\mathbb R}^{n+1})} S\cdot (\Phi_{\e}
\ast \nabla g)(x)\, dV(x,S)
\label{musm6}
\end{equation}
while by \eqref{defPhist}, 
\begin{equation}
\delta(\Phi_{\e}\ast V)(g)=\int_{{\bf G}_n({\mathbb R}^{n+1})} \Phi_{\e}\ast (S\cdot \nabla g)(x)\, dV(x,S).
\label{musm7}
\end{equation}
Since $\Phi_{\e}\ast$ commutes with $S\,\cdot$, \eqref{musm6} and \eqref{musm7} prove \eqref{musm4}.
\hfill{$\Box$}

The following is
used when we need to deal with error terms in the next section. 
\begin{lemma}
\label{toma}
For $V\in{\bf V}_n({\mathbb R}^{n+1})$ with $\|V\|(\Omega)<\infty$ 
and for all $x\in {\mathbb R}^{n+1}$ and $r>0$, we have
\begin{equation}
\Omega(x)\|V\|(B_r(x))\leq e^{\Cr{c_1} r}\|V\|(\Omega),
\label{toma1}
\end{equation}
\begin{equation}
\int_{{\mathbb R}^{n+1}} \Omega(x)\|V\|(B_r(x))\, dx\leq \omega_{n+1}
e^{\Cr{c_1}r}r^{n+1} \|V\|(\Omega).
\label{toma2}
\end{equation}
\end{lemma}
{\it Proof}. 
By \eqref{omegaproeq}, for $y\in B_r(x)$, we have $\Omega(x)\leq \Omega(y)e^{\Cr{c_1}r}$,
thus
\begin{equation*}
\Omega(x)\|V\|(B_r(x))\leq \int_{B_r(x)}\Omega(y)e^{\Cr{c_1}r}\, d\|V\|(y)\leq e^{\Cr{c_1}r}\|V\|(\Omega),
\end{equation*}
proving \eqref{toma1}. Similarly, since $\chi_{B_r(x)}(y)=\chi_{B_r(y)}(x)$,
\begin{equation*}
\begin{split}
\int_{\mathbb R^{n+1}} \Omega(x)\|V\|(B_r(x))\, dx&=\int_{{\mathbb R}^{n+1}}\int_{{\mathbb R}^{n+1}}\Omega(x)\chi_{B_r(x)}(y)\,dx\,d\|V\|(y)\\
&=\int_{{\mathbb R}^{n+1}}\int_{ B_r(y)}\Omega(x)\,dx\,d\|V\|(y)\\
&\leq \omega_{n+1} e^{\Cr{c_1}r}r^{n+1}  \int_{\mathbb R^{n+1}} \Omega(y)\, d\|V\|(y)=\omega_{n+1}
e^{\Cr{c_1}r}r^{n+1} \|V\|(\Omega),
\end{split}
\end{equation*}
proving \eqref{toma2}. 
\hfill{$\Box$} 
\section{Smoothed mean curvature vector $h_{\e}(\cdot,V)$}
\label{SMC}
Given $V\in {\bf V}_n({\mathbb R}^{n+1})$, if the first variation $\delta V$ is bounded and 
absolutely continuous with respect to $\|V\|$, the Radon-Nikodym derivative
$h(\cdot,V)=-\delta V/\|V\|$ defines the generalized mean curvature vector of $V$
as in \eqref{fvf}. 
Here, even for $V$ with unbounded first variation, 
we want to have a smooth analogue of $h(\cdot,V)$ to construct an 
approximate mean curvature flow. Thus we define a {\it smoothed mean curvature vector}
$h_{\e}(\cdot,V)$ for $\e\in (0,1)$ by
\begin{equation}
h_{\e}(\cdot,V):=-\Phi_{\e}\ast \big(\frac{\Phi_{\e}\ast \delta V}{\Phi_{\e}\ast \|V\|
+\e\Omega^{-1}}\big).
\label{curvature}
\end{equation}
We may often write $h_{\e}(\cdot,V)$ as $h_{\e}$ for simplicity. 
Note that this is a well-defined smooth vector field; 
since $\Omega^{-1}\geq 1$ by \eqref{omegacon1},
the denominator is strictly positive. 
Formally, as $\e\rightarrow 0+$, $h_{\e}$ will be more and more concentrated around
${\rm spt}\,\|V\|$ and 
we expect that $h_{\e}(\cdot,V)$ converges in a suitable sense to $h(\cdot,V)$, 
as long as there are some suitable bounds. 
The term ``smoothed mean curvature vector'' is used in \cite{Brakke}, but we should 
warn the reader that it may happen that the generalized mean curvature 
$h(\cdot,V)$ may not exist in general while $h_{\e}(\cdot,V)$ is always well-defined.
We also point out that there is a difference from \cite{Brakke} that we 
have the extra $\e\Omega^{-1}$ term to avoid division by 0 (see \cite[p.39]{Brakke}). In \cite{Brakke}, 
$\Phi_{\e}\ast\|V\|$ (with a different and more complicated $\Phi_{\e}$,
see \cite[p.37]{Brakke}) is prepared so that it is everywhere positive 
on ${\mathbb R}^{n+1}$ unless $\|V\|(\Omega)=0$. Though it is a simple 
modification, various computations are clearly tractable compared to \cite{Brakke}.
After some reading, one must admit that the corresponding computations in \cite{Brakke} 
are discouragingly difficult to follow in the original form. 
In the following, we also use the notation
\begin{equation}
\tilde{h}_{\e}:=-\frac{\Phi_{\e}\ast\delta V}{\Phi_{\e}\ast\|V\|+\e\Omega^{-1}}
\label{curvature2}
\end{equation}
for simplicity and note that $h_{\e}=\Phi_{\e}\ast{\tilde h}_{\e}$.
\subsection{Rough pointwise estimates on $h_{\e}(\cdot,V)$}
\begin{lemma}
There exists a
constant $\Cl[eps]{e_1}\in(0,1)$ depending only on $n$, $\Cr{c_1}$ and $M$ with the following property. 
Suppose $V\in {\bf V}_n({\mathbb R}^{n+1})$ with $\|V\|(\Omega)\leq M$ and $\e\in (0,\Cr{e_1})$.
Then, for all $x\in {\mathbb R}^{n+1}$, we have
\begin{equation}
|{\tilde h}_{\e}(x,V)|\leq 2\e^{-2},\hspace{.5cm}|h_{\e}(x,V)|\leq 2\e^{-2},
\label{hestimate1}
\end{equation}
\begin{equation}
\|\nabla h_{\e}(x,V)\|\leq 2\e^{-4},
\label{hestimate2}
\end{equation}
\begin{equation}
\|\nabla^2 h_{\e}(x,V)\|\leq 2\e^{-6}.
\label{hestimate3}
\end{equation}
\end{lemma}
{\it Proof}. First by \eqref{musm2} and \eqref{Phies1}, we have
\begin{equation}
\begin{split}
|(\Phi_{\e}\ast \delta V)(x)|&\leq \int_{B_1(x)} \frac{|y-x|}{\e^2}\Phi_{\e}(y-x)+c(n) \exp(-\e^{-1}) \, d\|V\|(y) \\
&\leq \e^{-2}(\Phi_{\e}\ast\|V\|)(x)+c(n)\exp(-\e^{-1}) \|V\|(B_1(x)),
\end{split}
\label{hestimate4}
\end{equation}
where $c(n)$ is as in Lemma~\ref{lemma1.2}. 
Combining \eqref{hestimate4} and \eqref{toma1}, we obtain
\begin{equation}
\frac{|\Phi_{\e}\ast\delta V|}{\Phi_{\e}\ast\|V\|+\e\Omega^{-1}}
\leq \e^{-2}+c(n)M\e^{-1} \exp(\Cr{c_1}-\e^{-1}).
\label{hestimate6}
\end{equation}
Choose $\Cr{e_1}$ so that $c(n)M\e\exp(\Cr{c_1}-\e^{-1})\leq 1$ if 
$\e\in (0,\Cr{e_1})$. 
Now recalling $\Phi_{\e}\ast 1=1$ and \eqref{curvature}, we obtain \eqref{hestimate1} 
from \eqref{hestimate6}. For \eqref{hestimate2}, we note that $|\nabla\Phi_{\e}|\ast 1
\leq \e^{-2}+c(n)\exp(-\e^{-1})\omega_n$ by \eqref{Phies1}. Thus using \eqref{hestimate6}
and choosing an appropriate $\Cr{e_1}$, we obtain \eqref{hestimate2}. Using \eqref{Phies2}, 
we similarly obtain \eqref{hestimate3}. 
\hfill{$\Box$}

The following quantity plays the role of $\Omega$-weighted ``approximate $L^2$-norm''
of smoothed mean curvature vector. The reason is that, roughly speaking, we expect
that
\begin{equation*}
\int |h_{\e}(\cdot, V)|^2\, d\|V\|\approx \int \frac{|\Phi_{\e}\ast\delta V|^2}{(\Phi_{\e}\ast\|V\|+\e\Omega^{-1})^2}\,
d(\Phi_{\e}\ast\|V\|)\approx \int \frac{|\Phi_{\e}\ast\delta V|^2}{\Phi_{\e}\ast\|V\|+\e\Omega^{-1}}\,dx.
\end{equation*}
\begin{lemma}
For $V\in {\bf V}_n({\mathbb R}^{n+1})$ with $\|V\|(\Omega)<\infty$ and $\e\in (0,\epsilon_1)$,
\begin{equation*}
\int_{{\mathbb R}^{n+1}} \frac{|\Phi_{\e}\ast \delta V|^2\Omega}{\Phi_{\e}
\ast\|V\|+\e\Omega^{-1}}\, dx<\infty.
\end{equation*}
\end{lemma}
{\it Proof}. The claim follows from \eqref{hestimate6}, 
\eqref{hestimate4}, \eqref{toma2}, \eqref{toma3} and \eqref{musm3}. 
\hfill{$\Box$}
\subsection{$L^2$ approximations}
This subsection establishes various error estimates of approximations. 
\begin{prop}
There exists a constant $\Cl[eps]{e_2}\in (0,1)$ depending only on $n$, $\Cr{c_1}$ and
$M$ such that, for any $g\in {\mathcal B}_j$, 
$V\in {\bf V}_n({\mathbb R}^{n+1})$ with $\|V\|(\Omega)\leq M$, $j\in {\mathbb N}$, $\e\in (0,\Cr{e_2})$ with 
\begin{equation}
j\leq \frac12 \e^{-\frac{1}{6}},
\label{jere}
\end{equation}
we have
\begin{equation}
\label{curmain}
\Big|\int_{{\mathbb R}^{n+1}} h_{\e}\cdot g\, d\|V\|+\int_{{\mathbb R}^{n+1}} (\Phi_{\e}\ast\delta V)
\cdot g\, dy\Big| \leq \e^{\frac14}
\big(\int_{{\mathbb R}^{n+1}} 
\frac{|\Phi_{\e}\ast\delta V|^2\Omega}{\Phi_{\e}\ast\|V\|+\e\Omega^{-1}}\, dy\big)^{
\frac12}.
\end{equation}.
\label{prop_mc}
\end{prop}
Note that one can draw an analogy between \eqref{curmain} and \eqref{fvf}.

{\it Proof}. By \eqref{curvature} and \eqref{curvature2}, we have
\begin{equation}
\int_{{\mathbb R}^{n+1}}h_\e\cdot g\, d\|V\|=\int_{{\mathbb R}^{n+1}}
(\Phi_{\e}\ast {\tilde h}_{\e})\cdot g\, d\|V\| 
=\int_{{\mathbb R}^{n+1}} {\tilde h}_{\e}(y)\cdot \int_{{\mathbb R}^{n+1}} \Phi_{\e}(\cdot-y)g(\cdot)\, d\|V\|dy.
\label{cur0}
\end{equation}
We may also rewrite using the notation \eqref{curvature2}
\begin{equation}
\int_{{\mathbb R}^{n+1}} (\Phi_{\e}\ast\delta V)\cdot g\, dy=-\int_{{\mathbb R}^{n+1}}
{\tilde h}_{\e}(\Phi_{\e}\ast\|V\|+\e\Omega^{-1})\cdot g\, dy.
\label{cur1}
\end{equation}
Summing \eqref{cur0} and \eqref{cur1}, we obtain
\begin{equation}
\begin{split}
\big|&\int_{{\mathbb R}^{n+1}} h_{\e}\cdot g\, d\|V\|+\int_{{\mathbb R}^{n+1}} (\Phi_{\e}\ast\delta V)
\cdot g\, dy\big| \leq \int_{{\mathbb R}^{n+1}} |g(y)||{\tilde h}_{\e}(y,V)|\e\Omega^{-1}(y)\, dy \\
&+\int_{{\mathbb R}^{n+1}} |{\tilde h}_{\e}(y,V)| \Big|\int_{{\mathbb R}^{n+1}}\Phi_{\e}(x-y)g(x)
\, d\|V\|(x)- (\Phi_{\e}\ast\|V\|)(y)g(y)\Big|\, dy=: I_1+I_2.
\end{split}
\label{cur2}
\end{equation}
By H\"{o}lder's inequality and \eqref{bdef},
\begin{equation}
I_1\leq\e  \big(\int_{{\mathbb R}^{n+1}} |g|^2\Omega^{-2}\, dy\big)^{\frac12}
\big(\int_{{\mathbb R}^{n+1}}|{\tilde h}_{\e}|^2 \, dy\big)^{\frac12}
\leq j\e\big(\int_{{\mathbb R}^{n+1}}|{\tilde h}_{\e}|^2 \, dy\big)^{\frac12}.
\label{cur5}
\end{equation}
Recalling \eqref{curvature2}, \eqref{cur5} in particular gives 
\begin{equation}
I_1\leq j\e^{\frac{1}{2}}\big(\int_{{\mathbb R}^{n+1}} \frac{|\Phi_{\e}\ast
\delta V|^2\Omega}{\Phi_{\e}\ast
\|V\|+\e\Omega^{-1}}\, dy\big)^{\frac12}.
\label{cur7}
\end{equation}
For $I_2$, using \eqref{phiprog} for $g\in {\mathcal B}_j$, 
\begin{equation}
\begin{split}
\Big|\int_{{\mathbb R}^{n+1}} \Phi_{\e}(\cdot-y)g(\cdot)\, d\|V\|-(\Phi_{\e}\ast\|V\|)g\Big|
&= \Big|\int_{{\mathbb R}^{n+1}} (g(x)-g(y))\Phi_{\e}(x-y)\, d\|V\|(x)\Big|\\
&\leq je^{\Cr{c_1}}\Omega(y)\int_{B_1(y)} |x-y|\Phi_{\e}(x-y)\, d\|V\|(x).
\end{split}
\label{cur3}
\end{equation}
Using the property of $\Phi_{\e}$ being exponentially small away from 
the origin, we have 
\begin{equation}
\sup_{x\in B_1(y)\setminus B_{\sqrt{\e}}(y)}|x-y|\Phi_{\e}(x-y)\leq c(n)\e^{-n-1}\exp(-(2\e)^{-1})
=: c_{\e}.
\label{cur6}
\end{equation} 
Thus \eqref{cur3} and \eqref{cur6} give 
\begin{equation}
I_2\leq je^{\Cr{c_1}}\e^{\frac12}\int_{{\mathbb R}^{n+1}}\Omega \,|{\tilde h}_{\e}|
(\Phi_{\e}\ast \|V\|)\, dy+je^{\Cr{c_1}}c_{\e}\int_{{\mathbb R}^{n+1}}\Omega\, |{\tilde h}_{\e}|\,\|V\|(B_1(y))\, dy=: I_{2,a}+I_{2,b}.
\label{cur4}
\end{equation}
For $I_{2,a}$, use H\"{o}lder's inequality to obtain
\begin{equation}
I_{2,a}\leq je^{\Cr{c_1}}\e^{\frac12}\big(\int_{{\mathbb R}^{n+1}} |{\tilde h}_{\e}|^2 
(\Phi_{\e}\ast\|V\|)\Omega\, dy\big)^{\frac12}\big( (\Phi_{\e}\ast\|V\|)(\Omega)\big)^{\frac12}.
\label{cur8}
\end{equation}
Substitution of \eqref{toma3} (with \eqref{musm3}) into \eqref{cur8} gives
\begin{equation}
I_{2,a}\leq je^{2\Cr{c_1}} \e^{\frac12}\big(\int_{{\mathbb R}^{n+1}} \frac{|\Phi_{\e}\ast \delta V|^2
\Omega}{
\Phi_{\e}\ast\|V\|+\e\Omega^{-1}}\, dy\big)^{\frac12}M^{\frac12}.
\label{cur10}
\end{equation}
For $I_{2,b}$, by H\"{o}lder's inequality, 
\begin{equation}
I_{2,b}\leq je^{\Cr{c_1}} c_{\e} \big(\int_{\mathbb R^{n+1}} 
\frac{|\Phi_{\e}\ast\delta V|^2\Omega}{\Phi_{\e}\ast\|V\|+\e\Omega^{-1}}\, dy\big)^{\frac12}
\big(\int_{\mathbb R^{n+1}} \frac{\|V\|(B_1(y))^2\Omega}{\Phi_{\e}\ast\|V\|+\e
\Omega^{-1}}\, dy\big)^{\frac12}.
\label{cur11}
\end{equation}
Using \eqref{toma1}, we have
\begin{equation}
\label{cur12}
\int_{\mathbb R^{n+1}}\frac{\|V\|(B_1(y))^2\Omega}{\Phi_{\e}\ast\|V\|+\e\Omega^{-1}}\,
dy\leq \e^{-1} e^{\Cr{c_1}} M\int_{\mathbb R^{n+1}}\|V\|(B_1(y))\Omega\, dy.
\end{equation}
Then \eqref{cur11}, \eqref{cur12} and \eqref{toma2} prove
\begin{equation}
\label{cur14}
I_{2,b}\leq je^{2\Cr{c_1}}c_{\e}\e^{-\frac12}\omega_{n+1}^{\frac12} M
\big(\int_{\mathbb R^{n+1}}\frac{|\Phi_{\e}\ast\delta V|^2\Omega}{\Phi_{\e}\ast
\|V\|+\e\Omega^{-1}}\, dy\big)^{\frac12}.
\end{equation}
Combining \eqref{cur2}, \eqref{cur7}, \eqref{cur4}, \eqref{cur10}, \eqref{cur14}, \eqref{jere} and 
choosing $\Cr{e_2}$ appropriately depending only on $n, \Cr{c_1}$ and $M$, we obtain \eqref{curmain}.
\hfill{$\Box$}
\begin{prop}
\label{prop_mc2}
There exists a constant $\Cl[eps]{e_3}\in (0,1)$ depending only on $n$, $\Cr{c_1}$ and $M$ with
the following property. 
For $V\in {\bf V}_n({\mathbb R}^{n+1})$ with $\|V\|(\Omega)\leq M$, $j\in {\mathbb N}$, $\phi\in {\mathcal A}_j$
and $\e\in (0,\Cr{e_3})$ with \eqref{jere},
we have
\begin{equation}
\label{mc0}
\Big|\delta V(\phi h_{\e})+\int_{{\mathbb R}^{n+1}}
\frac{\phi|\Phi_{\e}\ast \delta V|^2}{\Phi_{\e}\ast \|V\|+\e\Omega^{-1}}\,
dx\Big|
\leq \e^{\frac14}\big(\int_{{\mathbb R}^{n+1}}
\frac{\phi|\Phi_{\e}\ast \delta V|^2}{\Phi_{\e}\ast \|V\|+\e\Omega^{-1}}\,dx
+1
\big)
\end{equation}
and 
\begin{equation}
\label{mc01}
\int_{{\mathbb R}^{n+1}} |h_{\e}|^2 \phi\, d\|V\| \leq \int_{{\mathbb R}^{n+1}}
\frac{\phi|\Phi_{\e}\ast\delta V|^2}{\Phi_{\e}\ast\|V\|+\e\Omega^{-1}}\,(1+
\e^{\frac14})\, dx+\e^{\frac14}.
\end{equation}
\end{prop}
Note that \eqref{mc0} measures a deviation from $\delta V(\phi h)=-\int \phi |h|^2
\, d\|V\|$, which is \eqref{fvf} with $g=\phi h$ if all quantities are well-defined. 
We use \eqref{mc01} when we prove the lower semicontinuity of $L^2$-norm of
mean curvature vector.

{\it Proof}. From the definition of the first variation, we have
\begin{equation}
\label{mc1}
\begin{split}
\delta V(\phi h_{\e})&=\int_{{\bf G}_n({\mathbb R}^{n+1})} \nabla(\phi h_{\e})\cdot S\, dV(\cdot,S)
=\int_{{\bf G}_n({\mathbb R}^{n+1})} (\phi \nabla h_{\e}
+\nabla\phi\otimes h_{\e})\cdot S\, dV(\cdot,S) \\
&=\int_{{\bf G}_n({\mathbb R}^{n+1})}\int_{{\mathbb R}^{n+1}} (\phi(x)\nabla\Phi_{\e}(x-y)
+\nabla\phi(x)\Phi_{\e}(x-y))\otimes {\tilde h}_{\e}(y)\cdot S\, dydV(x,S)
\end{split}
\end{equation}
and by \eqref{musm2}, 
\begin{equation}
\label{mc2}
\begin{split}
\int_{{\mathbb R}^{n+1}}
&\frac{\phi|\Phi_{\e}\ast \delta V|^2}{\Phi_{\e}\ast \|V\|+\e\Omega^{-1}}\,
dx=-\int_{{\mathbb R}^{n+1}} \phi {\tilde h}_{\e}\cdot (\Phi_{\e}\ast \delta V)\, dy \\
&=-\int_{{\mathbb R}^{n+1}}\int_{{\bf G}_n({\mathbb R}^{n+1})}
\phi(y) S(\nabla\Phi_{\e}(x-y))\cdot {\tilde h}_{\e}(y)\, dV(x,S)dy.
\end{split}
\end{equation}
By summing \eqref{mc1} and \eqref{mc2}, we obtain
\begin{equation}
\begin{split}
&\delta V(\phi h_{\e})+\int_{{\mathbb R}^{n+1}}\frac{\phi|\Phi_{\e}\ast \delta V|^2}{\Phi_{\e}\ast \|V\|+\e\Omega^{-1}}\,
dx \\&
=\int_{{\mathbb R}^{n+1}}\int_{{\bf G}_n({\mathbb R}^{n+1})}
\big((\phi(x)-\phi(y))S(\nabla\Phi_{\e}(x-y)) +\Phi_{\e}(x-y)S(\nabla\phi(x))
\big)
\, dV(x,S)\cdot {\tilde h}_{\e}(y)\, dy.
\end{split}
\label{mc3}
\end{equation}
To continue, we carry out a second order approximation of $\phi$ and interpolate 
the right-hand side of \eqref{mc3} by defining (all integrations are over ${\mathbb R}^{n+1}
\times {\bf G}_n({\mathbb R}^{n+1})$)
\begin{equation}
\begin{split}
&I_1:= \iint\big(\phi(x)-\phi(y)-\nabla\phi(y)
\cdot(x-y))S(\nabla\Phi_{\e}(x-y)\big)\, dV(x,S)\cdot {\tilde h}_{\e}(y)\, dy, \\
&I_2:= \iint
\Phi_{\e}(x-y)S(\nabla\phi(x)-\nabla\phi(y))\, dV(x,S)\cdot {\tilde h}_{\e}(y)\, dy, \\
&I_3:=\iint \nabla\phi(y)\cdot(x-y)S(\nabla \Phi_{\e}(x-y))+\Phi_{\e}(x-y)S(\nabla\phi(y))\, dV(x,S)
\cdot {\tilde h}_{\e}(y)\, dy
\end{split}
\label{mc4}
\end{equation}
so that $I_1+I_2+I_3$ equals to \eqref{mc3}. In addition, we define
\begin{equation}
I_4:=-\e^2 \iint S[\nabla_x(\nabla\phi(y)\cdot\nabla\Phi_{\e}(x-y))]\, dV(x,S)\cdot
{\tilde h}_{\e}(y)\, dy,
\label{mc4.5}
\end{equation}
where $\nabla_x$ indicates (for clarity) that the differentiation is with respect to
$x$ variables. In the following, we estimate $I_1$, $I_2$, $I_3-I_4$ and $I_4$. 
\newline
{\it Estimate of $I_1$}. We use \eqref{phipro2} to squeeze out a $|x-y|^2$ term 
to deal with $\e^{-2}$ term coming from $\nabla \Phi_{\e}$. Then we separate the
domain of integration to $B_{\e^{\frac56}}(y)$ and the complement. On the latter, 
$\Phi_{\e}(\cdot-y)$ is exponentially small with respect to $\e$. 
With this in mind, we have by \eqref{phipro2} and \eqref{Phies1} that
\begin{equation}
\label{mc5}
\begin{split}
|I_1|&\leq j\int(|{\tilde h}_{\e}|\phi)(y)\int e^{j|\cdot-y|} |\cdot-y|^2\big(\frac{|\cdot-y|}{\e^{2}}\Phi_{\e}(\cdot-y)+c(n) e^{-\e^{-1}}\chi_{B_1(y)}\big)\, d\|V\|dy  \\
&\leq je^{j\e^{\frac56}}\e^{\frac12} \int(|{\tilde h}_{\e}|\phi)(y)
\int \Phi_{\e}(\cdot-y)\, d\|V\|dy
\,\, (\frac{|x-y|^3}{\e^2}\leq \e^{\frac12}\mbox{ on }
B_{\e^{\frac56}}(y)\mbox{ is used}) \\
& + je^{j}c(n) \e^{-n-3}e^{-\frac{\e^{-\frac13}}{2}}\int_{{\mathbb R}^{n+1}}
\|V\|(B_1(y))|{\tilde h}_{\e}(y)|\Omega(y)\, dy\\
&+ je^{j} c(n)e^{-\e^{-1}} \int_{{\mathbb R}^{n+1}} \|V\|(B_1(y)) |{\tilde h}_{\e}(y)|
\Omega(y)\, dy.
\end{split}
\end{equation}
The integration of the first term of \eqref{mc5} may be estimated as
\begin{equation}
\label{mc6}
\begin{split}
\int |{\tilde h}_{\e}|\phi\int\Phi_{\e}(\cdot-y)\, d\|V\|dy&=\int_{{\mathbb R}^{n+1}}
(\Phi_{\e}\ast \|V\|)|{\tilde h}_{\e}|\phi\, dy\\
&\leq \big((\Phi_{\e}\ast\|V\|) (\Omega)\big)^{\frac12}
\big(\int_{{\mathbb R}^{n+1}} (\Phi_{\e}\ast\|V\|)|{\tilde h}_{\e}|^2\phi\, dy\big)^{\frac12}\\
&\leq \big(e^{\Cr{c_1}}M\big)^{\frac12}
\big(\int_{{\mathbb R}^{n+1}} (\Phi_{\e}\ast\|V\|)|{\tilde h}_{\e}|^2\phi\, dy\big)^{\frac12}
\end{split}
\end{equation}
where we used \eqref{toma3} and \eqref{musm3}. 
Use \eqref{hestimate1} and \eqref{toma2} for the second and third 
terms of \eqref{mc5}. Combined with \eqref{mc6}, then, we have some 
$c$ depending only on $\Cr{c_1}$, $M$ and $n$
such that
\begin{equation}
\label{mc9}
\begin{split}
|I_1|&\leq je^{j\e^{\frac56}} \e^{\frac12} (e^{\Cr{c_1}}M)^{\frac12} \big(\int_{{\mathbb R}^{n+1}}
(\Phi_{\e}\ast\|V\|)|{\tilde h}_{\e}|^2\phi\, dy\big)^{\frac12}
+jc e^{j-\e^{-\frac16}} \\
&\leq j \e^{\frac12} \int_{{\mathbb R}^{n+1}}\frac{\phi |\Phi_{\e}\ast
\delta V|^2}{\Phi_{\e}\ast \|V\|+\e\Omega^{-1}}\,dy+jc\e^{\frac12}+jc e^{-\frac12
\e^{-\frac16}},
\end{split}
\end{equation}
where we also used \eqref{jere}. 
\newline
{\it Estimate of $I_2$}.
By the similar manner, we estimate $I_2$. Note that $\nabla\Phi_{\e}$ is not present
while we have only $|\nabla\phi(x)-\nabla\phi(y)|\leq j|x-y|\phi(x)e^{j|x-y|}$ this time. 
We separate the domain of integration to $B_{\e^{\frac12}}(y)$ and the complement, 
and estimate just like $I_1$ to obtain \eqref{mc9} for $I_2$ in place of $I_1$. We 
omit the detail since it is repetitive. 
\newline
{\it Estimate of $I_3-I_4$}. The first point is that the integrand
with respect to $V$ of $I_3$
can be expressed as
\begin{equation}
\label{mc10}
\begin{split}
\nabla\phi(y)&\cdot(x-y)S(\nabla\Phi_{\e}(x-y))+\Phi_{\e}(x-y)S(\nabla \phi(y))\\
&=S[\nabla\phi(y)\Phi_{\e}(x-y)+\nabla\phi(y)\cdot (x-y)\nabla\Phi_{\e}(x-y)] \\
&=S[\nabla_{x}((x-y)\cdot\nabla\phi(y)\Phi_{\e}(x-y))].
\end{split}
\end{equation}
The function $(x-y)\Phi_{\e}(x-y)$ may be replaced by
$-\e^2\nabla\Phi_{\e}(x-y)$ with exponentially small error due to \eqref{Phies3}.
So we first check that this replacement produces small error indeed. By \eqref{mc10},
\begin{equation}
\label{mc11}
\begin{split}
I_3-I_4= \iint S[\nabla_x (\nabla\phi(y)\cdot c(\e)\e^2\nabla\psi(x-y){\hat\Phi}_{\e}(x-y))]\, dV(x,S)
\cdot{\tilde h}_{\e}(y)\, dy.
\end{split}
\end{equation}
On the support of $\nabla\psi$, ${\hat \Phi}_{\e}$ is of the order of $e^{-\e^{-2}}$,
thus estimating as in the second and third terms of \eqref{mc5}, we obtain from \eqref{mc11} and \eqref{adef} that
\begin{equation}
\label{mc12a}
\big|I_3-I_4\big|\leq j c(n,\Cr{c_1},M) e^{-\e^{-1}} . 
\end{equation}
{\it Estimate of $I_4$}. 
To be clear about the indices, the $i$-th component of the integrand of $I_4$ 
with respect to $V$ is (the same indices imply summation over $1$ to $n+1$)
\begin{equation}
S_{ij}\nabla_{x_j}(\nabla_{y_l}\phi(y)\nabla_{x_l}\Phi_{\e}(x-y))
=-\nabla_{y_l}\phi(y) \nabla_{y_l}(S_{ij}\nabla_{x_j}\Phi_{\e}(x-y)).
\label{mc12}
\end{equation}
Recalling \eqref{musm2} and writing the $i$-th component of $\Phi_{\e}
\ast\delta V$ as $(\Phi_{\e}\ast\delta V)_{i}$, \eqref{mc12} shows
\begin{equation}
\label{mc13}
I_4=\e^2 \int_{{\mathbb R}^{n+1}} \nabla\phi \cdot \nabla(\Phi_{\e}\ast \delta V)_i
({\tilde h}_{\e})_i\, dy=-\frac{\e^2}{2}\int_{{\mathbb R}^{n+1}} 
\frac{\nabla\phi\cdot \nabla |\Phi_{\e}\ast \delta V|^2}{\Phi_{\e}\ast\|V\|+\e
\Omega^{-1}}\, dy.
\end{equation}
Here, we want to carry out  one integration by parts for $I_4$. 
Let $\psi_r$ be a cut-off function such that $\psi_r(x)=1$ for 
$x\in B_{r/2}$, $\psi_r(x)=0$ for $x\in \mathbb R^{n+1}\setminus B_r$ and $|\nabla\psi_r(x)|\leq 3/r$. For example, with $\psi$ defined in \eqref{defpsi}, we may
set $\psi_r(x):=\psi(x/r)$. Then we have
\begin{equation}
\begin{split}
&\int_{\mathbb R^{n+1}} \frac{\nabla\phi\cdot\nabla|\Phi_{\e}\ast\delta V|^2}{\Phi_{\e}\ast\|V\|+\e\Omega^{-1}}
\, dy=\lim_{r\rightarrow\infty} \int_{\mathbb R^{n+1}} \psi_{r} \frac{\nabla\phi\cdot\nabla|\Phi_{\e}\ast\delta V|^2}{\Phi_{\e}\ast\|V\|+\e\Omega^{-1}}
\, dy \\
& =-\int_{\mathbb R^{n+1}}\nabla\cdot\big(\frac{\nabla\phi}{\Phi_{\e}\ast
\|V\|+\e\Omega^{-1}}\big)|\Phi_{\e}\ast\delta V|^2\, dy -\lim_{r\rightarrow\infty} \int_{\mathbb R^{n+1}} 
\frac{(\nabla\psi_r\cdot\nabla\phi)|\Phi_{\e}\ast\delta V|^2}{
\Phi_{\e}\ast\|V\|+\e\Omega^{-1}}\, dy.
\end{split}
\label{mc13s1}
\end{equation}
For the second term of \eqref{mc13s1}, we use \eqref{curvature2}, \eqref{hestimate1} and \eqref{adef} to 
obtain 
\begin{equation}
\label{mc13s2}
\Big|\int_{\mathbb R^{n+1}} \frac{(\nabla\psi_r\cdot\nabla\phi)|\Phi_{\e}\ast\delta V|^2}{
\Phi_{\e}\ast\|V\|+\e\Omega^{-1}}\, dy\Big|\leq 2j \e^{-2} \int_{\mathbb R^{n+1}} |\nabla\psi_r||\Phi_{\e}
\ast\delta V|\Omega\, dy.
\end{equation}
By \eqref{hestimate4} and also noticing $(\Phi_{\e}\ast\|V\|)(x)\leq c(n,\e) \|V\|(B_1(x))$, with a suitable
constant $c(n,\e)$, we have 
\begin{equation}
\label{mc13s3}
\int_{\mathbb R^{n+1}} |\nabla\psi_r||\Phi_{\e}\ast\delta V|\Omega\, dy
\leq \frac{c(n,\e)}{r} \int_{B_{r}\setminus B_{r/2}} \|V\|(B_1(x)) \Omega(x)\, dx.
\end{equation}
By \eqref{mc13s1}-\eqref{mc13s3} and \eqref{toma2}, we may justify the integration by parts for
$I_4$ on $\mathbb R^{n+1}$. Hence,
\begin{equation}
\begin{split}
|I_4|&=\big|\frac{\e^2}{2}\int_{{\mathbb R}^{n+1}} \nabla\cdot\big(\frac{\nabla\phi}{\Phi_{\e}\ast
\|V\|+\e\Omega^{-1}}\big)|\Phi_{\e}\ast\delta V|^2\, dy \big|\\
& \leq \frac{\e^2}{2}\int_{{\mathbb R}^{n+1}} \big(\frac{((n+1)j+\Cr{c_1}j)\phi}
{\Phi_{\e}\ast\|V\| +\e\Omega^{-1}} +\frac{j\phi |\nabla\Phi_{\e}
\ast\|V\||}{(\Phi_{\e}\ast\|V\|+\e\Omega^{-1})^2}\big)|\Phi_{\e}\ast\delta
V|^2
\, dy,
\end{split}
\label{mc14}
\end{equation}
where we also used $|\Delta \phi|\leq (n+1)j\phi$ and $\e\Omega^{-2}|\nabla\phi
\cdot \nabla\Omega|(\Phi_{\e}\ast\|V\|+\e\Omega^{-1})^{-1}\leq \Cr{c_1}j
\phi$ due to \eqref{omegacon1} and \eqref{adef}. To estimate the second
term of \eqref{mc14}, we have
\begin{equation}
\begin{split}
|\nabla\Phi_{\e}\ast\|V\| (y)|&\leq \int_{{\mathbb R}^{n+1}} |\nabla\Phi_{\e}(x-y)|\,
d\|V\|(x)\\
&\leq \int_{{\mathbb R}^{n+1}} \frac{|x-y|}{\e^2}\Phi_{\e}(x-y)\, d\|V\|(x)
+ce^{-\e^{-1}}\|V\|(B_1(y))\mbox{ (by \eqref{Phies1})} \\
&\leq \e^{-\frac32} \Phi_{\e}\ast\|V\|(y)+ce^{-\e^{-\frac12}}\|V\|(B_1(y))
\end{split}
\label{mc15}
\end{equation}
where we split the integration of the first term into $B_{\e^\frac12}(y)$ and
the complement as in the case of $I_1$. By substituting \eqref{mc15}
into \eqref{mc14} and recalling estimates 
\eqref{toma2} and \eqref{hestimate6}, with a suitable constant $c$
depending only on $\Cr{c_1}$, $M$ and $n$, we obtain
\begin{equation}
|I_4|\leq cj\e^{\frac12}\int_{{\mathbb R}^{n+1}} 
\frac{\phi\,|\Phi_{\e}\ast\delta V|^2}{\Phi_{\e}\ast\|V\|+\e\Omega^{-1}}
\, dy+cje^{-\e^{-\frac16}}.
\label{mc16}
\end{equation}
Combining \eqref{mc9}, remark for the estimate of $I_2$, \eqref{mc12a},
\eqref{mc16} and \eqref{jere}, we obtain \eqref{mc0}  by suitably restricting $\epsilon_3$. 

For the proof of \eqref{mc01}, by \eqref{curvature2} and $h_{\e}=\Phi_{\e}
\ast{\tilde h}_{\e}$, we have
\begin{equation}
\label{mc17}
\begin{split}
\int_{{\mathbb R}^{n+1}} |h_{\e}|^2\phi\, d\|V\|&=\int_{{\mathbb R}^{n+1}} |\Phi_{\e}\ast
{\tilde h}_{\e}|^2\phi\, d\|V\|\leq \int_{{\mathbb R}^{n+1}} \phi\, (\Phi_{\e}\ast
 |{\tilde h}_{\e}|^2)\, d\|V\| \\
&=\int_{{\mathbb R}^{n+1}} |{\tilde h}_{\e}(y)|^2\int_{{\mathbb R}^{n+1}} \phi(x)\Phi_{\e}(x-y)\, d\|V\|(x)dy.
\end{split}
\end{equation}
We then use \eqref{phipro} to conclude
\begin{equation}
\begin{split}
&\int_{{\mathbb R}^{n+1}} \phi(x)\Phi_{\e}(x-y)\, d\|V\|(x)\\
&\leq \phi(y)(\Phi_{\e}\ast\|V\|)(y)+j\phi(y)\int_{{\mathbb R}^{n+1}} e^{j|x-y|}
|x-y|\Phi_{\e}(x-y)\, d\|V\|(x)
\end{split}
\label{mc18}
\end{equation}
while the last term of \eqref{mc18} may be estimated by separating the 
integration over $B_{\e^{\frac12}}(y)$ and the complement as
\begin{equation}
\int_{{\mathbb R}^{n+1}}e^{j|x-y|}|x-y|\Phi_{\e}(x-y)\, d\|V\|(x)\leq \e^{\frac12} e^{j\e^{\frac12}}
(\Phi_{\e}\ast\|V\|)(y)+c(n) e^{j-\e^{-\frac12}}\|V\|(B_1(y)).
\label{mc19}
\end{equation}
Substitutions of \eqref{mc18} and \eqref{mc19} into \eqref{mc17} (and 
use \eqref{adef} and \eqref{jere}) give
\begin{equation}
\int_{{\mathbb R}^{n+1}} |h_{\e}|^2\phi\, d\|V\|\leq \int_{{\mathbb R}^{n+1}}
|{\tilde h}_{\e}|^2 \big\{(\Phi_{\e}\ast\|V\|)\phi(1+je \e^{\frac12} )
+j c(n)e^{-\frac12 \e^{-\frac12}}\Omega(y)\|V\|(B_1(y))\big\} \,dy.
\label{mc20}
\end{equation}
Since $|{\tilde h}_{\e}|^2\leq 4\e^{-4}$ by \eqref{hestimate1}, the last term of \eqref{mc20}
may be bounded by $j c(n,\Cr{c_1},M) \e^{-4} e^{-\frac12\e^{-\frac12}}$, also using \eqref{toma2}. 
 By choosing an appropriate $\Cr{e_3}$ depending only on $n, \Cr{c_1}$ and $M$, 
 and again using \eqref{jere}, we 
 obtain \eqref{mc01}. 
\hfill{$\Box$}
\begin{prop}
There exists $\Cl[eps]{e_4}\in (0,1)$ depending only on $n$, $\Cr{c_1}$ and $M$ with 
the following property. Suppose $V\in {\bf V}_n({\mathbb R}^{n+1})$ with $\| V\|(\Omega)\leq M$,
$\e\in (0,\Cr{e_4})$, $g\in {\mathcal B}_j$ and $j\in {\mathbb N}$ satisfying
\eqref{jere}. Then we have
\begin{equation}
\left|\int_{{\mathbb R}^{n+1}} h_{\e}\cdot g\, d\|V\|+\delta V(g)\right|\leq \e^{\frac14}
+\e^{\frac14}\big(\int_{{\mathbb R}^{n+1}} \frac{|\Phi_{\e}\ast\delta V|^2\Omega}{
\Phi_{\e}\ast\|V\|+\e\Omega^{-1}}\, dx\big)^{\frac12}.
\label{mc21}
\end{equation}
\label{cl5}
\end{prop}
{\it Proof}. By \eqref{musm1} and a similar estimate as \eqref{phiprog} 
for $\nabla g$, we have
\begin{equation}
\begin{split}
&\Big|\int_{{\mathbb R}^{n+1}} (\Phi_{\e}\ast\delta V)\cdot g\, dy-\delta V(g)\Big|
=|\delta V(\Phi_{\e}\ast g)-\delta V(g)| \\
&\leq \int_{{\mathbb R}^{n+1}}
|\nabla(\Phi_{\e}\ast g)-\nabla g|\, d\|V\|\leq cj\int_{{\mathbb R}^{n+1}}
\int_{{\mathbb R}^{n+1}} |x-y|\Phi_{\e}(x-y)\Omega(x)\, d\|V\|(x)dy \\
&\leq cj \e^{\frac12}\|V\|(\Omega),
\end{split}
\label{mc23}
\end{equation}
where we estimated as in \eqref{cur4} and $c$ is a constant depending only 
on $n$ and $\Cr{c_1}$.  Combining \eqref{curmain}, \eqref{mc23}, \eqref{jere}
and restricting $\Cr{e_4}\leq \Cr{e_2}$ depending only on $n, \Cr{c_1}$ and $M$ further, 
we obtain \eqref{mc21}.
\hfill{$\Box$}
\label{eltwo}
\subsection{Curvature of limit}
By the estimates in the previous subsection, we obtain the following
\begin{prop} 
\label{col}
Suppose that we have $\{V_j\}_{j=1}^{\infty}\subset {\bf V}_{n}
(\mathbb R^{n+1})$ with 
\begin{enumerate}
\item $\sup_{j} \|V_j\|(\Omega)<\infty$,
\item
$\liminf_{j\rightarrow\infty} \int_{\mathbb R^{n+1}}
\frac{|\Phi_{\e_j}\ast\delta V_j|^2\Omega}{\Phi_{\e_j}\ast\|V_j\|+\e_j
\Omega^{-1}}\, dx<\infty$,
\item
$\lim_{j\rightarrow\infty} \e_j=0$.
\end{enumerate}
Then there exists a converging subsequence $\{V_{j_l}\}_{l=1}^{\infty}$, and 
the limit $V\in {\bf V}_n(\mathbb R^{n+1})$ has a generalized
mean curvature $h(\cdot,V)$ with
\begin{equation}
\int_{\mathbb R^{n+1}} |h(\cdot,V)|^2\phi\, d\|V\|
\leq \liminf_{l\rightarrow\infty}\int_{\mathbb R^{n+1}}
\frac{|\Phi_{\e_{j_l}}\ast\delta V_{j_l}|^2\phi}{\Phi_{\e_{j_l}}\ast\|V_{j_l}\|+\e_{j_l}
\Omega^{-1}}\, dx
\label{cl3}
\end{equation}
for any $\phi\in \cup_{i\in \mathbb N} \mathcal A_i$.
\end{prop}
{\it Proof}. 
By (1), we may choose a subsequence $\{V_{j_l}\}_{l=1}^{\infty}$
converging to a limit $V\in {\bf V}_n(\mathbb R^{n+1})$ and so that 
the integrals in (2) are uniformly bounded for this subsequence as well. 
Fix $\phi\in \mathcal A_i$ and consider a Hilbert space
\begin{equation*}
X_\phi: =\big\{g=(g_1,\ldots,g_{n+1}) ; g \in L_{loc}^2 (\|V\|),\, \int_{\mathbb R^{n+1}}|g|^2 \phi^{-1}\, d\|V\|
<\infty\big\}
\end{equation*}
equipped with inner product $(f,g)_{X_{\phi}}:=\int_{\mathbb R^{n+1}}
f\cdot g\, \phi^{-1}\, d\|V\|$. Recall that $\phi>0$ on $\mathbb R^{n+1}$, and
$C_c^{\infty}(\mathbb R^{n+1}; \mathbb R^{n+1})$ is a dense subspace in $X_{\phi}$. 
Fix arbitrary $g\in C_c^{\infty}(\mathbb R^{n+1};\mathbb R^{n+1})$. Corresponding
to $g$, there exists $j'\in \mathbb N$ such that $g\in \mathcal B_{j'}$. By Proposition \ref{cl5} with $j=j'$ and 
combined with (1) and (2), we have
\begin{equation}
\lim_{l\rightarrow\infty} \delta V_{j_l}(g)=- \lim_{l\rightarrow\infty}
\int_{\mathbb R^{n+1}} h_{\e_{j_l}}(\cdot,V_{j_l})\cdot g\, d\|V_{j_l}\|.
\label{cl6}
\end{equation}
The left-hand side is equal to $\delta V(g)$ by the varifold convergence.
For $\phi\in \mathcal A_i$, we have by \eqref{mc01} (with $j=i$) and (2)
that, writing $h_{\e_{j_l}}=h_{\e_{j_l}}(\cdot,V_{j_l})$,
\begin{equation}
\begin{split}
-\lim_{l\rightarrow\infty}\int_{\mathbb R^{n+1}} & h_{\e_{j_l}}\cdot g\, d\|V_{j_l}\|
\leq \liminf_{l\rightarrow\infty}\Big(\int_{\mathbb R^{n+1}}|h_{\e_{j_l}}|^2\phi\, d\|V_{j_l}\|\Big)^{\frac12}
\Big(\int_{\mathbb R^{n+1}} |g|^2\phi^{-1}\, d\|V_{j_l}\|\Big)^{\frac12} \\
&\leq \Big(\liminf_{l\rightarrow\infty} \int_{\mathbb R^{n+1}} \frac{\phi|\Phi_{\e_{j_l}}
\ast\delta V_{j_l}|^2}{\Phi_{\e_{j_l}}\ast\|V_{j_l}\|+\e_{j_l}\Omega^{-1}}\, dx
\Big)^{\frac12}\Big(\int_{\mathbb R^{n+1}} |g|^2\phi^{-1}\, d\|V\|\Big)^{\frac12}.
\end{split}
\label{cl7}
\end{equation}
Writing the first term on the right-hand side of \eqref{cl7} as $C_0$, 
\eqref{cl6} and \eqref{cl7} show
\begin{equation}
\delta V(g)\leq C_0 \|g\|_{X_\phi}
\label{cl8}
\end{equation}
for any $g\in C_c^{\infty}(\mathbb R^{n+1};\mathbb R^{n+1})$. By a density argument,
$\delta V$ may be uniquely extended as 
a bounded linear functional on $X_\phi$. By the Riesz
representation theorem, there exists a unique $f\in X_\phi$ with $\|f\|_{X_\phi}\leq C_0$ such that 
$\delta V(g)=(f,g)_{X_\phi}$ for all 
$g\in X_\phi$. 
Then, note that $-f\phi^{-1}$ is the generalized mean curvature 
$h(\cdot,V)$, and \eqref{cl3} is equivalent to $\|f\|_{X_\phi}\leq C_0$.
\hfill{$\Box$}

\subsection{Motion by smoothed mean curvature}
This subsection establishes an approximate motion law when a varifold
is moved by the smoothed mean curvature vector. 
\begin{prop} 
\label{propaf0}
There exists $\Cl[eps]{e_5}\in (0,1)$ depending only on $n$, $\Cr{c_1}$ and $M$ with
the following. Suppose $V\in {\bf V}_n({\mathbb R}^{n+1})$ with $\|V\|(\Omega)\leq M$,
$j\in {\mathbb N}$, 
$\phi\in {\mathcal A}_j$, $\e\in (0,\Cr{e_5})$ with \eqref{jere}, $\Delta t\in(2^{-1}\e^{\Cl[c]{c_a}}, \e^{\Cr{c_a}}]$, where we set 
\begin{equation}
\label{defca}
\Cr{c_a}:=3n+20.
\end{equation}
 Define 
\begin{equation*}
f(x):=x+h_{\e}(x,V)\, \Delta t .
\end{equation*}
Then we have
\begin{equation}
\Big|\frac{\|f_{\sharp}V\|(\phi)-\|V\|(\phi)}{\Delta t}-\delta(V,\phi)(h_{\e}(\cdot,V))\Big|
\leq \e^{\Cr{c_a}-10},
\label{af0}
\end{equation}
\begin{equation}
\label{af1}
\frac{\|f_{\sharp}V\|(\Omega)-\|V\|(\Omega)}{\Delta t}
+\frac14 \int_{{\mathbb R}^{n+1}} \frac{|\Phi_{\e}\ast\delta V|^2\, \Omega}{\Phi_{\e}
\ast\|V\| + \e\Omega^{-1}}\,dx\leq 3\e^{\frac14}+\frac{\Cr{c_1}^2}{2}
\|V\|(\Omega).
\end{equation}
Moreover, if $\|f_{\sharp} V\|(\Omega)\leq M$, then we have
\begin{equation}
\label{af2}
\big| \delta (V,\phi)(h_{\e}(\cdot,V)) -\delta(f_{\sharp}V,\phi)(h_{\e}(\cdot,f_{\sharp}V))\big|
\leq \e^{\Cr{c_a}-2n-19},
\end{equation}
\begin{equation}
\label{af3}
\Big|\int_{{\mathbb R}^{n+1}} \frac{|\Phi_{\e}\ast\delta V|^2\Omega}{\Phi_{\e}\ast
\|V\|+\e\Omega^{-1}}\, dx-\int_{{\mathbb R}^{n+1}} \frac{|\Phi_{\e}\ast \delta(f_{\sharp}
V)|^2\Omega}{\Phi_{\e}\ast\| f_{\sharp}V\|+\e\Omega^{-1}}\, dx\Big|
\leq \e^{\Cr{c_a}-3n-18}.
\end{equation}
\end{prop}
{\it Proof}.
For simplicity, write
$F(x):=  f(x)-x=h_{\e}(x,V)\Delta t$. We have
\begin{equation}
\label{est_F2}
|F(x)|=|h_{\e}(x,V)|\Delta t\leq  2 \e^{\Cr{c_a}-2}
\end{equation}
by \eqref{hestimate1},
\begin{equation}
\label{est_gradF2}
\|\nabla F(x)\|=\Delta t \|\nabla h_{\e}(x,V)\|\leq 2\e^{\Cr{c_a}-4}
\end{equation}
by \eqref{hestimate2},
\begin{equation}
\label{est_a}
|\phi(f(x))-\phi(x)|\leq j\Omega(x)\exp(j |F(x)|) |F(x)|
\leq  \e^{\Cr{c_a}-3}\Omega(x)
\end{equation}
by \eqref{phipro}, \eqref{adef}, \eqref{est_F2}, \eqref{jere} and restricting $\e$,
\begin{equation}
\label{est_b}
\big||\Lambda_n \nabla f(x)\circ S|-1\big|\leq c(n)\|\nabla F(x)\|
\leq \frac12 \e^{\Cr{c_a}-5}\leq \e^{-5}\Delta t
\end{equation}
by \eqref{est_gradF2} and restricting $\e$ depending only on $n$,
\begin{equation}
\label{est_c}
\begin{split}
|\phi(f(x))-\phi(x)-F(x)\cdot \nabla \phi(x)| &
\leq  j|F(x)|^2\Omega(x)\exp(j|F(x)|) 
\leq \frac12 \e^{2\Cr{c_a}-5}  \Omega(x) \\ & \leq \e^{\Cr{c_a}-5}\Omega(x)\Delta t
\end{split}
\end{equation}
by \eqref{phipro2}, \eqref{est_F2}, \eqref{jere} 
and by restricting $\e$,
\begin{equation}
\label{est_d}
\big||\Lambda_n \nabla f(x)\circ S|-1-\nabla F(x)\cdot S\big|
\leq c(n)\|\nabla F(x)\|^2\leq 4c(n)\e^{2\Cr{c_a}-8}\leq \e^{\Cr{c_a}-9}\Delta t
\end{equation}
by \eqref{est_gradF2} and restricting $\e$ depending only on $n$. 
Now recalling the definition of push-forward of varifold and \eqref{Bdef},
we have
\begin{equation}
\label{estvel1}
\begin{split}
\|f_{\sharp}V\|(\phi)&-\|V\|(\phi)-\delta (V,\phi)(h_{\e}(\cdot,V))\Delta t
= \|f_{\sharp} V\|(\phi)-\|V\|(\phi)-\delta(V,\phi)(F) \\
&=\int_{{\bf G}_n({\mathbb R}^{n+1})}( \phi(f(x))|\Lambda_n
\nabla f(x)\circ S|-\phi(x))\, dV(x,S)\\ &-\int_{{\bf G}_n({\mathbb R}^{n+1})} (\nabla F(x)\cdot S\, \phi(x)+F(x)\cdot\nabla
\phi(x))\, dV(x,S).
\end{split}
\end{equation}
We then interpolate \eqref{estvel1} and use \eqref{est_a}-\eqref{est_d} as
\begin{equation}
\label{estvel2}
\begin{split}
|\eqref{estvel1}|&\leq \int_{{\bf G}_n({\mathbb R}^{n+1})}\big|(\phi(f(x))-\phi(x))
|\Lambda_n\nabla f(x)\circ S|+(|\Lambda_n\nabla f(x)\circ S|-1)\phi(x)\\
& \hspace{1cm} -\nabla F(x)\cdot S\, \phi(x)-F(x)\cdot\nabla\phi(x)\big|\, dV(x,S)\\
&=\int_{{\bf G}_n({\mathbb R}^{n+1})} \big|(\phi(f(x))-\phi(x))(|\Lambda_n\nabla
f(x)\circ S|-1)+(\phi(f(x))-\phi(x) \\
& \hspace{1cm} -F(x)\cdot\nabla\phi(x)) +(|\Lambda_n\nabla f(x)\circ S|-1
-\nabla F(x)\cdot S)\phi(x)\big|\, dV(x,S)\\
&\leq (\e^{\Cr{c_a}-8}+\e^{\Cr{c_a}-5}+\e^{\Cr{c_a}-9})\|V\|(\Omega)\Delta t
\end{split}
\end{equation}
where we also used $\phi\leq \Omega$ for the last step. 
By restricting $\e$ so that $3\e M\leq 1$, we obtain \eqref{af0}. For \eqref{af1}, 
using \eqref{mc0} and \eqref{mc01} with $\phi=\Omega$, 
$j\in[\Cr{c_1}+1,\Cr{c_1}+2)$ and restricting $\e$ depending on $\Cr{c_1}$, 
we have
\begin{equation}
\label{estvel3}
\begin{split}
\delta(V,\Omega)(h_{\e})&=\delta V(\Omega h_{\e})+\int_{{\bf G}_n({\mathbb R}^{n+1})}
h_{\e}\cdot S^{\perp}(\nabla\Omega)\,dV(\cdot,S) \\
&\leq \delta V(\Omega h_{\e})+\frac12\int_{{\mathbb R}^{n+1}} |h_{\e}|^2\Omega
+|\nabla\Omega|^2 \Omega^{-1}\, d\|V\| \\
&\leq -\frac12(1-3\e^{\frac14})\int_{{\mathbb R}^{n+1}} \frac{|\Phi_{\e}\ast\delta V|^2
\Omega}{\Phi_{\e}\ast\|V\|+\e\Omega^{-1}}\, dx+2\e^{\frac14}+\frac{\Cr{c_1}^2}{2}
\|V\|(\Omega)
\end{split}
\end{equation}
where we also used \eqref{omegacon1}. Restrict $\Cr{e_5}$ so that 
$1-3\e^{\frac14}>\frac12$.  Then \eqref{estvel3} and \eqref{af0} give \eqref{af1}.

For \eqref{af2} and \eqref{af3}, for short, write ${\hat V}:=f_{\sharp} V$. 
Due to the assumption that
$\|f_{\sharp}V\|(\Omega)=\|\hat V\|(\Omega)\leq M$, we have \eqref{hestimate1}-\eqref{hestimate3} 
for $h_{\e}(\cdot,\hat V)$ as well. We first estimate 
$\Phi_{\e}\ast\|\hat V\|-\Phi_{\e}\ast\|V\|$ and 
$\Phi_{\e}\ast\delta \hat V-\Phi_{\e}\ast\delta V$, which lead to estimates
of $h_{\e}(\cdot,V)- h_{\e}(\cdot,\hat V)$. 
We have
\begin{equation}
\begin{split}
\big| \Phi_{\e}\ast \|\hat V\|(x)- &\Phi_{\e}\ast \|V\|(x)\big|=
\Big|\int\Phi_{\e}(z-x)\, d\|\hat V\|(z)-\int \Phi_{\e}(y-x)\, d\|V\|(y)\Big| \\
&=\Big|\int \Phi_{\e}(f(y)-x)|\Lambda_n \nabla f(y)\circ S|-\Phi_{\e}(y-x)\, dV(y,S)\Big| \\
&\leq \int |\Phi_{\e}(f(y)-x)-\Phi_{\e}(y-x)||\Lambda_n \nabla f(y)\circ S|\, dV(y,S) \\ 
&\hspace{.5cm}+\int \Phi_{\e}(y-x)\big||\Lambda_n \nabla f(y)\circ S|-1\big|\, dV(y,S).
\end{split}
\label{estvel14}
\end{equation}
By \eqref{est_F2} and \eqref{Phies1}, for some $\hat y$ lying on the line segment
connecting $y-x$ and $f(y)-x$, 
\begin{equation}
|\Phi_{\e}(f(y)-x)-\Phi_{\e}(y-x)|\leq |F(y)| |\nabla \Phi_{\e}(\hat y)| \leq c(n)\e^{\Cr{c_a}-n-5}
\chi_{B_2(x)}(y).
\label{estvel15}
\end{equation}
By \eqref{est_b}, 
\begin{equation}
\Phi_{\e}(y-x)\big||\Lambda_n\nabla f(y)\circ S|-1\big|\leq \e^{\Cr{c_a}-n-6} \chi_{B_1(x)}(y).
\label{estvel16}
\end{equation}
Combining \eqref{estvel14}-\eqref{estvel16}, we obtain
\begin{equation}
\big| \Phi_{\e}\ast \|\hat V\|(x)- \Phi_{\e}\ast \|V\|(x)\big|\leq \e^{\Cr{c_a}-n-7}\|V\|(B_2(x)).
\label{estvel17}
\end{equation}
Next, by \eqref{musm2},
\begin{equation}
\begin{split}
&\big|\Phi_{\e}\ast \delta\hat V (x)-\Phi_{\e}\ast \delta V (x)\big| \\ &=\Big|\int T
\big(\nabla\Phi_{\e}(z-x)\big)\, d\hat V(z,T)-\int S\big(\nabla\Phi_{\e}(y-x)\big)\, dV(y,S)
\Big| \\
&=\Big|\int \big\{(\nabla f(y)\circ S)(\nabla\Phi_{\e}(f(y)-x))|\Lambda_n \nabla f(y)\circ S|
-S(\nabla\Phi_{\e}(y-x))\big\}\, dV(y,S)\Big|.
\end{split}
\label{estvel18}
\end{equation}
By estimating $\nabla f(y)-I$ using \eqref{est_gradF2}
and using similar estimates as in \eqref{estvel15} and \eqref{estvel16} (where
$\Phi_{\e}$ is replaced by $\nabla\Phi_{\e}$, causing a multiplication by $\e^{-2}$), we obtain
\begin{equation}
\big|\Phi_{\e}\ast \delta\hat V (x)-\Phi_{\e}\ast \delta V (x)\big|
\leq \e^{\Cr{c_a}-n-9} \|V\|(B_2(x))
\label{estvel19}
\end{equation}
from \eqref{estvel18} by the similar interpolations. 
We also have rough estimates of
\begin{equation}
|\Phi_{\e}\ast \delta V(x)|,\,|\Phi_{\e}\ast \delta \hat V (x)|\leq \e^{-n-4}\|V\|(B_2(x)).
\label{estvel20}
\end{equation}
Using \eqref{estvel17}, \eqref{estvel19} and \eqref{estvel20}, we have
\begin{equation}
\begin{split}
&\Big|\frac{\Phi_{\e}\ast\delta\hat V}{\Phi_{\e}\ast\|\hat V\|+\e\Omega^{-1}}
-\frac{\Phi_{\e}\ast\delta V }{\Phi_{\e}\ast\| V\|+\e\Omega^{-1}}\Big| \\
&\leq \frac{|\Phi_{\e}\ast\delta\hat V -\Phi_{\e}\ast\delta V |}{\e\Omega^{-1}}
+\frac{|\Phi_{\e}\ast\delta V||\Phi_{\e}\ast\|\hat V\|-\Phi_{\e}\ast\| V\||}{
\e^2\Omega^{-2}} \\
&\leq \e^{\Cr{c_a}-n-10}\Omega(x)\|V\|(B_2(x))+\e^{\Cr{c_a}-2n-13}\Omega(x)^2 \|V\|(B_2(x))^2
\end{split}
\label{estvel21}
\end{equation}
and similarly
\begin{equation}
\begin{split}
&\Big|\frac{|\Phi_{\e}\ast\delta\hat V|^2\Omega}{\Phi_{\e}\ast\|\hat V\|+\e\Omega^{-1}}
-\frac{|\Phi_{\e}\ast\delta V|^2\Omega }{\Phi_{\e}\ast\| V\|+\e\Omega^{-1}}\Big| \\
&\leq \e^{\Cr{c_a}-2n-15}\Omega(x)^2\|V\|(B_2(x))^2+\e^{\Cr{c_a}-3n-17}\Omega(x)^3 \|V\|(B_2(x))^3.
\end{split}
\label{estvel22}
\end{equation}
By Lemma \ref{toma} with $r=2$, 
we obtain \eqref{af3} from \eqref{estvel22}.
Recalling the definition \eqref{curvature}, from \eqref{estvel21} and
with \eqref{toma1}, we obtain (writing $h_{\e}(\cdot, V)$ as $h_{\e}(V)$)
\begin{equation}
|h_{\e}(V)-h_{\e}(\hat V)|\leq \e^{\Cr{c_a}-2n-14} (M+M^2),\,\,
\|\nabla^l h_{\e}(V)-\nabla^l h_{\e}(\hat V )\|\leq \e^{\Cr{c_a}-2n-14-2l}(M+M^2)
\label{estvel23}
\end{equation}
for $l=1,2$.
Finally, we have
\begin{equation}
\begin{split}
&|\delta(V,\phi)(h_{\e}(V))-\delta(\hat V, \phi)(h_{\e}(\hat V))| =\Big|\int (\nabla h_{\e}(V)\cdot S\phi+h_{\e}(V)\cdot\nabla \phi)\, dV \\ &
-\int \big\{(\nabla h_{\e}(\hat V)\circ f)\cdot (\nabla f\circ S)(\phi\circ f) +( h_{\e}(\hat V)
\circ f)\cdot( \nabla \phi\circ f)\big\}|\Lambda_n \nabla f\circ S|\, dV\Big|.
\end{split}
\label{estvel24}
\end{equation}
Using \eqref{estvel23} as well as \eqref{est_F2}-\eqref{est_b} and \eqref{adef}, estimates by interpolations
on \eqref{estvel24} give \eqref{af2}.
\hfill{$\Box$}
\section{Existence of limit measures}
\begin{prop}
Given any ${\mathcal E}_0\in \mathcal{OP}_{\Omega}^N$ 
and $j\in {\mathbb N}$ with $j\geq \max\{1,\Cr{c_1}\}$,
there exist $\e_j\in (0,j^{-6})$, $p_j\in {\mathbb N}$, a family
$\mathcal E_{j,l}\in \mathcal{OP}_{\Omega}^N$ ($l=0,1,2,\ldots, j \,2^{p_j}$) with
the following property.
\begin{equation}
{\mathcal E}_{j,0}={\mathcal E}_0\,\,\,\mbox{for all }j\in {\mathbb N}
\label{exapp1}
\end{equation}
and with the notation of
\begin{equation}
\Delta t_j:=\frac{1}{2^{p_j}},
\label{exapp0}
\end{equation} 
we have
\begin{equation}
\|\partial \mathcal E_{j,l}\|(\Omega)\leq \|\partial \mathcal E_{0}\|(\Omega)
\exp\big(\frac{\Cr{c_1}^2 l}{2} \Delta t_j\big)+\frac{2\e_j^{\frac18}}{\Cr{c_1}^2}
\big(\exp\big(\frac{\Cr{c_1}^2 l}{2}\Delta t_j\big)-1\big),
\label{exapp2}
\end{equation}
\begin{equation}
\begin{split}
&\frac{\|\partial\mathcal E_{j,l}\|(\Omega)-\|\partial \mathcal E_{j,l-1}\|(\Omega)}{\Delta t_j}+\frac14 \int_{{\mathbb R}^{n+1}} \frac{|\Phi_{\e_j}\ast\delta(\partial \mathcal E_{j,l})
|^2\Omega}{\Phi_{\e_j}\ast \|\partial \mathcal E_{j,l}\|+\e_j\Omega^{-1}}\, dx \\
& -\frac{(1-j^{-5})}{\Delta t_j}\Delta_{j}\|\partial \mathcal E_{j,l-1}\|(\Omega)
\leq \e_j^{\frac18}+\frac{\Cr{c_1}^2}{2}\|\partial\mathcal E_{j,l-1}\|(\Omega) ,
\end{split}
\label{exapp3}
\end{equation}
\begin{equation}
\frac{\|\partial \mathcal E_{j,l}\|(\phi)-\|\partial\mathcal E_{j,l-1}\|(\phi)}{\Delta t_j}
\leq \delta(\partial \mathcal E_{j,l},\phi)(h_{\e_j}(\cdot,\partial\mathcal E_{j,l}))+\e_j^{\frac18}
\label{exapp4}
\end{equation}
for $l=1,2,\ldots,j \,2^{p_j}$ and $\phi\in \mathcal A_j$.
When $\Cr{c_1}=0$, the right-hand side of \eqref{exapp2} should
be understood as the limit $\Cr{c_1}\rightarrow 0+$.
\label{exapp}
\end{prop}
{\it Proof}. 
Given $\mathcal E_0\in \mathcal{OP}_{\Omega}^N$ and $j\in \mathbb N$
with $j\geq \max\{1,\Cr{c_1}\}$, define
\begin{equation}
M_j:=\|\partial\mathcal E_0\|(\Omega)\exp\big(\frac{\Cr{c_1}^2 j}{2}\big)+1.
\label{Mdef0}
\end{equation}
Let $\Cr{e_1},\ldots,\Cr{e_5}$ be chosen in the previous section
corresponding to $M_j$ as $M$, 
then we choose $\e_j$ so that 
$\e_j\leq \min\{\Cr{e_1},\ldots,\Cr{e_5}\}$,
\begin{equation}
\label{ejex}
\frac{2\e_j^{\frac18}}{\Cr{c_1}^2} \big(\exp\big(\frac{\Cr{c_1}^2 j}{2}\big)-1\big)<1,
\hspace{.5cm}
3\e_j^{\frac14}+\e_j^{\Cr{c_a}-3n-18}<\e_j^{\frac18}
\end{equation}
and \eqref{jere} hold. Let
$\Cr{c_a}$ be as in \eqref{defca},
and choose $p_j\in \mathbb N$ so that 
\begin{equation}
\frac{1}{2^{p_j}}\in (2^{-1} \e_j^{\Cr{c_a}},\e_j^{\Cr{c_a}}].
\label{gridef}
\end{equation}
Define $\Delta t_j$ as in \eqref{exapp0}. 
We proceed with inductive
argument. Set $\mathcal E_{j,0}=\mathcal E_0$. 
Assume that up to $k=l\in \{0,1,\ldots,j\,2^{p_j}-1\}$, $\mathcal E_{j,k}$ is determined with the
estimates \eqref{exapp2}-\eqref{exapp4}. We will define $\mathcal E_{j,l+1}$
satisfying the estimates. 
Choose $f_1\in {\bf E}({\mathcal E}_{j,l},j)$ (cf. Definition \ref{mrld}) such that
 \begin{equation}
 \label{ap1}
 \|\partial (f_1)_{\star} {\mathcal E}_{j,l} \|(\Omega)
 -\|\partial {\mathcal E}_{j,l}\|(\Omega)\leq (1-j^{-5})\Delta_{j}\| 
 \partial{ \mathcal E}_{j,l} \|(\Omega)
 \end{equation}
 and define
 \begin{equation}
 {\mathcal E}^*_{j,l+1}:=(f_1)_{\star}{\mathcal E}_{j,l}\in \mathcal{OP}_{\Omega}^N.
 \label{ap2}
 \end{equation}
 We note that
 \begin{equation}
 \|\partial {\mathcal E}^*_{j,l+1}\|(\Omega)\leq \|\partial{\mathcal E}_{j,l}\|(
 \Omega)\leq M_j
 \label{ap3}
 \end{equation}
 by \eqref{ap1}, \eqref{exapp2}, \eqref{ejex} and \eqref{Mdef0}. We next
define a smooth function $f_2:{\mathbb R}^{n+1}\rightarrow{\mathbb R}^{n+1}$ by
 \begin{equation}
 \label{ap4}
 f_2(x):=x+\Delta t_j\, h_{\e_j}(x,\partial {\mathcal E}^*_{j,l+1}).
 \end{equation}
 By the choice of $\e_j$ and $\Delta t_j$, and by \eqref{hestimate1}
 and \eqref{hestimate2}, we have
\begin{equation}
\label{ap5}
|\Delta t_j\,h_{\e_j}(x,\partial {\mathcal E}^*_{j,l+1})|\leq 2 \e_j^{\Cr{c_a}-2},
\hspace{.5cm} \|\nabla (\Delta t_j\,h_{\e_j}(x,\partial {\mathcal E}^*_{j,l+1}))\|
\leq 2\e_j^{\Cr{c_a}-4},
\end{equation}
thus $f_2$ is a diffeomorphism and ${\mathcal E}^*_{j,l+1}$-admissible
in particular. We then define
\begin{equation}
\label{ap6}
{\mathcal E}_{j,l+1}:=(f_2)_{\star}{\mathcal E}^*_{j,l+1}\in 
\mathcal{OP}_{\Omega}^N.
\end{equation} 
Note that, since $f_2$ is a diffeomorphism, 
if we write ${\mathcal E}^*_{j,l+1}=\{E_i\}_{i=1}^{N}$, then we have
${\mathcal E}_{j,l+1}=\{f_2(E_i)\}_{i=1}^{N}$. Furthermore, 
we have
\begin{equation}
(f_2)_{\sharp} \partial{\mathcal E}^*_{j,l+1}=(f_2)_{\sharp}|\cup_{i=1}^{N}\partial E_i|=
|\cup_{i=1}^{N}\partial (f_2(E_i))|=\partial {\mathcal E}_{j,l+1}. 
\label{ap7}
\end{equation}
To close the inductive 
argument, we need to check \eqref{exapp2}-\eqref{exapp4} with $l$ replaced by $l+1$.
To prove \eqref{exapp2}, we use \eqref{af1} with $M=M_j$,
$V=\partial {\mathcal E}^*_{j,l+1}$ as well as $3\e_j^{\frac14}<\e_j^{\frac18}$ 
of \eqref{ejex} to obtain
\begin{equation}
\begin{split}
\|(f_2)_{\sharp}\partial {\mathcal E}^*_{j,l+1}\|(\Omega) &
\leq
\|\partial {\mathcal E}^*_{j,l+1}
\|(\Omega)+\Delta t_j(\e_j^{\frac18} +\frac{\Cr{c_1}^2}{2}\|\partial {\mathcal E}^*_{j,l+1}
\|(\Omega)) \\
& \leq
\|\partial {\mathcal E}_{j,l}
\|(\Omega)+\Delta t_j(\e_j^{\frac18} +\frac{\Cr{c_1}^2}{2}\|\partial {\mathcal E}_{j,l}
\|(\Omega)),
\end{split}
\label{ap9}
\end{equation}
the last inequality due to \eqref{ap3}. 
By \eqref{ap9} and \eqref{exapp2}, a direct computation using 
$e^{(x+s)}\geq (1+s)e^x$ for $s\geq 0$ proves
\eqref{exapp2} with $l$ replaced by $l+1$.  In particular, this proves
that $\|\partial\mathcal E_{j,l+1}\|(\Omega)\leq M_j$, giving
the validity of \eqref{af2} and \eqref{af3} for the pair $V=
\partial\mathcal E^*_{j,l+1}$ and $f_{\sharp} V=\partial\mathcal E_{j,l+1}$. From \eqref{af1}, \eqref{ap3}, \eqref{af3}, \eqref{ap1} and \eqref{ejex}, 
we obtain \eqref{exapp3} for $l+1$ in place
of $l$. From \eqref{af0}, \eqref{af2}, \eqref{ejex} and $f_1\in {\bf E}(\mathcal E_{j,l}
,j)$, we obtain \eqref{exapp4}
for $l+1$ in place of $l$. This closes the inductive step, showing 
\eqref{exapp2}-\eqref{exapp4} up to $l=j\,2^{p_j}$. 
\hfill{$\Box$}
\begin{rem} Due to the choice of $\e_j$, each $\partial\mathcal E_{j,l}$
satisfies various estimates obtained in Section \ref{SMC} with
$V=\partial \mathcal E_{j,l}$, $\e=\e_j$. 
\end{rem}
\begin{rem}
It is convenient to define approximate solutions for all $t\geq 0$
instead of discrete times. 
For each $j\in {\mathbb N}$ with $j\geq \max\{1,\Cr{c_1}\}$, 
define a family $\mathcal E_j (t)\in
\mathcal{OP}_{\Omega}^N$ for $t\in [0,j]$ by
\begin{equation}
\mathcal E_j(t):=\mathcal E_{j,l} \mbox{ if }t\in ((l-1)\Delta t_j,
l\Delta t_j].
\label{contime1}
\end{equation}
\end{rem}
\begin{prop}
Under the assumptions of Proposition \ref{exapp},
there exist a subsequence $\{j_l\}_{l=1}^{\infty}$ and a family 
of Radon measures $\{\mu_t\}_{t\in {\mathbb R}^+}$ on ${\mathbb
R}^{n+1}$ such that
\begin{equation}
\lim_{l\rightarrow\infty} \|\partial\mathcal E_{j_l}(t)\|(\phi)
=\mu_t(\phi)
\label{meconeq}
\end{equation}
for all $\phi\in C_c({\mathbb R}^{n+1})$ and for all $t\in {\mathbb R}^+$. For 
all $T<\infty$, we have 
\begin{equation}
\label{mecones}
\limsup_{l\rightarrow\infty} \int_0^T\Big(\int_{\mathbb R^{n+1}} 
\frac{|\Phi_{\e_{j_l}} \ast \delta (\partial \mathcal E_{j_l}(t))|^2
\Omega}{\Phi_{\e_{j_l}}\ast\|\partial \mathcal E_{j_l}(t)\|+\e_{j_l}
\Omega^{-1}}\,dx-\frac{1}{\Delta t_{j_l}} \Delta_{j_l} \|\partial
\mathcal E_{j_l}(t)\|(\Omega)\Big)\,dt<\infty,
\end{equation}
and for a.e.~$t\in \mathbb R^+$, we have
\begin{equation}
\label{mecomake}
\lim_{l\rightarrow\infty}j_l^{2(n+1)}\Delta_{j_l}\|\partial\mathcal E_{j_l}(t)\|(\Omega)=0.
\end{equation}
\label{promecone}
\end{prop}
{\it Proof}. Let $2_{\mathbb Q}$  be the set of all non-negative
numbers of the form $\frac{i}{2^j}$ for some $i,j\in{\mathbb N}\cup
\{0\}$. $2_{\mathbb Q}$ is dense in ${\mathbb R}^+$ and countable. 
For each fixed $J\in {\mathbb N}$, $\limsup_{j\rightarrow\infty}
(\sup_{t\in [0,J]}\|\partial {\mathcal E}_j(t)\|(\Omega))\leq \|\partial \mathcal E_0\|(\Omega)\exp(\Cr{c_1}^2J/2)$ by \eqref{exapp2}. 
Thus, by diagonal argument, we may choose a subsequence and a family of Radon measures 
$\{\mu_t\}_{t\in 2_{\mathbb Q}}$ on ${\mathbb R}^{n+1}$ such that
\begin{equation}
\label{mecon}
\lim_{l\rightarrow\infty} \|\partial {\mathcal E}_{j_l}(t)\|(\phi)=\mu_{t}(\phi)
\end{equation}
for all $\phi\in C_c({\mathbb R}^{n+1})$ and $t\in 2_{\mathbb Q}$. 
We also have
\begin{equation}
\mu_t(\Omega)\leq \|\partial \mathcal E_0\|(\Omega) \exp(\Cr{c_1}^2 t/2)
\label{mecon2}
\end{equation}
for all $t\in 2_{\mathbb Q}$. 
Next, let $Z:=\{\phi_{q}\}_{q\in {\mathbb N}}$ be a countable subset of $C^2_c({\mathbb R}^{n+1};{\mathbb R}^+)$ which is dense in $C_c({\mathbb R}^{n+1};{\mathbb R}^+)$ with respect to the supremum norm. We claim that, for any 
given $J\in {\mathbb N}$, 
\begin{equation}
g_{q,J}(t):=\mu_t(\phi_q)-2t\|\nabla^2 \phi_q\|_{\infty} (\min_{x\in {\rm spt}\,
\phi_q}\Omega(x))^{-1} \|\partial\mathcal E_0\|(\Omega)\exp(\Cr{c_1}^2 J/2)
\label{mecon4}
\end{equation}
is a monotone decreasing function of $t\in [0,J]\cap 2_{\mathbb Q}$. Since 
$\phi_q$ has a compact support and due to the linear dependence of 
\eqref{mecon4} on $\phi_q$, we may assume $\phi_q<\Omega$ without loss of generality. 
To prove \eqref{mecon4}, just like \eqref{estvel3}, using \eqref{mc0} and \eqref{mc01}, 
we have
\begin{equation}
\label{mecon3}
\delta(\partial\mathcal E_j(t),\phi)(h_{\e_j}(\cdot,\partial\mathcal E_j(t)))\leq 
2\e_j^{\frac14} +\frac12 \int_{{\mathbb R}^{n+1}} 
\frac{|\nabla\phi|^2}{\phi}\, d\|\partial\mathcal E_j(t)\|
\end{equation}
for $\phi\in {\mathcal A}_j$ and $t\in [0,j]$. For any $\phi_q\in Z$ with $\phi_q<\Omega$
and sufficiently large $i\in {\mathbb N}$, choose $j_0\in {\mathbb N}$ so that 
$\phi_q+i^{-1}\Omega\in {\mathcal A}_{j_0}$ holds and $j_0\geq J$.
 For any $t_1,t_2\in [0,J]\cap 2_{\mathbb Q}$ with $t_2>t_1$ fixed,
 choose a larger $j_0$ so that $t_1$ and $t_2$ are 
 integer-multiples of $1/2^{j_0}$. Then, by \eqref{exapp4}
 and \eqref{mecon3}, we have
 \begin{equation}
 \label{mecon5}
 \begin{split}
 &\|\partial\mathcal E_{j_l}(t_2)\|(\phi_q+i^{-1}\Omega)
 -\|\partial\mathcal E_{j_l}(t_1)\|(\phi_q+i^{-1}\Omega) \\ 
& \leq (\e_{j_l}^{\frac18} +2\e_{j_l}^{\frac14})(t_2-t_1)+\frac12
 \int_{t_1}^{t_2}\int_{\mathbb R^{n+1}}\frac{|\nabla(\phi_q+i^{-1}\Omega)|^2}{\phi_q+i^{-1}\Omega}\, d\|\partial\mathcal E_{j_l}(t)\|dt
 \end{split}
 \end{equation}
 for all $j_l\geq j_0$. As $l\rightarrow \infty$, the left-hand side of \eqref{mecon5} may be
 bounded from below using \eqref{mecon} and \eqref{exapp2} as
 \begin{equation}
 \geq \mu_{t_2}(\phi_q)-\mu_{t_1}(\phi_q) - i^{-1}\|\partial\mathcal E_0\|
 (\Omega)\exp(\Cr{c_1}^2 J/2).
 \label{mecon6}
 \end{equation}
 To estimate the right-hand side of \eqref{mecon5}, note that
 \begin{equation}
 \label{mecon7}
 \frac{|\nabla(\phi_q+i^{-1}\Omega)|^2}{\phi_q+i^{-1}\Omega}
 \leq 2\frac{|\nabla\phi_q|^2}{\phi_q}
 +2i^{-1}\frac{|\nabla\Omega|^2}{\Omega}\leq 4\|\nabla^2\phi_q\|_{\infty}
 (\min_{x\in {\rm spt}\,\phi_q} \Omega(x))^{-1}\Omega+
 2i^{-1}\Cr{c_1}^2 \Omega.
 \end{equation}
Now, using \eqref{mecon5}-\eqref{mecon7}, and then letting
$i\rightarrow\infty$, we obtain 
\begin{equation}
\mu_{t_2}(\phi_q)-\mu_{t_1}(\phi_q)\leq 
2\|\nabla^2\phi_q\|_{\infty}(\min_{x\in {\rm spt}\,\phi_q}
\Omega(x))^{-1} \|\partial\mathcal E_0\|(\Omega) \exp(\Cr{c_1}^2
J/2)(t_2-t_1).
\label{mecon8}
\end{equation}
 Then \eqref{mecon8} proves that $g_{q,J}(t)$ defined in \eqref{mecon4} is monotone decreasing. Define
 \begin{equation}
 D:=\cup_{J\in {\mathbb N}}
 \{t\in (0,J) : \mbox{for some $q\in \mathbb N$, 
 $\lim_{s\rightarrow t-}g_{q,J}(s)>\lim_{s\rightarrow t+}g_{q,J}(s)$}\}.
 \label{mecon9}
 \end{equation}
 By the monotone property of $g_{q,J}$, $D$ is a countable set
 on $\mathbb R^+$, and $\mu_t(\phi_q)$ may be defined 
 continuously on the
 complement of $D$ uniquely from the values on $2_{\mathbb Q}$. 
 For any $t\in {\mathbb R}^+\setminus 
 (D\cup 2_{\mathbb Q})$ and $\phi_q\in Z$, we claim that
 \begin{equation}
 \lim_{l\rightarrow\infty}\|\partial \mathcal E_{j_l}(t)\|(\phi_q)=
 \mu_t(\phi_q).
 \label{mecon10}
 \end{equation}
Due to the definition of $\partial\mathcal E_{j_l}(t)$, there exists
a sequence $\{t_l\in 2_{\mathbb Q}\}_{l=1}^{\infty}$ such that
$\partial \mathcal E_{j_l}(t_l)=\partial\mathcal E_{j_l}(t)$ and that
$\lim_{l\rightarrow\infty} t_l=t+$. For any $s>t$ with 
$s\in 2_{\mathbb Q}$ and for all sufficiently 
large $l$, \eqref{mecon5} shows 
\begin{equation}
\|\partial\mathcal E_{j_l}(s)\|(\phi_q+i^{-1}\Omega)\leq \|\partial \mathcal E_{j_l}
(t_l)\|(\phi_q+i^{-1}\Omega)+ O(s-t).
\label{mecon11}
\end{equation}
Taking $\liminf_{l\rightarrow\infty}$ and taking $i\rightarrow\infty$
on both sides of \eqref{mecon11}, we have
\begin{equation}
\mu_s(\phi_q)\leq \liminf_{l\rightarrow\infty}\|\partial\mathcal E_{j_l}(t_l)\|
(\phi_q)+O(s-t).
\label{mecon12}
\end{equation}
By letting $s\rightarrow t+$, $\partial\mathcal E_{j_l}(t_l)=\partial
\mathcal E_{j_l}(t)$, \eqref{mecon12} and the continuity of
$\mu_s(\phi_q)$ at $s=t$ imply
\begin{equation}
\mu_t(\phi_q)\leq \liminf_{l\rightarrow\infty}\|\partial\mathcal E_{j_l}(t)\|
(\phi_q).
\label{mecon13}
\end{equation}
For any $s<t$ with $s\in 2_{\mathbb Q}$, we also have
\begin{equation}
\|\partial\mathcal E_{j_l}(t_l)\|(\phi_q+i^{-1}\Omega)\leq 
\|\partial\mathcal E_{j_l}(s)\|(\phi_q+i^{-1}\Omega)+O(t_l-s).
\label{mecon14}
\end{equation}
Take $\limsup_{l\rightarrow\infty}$, then let $i\rightarrow\infty$
to obtain from \eqref{mecon14}
\begin{equation}
\limsup_{l\rightarrow\infty}\|\partial\mathcal E_{j_l}(t)\|(\phi_q)
\leq \mu_s(\phi_q)+O(t-s).
\label{mecon15}
\end{equation}
By letting $s\rightarrow t-$ and by the continuity of $\mu_s(\phi_q)$, we have
\begin{equation}
\limsup_{l\rightarrow\infty}\|\partial\mathcal E_{j_l}(t)\|(\phi_q)
\leq \mu_t(\phi_q).
\label{mecon16}
\end{equation}
\eqref{mecon13} and \eqref{mecon16} prove \eqref{mecon10}
for all $\phi_q\in Z$. Since $Z$ is dense in $C_c({\mathbb R}^{n+1};
{\mathbb R}^+)$,
\eqref{mecon10} determines the limit measure uniquely and
the convergence also holds in general 
for $\phi\in C_c({\mathbb R}^{n+1})$. For $t\in D$, since $D$ is 
countable, we may choose a further subsequence by a diagonal
argument so that a further subsequence of 
$\{\|\partial\mathcal E_{j_l}(t)\|\}_{l=1}^{\infty}$ converges
for all $t\in {\mathbb R}^+$ to a Radon measure $\mu_t$.
Finally \eqref{mecones} follows from \eqref{exapp3}. 
Since $\Delta t_{j_l} \leq \e_{j_l}^{\Cr{c_a}}
\ll j_l^{-2(n+1)}$ by \eqref{gridef}, \eqref{defca} and \eqref{jere},
we have $\lim_{l\rightarrow\infty} 
\int_0^T -j_l^{2(n+1)}\Delta_{j_l}\|\partial\mathcal E_{j_l}(t)\|(\Omega)\,dt
\leq \lim_{l\rightarrow\infty}\Delta t_{j_l} j_l^{2(n+1)}=0$. Thus there exists a further subsequence such that
the integrand converges pointwise to 0 for a.e.~on $[0,T]$. As $T\rightarrow\infty$ and
carrying out a diagonal argument, we may conclude \eqref{mecomake} holds
for a.e.~$t\in \mathbb R^+$ for a subsequence. 
\hfill{$\Box$}
\begin{rem}
In \eqref{ap1}, we choose $f_1\in {\bf E}(\mathcal E_{j,l},j)$ so that $f_1$ 
nearly achieves $\inf$ among ${\bf E}(\mathcal E_{j,l},j)$. The choice of
factor $1-j^{-5}$ can be different, on the other hand. In fact, all we need is 
\eqref{mecomake} (which is needed to obtain integrality later) and we may replace
$1-j^{-5}$ by any fixed number in $(0,1)$, or even a sequence of numbers $\alpha_j$ 
as long as $\lim_{j\rightarrow\infty} j^{2(n+1)}\alpha^{-1}_j \Delta t_j=0$ is
satisfied. Such choice would give a different estimate in 
\eqref{exapp3} with different factor instead of $1-j^{-5}$ but otherwise, the 
proof is identical. Since $\Delta t_j$ goes to 0 very fast ($\Delta t_j\leq \e_j^{\Cr{c_a}}=\e_j^{3n+20}$ and $\e_j<j^{-6}$), we may make 
a choice so that $\alpha_j$ goes to 0 very fast. 
This means that, if we wish, we may choose $f_1\in {\bf E}(\mathcal E_{j,l},j)$
which only achieves a ``tiny fraction'' of $\inf$ in 
$\Delta_j\|\partial\mathcal E_{j,l}\|(\Omega)$, and asymptotically doing 
almost no apparent area reducing as $j\rightarrow\infty$. The choice should be reflected upon 
the singularities of the limiting $V_t$ but we do not know how to characterize this aspect.
\label{choicef}
\end{rem}
\section{Rectifiability theorem}
\label{secrec}
The main result of this section is Theorem \ref{rect1}, which is analogous to 
Allard's rectifiability theorem \cite[5.5(1)]{Allard} but with an added difficulty of 
having only a control of smoothed mean curvature vector up to the length scale
of $O(1/j^{2})$ and a certain area minimizing property in a smaller length scale.
Except for using the notions introduced in Section 4 such as $\mathcal E$-admissible functions and $\Delta_j\|\partial\mathcal E\|(\Omega)$, the content of Section 7 and 8 are more
or less independent of Section 5 and 6, and they can be of
independent interests. 

We first recall a formula usually referred to as the monotonicity formula from 
\cite[5.1(3)]{Allard}:
\begin{lemma}
\label{red_den}
Suppose $V\in {\bf V}_n(\mathbb R^{n+1})$, $0<r_1<r_2<\infty$, $x\in 
\mathbb R^{n+1}$, and for $0\leq s<\infty$,
\begin{equation}
\label{fmono1}
\|\delta V\|(B_r(x))\leq s\|V\|(B_r(x))
\end{equation}
whenever $r_1< r< r_2$. Then
\begin{equation}
\label{fmono2}
(\exp (sr))r^{-n} \|V\|(B_r(x))
\end{equation}
is nondecreasing in $r$ for $r_1< r <r_2$. 
\end{lemma}
The following Proposition~ \ref{rect0} is essential to prove the rectifiability of
the limit measure. For the similar purpose in \cite{Brakke}, 
Brakke cites a result in \cite{Almgren} of Almgren. 
The proof by Almgren requires extensive tools involving varifold slicing and 
piecewise smooth Lipschitz
deformation to cubical complexes. On the other hand, his proof does not provide
a deformation with $\mathcal E$-admissibility or volume estimate (Proposition \ref{rect0} (4))
which are essential in our proof. For codimension 1 case, we provide 
a more direct proof using radial projection as follows. 
\begin{prop}
There exist $\Cl[c]{c_2}, \Cl[c]{c_3}\in (0,\infty)$ depending only on $n$ 
with
the following property. 
For $\mathcal E=\{E_i\}_{i=1}^N\in \mathcal{OP}_{\Omega}^N$, suppose 
$0\in {\rm spt}\,\|\partial \mathcal E\|$ and $\|\partial\mathcal E\|(B_R)
\leq \Cr{c_2}R^{n}$. Then there exist a $\mathcal E$-admissible function $f$
and $r\in [\frac{R}{2},R]$ such that 
\begin{itemize}
\item[(1)] $f(x)=x$ for $x\in \mathbb R^{n+1}\setminus U_r$,
\item[(2)] $f(x)\in B_r$ for $x\in B_r$,
\item[(3)] $\|\partial f_{\star}\mathcal E\|(B_r)\leq \frac12
\|\partial\mathcal E\|(B_r)$,
\item[(4)] $\mathcal L^{n+1}(E_i\triangle \tilde E_i)\leq \Cr{c_3}
(\|\partial \mathcal E\|(B_r))^{\frac{n+1}{n}}$ for all $i$, where $\{\tilde E_i\}_{
i=1}^{N}=f_{\star}\mathcal E$.
\end{itemize}
\label{rect0}
\end{prop}
{\it Proof}. For $r>0$ let $\nu(r):=\|\partial \mathcal E\|(B_r)=\mathcal H^{n}
(B_r\cap \cup_{i=1}^N \partial E_i)$. Since $0\in {\rm spt}\,
\|\partial \mathcal E\|$, we have $\nu(r)>0$ for $r>0$ and 
$\nu(r)$ is a monotone increasing function which is
differentiable a.e.. We also have  
\begin{equation}
\mathcal H^{n-1}({\partial B_r\cap \cup_{i=1}^N \partial E_i})\leq \nu'(r)<\infty
\label{rlm1}
\end{equation}
whenever $\nu$ is differentiable. By the 
relative isoperimetric inequality \cite[p.152]{AFP}, there exists 
$\Cr{c_3}$ depending only on $n$ such that
\begin{equation}
\label{rlm1.1}
\min\{\mathcal L^{n+1}(U_R \cap E_i),\mathcal L^{n+1}(U_R\setminus E_i)\}
\leq \Cr{c_3}( \mathcal H^{n}(U_R \cap \partial E_i))^{\frac{n+1}{n}}.
\end{equation}
We assume 
\begin{equation}
\nu(R)\leq \Big(\frac{
\mathcal L^{n+1}(U_R) }{2^{n+2}\Cr{c_3}} \Big)^{\frac{n}{n+1}},
\label{rlm1.2}
\end{equation}
and we further restrict $\nu(R)$ in the following. 
Since $\mathcal H^{n}(U_R\cap \partial E_i)\leq \nu(R)$, 
\eqref{rlm1.1} and \eqref{rlm1.2} imply that there is a unique $i_0
\in \{1,\ldots,N\}$ such that
\begin{equation}
\mathcal L^{n+1}(U_R\setminus E_{i_0})\leq \Cr{c_3}(\nu(R))^{\frac{n+1}{n}}
\leq \frac{1}{2^{n+2}}\mathcal L^{n+1}(U_R),
\label{rlm1.3}
\end{equation}
i.e., $E_{i_0}$ takes up a major part of $U_R$. The reason for the existence of such $i_0$ is as follows. Otherwise, all $E_i$ would  have a small
measure in $U_R$. Since $U_R\cap \cup_{i=1}^N E_i$ is a full measure set,
there exists a combination $E_{i_1},\ldots, E_{i_J}$ such that $(\mathcal L^{n+1}(U_R))^{-1} \mathcal L^{n+1}
(\cup_{k=1}^J E_{i_k})\in (1/4,3/4)$. The relative isoperimetric inequality applied to 
$\hat E:=\cup_{k=1}^J E_{i_k}$ gives a lower bound $\Cr{c_3}(\|\nabla\chi_{\hat E}\|(U_R))^{\frac{n+1}{n}}\geq \mathcal L^{n+1}(U_R)/4$
while we have $\|\nabla \chi_{\hat E}\|(U_R)\leq \mathcal H^{n}(U_R\cap\cup_{i=1}^N 
\partial E_i)$. This gives a contradiction to \eqref{rlm1.2}.

For all $r\in [\frac{R}{2},R]$, \eqref{rlm1.3} also gives
$\mathcal L^{n+1}(U_r\setminus E_{i_0})\leq \frac12 \mathcal L^{n+1}(U_r)$,
thus \eqref{rlm1.1} with $R$ replaced by $r$ shows 
\begin{equation}
\mathcal L^{n+1}( U_r\setminus  E_{i_0})\leq \Cr{c_3} (\mathcal H^{n}
(U_r\cap \partial E_{i_0}))^{\frac{n+1}{n}}
\label{rlm1.4}
\end{equation}
for all $r\in [\frac{R}{2},R]$.
Next, let $\tilde A:=\{ r\in [\frac{R}{2},R] : \mathcal H^{n}(\partial B_r
\setminus E_{i_0})> \frac12 \mathcal H^{n}(\partial B_r)\}$
and $A:=[\frac{R}{2},R]\setminus \tilde A$. 
Since
\begin{equation}
\label{rlm1.6}
\mathcal L^{n+1} ((U_R\setminus B_{\frac{R}{2}})\setminus E_{i_0})
=\int_{\frac{R}{2}}^R \mathcal H^{n}(\partial B_r\setminus E_{i_0})\,
dr\geq \frac12 \mathcal L^1(\tilde A) 
\mathcal H^{n}(\partial B_{\frac{R}{2}}),
\end{equation}
\eqref{rlm1.3} and \eqref{rlm1.6} show
\begin{equation}
\mathcal L^1(\tilde A)\leq \frac{R}{2(n+1)}\,\, \mbox{ and }\,\,
\mathcal L^1 (A)\geq (\frac12-\frac{1}{2(n+1)})R
\geq \frac{R}{4}.
\label{rlm1.7}
\end{equation}
In particular, \eqref{rlm1.7} proves that 
\begin{equation}
\mathcal H^{n}(\partial B_r
\setminus E_{i_0})\leq  \frac12 \mathcal H^{n}(\partial B_r)
\mbox{ for
$r\in A\subset[\frac{R}{2},R]$ with $\mathcal L^1(A)\geq \frac{R}{4}$. }
\label{rlm1.8}
\end{equation}

Next, fix arbitrary $r\in A$ which also satisfies \eqref{rlm1}, and let $G_i:=E_i\cap
\partial B_r$. Each $G_i$ is open with respect to the 
topology on $\partial B_r$ and $\partial G_i\subset
\partial B_r\cap \partial E_i$. Note also that $\partial B_r
\setminus E_i=\partial B_r \setminus G_i$. By the relative
isoperimetric inequality on $\partial B_r$ and \eqref{rlm1.8}, there exists $\Cl[c]{c_4}$
depending only on $n$ such that
\begin{equation}
\label{rlm4}
\mathcal H^{n}(\partial B_r\setminus G_{i_0})=\min\{ \mathcal H^{n}(G_{i_0}),\mathcal H^{n}(\partial B_r\setminus G_{i_0})\}
\leq \Cr{c_4}(\mathcal H^{n-1}(\partial G_{i_0}))^{\frac{n}{n-1}}.
\end{equation}

Now we choose $B_{2r_0}(x_0)\subset U_r\cap E_{i_0}$ 
and choose a Lipschitz map $f$ as follows. $f(x)=x$ if 
$x\in \mathbb R^{n+1}\setminus U_r$, $f$ maps $B_{r_0}(x_0)$
to $B_r$ bijectively, and $B_r\setminus U_{r_0}(x_0)$ onto
$\partial B_r$ by radial projection centered at $x_0$.  See Figure \ref{fig_4} for a general idea of the map.
We claim that such $f$ is $\mathcal E$-admissible.
Let $\tilde E_i:={\rm int}(f(E_i))$. For $i\neq i_0$, $\tilde E_{i}=E_i
\setminus B_r$, because $f$ is identity on $\mathbb R^{n+1}
\setminus B_r$ and $f(E_i
\cap B_r)\subset \partial B_r$. On the other hand, $\tilde E_{i_0}=E_{i_0}\cup U_r$ since $U_r=f(U_{r_0}(x_0))$ and  $U_{r_0}(x_0)\subset
E_{i_0}$, and any $x\in \partial B_r\cap E_{i_0}$ is in $E_{i_0}\cup U_r$. 
For two open sets $A$ and $B$, we have $\partial (A\cap B)\subset(\partial A\cap
{\rm clos}\,B)\cup (\partial B\cap A)$ and $\partial (A\cup B)\subset
(\partial A\setminus {\rm clos}\,B)\cup (\partial B\setminus A)$. So 
\begin{equation}
\partial \tilde E_{i}=\partial (E_i\cap(\mathbb R^{n+1}\setminus
B_r))\subset(\partial E_i\cap{\rm clos}\,(\mathbb R^{n+1}\setminus B_r) )
\cup (\partial B_r\cap E_i)=(\partial E_i\setminus U_r)
\cup G_i
\label{rlm6}
\end{equation}
for $i\neq i_0$ while 
\begin{equation}
\partial \tilde E_{i_0}=\partial(E_{i_0}\cup U_r)\subset(\partial E_{i_0}
\setminus B_r)\cup (\partial B_r\setminus E_{i_0})=(\partial E_{i_0}
\setminus B_r)\cup (\partial B_r\setminus G_{i_0}).
\label{rlm7}
\end{equation}
We need to check $\mathbb R^{n+1}\setminus \cup_{i=1}^N \tilde E_i\subset f(\cup_{i=1}^N
\partial E_i)$. Since $\mathbb R^{n+1}\setminus \cup_{i=1}^N\tilde E_i$ does
not have any interior point, it is enough to prove $\cup_{i=1}^N \partial \tilde E_i
\subset f(\cup_{i=1}^N \partial E_i)$. For $i\neq i_0$, $\partial E_i\setminus
U_r\subset f(\partial E_i)$ since $f$ is identity on $\mathbb R^{n+1}\setminus U_r$.
For any $x\in G_i$, consider a line segment $I$ with two ends, $x_0$ and $x$.  
Since $x\in G_i=\partial B_r\cap E_i$, there is some neighborhood of $x$ of $I$
belonging to $E_i$. On the other hand, we have $B_{r_0}(x_0)\subset E_{i_0}$, 
thus there must be 
some point $\hat x\in I\cap \partial E_{i_0}$. Since $f$ on $B_r\setminus 
B_{r_0}(x_0)$ is a radial projection to $\partial B_r$, $f(\hat x)=x$. This proves
that $G_i\subset f(\partial E_{i_0})$. Then \eqref{rlm6} shows $\partial
\tilde E_i\subset f(\partial E_i\cup \partial E_{i_0})$ for $i\neq i_0$. 
For $i=i_0$, $\partial
E_{i_0}\setminus B_r=f(\partial E_{i_0}\setminus B_r)$
since $f$ is identity there. For any $x\in \partial B_r \setminus G_{i_0}=\partial
B_r\setminus E_{i_0}$, either
$x\in \partial E_i$ for some $i$ (including $i=i_0$), or $x\in E_i$ for some 
$i\neq i_0$. In the former case, since $f$ is identity on $\partial B_r$, 
$x\in f(\partial E_i)$. In the latter case, the line segment connecting $x_0$
and $x$ contains $\hat x\in \partial E_{i_0}$ just as before, hence $x\in f(\partial
E_{i_0})$. Thus by \eqref{rlm7}, we have $\partial\tilde E_{i_0}\subset 
f(\cup_{i=1}^N\partial E_{i})$. In all, we have proved that $\cup_{i=1}^N
\partial\tilde E_i\subset f(\cup_{i=1}^N \partial E_i)$, and this proves that
$f$ is $\mathcal E$-admissible. With 
$\tilde{\mathcal E} =f_{\star}\mathcal E=\{\tilde E_i\}_{i=1}^N$, 
we have from \eqref{rlm6}, \eqref{rlm7} and $\cup_{i\neq i_0}G_i\subset \partial
B_r\setminus G_{i_0}$ that
\begin{equation}
\begin{split}
\|\partial \tilde{\mathcal E}\|(B_r)&=\mathcal H^{n}(
\cup_{i=1}^N\partial \tilde E_i\cap B_r)\leq 
\mathcal H^{n}(\partial B_r\setminus G_{i_0})+\sum_{i\neq i_0}
\mathcal H^{n}(\partial E_i\cap \partial B_r) \\
&=\mathcal H^{n}(\partial B_r\setminus G_{i_0}),
\end{split}
\label{rlm8}
\end{equation}
the last equality due to \eqref{rlm1}. We next note that $E_i\triangle
\tilde E_i=E_i\cap B_r$ for $i\neq i_0$ and $=U_r\setminus E_{i_0}$ 
for $i=i_0$. Since both are included in $B_r\setminus E_{i_0}$, 
\eqref{rlm1.4} shows that the condition (4) is satisfied with  this 
$\Cr{c_3}$. Thus we conclude that $\mathcal E$-admissible 
function $f$ satisfies conditions (1), (2), (4) so far. 

If the conclusion were not true, then, we must have $\|\partial\tilde{\mathcal E}
\|(B_r)> \frac12 \|\partial \mathcal E\|(B_r)=\frac12\nu(r)$ if $r\in A$
with \eqref{rlm1}. 
Combining \eqref{rlm8}, \eqref{rlm4} and \eqref{rlm1}, we obtain
\begin{equation}
\frac12 \nu(r)\leq \Cr{c_4} (\nu'(r))^{\frac{n}{n-1}}.
\label{rlm9}
\end{equation} 
Since we have $\mathcal L^1(A)\geq
\frac{R}{4}$ by \eqref{rlm1.8}, 
\begin{equation}
\nu^{\frac{1}{n}}(R)\geq \int_A (\nu^{\frac{1}{n}}(r))'\,dr\geq n^{-1}(2\Cr{c_4})^{\frac{1-n}{n}} \frac{R}{4}.
\label{rlm12}
\end{equation}
We would obtain 
a contradiction to $\|\partial\mathcal E\|(B_R)=\nu(R)\leq \Cr{c_2}R^{n}$ by
choosing an appropriately small $\Cr{c_2}$ depending only on $n$.
\hfill{$\Box$}
\begin{thm} (cf. \cite[p.78]{Brakke})
\label{rect1}
Suppose that $\{\mathcal E_j\}_{j=1}^{\infty}\subset\mathcal{OP}_\Omega^N$ and $\{\e_j\}_{j=1}^{\infty}\subset (0,1)$
satisfy
\begin{enumerate}
\item
$\lim_{j\rightarrow\infty}j^4 \e_j =0$,
\item
$\sup_j \|\partial\mathcal E_j\|(\Omega)<\infty$,
\item
$\liminf_{j\rightarrow\infty}
\int_{\mathbb R^{n+1}} \frac{|\Phi_{\e_j}\ast\delta(\partial\mathcal E_j)|^2
\Omega}{\Phi_{\e_j}\ast\|\partial\mathcal E_j\|+\e_j\Omega^{-1}}\, dx
<\infty$,
\item
$\lim_{j\rightarrow\infty} \Delta_j \|\partial\mathcal E_j\|(\Omega)=0$.
\end{enumerate}
Then there exists a converging subsequence $\{\partial \mathcal E_{j_l}\}_{l=1}^{\infty}$
whose limit $V\in{\bf V}_{n}(\mathbb R^{n+1})$ satisfies
\begin{equation}
\label{rect3.5}
\theta^{*n}(\|V\|,x)\geq \frac{\Cr{c_2}}{16\omega_{n}} \mbox{ for 
$\|V\|$ a.e.~$x$}.
\end{equation}
Furthermore, $V\in {\bf RV}_{n}(\mathbb R^{n+1})$. 
\end{thm}
{\it Proof}. The existence of converging subsequence $\{\partial \mathcal E_{j_l}\}_{l=1}^{\infty}$ and the limit 
$V$ with 
\begin{equation}
\|V\|(\Omega)\leq \sup_l \|\partial \mathcal E_{j_l}\|(\Omega)\leq M
\label{rect6}
\end{equation}
for some $M\in (0,\infty)$ follows from the compactness of Radon measures. We may also assume
that the quantities in (3) are uniformly bounded also by $M$ for this subsequence. 
Fix $R\in(0,1)$ and $x_0\in \mathbb R^{n+1}$ and define
\begin{equation}
F_R:=\{x\in B_1(x_0) : R^{-n}\|V\|(B_R (x))<\Cr{c_2}/16\},
\label{rect6.5}
\end{equation} 
where $\Cr{c_2}$ is the constant given by Proposition \ref{rect0}. 
We will prove that $\lim_{R\rightarrow 0}\|V\|(F_R)=0$ which proves 
\eqref{rect3.5} in $B_1(x_0)$. Since $x_0$ is arbitrary, we have 
\eqref{rect3.5} on $\mathbb R^{n+1}$. 

For $x\in F_R$, we may choose $\phi\in C_c^{\infty}(\mathbb R^{n+1})$ approximating
$\chi_{B_R(x)}$ such that $\phi=1$ on
$B_R(x)$, $\phi=0$ outside $B_{2R}(x)$ and $0\leq \phi\leq 1$ with 
$R^{-n}\|V\|(\phi)<\Cr{c_2}/16$. Since $\lim_{l\rightarrow\infty}
\|\partial\mathcal E_{j_l}\|=\|V\|$, for all sufficiently large $l$ depending on $x$, 
we have 
\begin{equation}
R^{-n}\|\partial\mathcal E_{j_l}\|(\phi)<\Cr{c_2}/16.
\label{rect8}
\end{equation}
Since $\Phi_{\e_{j_l}}\ast\phi$ converges uniformly to $\phi$ on 
$B_{2R+1}(x)$ by (1) and is equal to $0$ outside,
\begin{equation}
\big|\|\Phi_{\e_{j_l}}\ast\partial\mathcal E_{j_l}\|(\phi)-\|\partial\mathcal E_{j_l}\|(\phi)\big|=\big|\|\partial\mathcal E_{j_l}\|
(\Phi_{\e_{j_l}}\ast\phi-\phi)\big|\leq \sup_{B_{2R+1}(x)}(|\Phi_{\e_{j_l}}\ast\phi-\phi|\Omega^{-1}) M
\label{rect7}
\end{equation}
converges to 0. Thus, by \eqref{rect8} and \eqref{rect7}, for $x\in F_R$ there exists $m_x
\in \mathbb N$ such that 
$
R^{-n}\|\Phi_{\e_{j_l}}\ast\partial\mathcal E_{j_l}\|(B_R(x))
<\Cr{c_2}/16$
for all $l\geq m_x$. 
Thus, if we define 
\begin{equation}
F_{R,m}:=\{x\in F_R : R^{-n}\|\Phi_{\e_{j_l}}\ast\partial\mathcal E_{j_l}\|
(B_R (x))<\Cr{c_2}/16\mbox{ for all }l\geq m\},
\label{rect9}
\end{equation}
$F_{R,m}
\subset F_{R,m+1}$ for all $m\in \mathbb N$ with $\cup_{m\in\mathbb N}
F_{R,m}=F_R$.
Hence we may choose $m_1\in \mathbb N$ with
\begin{equation}
\label{rect10}
\|V\|\lfloor_\Omega (F_{R,m_1})\geq \frac12 \|V\|\lfloor_{\Omega}(F_R).
\end{equation}
Next, define
\begin{equation}
\label{rect11}
G_R:=\{x\in \mathbb R^{n+1} : {\rm dist}\,(x, F_{R,m_1})<(1-2^{-\frac{1}{n}})
R\}.
\end{equation}
By definition, $G_R$ is open, and for any $x\in G_R$, there exists 
$y\in F_{R,m_1}$ with $|x-y|<(1-2^{-\frac{1}{n}})R$. By \eqref{rect9}, 
\begin{equation}
\label{rect12}
(2^{-\frac{1}{n}} R)^{-n} \|\Phi_{\e_{j_l}}\ast\partial \mathcal E_{j_l}\|
(B_{2^{-\frac{1}{n}}R}(x))\leq 2 R^{-n} \|\Phi_{\e_{j_l}}\ast
\partial \mathcal E_{j_l}\|(B_R(y))<\Cr{c_2}/8
\end{equation}
for all $l\geq m_1$ and $x\in G_R$. Since $G_R$ is open, we may choose
$m_2\in \mathbb N$ with $m_2\geq m_1$ such that 
\begin{equation}
\|\partial \mathcal E_{j_l}\|\lfloor_\Omega (G_R)\geq \frac12 \|V\|\lfloor_\Omega (
G_R)
\label{rect13}
\end{equation}
for all $l\geq m_2$. Since $F_{R,m_1}\subset G_R$, \eqref{rect13} and 
\eqref{rect10} show
\begin{equation}
\|\partial \mathcal E_{j_l}\|\lfloor_\Omega ( G_R)\geq \frac14 \|V\|\lfloor_\Omega ( F_R)
\label{rect14}
\end{equation}
for all $l\geq m_2$. Choose $m_3\in \mathbb N$ such that $m_3\geq m_2$
and 
\begin{equation}
\label{rect14.5}
\frac{1}{2j_{m_3}^{2}}<\frac{R}{2}.
\end{equation}
Define
\begin{equation}
\label{rect15}
G_{R,j_l,1}:=\{x\in G_R : \theta^{n}(\|\partial \mathcal E_{j_l}\|,x)= 1
\mbox{ and }(2j_l^2)^{n}\|\Phi_{\e_{j_l}}\ast \partial \mathcal E_{j_l}\|
(B_{\frac{1}{2j_l^2}}(x))>\Cr{c_2}/4\}
\end{equation}
and
\begin{equation}
\label{rect16}
G_{R,j_l,2}:=\{x\in G_R : \theta^{n}(\|\partial \mathcal E_{j_l}\|,x)= 1
\mbox{ and }(2j_l^2)^{n}\|\Phi_{\e_{j_l}}\ast \partial \mathcal E_{j_l}\|
(B_{\frac{1}{2j_l^2}}(x))\leq \Cr{c_2}/4\}.
\end{equation}
Since $\theta^{n}(\|\partial\mathcal E_{j_l}\|,x)=1$ for
$\|\partial\mathcal E_{j_l}\|$ a.e.~$x$, we have
\begin{equation}
\label{rect17}
\|\partial\mathcal E_{j_l}\|\lfloor_\Omega ( G_{R,j_l,1}\cup G_{R,j_l,2})=
\|\partial\mathcal E_{j_l}\|\lfloor_\Omega ( G_R).
\end{equation}
First we consider the case $x\in G_{R,j_l,1}$ with $l\geq m_3$. We use
$r_1=\frac{1}{2j_l^2}<2^{-\frac{1}{n}} R=r_2$ in Lemma \ref{red_den}.
Here, the inequality follows from \eqref{rect14.5}. If \eqref{fmono1}
holds with $s:=(2^{-\frac{1}{n}}
R-\frac{1}{2j_l^2})^{-1}(\ln 2)$, then we would have a contradiction to 
\eqref{rect12} and \eqref{rect15}. Thus there exists $\frac{1}{2j_l^2}<
r_x<2^{-\frac{1}{n}} R$ such that \eqref{fmono1} does not hold, i.e.,
\begin{equation}
\label{rect18}
\|\delta(\Phi_{\e_{j_l}}\ast\partial\mathcal E_{j_l})\|(B_{r_x}(x))
>s \|\Phi_{\e_{j_l}}
\ast \partial \mathcal E_{j_l}\|(B_{r_x}(x))\geq \frac{1}{2R} 
\|\Phi_{\e_{j_l}}
\ast \partial \mathcal E_{j_l}\|(B_{r_x}(x)),
\end{equation}
where the last inequality holds from the definition of $s$. Since 
$\e_{j_l}\leq   j_l^{-4}<j_l^{-2}<2r_x$ by (1) for all large $l$, 
$\Phi_{\e_{j_l}}\ast \chi_{B_{r_x}(x)}\geq 
\frac14$ on $B_{r_x}(x)$. Thus we have
\begin{equation}
\label{rect19}
\|\Phi_{\e_{j_l}}\ast\partial \mathcal E_{j_l}\|(B_{r_x}(x))
=\|\partial\mathcal E_{j_l}\|(\Phi_{\e_{j_l}}\ast\chi_{B_{r_x}(x)})
\geq \frac14 \|\partial \mathcal E_{j_l}\|(B_{r_x}(x)).
\end{equation}
By \eqref{musm4}, \eqref{omegaproeq}, \eqref{rect18} and \eqref{rect19},
we have
\begin{equation}
\label{rect20}
\begin{split}
\|\Phi_{\e_{j_l}}\ast\delta(\partial\mathcal E_{j_l})\|\lfloor_\Omega (
B_{r_x}(x))& = \|\delta(\Phi_{\e_{j_l}}\ast \partial\mathcal E_{j_l})\|
\lfloor_\Omega ( B_{r_x}(x)) \\
& \geq \Omega(x)\exp(-2\Cr{c_1} R) \|\delta(\Phi_{\e_{j_l}}\ast \partial 
\mathcal E_{j_l})\|(B_{r_x}(x)) \\
&\geq \frac{1}{8R} \Omega(x)\exp(-2\Cr{c_1}R)\|\partial\mathcal E_{j_l}\|(B_{r_x}(x)) \\
&\geq \frac{1}{8R} \exp(-4\Cr{c_1}R) \|\partial\mathcal E_{j_l}\|\lfloor_\Omega (
B_{r_x}(x)).
\end{split}
\end{equation}
Let $\mathcal C:=\{B_{r_x}(x) : x\in G_{R,j_l,1}\}$, where $r_x$ is 
as above. By the Besicovitch covering theorem, there exists a 
collection of subfamilies $\mathcal C_1,\ldots,\mathcal C_{{\bf B}_{n+1}}$,
each of them consisting of mutually disjoint balls and such that
\begin{equation}
G_{R,j_l,1}\subset \cup_{i=1}^{{\bf B}_{n+1}} \cup_{B_{r_x}(x)\in \mathcal C_i}
B_{r_x}(x).
\label{rect21}
\end{equation} 
Then for some $i_0\in \{1,\ldots, {\bf B}_{n+1}\}$, we have
\begin{equation}
\label{rect22}
\begin{split}
\|\partial\mathcal E_{j_l}\|\lfloor_{\Omega}( G_{R,j_l,1})& \leq {\bf B}_{n+1}
\sum_{B_{r_x}(x)\in \mathcal C_{i_0}} \|\partial\mathcal E_{j_l}\|\lfloor_{\Omega}(
B_{r_x}(x)) \\
&\leq 8R \exp(4\Cr{c_1} R) {\bf B}_{n+1} \sum_{B_{r_x}(x)\in \mathcal C_{i_0}}  
\|\Phi_{\e_{j_l}}\ast\delta(\partial\mathcal E_{j_l})\|\lfloor_\Omega(
B_{r_x}(x)) \\
& \leq 8R \exp(4\Cr{c_1} R) {\bf B}_{n+1} \|\Phi_{\e_{j_l}}\ast\delta(\partial\mathcal E_{j_l})\|\lfloor_\Omega( B_{1+2R}(x_0))
\end{split}
\end{equation}
by \eqref{rect20} and $G_R\subset B_{1+R}(x_0)$. In addition, by 
\eqref{musm2} and the Cauchy-Schwarz inequality, we obtain
\begin{equation}
\label{rect23}
\begin{split}
&\|\Phi_{\e_{j_l}}\ast\delta(\partial\mathcal E_{j_l})\|\lfloor_\Omega (B_{1+2R}(x_0))
= \int_{B_{1+2R}(x_0)} \Omega |\Phi_{\e_{j_l}}\ast \delta(\partial
\mathcal E_{j_l})|\, dx \\
&\leq \Big(\int_{\mathbb R^{n+1}} \frac{\Omega |\Phi_{\e_{j_l}}\ast
\delta(\partial\mathcal E_{j_l})|^2}{\Phi_{\e_{j_l}}\ast\|\partial\mathcal E_{j_l}
\|+\e_{j_l}\Omega^{-1}}\Big)^{\frac12}\Big(\int_{B_{1+2R}(x_0)}
\Omega(\Phi_{\e_{j_l}}\ast\|\partial\mathcal E_{j_l}\|+\e_{j_l}\Omega^{-1})
\Big)^{\frac12} \\
&\leq M^{\frac12}(M+c(n)\e_{j_l})^{\frac12}.
\end{split}
\end{equation}
\eqref{rect22} and \eqref{rect23} prove that, for all fixed $0<R<1$, 
\begin{equation}
\limsup_{l\rightarrow\infty} \|\partial\mathcal E_{j_l}\|\lfloor_\Omega ( G_{R,j_l,1})
\leq 8R \exp(4\Cr{c_1} R){\bf B}_{n+1}M.
\label{rect24}
\end{equation}
Next, suppose that $x\in G_{R,j_l,2}$. From \eqref{rect16} and 
\eqref{rect19} (where $r_x$ may be replaced by $(2j_l^2)^{-1}$
for the same reason), we have
\begin{equation}
\label{rect25}
(2j_l^2)^{n}\|\partial \mathcal E_{j_l}\|(B_{\frac{1}{2j_l^2}}(x))
\leq \Cr{c_2}.
\end{equation}
Note that $x\in {\rm spt}\,\|\partial \mathcal E_{j_l}\|$. Then, 
Proposition \ref{rect0} shows the existence of $r_x\in [\frac{1}{4j_l^2},
\frac{1}{2j_l^2}]$ and a $\mathcal E_{j_l}$-admissible function 
$f_x$ such that
\begin{enumerate}
\item[(i)] $f_x (y)=y$ for $y\in \mathbb R^{n+1}\setminus U_{r_x}(x)$,
\item[(ii)] $f_x(y)\in B_{r_x}(x)$ for $y\in B_{r_x}(x)$,
\item[(iii)] $\|\partial (f_x)_\star \mathcal E_{j_l}\|(B_{r_x}(x))
\leq \frac12 \|\partial \mathcal E_{j_l}\|(B_{r_x}(x)),$
\item[(iv)] $\mathcal L^{n+1} (E_{i}\triangle \tilde E_{x,i})\leq \Cr{c_3} 
(\|\partial \mathcal E_{j_l}\|(B_{r_x}(x)))^{\frac{n+1}{n}}$
for all $i=1,\ldots, N$, 
\end{enumerate}
where $\{E_i\}_{i=1}^{N}
=\mathcal E_{j_l}$ and  $\{\tilde E_{x,i}\}_{i=1}^N=
(f_x)_{\star}\mathcal E_{j_l}.$
By \eqref{omegaproeq}, (iii) may be replaced by
\begin{equation}
\label{rect26}
\|\partial (f_x)_\star \mathcal E_{j_l}\|\lfloor_\Omega ( B_{r_x}(x))
\leq 2^{-\frac12} \|\partial \mathcal E_{j_l}\|\lfloor_\Omega (B_{r_x}(x))
\end{equation}
for all sufficiently large $l$ depending only on $\Cr{c_1}$. 
Applying the Besicovitch covering theorem to the family 
$\{B_{r_x}(x)\}_{x\in G_{R,j_l,2}}$, we have a finite set $\{x_k\}_{k=1}^{
\Lambda}$ such that $\{B_{r_{x_k}}(x_k)\}_{k=1}^{\Lambda}$ is 
mutually disjoint and (writing $B_{r_{x_k}}(x_k)$ as $B(k)$)
\begin{equation}
\|\partial\mathcal E_{j_l}\|\lfloor_\Omega ( \cup_{k=1}^{\Lambda} B(k))\geq {\bf B}_{n+1}^{-1} \|\partial \mathcal E_{j_l}\|\lfloor_\Omega (
G_{R,j_l,2}).
\label{rect27}
\end{equation}
Note that the finiteness of $\Lambda$ follows from $r_x\geq \frac{1}{
4j_l^2}$ and $G_R\subset B_{1+R}(x_0)$. With this choice, define
$f:\mathbb R^{n+1}\rightarrow\mathbb R^{n+1}$ by
\begin{equation}
f(x):=\left\{\begin{array}{ll} f_{x_k}(x) & \mbox{if $x\in B(k)$
for some $k\in \{1,\ldots,\Lambda\}$}, \\
x & \mbox{otherwise}.
\end{array} \right.
\label{rect28}
\end{equation}
Since $f_{x_k}$ is $\mathcal E_{j_l}$-admissible, due to the disjointness of $\{B(k)\}_{k=1}^{\Lambda}$, 
so is $f$. In addition, $f$ belongs to 
${\bf E}(\mathcal E_{j_l},j_l)$. For this, we need to check the
conditions of Definition \ref{mrld} (a)-(c). (a) is satisfied since $\max |f(x)-x|
\leq\max ( {\rm diam}\,B(k))\leq \frac{1}{j_l^2}$. For (b), 
write $f_{\star}\mathcal E_{j_l}=:\{\tilde E_i\}_{i=1}^N$.
Then we have
$E_i\triangle \tilde E_i=\cup_{k=1}^{\Lambda} E_i\triangle
\tilde E_{x_ k,i}$ and (iv) and \eqref{rect25} give
\begin{equation}
\label{rect29}
\begin{split}
\mathcal L^{n+1}(E_i\triangle\tilde E_i)&\leq\Cr{c_3} \sum_{k=1}^{\Lambda}
(\|\partial \mathcal E_{j_l}\|(B(k)))^{\frac{n+1}{n}} \leq 
\frac{\Cr{c_3}\Cr{c_2}^{\frac{1}{n}}}{ 2j_l^2}
\|\partial \mathcal E_{j_l}\|
(\cup_{k=1}^{\Lambda} B(k)) \\
& \leq c(n)(\min_{B_3(x_0)}\Omega)^{-1}\frac{M}{j_l^2}.
\end{split}
\end{equation}
Thus, for all sufficiently large $l$, we have (b). For (c), using ${\rm diam}\,
B(k)\leq 1/j_l^2$ and arguing as in the proof of Lemma \ref{eslip}
with (iii), we may prove 
\begin{equation}
\|\partial f_{\star}\mathcal E_{j_l}\|(\phi)
-\|\partial \mathcal E_{j_l}\|(\phi)=\sum_{k=1}^{\Lambda}
\big(\|\partial (f_{x_k})_{\star}\mathcal E_{j_l}\|\lfloor_\phi ( B(k))
-\|\partial\mathcal E_{j_l}\|\lfloor_\phi ( B(k))\big)\leq 0
\label{rect30}
\end{equation}
for $\phi\in \mathcal A_{j_l}$
for all sufficiently large $l$. Thus we proved $f\in {\bf E}(\mathcal E_{j_l},j_l)$. By \eqref{mrld2}, \eqref{rect26} and \eqref{rect27}, then, we have
\begin{equation}
\label{rect31}
\begin{split}
\Delta_{j_l}\|\partial\mathcal E_{j_l}\|(\Omega) & \leq \|\partial f_{\star}
\mathcal E_{j_l}\|(\Omega)-\|\partial\mathcal E_{j_l}\|(\Omega) \\
&=\sum_{k=1}^{\Lambda} \Big(\|\partial(f_{x_k})_{\star}\mathcal E_{j_l}
\|\lfloor_\Omega ( B(k))-\|\partial \mathcal E_{j_l}\|\lfloor_\Omega ( B(k))\Big)\\
& \leq (2^{-\frac12}-1)\sum_{k=1}^{\Lambda} \|\partial \mathcal E_{j_l}\|
\lfloor_\Omega (B(k)) \leq (2^{-\frac12}-1){\bf B}_{n+1}^{-1} \|\partial 
\mathcal E_{j_l}\|\lfloor_\Omega ( G_{R,j_l,2}).
\end{split}
\end{equation}
\eqref{rect31} and (4) prove 
\begin{equation}
\lim_{l\rightarrow\infty} \|\partial 
\mathcal E_{j_l}\|\lfloor_\Omega ( G_{R,j_l,2})=0,
\label{rect32}
\end{equation}
and \eqref{rect17}, \eqref{rect24} and \eqref{rect32} prove
\begin{equation}
\limsup_{l\rightarrow\infty} \|\partial\mathcal E_{j_l}\|\lfloor_\Omega (
G_R)\leq 8R\exp(4\Cr{c_1}R){\bf B}_{n+1} M.
\label{rect33}
\end{equation}
Recalling \eqref{rect14}, \eqref{rect33} proves $\lim_{R\rightarrow 0}
\|V\|\lfloor_\Omega (F_R)=0$, which proves \eqref{rect3.5}. 
From Proposition \ref{col}, $\|\delta V\|$ is a Radon measure and 
applying Allard's rectifiability theorem \cite[5.5(1)]{Allard}, $V$ is rectifiable. 
\hfill{$\Box$}

\section{Integrality theorem}
\label{secint}
In the following, we write $T\in {\bf G}(n+1,n)$ as the subspace 
corresponding to $\{x_{n+1}=0\}$ and $T^{\perp}\in {\bf G}(n+1,1)$ as the
orthogonal complement $\{x_1=\cdots=x_n=0\}$. As usual, they are identified with the 
$(n+1)\times (n+1)$ matrices 
representing the orthogonal projections to these subspaces. 
Given a set $Y\subset T^{\perp}$ and $r_1,r_2\in (0,\infty)$,
define a closed set
\begin{equation}
E(r_1,r_2):=\{x\in \mathbb R^{n+1} : |T(x)|\leq r_1,\, {\rm dist}\,(T^{\perp}
(x),Y)\leq r_2\}.
\label{itg1}
\end{equation}
\begin{lemma} (\cite[4.20]{Brakke}) 
Corresponding to $n,\nu\in \mathbb N$, $\alpha\in(0,1)$ and $\zeta\in (0,1)$, there exist $\gamma\in (0,1)$ and $j_0\in \mathbb N$ with
the following property.
Assume 
\begin{enumerate}
\item
$\mathcal E=\{E_i\}_{i=1}^N \in \mathcal{OP}_{\Omega}^N$, $j\in 
\mathbb N$ with $j\geq j_0$, 
$R\in (0,\frac12 j^{-2})$,  $\rho\in (0,\frac12 j^{-2})$,
\item
$\rho\geq \alpha R$,
\item 
$Y\subset T^{\perp}$ has no more than $\nu$ elements, ${\rm diam}\, Y<j^{-2}$
and $\theta^{n}(\|\partial \mathcal E\|,y)=1$ for all $y\in Y$,
\newline
and writing $E^*(r):=E(r,(1+R^{-1}r)\rho)$ for short, assume further that
\item
$\int_{{\bf G}_{n}(E^*(r))}\|S-T\|\, d(\partial \mathcal E)
(x,S)\leq \gamma \|\partial\mathcal E\|(E^*(r))$ for all $r\in (0,R)$,
\item 
$\Delta_j \|\partial \mathcal E\|(E^*(r))
\geq -\gamma \|\partial\mathcal E\|(E^*(r))$ for all $r\in (0,R)$.
\end{enumerate}
Then we have
\begin{equation}
\|\partial \mathcal E\|(E(R,2\rho))\geq (\mathcal H^0(Y)-\zeta)\omega_{n} R^{n}.
\label{itg2}
\end{equation}
\label{itglem}
\end{lemma}
\begin{rem} We note that conditions (3), (4) and (5) are different from Brakke's.
The differences are essential to complete the proof of integrality. 
\end{rem}
{\it Proof}. We may assume that 
\begin{equation}
\label{itg4.0}
\mathcal H^0(Y)=\nu
\end{equation}
since the lesser cases
$\mathcal H^0(Y)\in \{1,\ldots,\nu-1\}$ can be equally proved and we may simply 
choose the smallest $\gamma$ and the largest $j_0$ among them. 
We choose and fix a large $j_0\in \mathbb N$ so that
\begin{equation}
\big(\nu-2^{-1}(1+\zeta)\big)(\nu-\zeta)^{-1}<\exp(-4j_0 ^{-1})
\label{itg0}
\end{equation}
which depends only on $\nu$ and $\zeta$. In the following, 
we assume that $\mathcal E$, $j$, $R$, $\rho$ and $Y$ 
satisfy (1)-(5). 
Next we set
\begin{equation}
\label{itg3}
r_1:=\inf\{\lambda>0 : \|\partial \mathcal E\|(E(\lambda,(1+R^{-1}\lambda)\rho))
\leq (\nu-\zeta)\omega_{n} \lambda^{n}\}.
\end{equation}
Since $\cup_{y\in Y} U_\lambda (y)\subset E(\lambda,(1+R^{-1}\lambda)\rho)$ for
$\lambda<\rho$, 
\begin{equation}
\liminf_{\lambda\rightarrow 0} (\omega_{n} \lambda^{n})^{-1}
\|\partial \mathcal E\|(E(\lambda,(1+R^{-1}\lambda)\rho))\geq \sum_{y\in Y}
\theta^{n}(\|\partial \mathcal E\|,y)=\nu
\label{itg4}
\end{equation}
by (3) and \eqref{itg4.0}. Thus, \eqref{itg4} shows $r_1>0$. If $r_1\geq R$, then, 
we would have the opposite inequality in \eqref{itg3} for all $\lambda<R$. By letting 
$\lambda\nearrow R$, we would obtain \eqref{itg2}. In the following, we assume
that $r_1<R$, and look for a contradiction to (5), with an 
appropriate choice of $\gamma$. For the repeated use, we define
\begin{equation}
\rho_1:=(1+R^{-1}r_1)\rho
\label{rhoitg}
\end{equation} 
and note that
\begin{equation}
\|\partial\mathcal E\|(E(r_1,\rho_1))=(\nu-\zeta)\omega_{n} r_1^{n}.
\label{itg5}
\end{equation}
This is because, considering the inequality for $\lambda<r_1$
and letting $\lambda\nearrow r_1$, we have $\geq$. On the other hand, 
there exists a sequence $\lambda_i\geq r_1$ satisfying the inequality in 
\eqref{itg3} and 
letting $i\rightarrow\infty$, we obtain $\leq$. Combined with (4) and (5), 
\eqref{itg5} gives
\begin{equation}
\label{itg5n1}
\int_{{\bf G}_{n}(E(r_1,\rho_1))} \|S-T\|\, d(\partial\mathcal E)(x,S)\leq \gamma (\nu-\zeta)
\omega_{n}r_1^{n}
\end{equation}
and 
\begin{equation}
\Delta_j \|\partial\mathcal E\|(E(r_1,\rho_1))\geq -\gamma(\nu-\zeta)
\omega_{n}r_1^{n}.
\label{itg5n2}
\end{equation}
Next, define
\begin{equation}
V:=\partial \mathcal E\lfloor_{{\bf G}_{n}(E(r_1,\rho_1))}
(=|E(r_1,\rho_1)\cap \cup_{i=1}^N \partial E_i|)
\label{itg6}
\end{equation}
and consider $T_{\sharp}V$, 
the usual push-forward of varifold counting
multiplicities. One notes that 
\begin{equation}
T_{\sharp} V(\phi)
=\int_T \phi(x,T)\mathcal H^0(T^{-1}(x)\cap \big(\cup_{i=1}^N \partial
E_i\big)\cap E(r_1,\rho_1))\, d\mathcal H^{n}(x)
\label{itg6.5}
\end{equation}
for $\phi\in C_c({\bf G}_{n}(\mathbb R^{n+1}))$ and $\theta^{n}(\|T_{\sharp}
V\|,x)=\mathcal H^0(T^{-1}(x)\cap \big(\cup_{i=1}^N \partial
E_i\big)\cap E(r_1,\rho_1))$ for $\mathcal H^{n}$ a.e.~$x\in 
T$. Define
\begin{equation}
A_0:=\{x\in U_{r_1}^{n}, \, \theta^{n}(\|T_{\sharp} V\|,x)\leq \nu-1\}.
\label{itg7}
\end{equation}
For $\mathcal H^{n}$ a.e.~ $x\in U^{n}_{r_1}\setminus A_0$, we
have $\theta^{n}(\|T_{\sharp}V\|,x)\geq \nu$. Thus, 
\begin{equation}
\begin{split}
\nu(\omega_{n}r_1^{n} -\mathcal H^{n}(A_0))&\leq \|T_{\sharp} V\|
(U^{n}_{r_1})=\int_{{\bf G}_{n}(E(r_1,\rho_1))} |\Lambda_{n} T\circ S|\, dV(x,S) \\
&\leq \|V\|(E(r_1,\rho_1))=(\nu-\zeta)\omega_{n} r_1^{n},
\end{split}
\label{itg8}
\end{equation}
where we used \eqref{itg5} and \eqref{itg6} in the last line. By \eqref{itg8} 
we obtain
\begin{equation}
\mathcal H^{n}(A_0)\geq \nu^{-1}\zeta\omega_{n}r_1^{n}.
\label{itg9}
\end{equation}
We next set 
\begin{equation}
\label{itg9.0}
\eta:=\frac{1-\zeta}{8}
\end{equation}
in the following. 
We then choose $s\in (0,1/2)$ depending only on $\nu,\zeta$ and $n$ so that
$
\mathcal H^{n}(U^{n}_1\setminus U^{n}_{1-2s})\leq \eta (2\nu)^{-1}\zeta 
\omega_{n}.$
This implies from \eqref{itg9} that
\begin{equation}
\label{itg9.1}
\mathcal H^{n}(A_0 \cap U^{n}_{r_1(1-2s)})\geq (1-2^{-1}\eta)\nu^{-1}\zeta 
\omega_{n} r_1^{n}.
\end{equation}
We then claim that 
there exist
\begin{equation}
\delta\in (0,s r_1 ) \mbox{ and }A\subset A_0
\label{itg9.2}
\end{equation}
such that
\begin{equation}
\label{itg10}
A\subset U_{r_1(1-2s)}^{n}
\mbox{ and }
\mathcal H^{n}(A)\geq (1-\eta)\nu^{-1}\zeta\omega_{n} r_1^{n},
\end{equation}
and for each $a\in A$, we have 
\begin{equation}
\label{itg11}
\int_{{\bf G}_{n}(C(T,a,\delta))} |\Lambda_{n} T\circ S|\, dV(x,S)\leq
(\nu-1+\eta)\omega_{n} \delta^{n}
\end{equation}
and 
\begin{equation}
\|V\|(C(T,a,\delta))\leq \eta\omega_{n} \delta^{n-1} r_1.
\label{itg12}
\end{equation}
Here, $C(T,a,\delta):=\{x\in \mathbb R^{n+1} : |T(x)-a|\leq \delta\}$. 
The reason for the existence of $A$ and $\delta$ is as follows. 
Since $\theta^{n}(\|T_{\sharp} V\|,\cdot)\leq \nu-1$ on $A_0$, we have
\begin{equation}
\label{itg13}
\lim_{r\rightarrow 0}\frac{1}{\omega_{n}r^{n}}\int_{B_r^{n}(x)}
\theta^{n}(\|T_{\sharp} V\|,y)\, d\mathcal H^{n}(y)
\leq \nu-1
\end{equation}
for a.e.~$x\in A_0$ by the Lebesgue theorem. On the other hand,
\begin{equation}
\label{itg14}
\int_{B_r^{n}(x)}\theta^{n}(\|T_{\sharp}V\|,y)\, d\mathcal H^{n}(y) =
\|T_{\sharp} V\|(B^{n}_r(x)) =\int_{{\bf G}_{n}(C(T,x,r))} |\Lambda_{n} T\circ S|\, dV(y,S).
\end{equation}
Combining \eqref{itg13} and \eqref{itg14}, one may argue that 
for sufficiently small $\delta$, \eqref{itg11} is
satisfied for a set in $A_0$ whose complement can be arbitrarily small 
in measure. For \eqref{itg12}, define $A_{0,\delta}:=\{a\in A_0 : 
\|V\|(C(T,a,\delta))\geq \eta\omega_{n} \delta^{n-1}r_1\}$. By the 
Besicovitch covering theorem, there exists a disjoint family 
$\{B_\delta^{n} (x_i)\}_{i=1}^m$ such that 
\begin{equation}
\label{itg15}
\begin{split}
\mathcal H^{n}
(A_{0,\delta})&\leq {\bf B}_{n} m\omega_{n} \delta^{n} 
\leq {\bf B}_{n} \delta (\eta r_1)^{-1} \sum_{i=1}^m \|V\|(C(T,x_i,\delta))
\\ &
\leq {\bf B}_{n} \delta (\eta r_1)^{-1} (\nu-\zeta)\omega_{n}r_1^{n},
\end{split}
\end{equation}
where we also used \eqref{itg5} and \eqref{itg6}. 
Thus \eqref{itg15} shows that we may choose $\delta$ sufficiently 
small so that the measure of $A_{0,\delta}$ is small. On the 
complement of $A_{0,\delta}$, we have \eqref{itg12}. Comparing 
\eqref{itg9}, \eqref{itg9.1} and \eqref{itg10}, we may thus choose $\delta$ 
and $A\subset A_0$ so that
\eqref{itg10}-\eqref{itg12} are satisfied. We should emphasize that the choice
of $s$ is solely determined by $\zeta, \nu$ and $n$ while
$\delta$ may depend additionally on other quantities. 

Let $\xi\in (0,\frac{\rho_1 r_1}{R})$ be arbitrary and for each $a
\in A$, define 
\begin{equation}
a^{*}:=\frac{r_1 a}{r_1-\delta} ,
\label{itg16}
\end{equation}
\begin{equation}
E_1(a):=\{x\in C(T,a,\delta) : |T(x)-a^{*}|\leq 2\delta\xi^{-1}
(\rho_1 -{\rm dist}\,(T^{\perp}(x),Y))\},
\label{itg17}
\end{equation}
\begin{equation}
E_2(a):=\{x\in C(T,0,r_1)\setminus E_1(a) : |T(x)-a^{*}|\leq 2r_1\xi^{-1}(\rho_1
-{\rm dist}\,(T^{\perp}(x),Y))\}.
\label{itg18}
\end{equation}
We give a few remarks on the definitions \eqref{itg16}-\eqref{itg18}. 
We have 
\begin{equation}
|a-a^{*}|=\frac{\delta}{r_1-\delta}|a|< \frac{\delta r_1}{r_1-\delta}(1-2s)<
\delta<r_1 s
\label{itg20.1}
\end{equation}
by $a\in A$, \eqref{itg9.2} and
\eqref{itg10}, 
so in particular 
\begin{equation}
\label{itg21}
a^*\in U_{r_1(1-s)}^{n}\cap U^{n}_{\delta}(a).
\end{equation}
The choice of $a^*$ is made
so that the radial expansion centered at $T^{-1}(a^*)$ by the factor
of $r_1/\delta$ maps $E_1(a)$ to $E_1(a)\cup E_2(a)$ one-to-one. 
More precisely, let $F_a:\mathbb R^{n+1}\rightarrow\mathbb R^{n+1}$ be defined by 
\begin{equation}
F_a(x):=T^{\perp}(x)+\frac{r_1}{\delta} (T(x)-a^*)+a^*.
\label{itg22}
\end{equation}
Then, one can check that $|T(x)-a|\leq \delta$ if and only if 
$|T(F_a(x))|\leq r_1$ using \eqref{itg16}. The 
latter conditions involving $|T(x)-a^*|$ on $E_1(a)$ and $E_2(a)$ are also equivalent
for $x$ and $F_a(x)$. Thus we have a one-to-one correspondence between $E_1(a)$ 
and $E_1(a)\cup E_2(a)$ by $F_a$, i.e., 
\begin{equation}
F_a(E_1(a))=E_1(a)\cup E_2(a).
\label{itg23}
\end{equation}
By the definition of $E(r_1,
\rho_1)$, one can check that $E_i(a)\subset E(r_1,\rho_1)$
for $i=1,2$. 

With these sets defined, let $f_a:\mathbb R^{n+1}\rightarrow\mathbb R^{n+1}$ be a Lipschitz map such that $f_a(x)=x$ on $\mathbb R^{n+1}
\setminus (E_1(a)\cup E_2(a))$, $f_a\lfloor_{E_1(a)}=F_a\lfloor_{E_1(a)}$, and $f_a$ radially projects $E_2(a)$ onto
$\partial (E_1(a)\cup E_2(a))$ from $T^{-1}(a^*)$. By \eqref{itg23}, 
$f_a$ expands $E_1(a)$ to $E_1(a)\cup E_2(a)$ and ``crushes''
$E_2(a)$ to the boundary $\partial(E_1(a)\cup E_2(a))$. 
It is not difficult to check that $f_a$ is $\mathcal E$-admissible. 
Write $\tilde E_i:={\rm int}\,(f_a(E_i))$. 
We need to check (a)-(c) of Definition \ref{LFOP3}. (c) is
trivial. (a) follows from the bijective nature between $E_1(a)$ 
and $E_1(a)\cup E_2(a)$. For (b), suppose $x\in \partial(E_1(a)
\cup E_2(a))\setminus\cup_{i=1}^N \tilde E_i$. If $x\in \partial E_i$
for some $i$, then $x\in f_a(\partial E_i)$ since $f_a$ is identity there.
If $x\notin \partial E_i$ for all $i$, then there exists some $i$ such that
$x\in E_i$ due to \eqref{siqe}. $f_a^{-1}(x)$ is a closed 
line segment or a point. If this set is all included in $E_i$, then, 
we would have $x\in {\rm int}(f_a(E_i))=\tilde E_i$, a contradiction. 
Thus there is some $y\in \partial E_i\cap f_a^{-1}(x)$ and this 
shows $x\in f_a(\partial E_i)$. Other case when $x\notin  \partial(E_1(a)
\cup E_2(a))$ is easily handled to conclude that (b) holds. Thus 
$f_a$ is $\mathcal E$-admissible. 

\begin{figure}[h]
\centering
\def\svgwidth{0.9\textwidth}
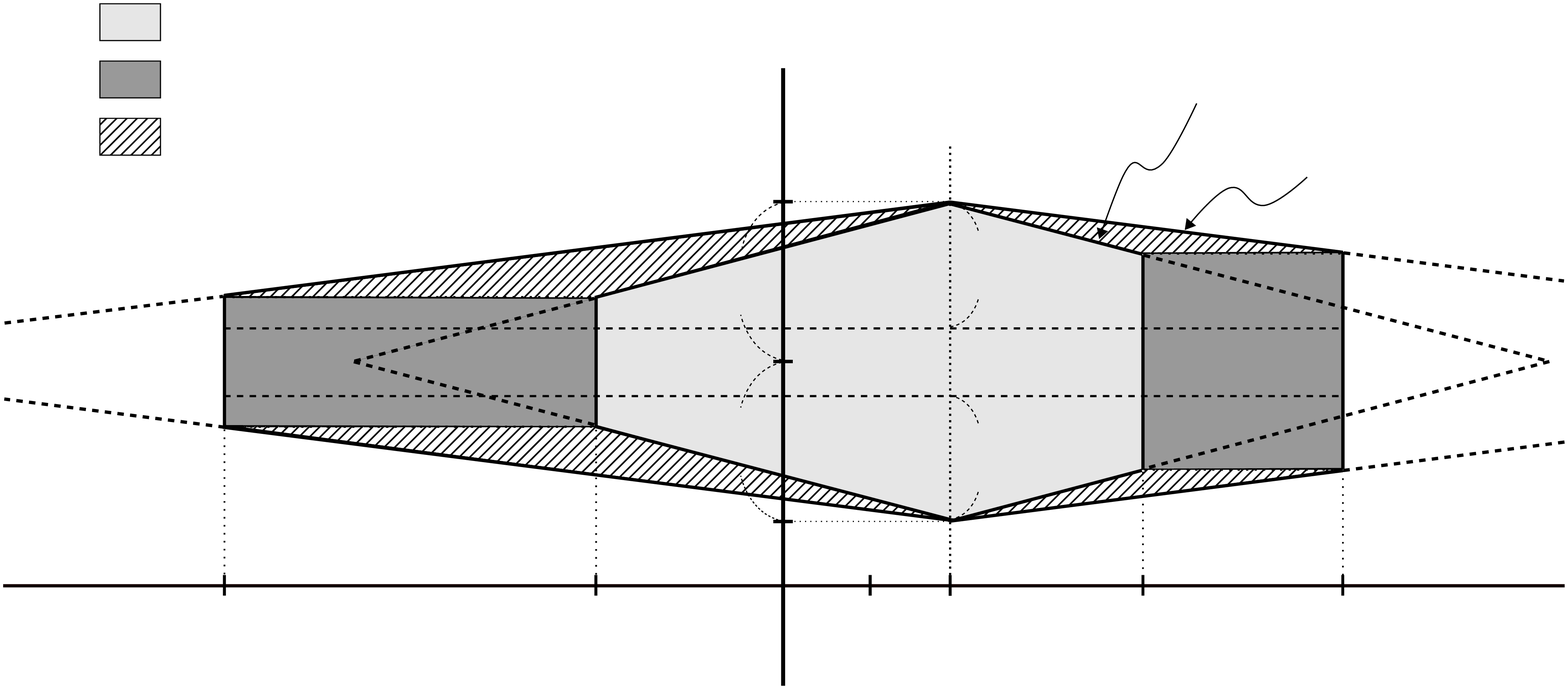
\caption{}
\label{fig_5}
\end{figure}

To separate $E_2(a)$ into two parts, we next define 
\begin{equation}
\label{itg19}
E_3(a):=\{x\in E_2(a) : f_a(x)\in \partial C(T,0,r_1)\},
\end{equation}
\begin{equation}
E_4(a):=E_2(a)\setminus E_3(a).
\label{itg20}
\end{equation}
Note that $\partial (E_1(a)\cup E_2(a))$ consists of the sets in a cylinder
$\partial C(T,0,r_1)$ and cones of type $\{x: |T(x)-a^*|=2r_1\xi^{-1}(\rho_1-{\rm 
dist}\,(T^{\perp}(x),Y))\}$ (see Figure \ref{fig_5} for $n=1$). 
The set $E_3(a)$ thus is the one mapped to the 
cylinder by $f_a$ and $E_4(a)$ is the one to the cones. 

We note that
\begin{equation}
E_4(a)\subset \{x\in E_2(a) : {\rm dist}\,(T^{\perp}(x),Y)\geq\rho_1-\xi\}
\label{itg24}
\end{equation}
and
\begin{equation}
\label{itg24.2}
E(r_1,\rho_1-\xi)\subset E_1(a)\cup E_2(a).
\end{equation}
To see these, for $x\in E_4(a)$, since $f_a(x)$ is a point on the 
cone, we have $|T(f_a(x))-a^*|=2r_1\xi^{-1}(\rho_1-{\rm dist}\,(T^{\perp}(f_a(x))
,Y))$. Since $f_a(x), a^*\in C(T,0,r_1)$, $|T(f_a(x))-a^*|
\leq 2r_1$. By the definition of $f_a$, we have 
$T^{\perp}(f_a(x))=T^{\perp}(x)$. These considerations show \eqref{itg24}.
If $x\in E(r_1,\rho_1-\xi)$, by \eqref{itg1}, $|T(x)|\leq r_1$ and 
${\rm dist}\,(T^{\perp}(x),Y)\leq \rho_1-\xi$. Then we have 
$|T(x)-a^*|\leq 2r_1\leq 2r_1\xi^{-1}(\rho_1-{\rm dist}\,(T^{\perp}(x),Y))$
and \eqref{itg24.2} follows. 

For a given $x\in E_1(a)\cup E_2(a)$, let $v_1,\ldots, v_{n+1}$ be a
set of orthonormal vectors such that $v_1=\frac{T(x)-a^*}{|T(x)-a^*|}$, 
$v_2,\ldots,v_{n}\in T$ and $v_{n+1}\in T^{\perp}$. 
Direct computations show
\begin{equation}
\nabla_{v_i}f_a(x)=\frac{r_1}{\delta} v_i \mbox{ if }1\leq i\leq n\mbox{ and }
\nabla_{v_{n+1}}f_a(x)=v_{n+1} \, \mbox{ on }E_1(a),
\label{nitg1}
\end{equation}
\begin{equation}
\nabla_{v_1}f_a(x)=0\mbox{ on } E_2(a),
\label{nitg2}
\end{equation}
\begin{equation}
\nabla_{v_{n+1}} f_a(x)=v_{n+1}\mbox{ on }E_3(a),
\label{nitg3}
\end{equation}
\begin{equation} \nabla_{v_i}f_a(x)\in T\mbox{ and }|\nabla_{v_i}f_a(x)|\leq \frac{4r_1}{|T(x)-a^*|\sqrt{s}}
\mbox{ if }2\leq i\leq n\mbox{ on }E_3(a),
\label{nitg4}
\end{equation}
\begin{equation}
\nabla_{v_{n+1}}f_a(x)=v_{n+1}\pm 2r_1\xi^{-1} v_1 \mbox{ on }E_4(a),
\label{nitg5}
\end{equation}
\begin{equation}
\label{nitg6}
\nabla_{v_i}f_a(x)\parallel v_i\mbox{ and }|\nabla_{v_i}f_a(x)|\leq\frac{2r_1}{|T(x)-a^*|}\mbox{ if }2\leq i\leq n\mbox{ on }E_4(a).
\end{equation}
Above computations are all valid whenever $\nabla_{v_i}f_a(x)$ is defined. 
On $E_1(a)$, \eqref{itg22} gives \eqref{nitg1}. On $E_2(a)$, since $f_a$ is 
a radial projection in the direction of $v_1$ 
to $\partial(E_1(a)\cup E_2(a))$, we have \eqref{nitg2}. 

For $x\in E_3(a)$ more precisely, $f_a$ is a radial projection from $T^{-1}(a^*)$ of 
$C(T,0,r_1)\setminus C(T,a,\delta)$ to $\partial C(T,0,r_1)$. 
Thus, it is clear that we have \eqref{nitg3}. 
One may express the formula of $f_a$ implicitly by introducing a ``stretching
factor'' $t=t(x)>0$ as
\begin{equation}
f_a(x)=T^{\perp}(x)+ t(T(x)-a^*)+a^*,\,\,
|t(T(x)-a^*)+a^*|^2=r_1^2.
\label{itg29}
\end{equation}
Differentiating both identities of \eqref{itg29} with respect to $v_i$ ($i=2,\ldots,n$), we obtain
\begin{equation*}
\nabla_{v_i} f_a(x)=\nabla_{v_i}t(T(x)-a^*)+tv_i,\,\, f_a(x)\cdot(\nabla_{v_i}t (T(x)-a^*)
+tv_i)=0
\end{equation*}
and
\begin{equation}
\nabla_{v_i}f_a(x)=tv_i-t \frac{f_a(x)\cdot v_i}{f_a(x)\cdot(T(x)-a^*)} (T(x)-a^*)
=tv_i-t\frac{f_a(x)\cdot v_i}{f_a(x)\cdot v_1} v_1.
\label{itg30}
\end{equation}
We need a lower bound of $|f_a(x)\cdot v_1 |$ to estimate \eqref{itg30}.
From \eqref{itg29} and by the definition of $v_2,\ldots,v_n$, we have $f_a(x)\cdot v_i
=a^*\cdot v_i$ for $i=2,\ldots,n$. Then, we have 
\begin{equation}
\label{itg30sup}
|f_a(x)\cdot v_1|^2
=|T(f_a(x))|^2-\sum_{i=2}^n |T(f_a(x))\cdot v_i|^2=r_1^2-|a^*|^2\geq r_1^2-(1-s)^2 r_1^2
\end{equation}
where we used $|T(f_a(x))|=r_1$ and 
$|a^*|<r_1(1-s)$ from 
\eqref{itg21}. Thus we have from \eqref{itg30} and \eqref{itg30sup}
that
\begin{equation}
|\nabla_{v_i} f_a(x)|\leq t\big(1+\frac{1}{\sqrt{2s-s^2}}\big)
\leq  \frac{4r_1}{|T(x)-a^*|\sqrt{s}}.
\label{itg31}
\end{equation}
The last inequality is due to 
$|t(T(x)-a^*)|\leq |t(T(x)-a^*)+a^*|+|a^*|\leq 2r_1$ and $s<1/2$. 
Combined with the expression of \eqref{itg30}, this proves \eqref{nitg4}.

For $x\in E_4(a)$, one can check that
\begin{equation}
f_a(x)=T^{\perp}(x)+\frac{T(x)-a^*}{|T(x)-a^*|}2r_1\xi^{-1}
(\rho_1-{\rm dist}\,(T^{\perp}(x),Y))+a^*.
\label{itg28}
\end{equation}
We have
$\nabla_{v_{n+1}}f_a(x)=v_{n+1}\pm \frac{T(x)-a^*}{|T(x)-a^*|}2r_1\xi^{-1}$, which gives
\eqref{nitg5}. 
For $i=2,\ldots,n$, we have $\nabla_{v_i}f_a(x)=\frac{2r_1\xi^{-1}
(\rho_1-{\rm dist}\,(T^{\perp}(x),Y))}{|T(x)-a^*|} v_i$ since $T(x)-a^*
\parallel v_1$ and $v_1\perp v_i$. Using 
\eqref{itg24}, we obtain \eqref{nitg6}. 

We next need to compute the Jacobian $|\Lambda_{n} \nabla f_a(x)\circ S|$ 
for arbitrary $S\in {\bf G}(n+1,n)$ to compute $\|(f_a)_{\sharp} V\|$. 
As we will check, we may estimate as
\begin{equation}
|\Lambda_{n} \nabla f_a(x)\circ S|\leq\big( \frac{r_1}{\delta}\big)^{n}
|\Lambda_{n} T\circ S|+\big(\frac{r_1}{\delta}\big)^{n-1}\mbox{ on }
E_1(a),
\label{nitg7}
\end{equation}
\begin{equation}
|\Lambda_{n} \nabla f_a(x)\circ S|\leq \|S-T\|\big(\frac{4r_1}{|T(x)-a^*|\sqrt{s}}\big)^{n-1}\mbox{ on }E_3(a),
\label{nitg8}
\end{equation}
\begin{equation}
\label{nitg9}
|\Lambda_{n} \nabla f_a(x)\circ S|\leq \|S-T\|\frac{\sqrt{4r_1^2+\xi^2}}{\xi}\big(\frac{2r_1}{|T(x)-a^*|}\big)^{n-1}
\mbox{ on }E_4(a).
\end{equation}
To see \eqref{nitg7}-\eqref{nitg9}, after an orthogonal rotation, we may 
consider that $v_1,\ldots, v_{n+1}$ are parallel to coordinate axis of $x_1,\ldots,
x_{n+1}$, respectively, and let 
$u_1=(u_{1,1},  \ldots,u_{n+1,1})^{\top}$,
$\ldots$, $u_{n+1}=(u_{1,n+1},\ldots,u_{n+1,n+1})^{\top}$
be a set of orthonormal vectors such 
that $u_1,\ldots,u_{n}$ span $S$ and $u_{n+1}\in S^{\perp}$. 
Then, $|\Lambda_{n} \nabla f_a(x)\circ S|$
is the volume of $n$-dimensional parallelepiped formed by 
$\nabla f_a(x)\circ u_1,\ldots,\nabla f_a(x)\circ u_{n}\in \mathbb R^{n+1}$. Let $L=(L_{i,j})$ 
be the
$(n+1)\times n$ matrix whose column vectors are formed by $\nabla f_a(x)\circ u_1,\ldots,\nabla f_a(x)\circ u_{n}$. Then we have by the Binet-Cauchy formula (\cite[Theorem 3.7]{Evans-Gariepy})
\begin{equation}
|\Lambda _{n} \nabla f_a(x)\circ S|^2={\rm det}( L^\top \circ L)
=\sum_{l=1}^{n+1} ( \det [(L_{i,j})_{i\neq l, 1\leq j\leq n}])^2.
\label{nitg9.5}
\end{equation}
\newline
{\it Computation for \eqref{nitg7}}.
On $E_1(a)$, due to \eqref{nitg1}, 
$\nabla f_a(x)$ is the $(n+1)\times (n+1)$ diagonal matrix whose first $n$ diagonal
elements are all $r_1/\delta$ and whose last diagonal element is $1$. 
Then, the minor formed by eliminating the last row of $L$ is 
$(u_{i,j})_{1\leq i,j\leq n}$ times $r_1/\delta$, and its determinant is 
$(r_1/\delta)^{n}$ times determinant of $(u_{i,j})_{1\leq i,j\leq n}$. 
Note that the determinant of  $(u_{i,j})_{1\leq i,j\leq n}$
is precisely $|\Lambda_{n} T\circ S|$ since $T$ now is the diagonal
matrix whose first $n$ diagonal elements are all $1$ and whose $n+1$-th
diagonal element is $0$. 
For a minor formed by eliminating the $l$-th row of $L$, $1\leq l\leq n$, 
the determinant is $(r_1/\delta)^{n-1}$ times the determinant of 
$(u_{i,j})_{i\neq l, 1\leq j\leq n}$. Considering the orthogonality of the matrix
$(u_{i,j})_{1\leq i,j\leq n+1}$ and the formula for the inverse matrix, the 
determinant is given by $(-1)^{l+n+1} u_{l,n+1}$.  Thus, from \eqref{nitg9.5}, 
we have
\begin{equation}
|\Lambda_{n} \nabla f_a(x)\circ S|^2=\big(\frac{r_1}{\delta}\big)^{2n}
|\Lambda_{n} T\circ S|^2 +\big(\frac{r_1}{\delta}\big)^{2(n-1)}
\sum_{l=1}^{n} (u_{l,n+1})^2.
\label{nitg10}
\end{equation}
Since $|u_{n+1}|=1$, \eqref{nitg10} gives \eqref{nitg7}. 
\newline
{\it Computation for \eqref{nitg8} and \eqref{nitg9}}. 
Here let us write $\nabla f_a(x)$ as $\nabla f$ for short and the $(i,j)$-element of $\nabla f$
as $\nabla f_{i,j}$. 
From \eqref{nitg2}, we have $\nabla f_{i,1}=0$ for all $1\leq i\leq n+1$. 
Then, from \eqref{nitg9.5}, we have
\begin{equation}
\begin{split}
|\Lambda_{n} \nabla f\circ S|^2&={\rm det}\,[ (u_1,\ldots,u_{n})^{\top}\circ
(\nabla f)^{\top}\circ \nabla f\circ(u_1,\ldots,u_{n})] \\
&=({\rm det}\,[(u_{i,j})_{2\leq i\leq n+1, 1\leq j\leq n}])^2 {\rm det}\,[((\nabla f)^{\top}\circ
\nabla f)_{2\leq i,j\leq n+1}].
\end{split}
\label{nitg11}
\end{equation}
By the orthogonality again, we have ${\rm det}\,[(u_{i,j})_{2\leq i\leq n+1,
1\leq j\leq n}]=(-1)^{n}u_{1,n+1}$. Note that $|u_{1,n+1}|
\leq \big(\sum_{i=1}^n (u_{i,n+1})^2\big)^{\frac12}\leq |(T-S)\circ u_{n+1}|$, 
so that $|u_{1,n+1}|\leq \|T-S\|$. Also, considering the fact that 
${\rm det}\,[((\nabla f)^{\top}\circ \nabla f)_{2\leq i,j\leq n+1}]$ is the square of 
$n$-dimensional volume of parallelepiped formed by vectors
$(\nabla f_{1,j},\ldots,\nabla f_{n+1,j})^{\top}$, $j=2,\ldots,n+1$, it is bounded by
$\prod_{j=2}^{n+1}|\nabla _{v_j} f|^2$. These considerations combined with \eqref{nitg11},
\eqref{nitg3} and \eqref{nitg4} give \eqref{nitg8}. 
Similarly using \eqref{nitg5}, \eqref{nitg6} and \eqref{nitg11}, we obtain \eqref{nitg9}. 

We next calculate the mass of $(f_a)_{\sharp}V$. For later use, we note the following.  Since
$\cup_{i=1}^N\partial \tilde E_i\subset f(\cup_{i=1}^N \partial E_i)$
and the varifold push-forward counts the multiplicities of the image, we have
\begin{equation}
\label{itg25}
\|\partial (f_a)_{\star} \mathcal E\| (E(r_1,\rho_1))
=\mathcal H^{n}(\cup_{i=1}^N \partial \tilde E_i\cap E(r_1,\rho_1))
\leq  \|(f_a)_{\sharp} V\|(E(r_1,\rho_1)).
\end{equation}
Now, using \eqref{nitg7}-\eqref{nitg9}, 
we have
\begin{equation}
\begin{split}
\|(f_a)_{\sharp} V\|(E(r_1,\rho_1))& = \int_{{\bf G}_{n}(E(r_1,\rho_1))}
|\Lambda_{n} \nabla f_a(x)\circ S|\, dV(x,S) \\
&\leq \int_{{\bf G}_{n}(E_1(a))} r_1^{n}\delta^{-n}|\Lambda_{n}
T\circ S|+r_1^{n-1}\delta^{1-n}\, dV(x,S) \\
& + \int_{{\bf G}_{n}(E_3(a))} \|S-T\|\big(\frac{4r_1}{|T(x)-a^*|\sqrt{s}}
\big)^{n-1}\, dV(x,S) \\
& + \int_{{\bf G}_{n}(E_4(a))} \|S-T\|\xi^{-1}\sqrt{4r_1^2+\xi^2}
\big(\frac{2r_1}{|T(x)-a^*|}\big)^{n-1}\, dV(x,S) \\
& + \|V\|(E(r_1,\rho_1)\setminus (E_1(a)\cup E_2(a)))=: I_1+\ldots+I_4.
\end{split}
\label{sitg1}
\end{equation}
Since $E_1(a)\subset C(T,a,\delta)$, and by \eqref{itg11} and \eqref{itg12}, 
we have
\begin{equation}
I_1\leq 
(\nu-1+\eta)\omega_{n}r_1^{n}+\eta\omega_{n}r_1^{n}
=(\nu-1+2\eta)\omega_{n}r_1^{n}.
\label{sitg2}
\end{equation}
By defining 
\begin{equation}
c(r_1\xi^{-1}):=\max\{4^{n} s^{\frac{1-n}{2}}, 2^{n}(2r_1\xi^{-1}+1)\},
\label{sitg3}
\end{equation}
we have
\begin{equation}
\label{sitg4}
I_2+I_3\leq c(r_1\xi^{-1})\int_{{\bf G}_{n}(E(r_1,\rho_1))} \|S-T\|\big(\frac{r_1}{|T(x)-a^*|}
\big)^{n-1}\, dV(x,S).
\end{equation}
By \eqref{itg24.2}, we have
\begin{equation}
\label{sitg5}
I_4\leq \|V\|(E(r_1,\rho_1)\setminus E(r_1,\rho_1-\xi))=\|V\|(E(r_1,\rho_1))
-\|V\|(E(r_1,\rho_1-\xi)).
\end{equation}
On the other hand, due to \eqref{itg3}, we have 
$\|V\|(E(\lambda,(1+\lambda R^{-1})\rho))>(\nu-\zeta)\omega_{n}\lambda^{n}$
for $\lambda<r_1$. Hence, for $\lambda_*:=r_1-R\xi\rho^{-1}$ which solves
$\rho_1-\xi=(1+\lambda_*R^{-1})\rho$, we have
\begin{equation}
\|V\|(E(r_1,\rho_1-\xi))\geq \|V\|(E(\lambda_*,(1+\lambda_*/R)\rho))
>(\nu-\zeta)\omega_{n}\lambda_*^{n}.
\label{sitg6}
\end{equation}
By Bernoulli's inequality, we have $\lambda_*^{n}=(r_1-R\xi\rho^{-1})^{n}
\geq r_1^{n}-nr_1^{n-1}R\xi\rho^{-1}$, and \eqref{sitg5}, \eqref{sitg6}
and \eqref{itg5} show
\begin{equation}
\label{sitg7}
I_4\leq (\nu-\zeta)\omega_{n}n r_1^{n-1}R\xi\rho^{-1}
\leq \nu n\omega_{n} \alpha^{-1} (\xi r_1^{-1}) r_1^{n},
\end{equation}
where we used (2) ($\rho\geq \alpha R$) in the last inequality. 
The estimates so far hold for any $a\in A$. To estimate $I_2+I_3$, 
we integrate the right-hand side of \eqref{sitg4} with respect to $a$. 
For any fixed $x\in E(r_1,\rho_1)$, by \eqref{itg16}, 
\begin{equation}
\begin{split}
\int_A 
\big(\frac{r_1}{|T(x)-a^*|}\big)^{n-1}\, & d\mathcal H^{n}(a)=
\big(\frac{r_1-\delta}{r_1}\big)^{n-1} \int_A \big(\frac{r_1}{|\frac{r_1-\delta}{r_1}T(x)-a|}
\big)^{n-1}\, d\mathcal H^{n}(a) \\
&\leq \int_{B_{2r_1}^{n}}\big(\frac{r_1}{|y|}\big)^{n-1}\, d\mathcal H^{n}(y)
= 2n\omega_{n} r_1^{n}\end{split}
\label{sitg7.5}
\end{equation}
after a change of variable $y=\frac{r_1-\delta}{r_1}T(x)-a$ and using 
$\{y : \frac{r_1
-\delta}{r_1}T(x)-y\in A\}\subset B_{2r_1}^{n}$. Then, by Fubini's theorem
and \eqref{sitg7.5}, 
\begin{equation}
\begin{split}
&\int_A\, d\mathcal H^{n}(a)\int_{{\bf G}_{n}(E(r_1,\rho_1))} \|S-T\|
\big(\frac{r_1}{|T(x)-a^*|}\big)^{n-1}\, dV(x,S) \\
&\leq 2n\omega_{n} r_1^{n} \int_{{\bf G}_{n}(E(r_1,\rho_1))}\|S-T\|\, dV(
x,S) \leq 2n\omega_{n}^2\nu  r_1^{2n} \gamma
\end{split}
\label{sitg8}
\end{equation}
where \eqref{itg5n1} is used. By \eqref{itg10} and \eqref{sitg8}, there exists $a\in
A$ such that we have
\begin{equation}
\int_{{\bf G}_{n}(E(r_1,\rho_1))} \|S-T\|
\big(\frac{r_1}{|T(x)-a^*|}\big)^{n-1}\, dV(x,S)\leq 2n(1-\eta)^{-1}\nu^2\zeta^{-1}
\omega_{n}\gamma r_1^{n}. 
\label{sitg9}
\end{equation}
With this choice of $a$, \eqref{sitg1}, \eqref{sitg2}, \eqref{sitg4}, \eqref{sitg7}
and \eqref{sitg9} show
\begin{equation}
\begin{split}
&\|(f_a)_{\sharp} V\|(E(r_1,\rho_1)) \\ &\leq \{\nu-1+2\eta+c(r_1\xi^{-1})2n(1-\eta
)^{-1}\nu^2\zeta^{-1}\gamma+\nu n\alpha^{-1}\xi r_1^{-1}\}\omega_{n}r_1^{n}.
\end{split}
\label{sitg10}
\end{equation}
Up to this point, $\xi\in (0,\frac{\rho_1 r_1}{R})$ is arbitrary. Fix $\xi$ 
so that $\nu n\alpha^{-1} \xi r_1^{-1}=\eta$. Since $\rho_1>\rho$ and
$\rho\geq \alpha R$, one can check that $\xi\in (0,\rho_1r_1/R)$. 
The choice of $\xi r_1^{-1}$
depends only on $\nu, \zeta, n, \alpha$. This fixes $c(r_1\xi^{-1})$ in \eqref{sitg3},
and $c(r_1\xi^{-1})$ depends only on $\nu,\zeta,n,\alpha$. We then 
restrict $\gamma$ so that $c(r_1\xi^{-1})2n(1-\eta)^{-1}\nu^2\zeta^{-1}
\gamma\leq \eta$, which again depends only on the same constants. 
Then we have from \eqref{sitg10} and \eqref{itg9.0} that 
\begin{equation}
\label{sitg11}
\|(f_a)_{\sharp} V\|(E(r_1,\rho_1))\leq (\nu-1+4\eta)\omega_{n} r_1^{n}
=\big(\nu-1+2^{-1}(1-\zeta)\big)\omega_{n}r_1^{n}.
\end{equation}
We next check that $f_a\in {\bf E}(\mathcal E,E(r_1,\rho_1), j)$ by using Lemma \ref{eslip}.
We have already seen that $f_a$ is $\mathcal E$-admissible. We may 
take $C=E(r_1,\rho_1)$ in Lemma \ref{eslip} and (a) is satisfied. 
Since $T^{\perp}(f_a(x) -x)=0$, $f_a(x)\in C(T,0,r_1)$ for $x\in 
E(r_1,\rho_1)$ and $r_1<R<\frac12
j^{-2}$ (by (1)), we have $|f_a(x)-x|\leq 2r_1<j^{-2}$ so we have (b) satisfied.
For (c), we have $\tilde E_i\triangle E_i\subset E(r_1,\rho_1)$ and 
due to (1) and (3), ${\rm diam}\, E(r_1,\rho_1)<4j^{-2}$ (note \eqref{rhoitg}).
Thus for suitably restricted $j$ depending on $n$, we have (c). 
For (d), by \eqref{itg25}, \eqref{sitg11}, \eqref{itg5} and \eqref{itg0} we have
\begin{equation}
\label{sitg12}
\|\partial (f_a)_{\star}\mathcal E\|(E(r_1,\rho_1))
\leq \exp(-4 j_0^{-1}) \|\partial \mathcal E\|(E(r_1,\rho_1)).
\end{equation}
Since ${\rm diam}\, E(r_1,\rho_1)< 4j^{-2}$, we have (d), and we have
$f_a\in {\bf E}(\mathcal E,E(r_1,\rho_1),j)$. Finally, consider $\Delta_j\|\partial \mathcal
E\|(E(r_1,\rho_1))$. By \eqref{itg5n2}, \eqref{mrld3}, 
\eqref{itg5} and \eqref{sitg11},
we have
\begin{equation}
\label{sitg13}
\begin{split}
-\gamma(\nu-\zeta)\omega_{n}r_1^{n}&\leq \Delta_j\|\partial\mathcal E\|(E(r_1,\rho_1)) \\
&\leq \|\partial (f_a)_{\star}\mathcal E\|(E(r_1,\rho_1))-\|\partial \mathcal E
\|(E(r_1,\rho_1)) \\
&\leq -2^{-1}(1-\zeta)\omega_{n}r_1^{n}.
\end{split}
\end{equation}
By restricting $\gamma$ further depending only on $\zeta$ and $\nu$, 
\eqref{sitg13}
is a contradiction. This concludes the proof. 
\hfill{$\Box$}

For large length scale ($\geq j^{-2}$), we use the following. 
\begin{lemma}(\cite[Section 4.21]{Brakke})
Suppose
\begin{enumerate}
\item
$\nu\in {\mathbb N}$, $\xi\in (0,1)$, $M\in (1,\infty)$, $0<r_0<R<\infty$, $T\in {\bf G}(n+1,n)$ and $V\in {\bf V}_n({\mathbb R}^{n+1})$,
\item
$Y\subset T^{\perp}$ has no more than $\nu+1$ elements,
\item
$(M+1){\rm diam}\, Y\leq R$,
\item
$r_0<(3\nu)^{-1}{\rm diam}\, Y$,
\item
$R\|\delta V\|(B_r(y))\leq \xi \|V\|(B_r(y))$ for all $y\in Y$ and $r\in (r_0,R)$,
\item
$\int_{{\bf G}_{n}(B_r(y))}\|S-T\|\, dV(x,S)\leq \xi \|V\|(B_r(y))$ for all $y\in Y$ and $r
\in (r_0,R)$.
\end{enumerate}
Then there are $V_1,\, V_2\in {\bf V}_n({\mathbb R}^{n+1})$ and a partition of $Y$ 
into subsets $Y_0,\, Y_1,\, Y_2$ such that 
\begin{equation}
V\geq V_1+V_2,
\label{4-21-7}
\end{equation}
\begin{equation}
\mbox{neither } Y_1\mbox{ nor } Y_2\mbox{ has more than }\nu\mbox{ elements},
\label{4-21-8}
\end{equation}
\begin{equation}
(M\,{\rm diam}\, Y)\|\delta V_j\|(B_r(y))\leq 2M(\nu+1)(3\nu M)^{n+1}(\exp \xi)\xi \|V_j\|(B_r(y))
\label{4-21-9}
\end{equation}
\centerline{for all $y\in Y_j$, $r\in (r_0,M\, {\rm diam}\, Y)$ and $j=1,2$,}
\begin{equation}
\int_{{\bf G}_n(B_r(y))}\|S-T\|\, dV_j(x,S)\leq M(3\nu M)^n (\exp\xi)\xi \|V_j\|(B_r(y))
\label{4-21-10}
\end{equation}
\centerline{for all $y\in Y_j$, $r\in (r_0,M\, {\rm diam}\, Y)$ and $j=1,2$,}
\begin{equation}
V_j\geq V\lfloor\{x\in {\mathbb R}^{n+1}:{\rm dist}\, (T^{\perp}(x),Y_j)\leq r_0\} \times{\bf G}(n+1,n)\,\, \mbox{ for } j=1,2,
\label{4-21-11}
\end{equation}
\begin{equation}
\begin{split}
&\{(1+1/M)^n+(\nu+1)/M\}  (\exp\xi)\frac{\|V\|(\{x:{\rm dist}\, (x,Y)\leq R\})}{\omega_n R^n} \\
&\geq \sum_{y\in Y_0} \frac{\|V\|(B_{r_0}(y))}{\omega_n r_0^n} 
+\sum_{j=1,2} \frac{\|V_j\|(\{x :{\rm dist}\,(x,Y_j)\leq M\, {\rm diam}\, Y\})}{\omega_n 
(M\, {\rm diam}\, Y)^n}.
\end{split}
\label{4-21-12}
\end{equation}
\label{lemma4.21}
\end{lemma}
The proof of Lemma \ref{lemma4.21} is the same as \cite[Lemma 6.1]{Allard} except that $r_0\rightarrow 0$ in
\cite{Allard} while it is stopped at a positive radius $r_0$ here. 
\begin{lemma}
Corresponding to $n,\nu\in {\mathbb N}$ and $\lambda\in (1,2)$, there exist $\gamma\in (0,1)$ 
and ${\tilde M}\in (1,\infty)$ 
with the following property. Suppose 
\begin{enumerate}
\item
$0<r_0<R<\infty,\, T\in {\bf G}(n+1,n),\, V\in {\bf V}_n({\mathbb R}^{n+1})$,
\item
$Y\subset T^{\perp}$ has  no more than $\nu+1$ elements,
\item
$\{(1+3\nu)^2+{\tilde M}^2\}^{\frac12} r_0 <R$,
\item
${\rm diam}\, Y\leq \gamma R$,
\item
$R\|\delta V\|(B_r(y))\leq \gamma \|V\|(B_r(y))$ for all $y\in Y$ and $r\in (r_0,R)$,
\item
$\int_{{\bf G}_n(B_r(y))}\|S-T\|\, dV(x,S)\leq \gamma \|V\|(B_r(y))$ for all $y\in Y$ and
$r\in (r_0,R)$.
\end{enumerate}
Then there exists a partition of $Y$ into subsets $Y_0,Y_1,\ldots, Y_J$ such that
\begin{equation}
{\rm diam}\, Y_j\leq 3\nu  r_0\,\mbox{ for all } j\in \{1,\ldots,J\},
\label{4-21-17}
\end{equation}
\begin{equation}
\begin{split}
\lambda &\frac{\|V\|(\{x : {\rm dist}\,(x,Y)\leq R\})}{\omega_n R^n}\geq
\sum_{y\in Y_0} \frac{\|V\|(B_{r_0}(y))}{\omega_n r_0^n} \\
&+\sum_{j=1}^J \frac{\|V\|(\{x: {\rm dist}\,(T^{\perp}(x),Y_j)\leq 
r_0,\, |T(x)|\leq {\tilde M}r_0\})}{\omega_n ({\tilde M}r_0)^n}.
\end{split}
\label{4-21-18}
\end{equation}
\label{thm4.21.1}
\end{lemma}
{\it Proof}. We use Lemma \ref{lemma4.21} to partition $Y$ into subsets whose diameters are all smaller than $3\nu r_0$. 
In the case $Y$ consists of only one element, we may take $Y_0:=Y$ and 
Lemma \ref{red_den} shows 
\eqref{4-21-18} by choosing an appropriately small $\gamma$ in (5)
depending only on $\lambda$. We do 
not need ${\tilde M}$ in this case. If $Y$ consists of more than one element, we apply Lemma \ref{lemma4.21}. We separate into two cases first.
\newline
{\it First inductive step : Case 1}. Suppose (4) of Lemma \ref{lemma4.21} is not satisfied, i.e., 
\begin{equation}
 {\rm diam}\, Y\leq 3\nu r_0.
\label{4-21-19}
\end{equation}
In this case, we set $J=1$, 
$Y_1:=Y$ and $Y_0=\emptyset$. 
For any $y\in Y$, we have by \eqref{4-21-19} 
\begin{equation}
\{x: {\rm dist}\,(T^{\perp}(x),Y_1)\leq r_0,\, |T(x)|\leq {\tilde M}r_0\}\subset B_{r_0((1+3\nu)^2+{\tilde M}^2)^\frac12}(y).
\label{4-21-19.2}
\end{equation}
We have
\begin{equation}
\frac{\|V\|(B_{ r_0((1+3\nu)^2+{\tilde M}^2)^{\frac12}}(y))}{\omega_n (r_0{\tilde M})^n}
= \frac{\|V\|(B_{r_0 ((1+3\nu)^2+{\tilde M}^2)^{\frac12}}(y))}{\omega_n ( r_0
((1+3\nu)^2+{\tilde M}^2)^{\frac{1}{2}})^n}\big(1+\frac{(1+3\nu)^2}{{\tilde M}^2}\big)^\frac{n}{2}.
\label{4-21-19.1}
\end{equation}
By Lemma \ref{red_den} with (5), \eqref{4-21-19}, (3)
and \eqref{4-21-19.1}, we have
\begin{equation}
\frac{\|V\|(B_{ r_0((1+3\nu)^2+{\tilde M}^2)^{\frac12}}(y))}{\omega_n (r_0{\tilde M})^n}
\leq (\exp \gamma)\big(1+\frac{(1+3\nu)^2}{{\tilde M}^2}\big)^\frac{n}{2} \frac{\|V\|(B_R(y))}{\omega_n R^n}.
\label{4-21-20}
\end{equation}
Since $B_R(y)\subset \{x:{\rm dist}\,(x,Y)\leq R\}$, combining \eqref{4-21-19.2}, \eqref{4-21-20},
we choose large ${\tilde M}$
and small $\gamma$ depending only on $n,\nu$ and $\lambda$ so that
\eqref{4-21-18} is satisfied. 
\newline
{\it First inductive step : Case 2}. Suppose
(4) of Lemma \ref{lemma4.21} is satisfied. With $M$
satisfying (3) of Lemma \ref{lemma4.21} and $\xi=\gamma$, we apply
Lemma \ref{lemma4.21}. Thus we have a partition of $Y$ into $Y_0,Y_1,Y_2$
with \eqref{4-21-7}-\eqref{4-21-12}. 
\newline
{\it Second inductive step for $Y_1$ and $Y_2$}. We next proceed just like before for
$Y_1$ and $Y_2$. That is, for each $j=1,2$, if $Y_j=\{y\}$, 
we use Lemma \ref{red_den} with \eqref{4-21-9} to derive
\begin{equation}
\frac{\|V\|(B_{r_0}(y))}{\omega_n r_0^n}\leq 
\exp \{2M(\nu+1)(3\nu M)^{n+1}(\exp \gamma)\gamma\}
\frac{\|V_j\|(B_{ M {\rm diam}\, Y}(y))}{\omega_n (M{\rm diam}\,Y)^n}
\label{4-21-21}
\end{equation}
where we have also used \eqref{4-21-11}. Note that the right-hand side
of \eqref{4-21-21} is bounded from above via \eqref{4-21-12}. We add this $Y_j$
to $Y_0$. 
Suppose $Y_j$ consists of more than
one point, and furthermore, \eqref{4-21-19} is satisfied with $Y_j$ in place of $Y$.
Note that \eqref{4-21-11} shows
\begin{equation}
\begin{split}
\|V\|(\{x: {\rm dist}&\,(T^{\perp}(x),Y_j)\leq r_0,\, |T(x)|\leq {\tilde M}r_0\}) \\
&\leq \|V_j\|(\{x: {\rm dist}\,(T^{\perp}(x),Y_j)\leq r_0,\, |T(x)|\leq {\tilde M}r_0\}).
\end{split}
\label{4-21-22}
\end{equation}
We then go through the same argument \eqref{4-21-19.2}-\eqref{4-21-20} 
for $V_j$ in place of $V$ and for $M {\rm diam}\, Y$ in place of $R$ there. 
Note that we may apply Lemma \ref{red_den} for $V_j$ due to \eqref{4-21-9}. 
For doing so, we may achieve $r_0((1+3\nu)^2+{\tilde M}^2)^{\frac12}< M\,{\rm diam}\, Y$ 
since ${\rm diam}\, Y>3\nu r_0$ holds and since we may choose $M$ greater than ${\tilde M}$
by a factor depending only on $\nu$. If $Y_j$ does not satisfy \eqref{4-21-19}, then we 
again apply Lemma \ref{lemma4.21} to $Y_j$ to obtain a partition. Since the number of
elements of $Y_j$ is strictly decreasing in each step, the process ends at most after $\nu$
times. Depending only on $n, \nu$ and $\lambda$, choose large ${\tilde M}$ and $M$,
and then small $\gamma$.  Note that we need not take the same $M$ in this inductive step.
If we take $M$ in the first step, we may take $M-1$ as $M$ in Lemma \ref{lemma4.21}
in the next step so that (3) of Lemma \ref{lemma4.21} is
automatically satisfied (since $((M-1)+1) {\rm diam}\, Y_1\leq M\,{\rm diam}\, Y$, for example).
\hfill{$\Box$}
\begin{lemma}
Corresponding to $n,\nu\in \mathbb N$ and $\lambda\in (1,2)$, there exist
$\gamma, \eta\in (0,1)$, $\tilde M\in (1,\infty)$ and $j_0\in \mathbb N$ with the following property. 
Suppose 
\begin{enumerate}
\item $\mathcal E\in \mathcal{OP}_{\Omega}^N$,
$j\in \mathbb N$ with $j\geq j_0$,
\item $\e\leq \gamma j^{-4}$,
\item $\eta j^{-2}<R$,
\item $Y\subset T^{\perp}$ has no more than $\nu$ elements and $\theta^{n}
(\|\partial\mathcal E\|,y)=1$ for each $y\in Y$,
\item ${\rm diam}\, Y\leq \gamma R$, 
\item $R\|\delta(\Phi_{\e}\ast \partial\mathcal E)\|(B_r(y))\leq \gamma \|\Phi_{\e}\ast \partial \mathcal E\|(B_r(y))
$ for all $y\in Y$ and $r\in (\eta^2 j^{-2},R)$, 
\item $\int_{{\bf G}_{n}(B_r(y))} \|S-T\|\, d(\Phi_{\e}\ast \partial \mathcal E)(x,S)\leq \gamma \|\Phi_{\e}
\ast \partial \mathcal E\|(B_r(y))$ for all $y\in Y$ and $r\in (\eta^2 j^{-2},R)$, 
\newline
and writing
\begin{enumerate}
\item $\tilde R_1:=\eta^2 j^{-2}\lambda^{-\frac{1}{4n}}$,
\item $\tilde R_2:=\tilde M \eta^2 j^{-2}\lambda^{-\frac{1}{4n}}$,
\item $\rho:=\frac12\eta^2j^{-2}(1-\lambda^{-\frac{1}{4n}})$,
\newline
and for any subset $Y'\subset Y$, define
\item $E^*_1(r,Y'):=\{x\in\mathbb R^{n+1} : |T(x)|\leq r, {\rm dist}(Y',T^{\perp}(x))
\leq (1+\tilde R_1^{-1}r)\rho\}$,
\item $E^*_2(r,Y'):=\{x\in\mathbb R^{n+1} : |T(x)|\leq r, {\rm dist}(Y',T^{\perp}(x))
\leq (1+\tilde R_2^{-1}r)\rho\}$, 
\end{enumerate}
and assume for all $Y'\subset Y$ with ${\rm diam}\,Y'<j^{-2}$, $i=1,2$
and 
$r\in (0,j^{-2})$ that
\item $\int_{{\bf G}_{n}(E^*_i(r,Y'))}
\|S-T\|\, d(\partial \mathcal E)(x,S)\leq \gamma \|\partial\mathcal E\|(E^*_i(r,Y'))$,
\item $\Delta_j \|\partial \mathcal E\|(E^*_i(r,Y'))
\geq -\gamma \|\partial\mathcal E\|(E^*_i(r,Y')) $.
\end{enumerate}
Then we have
\begin{equation}
\lambda \|\Phi_{\e}\ast \partial \mathcal E\|(\{x : {\rm dist}\,(x,Y)\leq R\})\geq 
\omega_{n} R^{n}\mathcal H^0(Y).
\label{afitg7}
\end{equation}
\label{lemafi}
\end{lemma}
{\it Proof}. Given $\lambda\in (1,2)$, 
we first use Lemma \ref{thm4.21.1} with $\lambda$ there
replaced by $\lambda^{\frac14}$ to obtain $\gamma_1
\in (0,1)$
and $\tilde M\in (1,\infty)$ depending only on $n,\nu$ and $\lambda$
with the stated property. Choose $\eta\in (0,1)$ depending only on 
$n,\nu$ and $\lambda$ so that
\begin{equation}
(2 \tilde M +3\nu)\eta <1.
\label{fitg0}
\end{equation}
By setting 
\begin{equation}
\alpha:=\frac{1}{2 \tilde M}\lambda^{\frac{1}{4n}}(1-\lambda^{-\frac{1}{4n}})
\in (0,1)
\label{fitg1}
\end{equation}
and fixing
\begin{equation}
\label{fitg2}
\zeta:=1-\lambda^{-\frac14}\in (0,1),
\end{equation}
we use Lemma \ref{itglem} to obtain $\gamma_2\in (0,1)$ and $j_0\in
\mathbb N$ depending only on $n,\nu$ and $\lambda$ 
with the stated property. We assume that $\gamma\leq \min\{\gamma_1,\gamma_2\}$ and assume that we have (1)-(9).  
We set 
\begin{equation}
r_0:=\eta^2 j^{-2}
\label{fitg3}
\end{equation}
in Lemma \ref{thm4.21.1}. We first check that
the assumptions of Lemma \ref{thm4.21.1} are satisfied, where 
$V$ is replaced by $\Phi_{\e}\ast \partial \mathcal E$. By (3), 
we have $r_0<R$. By \eqref{fitg3}, \eqref{fitg0}
and (3), we have $\{(1+3\nu)^2+\tilde M^2\}^{\frac12}r_0
\leq (2\tilde M+3\nu)\eta^2 j^{-2}<\eta j^{-2}<R$. Thus we have 
(3) of Lemma \ref{thm4.21.1}. Note that (2), (4)-(6) of Lemma \ref{thm4.21.1} 
follow from (4)-(7) of Lemma \ref{lemafi}. Thus all the assumptions of 
Lemma \ref{thm4.21.1} are
satisfied, and there exists a partition of $Y$ into $Y_0,Y_1,\ldots,Y_J$
such that
\begin{equation}
\label{fitg4}
{\rm diam}\,Y_l\leq 3\nu \eta^2 j^{-2}<j^{-2} \mbox{ for all }l\in \{1,\ldots,
J\},
\end{equation}
\begin{equation}
\begin{split}
&\lambda^{\frac14} \frac{\|\Phi_{\e}\ast \partial \mathcal E\|(\{x : {\rm dist}\,(x,Y)\leq R\})}{
\omega_{n}R^{n}}\geq \sum_{y\in Y_0}\frac{\|\Phi_{\e}\ast \partial \mathcal E\|(
B_{\eta^2 j^{-2}}(y))}{\omega_{n}(\eta^2 j^{-2})^{n}} \\
&+\sum_{l=1}^J \frac{\|\Phi_{\e}\ast \partial \mathcal E\|(\{x : {\rm dist}\,(T^{\perp}(x),Y_l)
\leq \eta^2 j^{-2}, |T(x)|\leq \tilde M \eta^2 j^{-2}\})}{\omega_{n} (\tilde M
\eta^2 j^{-2})^{n}} .
\end{split}
\label{fitg5}
\end{equation}
Depending only on $n,\nu$ and $\lambda$, there exists $\gamma_3\in (0,\eta^8)$ such 
that, if $\e<\gamma_3 j^{-4}$, 
\begin{equation}
\begin{split}
&\lambda^{\frac{1}{4}} \Phi_{\e}\ast \chi_{B_{\eta^2 j^{-2}}(y)}
\geq 1 \\ &\mbox{ on }S_0(y):=\{x : |T^{\perp}(x)-y|\leq \eta^2 j^{-2}(1-\lambda^{-
\frac{1}{4n}}),|T(x)|\leq \eta^2 j^{-2}\lambda^{-\frac{1}{4n}}\},
\end{split}
\label{fitg6}
\end{equation}
\begin{equation}
\label{fitg7}
\begin{split}
&\lambda^{\frac{1}{4}}\Phi_{\e}\ast \chi_{\{x \,:\, {\rm dist}\,(T^{\perp}(x),Y_l)
\leq \eta^2 j^{-2}, |T(x)|\leq \tilde M \eta^2 j^{-2}\}} \geq 1 \\
& \mbox{ on }S_l:=\{x : {\rm dist}\,(T^{\perp}(x),Y_l)
\leq \eta^2 j^{-2}(1-\lambda^{-\frac{1}{4n}}), |T(x)|\leq \tilde M \eta^2 j^{-2}
\lambda^{-\frac{1}{4n}}\}.
\end{split}
\end{equation}
Since $\|\Phi_{\e}\ast \partial \mathcal E\|(B_{\eta^2 j^{-2}}(y))=\|\partial \mathcal E\|(\Phi_{\e}\ast 
\chi_{B_{\eta^2 j^{-2}}(y)})$ and similarly for the other cases, 
\eqref{fitg5}-\eqref{fitg7} show
\begin{equation}
\label{fitg8}
\begin{split}
\lambda^{\frac34}& \frac{\|\Phi_{\e}\ast \partial \mathcal E\|(\{x : {\rm dist}\,(x,Y)\leq R\})}{
\omega_{n}R^{n}} \\ &\geq \sum_{y\in Y_0}\frac{\|\partial \mathcal E\|(
S_0(y))}{\omega_{n}(\eta^2 j^{-2} \lambda^{-\frac{1}{4n}})^{n}} 
+\sum_{l=1}^J \frac{\|\partial \mathcal E\|(S_l)}{\omega_{n} (\tilde M
\eta^2 j^{-2}\lambda^{-\frac{1}{4n}})^{n}} .
\end{split}
\end{equation}
We now use Lemma \ref{itglem}. For elements in $Y_0$, we let 
$R=\eta^2 j^{-2} \lambda^{-\frac{1}{4n}}$ (the reader should not
confuse this $R$ with $R$ in the statement of the present Lemma) and $\rho=
\frac12\eta^2 j^{-2}(1-\lambda^{-\frac{1}{4n}})$, and for $Y_1,\ldots,Y_J$,
we let 
$R=\tilde M\eta^2 j^{-2} \lambda^{-\frac{1}{4n}}$ and the same $\rho$.
Since they are similar, we only give the detail for $Y_l$, $l\in \{1,\ldots,J\}$. 
We check that the assumptions of Lemma \ref{itglem} are satisfied first.
We have already assumed $j\geq j_0$, and $\tilde M\eta^2 j^{-2} \lambda^{-\frac{1}{4n}}<\eta j^{-2}<\frac12 j^{-2}$ by \eqref{fitg0}. We also have
$\frac12\eta^2 j^{-2}(1-\lambda^{-\frac{1}{4n}})<\frac12 j^{-2}$, thus (1) 
is satisfied. For (2), note that $\frac12\eta^2 j^{-2}(1-\lambda^{-\frac{1}{4n}})/(\tilde M\eta^2 j^{-2} \lambda^{-\frac{1}{4n}})=\alpha$ by 
\eqref{fitg1}, thus we have (2). (3) is satisfied due to \eqref{fitg4}. 
(4) and (5) are satisfied respectively due to (8) and (9) of Lemma \ref{lemafi}. 
Thus the assumptions of Lemma \ref{itglem} are all satisfied for $Y_l$, and
we have
\begin{equation}
\label{fitg9}
\frac{\|\partial \mathcal E\|(S_l)}{\omega_{n}(\tilde M 
\eta^2 j^{-2}\lambda^{-\frac{1}{4n}})^{n}} \geq 
\mathcal H^0(Y_l)-\zeta\geq \lambda^{-\frac14} \mathcal H^0(Y_l)
\end{equation}
where we used \eqref{fitg2}. The similar formula holds for $y\in Y_0$, and
\eqref{fitg8} and \eqref{fitg9} show \eqref{afitg7}. Finally we let $\gamma$ be
re-defined as $\min\{\gamma_1,\gamma_2,\gamma_3\}$ if necessary.  
\hfill{$\Box$}

\begin{thm}
(\cite[4.24]{Brakke})
Suppose that $\{\mathcal E_j\}_{j=1}^{\infty}\subset 
\mathcal{OP}_{\Omega}^N$ and $\{\e_j\}_{j=1}^{\infty}\subset (0,1)$ satisfy
\begin{itemize}
\item[(1)] $\lim_{j\rightarrow\infty} j^4 \e_j=0$,
\item[(2)] $\sup_j \|\partial\mathcal E_j\|(\Omega)<\infty$,
\item[(3)] $\liminf_{j\rightarrow\infty} \int_{\mathbb R^{n+1}} \frac{|\Phi_{\e_j}\ast
\delta(\partial\mathcal E_j)|^2\Omega}{\Phi_{\e_j}\ast\|\partial\mathcal E_j\|+
\e_j \Omega^{-1}}\, dx<\infty$,
\item[(4)] $\lim_{j\rightarrow\infty} j^{2(n+1)}\Delta_j\|\partial\mathcal E_j\|(\Omega)=0$.
\end{itemize}
Then there exists a converging subsequence $\{\partial\mathcal E_{j_l}\}_{l=1}^{\infty}$ whose limit satisfies $V\in {\bf IV}_{n}(\mathbb R^{n+1})$. 
\label{mainit}
\end{thm}
{\it Proof}.
We may choose a 
subsequence $\{j_l\}_{l=1}^{\infty}$ such that the quantities in (2) and (3) are uniformly 
bounded by $M$ and the sequence $\{\partial\mathcal E_{j_l}\}_{l=1}^{\infty}$
converges to $V\in {\bf RV}_{n}(\mathbb R^{n+1})$ by Theorem \ref{rect1}.
Without loss of generality, it is 
enough to prove that $V$ is integral in $U_1$. 
For each pair of positive integers $j$ and $q$, let 
$A_{j,q}$ be a set consisting of all $x\in B_1$ such that
\begin{equation}
\|\delta(\Phi_{\e_j}\ast \partial\mathcal E_j)\|(B_r(x))\leq q \|\Phi_{\e_j}
\ast \partial \mathcal E_j\|(B_r(x))
\label{xitg1}
\end{equation}
for all $r\in (j^{-2},1)$.
For any $x\in B_1 \setminus A_{j,q}$, we have
\begin{equation}
\label{xitg3}
\|\delta(\Phi_{\e_j}\ast \partial\mathcal E_j)\|(B_r(x))> q \|\Phi_{\e_j}
\ast \partial \mathcal E_j\|(B_r(x))
\end{equation}
for some $r\in (j^{-2},1)$.
Since $\Phi_{\e_j}\ast
\chi_{B_r(x)}\geq \frac14 \chi_{B_r(x)}$ as long as $\e_j\ll r^2$, we have
\begin{equation}
\label{xitg5}
\|\delta(\Phi_{\e_j}\ast \partial\mathcal E_j)\|(B_r(x))> \frac{q}{4}\|\partial \mathcal E_j\|(B_r(x)). 
\end{equation}
For sufficiently large $j$, (1) and $r\in (j^{-2},1)$ guarantee
that $\e_j\ll r^2$. 
Applying the Besicovitch covering theorem to a collection of such balls 
covering $B_1\setminus A_{j,q}$, there exists a
family $\mathcal C$ of disjoint balls such that
\begin{equation}
\|\partial \mathcal E_j\| (B_1 \setminus A_{j,q})
\leq{\bf B}_{n+1} \sum_{B_r(x)\in \mathcal C} \|\partial\mathcal E_j\|
(B_{r}(x)).
\label{xitg1.1}
\end{equation}
Thus, with \eqref{xitg1.1} and \eqref{xitg5}, we obtain
\begin{equation}
\label{xitg13}
\|\partial \mathcal E_j\| (B_1 \setminus A_{j,q})
\leq \frac{4{\bf B}_{n+1}}{q}\|\delta(\Phi_{\e_j}\ast\partial\mathcal E_j)\|(B_2).
\end{equation}
By the Cauchy-Schwarz inequality and \eqref{musm4}, 
\begin{equation}
\label{xitg14}
\|\delta(\Phi_{\e_j}\ast\partial\mathcal E_j)\|(B_2) 
\leq \big(\int_{B_2} \frac{|\Phi_{\e_j}\ast\delta(\partial\mathcal E_j)|^2}{
\Phi_{\e_j}\ast\|\partial\mathcal E_j\|+\e_j\Omega^{-1}}\, dx\big)^{\frac12}
\big(\int_{B_2} \Phi_{\e_j}\ast\|\partial\mathcal E_j\|+\e_j\Omega^{-1}\, dx\big)^{\frac12}.
\end{equation}
The right-hand side of
\eqref{xitg14} for $j_l$ is bounded by $(\min_{B_3}\Omega)^{-1} M^{\frac12}
(M^{\frac12}+2^{n+1}\omega_{n+1})$ for all $l$. 
Then \eqref{xitg13} and \eqref{xitg14} 
show 
\begin{equation}
\label{xitg15}
 \|\partial\mathcal E_{j_l}\|(B_1\setminus A_{j_l,q})\leq\frac{c(n,\Omega,M)}{q}
\end{equation}
for all $l, q\in \mathbb N$. 
Now for each $q\in \mathbb N$, set 
\begin{equation}
A_q:=\{x\in B_1 : \mbox{there exist $x_l\in A_{j_l,q}$ for infinitely many $l$ with
$x_l\rightarrow x$}\}
\label{xitg16}
\end{equation}
and define
\begin{equation}
A:=\cup_{q=1}^{\infty} A_q.
\label{xitg17}
\end{equation}
Then we have 
\begin{equation}
\|V\|(U_1\setminus A)=0.
\label{xitg18}
\end{equation}
This can be seen as follows. Take arbitrary compact set 
$K\subset U_1\setminus A$. For any $q\in \mathbb N$ we have
$K\subset U_1\setminus A_q$ by \eqref{xitg17}. For each point $x\in K$, by 
\eqref{xitg16}, there exists a neighborhood of $x$ which does not intersect with 
$A_{j_l,q}$ for all sufficiently large $l$. Due to the compactness of $K$, there exist $l_0\in \mathbb N$ and an 
open set $O_q\subset U_1$ 
such that $K\subset O_q$ and $O_q\cap A_{j_l,q}=\emptyset$
for all $l\geq l_0$. Let $\phi_q\in C_c(O_q;\mathbb R^+)$ be such that $0\leq \phi_q
\leq 1$ and $\phi_q=1$ on $K$. Then 
\begin{equation}
\label{xitg19}
\begin{split}
\|V\|(K) & \leq \|V\|(\phi_q)=\lim_{l\rightarrow\infty}\|\partial \mathcal E_{j_l}\|
(\phi_q)=\lim_{l\rightarrow\infty}\|\partial\mathcal E_{j_l}\|\lfloor_{B_1\setminus
A_{j_l,q}}(\phi_q) \\ & \leq \liminf_{l\rightarrow\infty}\|\partial\mathcal E_{j_l}\|(
B_1\setminus A_{j_l,q})\leq \frac{c(n,\Omega,M)}{q}
\end{split}
\end{equation}
where we used \eqref{xitg15}. Since $q$ is arbitrary, \eqref{xitg19} gives
$\|V\|(K)=0$, proving \eqref{xitg18}. 

Let $A^*$ be a set of points in $U_1$ such that the approximate tangent 
space of $V$ exists, i.e., 
\begin{equation}
\label{xitg20}
\begin{split}
A^*:= & \{x\in U_1 : \theta^{n}(\|V\|,x)\in (0,\infty), {\rm Tan}^{n}(\|V\|,x)
\in {\bf G}(n+1,n), \\ & \lim_{r\rightarrow 0+}( f_{(r)}\circ\, \tau_{(-x)})_{\sharp} V
=\theta^{n}(\|V\|,x) |{\rm Tan}^{n}(\|V\|,x)|\}.
\end{split}
\end{equation}
Here, $f_{(r)}(y):=r^{-1} y$ and $\tau_{(-x)}(y)=y-x$ for 
$y\in \mathbb R^{n+1}$. Since $V\in {\bf RV}_{n}(\mathbb R^{n+1})$, we have
$\|V\|(U_1\setminus A^*)=0$. Thus, for $\|V\|$ a.e.~$x\in U_1$, we have
$x\in A^*\cap A$. In the following, we fix $x$ and prove that $\theta^{n}(\|V\|,x)\in \mathbb N$ 
for such $x$, which proves that $V\in {\bf IV}_{n}(\mathbb R^{n+1})$. 
For simplicity, we write
\begin{equation}
\label{xitg22.5}
d:=\theta^{n}(\|V\|,x),\, T:={\rm Tan}^{n}(\|V\|,x).
\end{equation}
By an appropriate change of variables, we may assume that $x=0$ and 
$T=\{x_{n+1}=0\}$, with the understanding that all the relevant quantities are
re-defined accordingly with no loss of generality.   
By \eqref{xitg17}, there exists $q\in \mathbb N$ such that
$x=0\in A^*\cap A_q$, hence there exists a further subsequence of $\{j_l\}_{l=1}^{\infty}$ (denoted by the same index) such that $x_{j_l}\in A_{j_l,q}$ with
$\lim_{l\rightarrow\infty}x_{j_l}=0$. Set $r_l:=l^{-1}$, and choose a further subsequence so that
\begin{equation}
\label{xitg21}
\lim_{l\rightarrow\infty} (f_{(r_l)})_{\sharp} \partial \mathcal E_{j_l}
=\lim_{l\rightarrow\infty} (f_{(r_l)})_{\sharp} (\Phi_{\e_{j_l}}\ast\partial \mathcal E_{j_l})
=d|T|,
\end{equation}
\begin{equation}
\lim_{l\rightarrow\infty} \frac{x_{j_l}}{r_l}=0
\label{xitg22}
\end{equation}
and
\begin{equation}
\lim_{l\rightarrow\infty}\frac{ j_l^{-1}}{ r_l}=\lim_{l\rightarrow\infty}
\frac{l}{j_l}=0.
\label{xitg23f}
\end{equation}
We define
\begin{equation}
V_{j_l}:=(f_{(r_l)})_{\sharp}\partial\mathcal E_{j_l}, \hspace{1cm} 
\tilde V_{j_l}:=(f_{(r_l)})_{\sharp}(\Phi_{\e_{j_l}}\ast\partial\mathcal E_{j_l})
\label{xitg24}
\end{equation}
for simplicity in the following. 

Suppose that $\nu$ is the smallest positive integer strictly 
greater than $d$, i.e., 
\begin{equation}
\nu\in \mathbb N\,\,\mbox{ and }\,\, \nu\in (d,d+1].
\label{eravi}
\end{equation} 
Choose $\lambda\in (1,2)$ such that 
\begin{equation}
\label{xitg23}
\lambda^{n+1} d<\nu.
\end{equation}
Corresponding to such $\lambda$ and $\nu$, we choose $\gamma,\eta\in (0,1)$,
$\tilde M\in (1,\infty)$ 
and $j_0\in \mathbb N$ using Lemma \ref{lemafi}. 
We use Lemma \ref{lemafi} 
with $R=r_l$ in the following. To do so, as a first step, 
we prove that the first variations of $\tilde V_{j_l}$ converge to 0, i.e., 
\begin{equation}
\lim_{l\rightarrow\infty}\|\delta \tilde V_{j_l}\|(B_s)=\lim_{l\rightarrow\infty}
r_l^{1-n}\|\delta (\Phi_{\e_{j_l}}\ast\partial\mathcal E_{j_l})
\|(B_{sr_l})=0
\label{xitg25}
\end{equation}
for all $s>0$. To see this, note that we have $x_{j_l}\in A_{j_l,q}$, so that
\begin{equation}
\|\delta(\Phi_{\e_{j_l}}\ast\partial\mathcal E_{j_l})\|(B_{sr_l}(x_{j_l}))
\leq q\|\Phi_{\e_{j_l}}\ast\partial\mathcal E_{j_l}\|(B_{sr_l}(x_{j_l}))
\label{xitg26}
\end{equation}
by \eqref{xitg1}, where we note that $sr_l\in (j_l^{-2},1)$ for all sufficiently large
$l$ due to \eqref{xitg23f}. 
One can check that \eqref{xitg26} is equivalent to
\begin{equation}
\|\delta \tilde V_{j_l}\|(B_s(r_l^{-1}x_{j_l}))\leq  r_l q \|\tilde V_{j_l}\|(B_s(r_l^{-1}x_{j_l})).
\label{xitg27}
\end{equation}
By \eqref{xitg22}, $r_l^{-1} x_{j_l}\rightarrow 0$, and by \eqref{xitg21}, 
$\|\tilde V_{j_l}\|\rightarrow \| d|T|\|$. Since $r_l=l^{-1}$, by letting 
$l\rightarrow\infty$, \eqref{xitg27} proves
\eqref{xitg25}. We also need 
\begin{equation}
\lim_{l\rightarrow\infty} \int_{{\bf G}_{n}(B_s)} \|S-T\|\, d\tilde V_{j_l}
=\lim_{l\rightarrow\infty} r_l^{-n} \int_{{\bf G}_{n}(B_{sr_l})}
\|S-T\|\, d(\Phi_{\e_{j_l}}\ast\partial\mathcal E_{j_l})=0
\label{xitg28}
\end{equation}
and
\begin{equation}
\lim_{l\rightarrow\infty} \int_{{\bf G}_{n}(B_s)} \|S-T\|\, d V_{j_l}
=\lim_{l\rightarrow\infty} r_l^{-n} \int_{{\bf G}_{n}(B_{sr_l})}
\|S-T\|\, d(\partial\mathcal E_{j_l})=0
\label{xitg28.5}
\end{equation}
for all $s>0$, but these follow directly from the varifold convergence of
\eqref{xitg21} to $d|T|$. 

For each $l\in \mathbb N$ define
\begin{equation}
\begin{split}
G_l:=&\big\{x\in B_{(\lambda -1)r_l} : r_l \|\delta(\Phi_{\e_{j_l}}\ast\partial
\mathcal E_{j_l})\|(B_s(x))\leq \gamma \|\Phi_{\e_{j_l}}\ast\partial \mathcal E_{j_l}
\|(B_s(x))\mbox{ and } \\
& \int_{{\bf G}_{n}(B_s(x))} \|S-T\|\, d(\Phi_{\e_{j_l}}\ast\partial\mathcal E_{j_l})
\leq \gamma \|\Phi_{\e_{j_l}}\ast\partial\mathcal E_{j_l}\|(B_s(x)) \\
& \hspace{5cm}\mbox{ for all }s\in (\eta^2 j_l^{-2},r_l)\big\}.
\end{split}
\label{xitg29}
\end{equation}
By exactly the same line of argument as in \eqref{xitg1}-\eqref{xitg13}, we have
\begin{equation}
\begin{split}
\|&\partial\mathcal E_{j_l}\|(B_{(\lambda-1)r_l}
\setminus G_l)\\&
\leq 4{\bf B}_{n+1}\gamma^{-1} \Big(r_l \|\delta(\Phi_{\e_{j_l}}\ast\partial
\mathcal E_{j_l})\|(B_{\lambda r_l})+\int_{{\bf G}_{n}(B_{\lambda r_l})}
\|S-T\|\, d(\Phi_{\e_{j_l}}\ast\partial\mathcal E_{j_l})\Big).
\end{split}
\label{xitg30}
\end{equation}
Then, \eqref{xitg25}, \eqref{xitg28} and \eqref{xitg30} show that
\begin{equation}
\label{xitg31}
\lim_{l\rightarrow\infty} r_l^{-n}\|\partial\mathcal E_{j_l}\|(B_{(\lambda-1)r_l}
\setminus G_l)=0.
\end{equation}
Define
\begin{equation}
G_l^*:=\{x\in G_l : \theta^{n}(\|\partial\mathcal E_{j_l}\|,x)=1\}.
\label{xitg31a}
\end{equation}
Since $\partial\mathcal E_{j_l}$ is a unit density varifold, 
\begin{equation}
\|\partial\mathcal E_{j_l}\|(G_l\setminus G_l^*)=0.
\label{xitg31b}
\end{equation}
We next define, as in Lemma \ref{lemafi} (a)-(c), 
\begin{equation}
\tilde R_{1,l}:=\eta^2 j_l^{-2}\lambda^{-\frac{1}{4n}}, 
\tilde R_{2,l}:=\tilde M\eta^2 j_l^{-2}\lambda^{-\frac{1}{4n}}, 
\rho_l:=\frac12 \eta^2 j_l^{-2}(1-\lambda^{-\frac{1}{4n}}).
\label{xitg31d}
\end{equation}
We wish to apply Lemma \ref{lemafi} and define $G^{**}_l \subset G_l^*$ as 
follows. For $x\in G_l^*$, take any arbitrary finite set $Y'=\{y_1,\ldots, y_m\}
\subset G_l^*$ with $y_1=x$, 
$T(x-y_i)=0$ for $i\in \{2,\ldots, m\}$ and ${\rm diam}\, Y'<j_l^{-2}$.
We do not exclude the 
possibility that $Y'=\{y_1\}=\{x\}$.
Define 
\begin{equation}
\label{xitg31z}
E^*_{i,l}(r,Y'):=\{ z\in \mathbb R^{n+1} : |T(z-x)|\leq r,\, {\rm dist}\,(T^{\perp}(Y')
,T^{\perp}(z))\leq  (1+\tilde R_{i,l}^{-1}
r)\rho_l\}
\end{equation}
for $i=1,2$. We define $G_l^{**}$ as a set of point $x\in G_l^*$ such that, 
for arbitrary such $Y'$ described above and for all $r\in (0,j_l^{-2})$ and $i=1,2$, 
we have
\begin{equation}
\label{xitg32}
\begin{split}
& \int_{{\bf G}_{n}(E^*_{i,l}(r,Y'))}
\|S-T\|\, d(\partial\mathcal E_{j_l})\leq \gamma\|\partial\mathcal E_{j_l}\|(E^*_{i,l}(r,Y'))
\mbox{ and } \\ & \Delta_{j_l}\|\partial\mathcal E_{j_l}\|(E^*_{i,l}(r,Y'))\geq -\gamma 
\|\partial\mathcal E_{j_l}\|(E^*_{i,l}(r,Y')).
\end{split}
\end{equation}
We wish to show that $\|\partial\mathcal E_{j_l}\|(G_l^*
\setminus G_l^{**})$, which is a missed mass we cannot apply Lemma
\ref{lemafi}, is small. 
Whenever $x\in G_l^*\setminus G_l^{**}$, there exist a finite set 
$Y'_x=\{y_1,\ldots,y_m\}\subset G_l^*$ with
\begin{equation}
y_1=x,\,
T(x-y_i)=0
\mbox{ for }i\in \{2,\ldots,m\}, {\rm diam}\, Y_x'<j_l^{-2},
\label{xitg32a0}
\end{equation}
and $r_x \in (0,j_l^{-2})$ such that
\begin{equation}
\begin{split}
& \int_{{\bf G}_{n}(E^*_{i,l}(r_x,Y'_x))}
\|S-T\|\, d(\partial\mathcal E_{j_l})> \gamma\|\partial\mathcal E_{j_l}\|(E^*_{i,l}
(r_x,Y'_x))
\mbox{ for }i=1\mbox{ or }i=2\mbox{ or } \\ &
\Delta_{j_l}\|\partial\mathcal E_{j_l}\|(E^*_{i,l}(r_x,Y'_x))< -\gamma 
\|\partial\mathcal E_{j_l}\|(E^*_{i,l}(r_x,Y'_x)) \mbox{ for }i=1\mbox{ or }i=2.
\end{split}
\label{xitg32a}
\end{equation}
We separate $G_l^*\setminus G_l^{**}$ into four sets depending on the
conditions in \eqref{xitg32a},
\begin{equation}
\label{xitg32b}
W_{i,l}:=\{ x\in G_l^*\setminus G_l^{**} :  \int_{{\bf G}_{n}(E^*_{i,l}(r_x,Y'_x))}
\|S-T\|\, d(\partial\mathcal E_{j_l})> \gamma\|\partial\mathcal E_{j_l}\|(E^*_{i,l}
(r_x,Y'_x))\}
\end{equation}
for $i=1,2$ and 
\begin{equation}
\label{xitg32c}
\tilde W_{i,l}:=\{ x\in G_l^*\setminus G_l^{**} : 
\Delta_{j_l}\|\partial\mathcal E_{j_l}\|(E^*_{i,l}(r_x,Y'_x))< -\gamma 
\|\partial\mathcal E_{j_l}\|(E^*_{i,l}(r_x,Y'_x))\}
\end{equation}
for $i=1,2$ so that
\begin{equation}
\label{xitg32d}
G_l^*\setminus G_l^{**}=\cup_{i=1}^2(W_{i,l}\cup \tilde W_{i,l}).
\end{equation}
Typically, we would use the Besicovitch covering theorem to estimate
the missed mass, but here, the elements of 
covering of $G_l^*\setminus G_l^{**}$
are $E_{i,l}^*(r_x,Y'_x)$, which are not closed balls. Thus, direct use of the
Besicovitch is not possible. On the other hand, 
note that at any point in $W_{i,l}$, the covering 
$E_{i,l}^*(r_x,Y'_x)$ has always ``height''  bigger than $\rho_l$ 
in $T^{\perp}$ direction, and $\rho_l$ is $O(j_l^{-2})$. 
We take advantage of this property in the following. 
We estimate $\|\partial\mathcal E_{j_l}\|( W_{i,l})$ for $i=1,2$ first.
We choose a finite 
set of points $\{w_{l,k}\}_{k=1}^{K_l}$ in $B_{(\lambda-1)r_l}$ so that
\begin{equation}
\label{xitg32fa}
B_{(\lambda-1)r_l}\subset \cup_{k=1}^{K_l} B_{j_l^{-2}}
(w_{l,k}) 
\end{equation}
and the number of intersection $\{k' : B_{4 j_l^{-2}}(w_{l,k'})\cap
B_{4 j_l^{-2}}(w_{l,k})\neq \emptyset\}$ for each $k$ is less than a constant $c(n)$ 
depending only on
$n$.  Such a set of points 
can be lattice points with width $j_l^{-2}$ in $B_{(\lambda-1)r_l}$,
for example. We then have
\begin{equation}
\sum_{k=1}^{K_l}\int_{{\bf G}_{n}(B_{4j_l^{-2}}(w_{l,k}))} \|S-T\|\, d(\partial\mathcal E_{j_l})(x,S)
\leq c(n)\int_{{\bf G}_{n}(B_{\lambda r_l})} \|S-T\|\, d(\partial\mathcal E_{j_l})(x,S).
\label{xitg32f}
\end{equation}
If we set for $k\in \{1,\ldots,K_l\}$
\begin{equation}
W_{i,l,k}:=W_{i,l}\cap B_{j_l^{-2}}(w_{l,k}),
\label{xitg32g}
\end{equation}
by \eqref{xitg32fa}, we have
\begin{equation}
\cup_{k=1}^{K_l} W_{i,l,k}=W_{i,l}.
\label{xitg32h}
\end{equation}
We next separate each $W_{i,l,k}$ into a stacked regions of width $\rho_l$ 
in $T^{\perp}$ direction. Define for $m\in \mathbb Z$ with $|m|< j_l^{-2}\rho_l^{-1}+1$ 
\begin{equation}
W_{i,l,k,m}:=W_{i,l,k}\cap \{x\in\mathbb R^{n+1} : m\rho_l<T^{\perp}(x-w_{l,k})\leq
(m+1)\rho_l\}.
\label{xitg32i}
\end{equation}
Since $W_{i,l,k}\subset B_{j_l^{-2}}(w_{l,k})$, we have
\begin{equation}
\label{xitg32j}
W_{i,l,k}=\cup_{|m|<j_l^{-2}\rho_l^{-1}+1} W_{i,l,k,m}
\end{equation}
and it is important to note that $j_l^{-2}\rho_l^{-1}+1$ is a constant
depending only on $\eta$ and $\lambda$, so ultimately only on 
$n,\nu$ and $\lambda$. 
For each $x\in W_{i,l,k,m}$, there exist $Y_x\subset G_l^*$ and $r_x\in
(0,j_l^{-2})$ with the inequality of \eqref{xitg32b}. Define
\begin{equation}
\mathcal C_{i,l,k,m}:=\{B_{r_x}^{n}(T(x))\subset\mathbb R^{n} : x\in W_{i,l,k,m}\}
\label{xitg32k}
\end{equation}
which is a covering of $T(W_{i,l,k,m})$. Observe that, if there is a subfamily 
$\hat{\mathcal C}_{i,l,k,m}\subset \mathcal C_{i,l,k,m}$ such that
$T(W_{i,l,k,m})\subset \cup_{C\in \hat{\mathcal C}_{i,l,k,m}} C$, we have 
\begin{equation}
W_{i,l,k,m}\subset \cup_{B_{r_x}^{n}(T(x))\in \hat{\mathcal C}_{i,l,k,m}}
\{y : |T(x-y)|\leq r_x, |T^{\perp}(x-y)|\leq \rho_l\}.
\label{xitg32l}
\end{equation}
This is because, for any $x'\in W_{i,l,k,m}$, we have some $B_{r_x}^{n}(T(x))
\in \hat{\mathcal C}_{i,l,k,m}$ with $T(x')\in B_{r_x}^{n}(T(x))$. Since $x',x\in 
W_{i,l,k,m}$, $|T^{\perp}(x'-x)|<\rho_l$, so $x'\in \{y : |T(x-y)|\leq r_x, |T^{\perp}
(x-y)|\leq \rho_l\}$, which proves \eqref{xitg32l}.
We apply the Besicovitch covering theorem to $\mathcal C_{i,l,k,m}$ and 
obtain a set of subfamilies $\mathcal C^{(1)}_{i,l,k,m},\ldots,\mathcal C^{(L_{i,l,k,m})}_{i,l,k,m}\subset \mathcal C_{i,l,k,m}$ such that 
\begin{equation}
\label{xitg32nn}
L_{i,l,k,m}\leq {\bf B}_{n},
\end{equation}
each $\mathcal C^{(h)}$ ($h=1,\ldots,L_{i,l,k,m}$) consists of disjoint sets and 
$T(W_{i,l,k,m})\subset \cup_{h=1}^{L_{i,l,k,m}} \cup_{C\in \mathcal C^{(h)}_{i,l,k,m}} C$. 
Then \eqref{xitg32l} shows that we have
\begin{equation}
\label{xitg32m}
W_{i,l,k,m}\subset 
\cup_{h=1}^{L_{i,l,k,m}} \cup_{B_{r_x}^{n}(x)\in \mathcal C^{(h)}_{i,l,k,m}} 
\{ y : |T(x-y)|\leq r_x, |T^{\perp}(x-y)|\leq \rho_l\}.
\end{equation}
For $x\in W_{i,l,k,m}$, 
\begin{equation}
\label{xitg32ns}
\{y : |T(x-y)|\leq r_x, |T^{\perp}(x-y)|\leq \rho_l\}
\subset E^*_{i,l}(r_x,Y'_x).
\end{equation}
We note that if $B_{r_x}^{n}(x)\cap B_{r_{x'}}^{n}
(x')=\emptyset$, then $E^*_{i,l}(r_x,Y'_x)\cap E^*_{i,l}(r_{x'},Y'_{x'})=\emptyset$ 
since their projections to $T$ is $B_{r_x}^{n}(x)\cap B_{r_{x'}}^{n}
(x')$. Also we note that 
\begin{equation}
\label{xitg32n}
E^*_{i,l}(r_x,Y'_x)\subset B_{4j_l^{-2}}(w_{l,k})
\end{equation}
since
$x\in B_{j_l^{-2}}(w_{l,k})$, $Y'_x\in T^{\perp}(B_{j_l^{-2}}(x))$ (by \eqref{xitg32a0}), $r_x\in (0,j_l^{-2})$, $(1+\tilde R_{i,l}^{-1} r_x)\rho_l\leq \rho_l
+\frac{r_x}{2}<j_l^{-2}$ (by \eqref{xitg31d}
and \eqref{xitg31z}). 
We have by \eqref{xitg32m}, \eqref{xitg32ns}, \eqref{xitg32b}, \eqref{xitg32nn}
and \eqref{xitg32n} that
\begin{equation}
\label{xitg32p}
\begin{split}
\|\partial\mathcal E_{j_l}\|(W_{i,l,k,m})&\leq \sum_{h=1}^{L_{i,l,k,m}}
\sum_{B_{r_x}^{n}(x)\in \mathcal C^{(h)}_{i,l,k,m}}
\|\partial\mathcal E_{j_l}\|(E_{i,l}^* (r_x,Y_x')) \\
&\leq \sum_{h=1}^{L_{i,l,k,m}}
\sum_{B_{r_x}^{n}(x)\in \mathcal C^{(h)}_{i,l,k,m}}
\gamma^{-1} \int_{{\bf G}_{n}(E_{i,l}^*(r_x,Y'_x))} \|S-T\|\,d(\partial
\mathcal E_{j_l})(x,S) \\
&\leq \gamma^{-1}{\bf B}_{n} \int_{{\bf G}_{n}(B_{4j_l^{-2}}(w_{l,k}))}
\|S-T\|\, d(\partial\mathcal E_{j_l})(x,S).
\end{split}
\end{equation}
Now summing \eqref{xitg32p} over $|m|<j_l^{-2}\rho_l^{-1}+1$ (note \eqref{xitg32j} and the following remark), we have
\begin{equation}
\|\partial\mathcal E_{j_l}\|(W_{i,l,k})\leq \gamma^{-1}c(n,\nu,\lambda)
\int_{{\bf G}_{n}(B_{4j_l^{-2}}(w_{l,k}))}
\|S-T\|\, d(\partial\mathcal E_{j_l})(x,S).
\label{xitg32q}
\end{equation}
Summing \eqref{xitg32q} over $k=1,\ldots,K_l$ and by \eqref{xitg32h}
and \eqref{xitg32f}, we obtain
\begin{equation}
\|\partial\mathcal E_{j_l}\|(W_{i,l})\leq \gamma^{-1}c(n,\nu,\lambda)
\int_{{\bf G}_{n}(B_{\lambda r_l})}\|S-T\|\, d(\partial\mathcal E_{j_l})(x,S).
\label{xitg32s}
\end{equation}
By \eqref{xitg28.5} and \eqref{xitg32s}, we obtain
\begin{equation}
\lim_{l\rightarrow\infty}r_l^{-n}\|\partial\mathcal E_{j_l}\|(W_{i,l})=0.
\label{gs7}
\end{equation}
Next we estimate $\|\partial\mathcal E_{j_l}\|(\tilde W_{i,l})$ for $i=1,2$.
The argument is identical up to the second line of \eqref{xitg32p}
except that we use the covering satisfying the inequality of \eqref{xitg32c}
in place of \eqref{xitg32b}. By using the same notation, we obtain
\begin{equation}
\label{gs1}
\|\partial\mathcal E_{j_l}\|(\tilde W_{i,l,k,m})\leq - \sum_{h=1}^{L_{i,l,k,m}}
\sum_{B_{r_x}^{n}(x)\in \mathcal C_{i,l,k,m}^{(h)}} \gamma^{-1}
\Delta_{j_l}\|\partial\mathcal E_{j_l}\|(E_{i,l}^*(r_x,Y'_x)).
\end{equation}
Recall that $\{E_{i,l}^*(r_x,Y'_x)\}_{B_{r_x}^{n}(x)\in \mathcal C^{(h)}_{i,l,k,m}}$
is disjoint and we have \eqref{xitg32n}. Since $\mathcal L^{n+1} (B_{4j_l^{-2}}(w_{l,k}))
<j_l^{-1}$ for large $l$, Lemma \ref{escov2} shows
\begin{equation}
\Delta_{j_l}\|\partial\mathcal E_{j_l}\|(B_{4j_l^{-2}}(w_{l,k}))\leq 
\sum_{B_{r_x}^{n}(x)\in \mathcal C^{(h)}_{i,l,k,m}}\Delta_{j_l}\|\partial\mathcal
E_{j_l}\|(E_{i,l}^*(r_x,Y'_x))
\label{gs2}
\end{equation}
for each $h$.  Hence \eqref{gs1}, \eqref{gs2} and \eqref{xitg32nn} show
\begin{equation}
\label{gs3}
\|\partial\mathcal E_{j_l}\|(\tilde W_{i,l,k,m})\leq -{\bf B}_{n} \gamma^{-1}
\Delta_{j_l}\|\partial\mathcal E_{j_l}\|(B_{4j_l^{-2}}(w_{l,k}))
\end{equation}
and summation over $|m|<j_l^{-2}\rho_l^{-1}+1$ gives
\begin{equation}
\label{gs4}
\|\partial\mathcal E_{j_l}\|(\tilde W_{i,l,k})\leq -\gamma^{-1}c(n,\nu,\lambda)
\Delta_{j_l}\|\partial\mathcal E_{j_l}\|(B_{4j_l^{-2}}(w_{l,k})).
\end{equation}
By Lemma \ref{escov}, we have
\begin{equation}
\label{gs4a}
\begin{split}
- \Delta_{j_l}\|\partial\mathcal E_{j_l}\|(B_{4j_l^{-2}}(w_{l,k}))
\leq& -(\max_{B_{4j_l^{-2}}(w_{l,k})} \Omega)^{-1} 
\Delta_{j_l}\|\partial\mathcal E_{j_l}\|(\Omega) \\ &+(1-e^{-4\Cr{c_1}j_l^{-2}})
\|\partial\mathcal E_{j_l}\|(B_{4j_l^{-2}}(w_{l,k})).
\end{split}
\end{equation}
Noticing that $K_l$ in \eqref{xitg32fa} satisfies $K_l\leq c(n)(r_l j_l^2)^{n+1}$,
summation over $k$ of \eqref{gs4} combined with \eqref{gs4a} gives
\begin{equation}
\label{gs5}
\begin{split}
\|\partial\mathcal E_{j_l}\|(\tilde W_{i,l})\leq &
 \gamma^{-1} c(n,\nu,\lambda,\Omega)
\{ -(r_l j_l^2)^{n+1}\Delta_{j_l}\|\partial
\mathcal E_{j_l}\|(\Omega) \\ &+(1-e^{-4\Cr{c_1}j_l^{-2}})\|\partial\mathcal E_{j_l}\|(B_{\lambda r_l})\}.
\end{split}
\end{equation}
With (4), 
\eqref{xitg21} and \eqref{gs5}, we conclude that
\begin{equation}
\lim_{l\rightarrow\infty} r_l^{-n}\|\partial\mathcal E_{j_l}\|(\tilde W_{i,l})=0.
\label{gs6}
\end{equation}
Now, by \eqref{xitg32d}, \eqref{gs7} and \eqref{gs6} we have 
\begin{equation}
\label{gs8}
\lim_{l\rightarrow\infty}
r_l^{-n}\|\partial\mathcal E_{j_l}\|(G_l^{*}\setminus G_l^{**})=0.
\end{equation}
Combining \eqref{xitg31}, \eqref{xitg31b} and \eqref{gs8}, we have
\begin{equation}
\label{gs9}
\lim_{l\rightarrow\infty} r_l^{-n}\|\partial\mathcal E_{j_l}\|(B_{(\lambda-1)r_l}
\setminus G_l^{**})=0.
\end{equation}
Since $G_l^{**}\subset G_l^*\subset G_l$, $x\in G_l^{**}$ satisfies
\eqref{xitg29}, \eqref{xitg31a} and \eqref{xitg32}. 
Given any $s\in (0,\frac14)$ and $x\in G_l^{**}$, we 
use Lemma \ref{lemafi} with $R=r_l s$ 
for $Y=\{T^{\perp}(x)\}$, a single element case.
For all sufficiently large $j_l$, assumptions of Lemma \ref{lemafi} are all satisfied: 
(1) is fine for large $j_l$, (2) from Theorem~\ref{mainit} (1) for large $j_l$, 
(3) from \eqref{xitg23f} for large $j_l$, (4) from $Y$ having single element and 
$x\in G_l^*$, (5) from ${\rm diam}\, Y=0$, (6) and (7) from \eqref{xitg29}, (8) and
(9) from \eqref{xitg32}. Thus we have \eqref{afitg7}, or
\begin{equation}
\label{gs10}
\lambda\|\Phi_{\e_{j_l}}\ast\partial\mathcal E_{j_l}\|(B_{r_l s}(x))\geq \omega_{n}
(r_l s)^{n}
\end{equation}
for all large $j_l$. \eqref{gs10} implies 
\begin{equation}
G_l^{**}\subset B_{(\lambda-1)r_l}\cap \{x : |T^{\perp}(x)|\leq 3 r_l s\}
\label{gs11}
\end{equation}
for all sufficiently large $j_l$. This is because, if \eqref{gs11} were not true, 
then there would exist a subsequence (denoted by the same index) 
$x_{j_l}\in G_l^{**}$ with $|T^{\perp}(x_{j_l})|>3r_l s$ and we may assume
that $r_l^{-1} x_{j_l}\in B_{\lambda-1}$ converges to $\bar x\in B_{\lambda-1}
\cap \{x : |T^{\perp}(x)|\geq 3s\}$. By \eqref{xitg21}, since
$B_{2s}(\bar x)\cap T=\emptyset$, we have
\begin{equation}
0=\lim_{l\rightarrow\infty}\|(f_{(r_l)})_{\sharp}(\Phi_{\e_{j_l}}\ast
\partial\mathcal E_{j_l})\|(B_{2s}(\bar x))
=\lim_{l\rightarrow\infty} r_l^{-n}\|\Phi_{\e_{j_l}}\ast\partial\mathcal E_{j_l}\|
(B_{2r_l s}(r_l \bar x)).
\label{gs12}
\end{equation}
Since $\lim_{l\rightarrow\infty} r_l^{-1} |r_l \bar x -x_{j_l}|=0$, for sufficiently large
$j_l$, we have $B_{r_l s}(x_{j_l})\subset B_{2r_l s}(r_l\bar x)$.
Hence, continuing from \eqref{gs12}, we have
\begin{equation}
\label{gs13}
\geq \limsup_{l\rightarrow\infty} r_l^{-n}\|\Phi_{\e_{j_l}}\ast\partial\mathcal E_{j_l}\|
(B_{r_l s}(x_{j_l}))\geq \lambda^{-1}\omega_{n}
s^{n}
\end{equation}
where \eqref{gs10} is used in the last step, and we have a contradiction. This 
proves \eqref{gs11}. We next show that, for all sufficiently large $j_l$,
\begin{equation}
\label{gs14}
\mathcal H^0(\{x\in G_l^{**}\, :\, T(x)=a\})\leq \nu-1
\end{equation}
for all $a\in B_{(\lambda-1)r_l}\cap T$. 
For a contradiction, suppose we had some $a_l\in B_{(\lambda-1)r_l}\cap T$ such
that \eqref{gs14} fails. Then there exists $Y_l\subset 
T^{-1}(\{x \in G_l^{**} : T(x)=a_l\})$ with $\mathcal H^0(Y_l)=\nu$. 
We use Lemma \ref{lemafi} to $Y_l$
and $R=r_l$. One can check that
the assumptions are all satisfied just as for the single element case, except for
(5), which was trivial before. This time, on the other hand, due to \eqref{gs11}, 
we have ${\rm diam}\,Y_l \leq \gamma r_l$ by choosing $s=\gamma/6$, 
so (5) is also satisfied. Thus we have
\begin{equation}
\lambda\|\Phi_{\e_{j_l}}\ast\partial\mathcal E_{j_l}\|(\{x : {\rm dist}\,(x,Y_l)
\leq r_l\})\geq \omega_{n} r_l^{n}\nu.
\label{gs15s}
\end{equation}
We may assume after choosing a subsequence
that $r_l^{-1}a_l$ converges to $\bar a\in B_{\lambda-1}\cap T$.
By \eqref{xitg21}, 
\begin{equation}
\label{gs16s}
 \lambda^{n}\omega_{n}d =\lim_{l\rightarrow\infty}\|(f_{(r_l)})_{\sharp}(\Phi_{\e_{j_l}}\ast\partial\mathcal
E_{j_l})\| (B_{\lambda }(\bar a))=
\lim_{l\rightarrow\infty}r_l^{-n}\|\Phi_{\e_{j_l}}\ast\partial\mathcal
E_{j_l}\|(B_{\lambda r_l}(r_l \bar a)).
\end{equation} 
For large $j_l$, by \eqref{gs11} taking $s=(\sqrt{\lambda}-1)/6$, 
$\{ x : {\rm dist}\, (x,Y_l)\leq r_l\}\subset B_{\sqrt{\lambda} r_l}(a_l)\subset B_{\lambda r_l}(r_l \bar a)$. Hence \eqref{gs15s} and \eqref{gs16s} show
$\lambda^{n+1} d\geq \nu$ which is a contradiction to \eqref{xitg23}. This proves
\eqref{gs14}. 
Finally, we note that
\begin{equation}
\label{gs15}
\lim_{l\rightarrow\infty}r_l^{-n}\|T_{\sharp}\partial\mathcal E_{j_l}\|(B_{(\lambda-1)r_l}\setminus G_l^{**})
\leq \lim_{l\rightarrow\infty} r^{-n}_l\|\partial\mathcal E_{j_l}\|(B_{(\lambda-1)r_l}\setminus G_l^{**})=0
\end{equation}
due to \eqref{gs9} while
\begin{equation}
\begin{split}
\|T_{\sharp}\partial\mathcal E_{j_l}\|(G_l^{**}) &=
\int_{B_{(\lambda-1)r_l}\cap T} \mathcal \sum_{\{x\in G_l^{**}\, : \,T(x)=a\}}
\theta^{n}(\|\partial\mathcal E_{j_l}\|,x)\, d\mathcal H^{n}(a) \\ &
\leq \omega_{n}((\lambda-1)r_l)^{n}(\nu-1)
\end{split}
\label{gs16}
\end{equation}
by \eqref{gs14} for all large $j_l$. By \eqref{xitg21}, 
\begin{equation}
\begin{split}
\lim_{l\rightarrow\infty} r_l^{-n}\|T_{\sharp}\partial\mathcal E_{j_l}\|(B_{(\lambda
-1)r_l})&=\lim_{l\rightarrow\infty} \|T_{\sharp}V_{j_l}\|(B_{\lambda-1})
=\|T_{\sharp}d|T|\|(B_{\lambda-1}) \\ &=\omega_{n}(\lambda-1)^{n}d
\end{split}
\label{gs17}
\end{equation}
and \eqref{gs15}-\eqref{gs17} show $d\leq \nu-1$. By \eqref{eravi},
this proves $d=\nu-1$.
\hfill{$\Box$}

\section{Proof of Brakke's inequality}
\label{secBra}
Here, the main objective is to prove the inequality \eqref{sthm2} usually 
referred as Brakke's inequality. We are interested in proving integral form 
instead of differential form as in \cite{Brakke}. The proof is different from 
\cite{Brakke} and we adopt the proof of \cite{Takasao} which we believe is 
more transparent. 
\begin{lemma}
\label{pbi1}
Let $\{\partial
\mathcal E_{j_l}(t)\}_{t\in \mathbb R^+}$ ($l\in\mathbb N$) 
and $\{\mu_t\}_{t\in \mathbb R^+}$ be
as in Proposition \ref{promecone} satisfying \eqref{meconeq}, 
\eqref{mecones} and \eqref{mecomake}. Then we have the following.
\begin{enumerate}
\item[(a)] For a.e.~$t\in \mathbb R^+$, $\mu_t$ is integral, i.e., there exists 
$V_t\in {\bf IV}_{n}(\mathbb R^{n+1})$ such that $\mu_t=\|V_t\|$. 
\item[(b)] For a.e.~$t\in \mathbb R^+$, if a subsequence $\{j'_l\}_{l=1}^{\infty}\subset \{j_l\}_{l=1}^{\infty}$ satisfies 
\begin{equation}
\label{pbi3}
\sup_{l\in \mathbb N} \int_{\mathbb R^{n+1}} \frac{|\Phi_{\e_{j'_l}}\ast\delta(
\partial\mathcal E_{j'_l}(t))|^2\Omega}{\Phi_{\e_{j'_l}}\ast\|\partial\mathcal E_{j'_l}
(t)\|+\e_{j'_l}\Omega^{-1}}\, dx<\infty,
\end{equation}
then we have 
$\lim_{l\rightarrow\infty}\partial\mathcal E_{j'_l}(t)=V_t\in {\bf IV}_{n}(\mathbb R^{n+1})$
as varifolds and $\mu_t=\|V_t\|$. 
\item[(c)] Furthermore, for a.e.~$t\in\mathbb R^+$, $V_t$ has a generalized mean curvature $h(
\cdot,V_t)$ which satisfies
\begin{equation}
\label{l2gmc}
\int_{\mathbb R^{n+1}}  |h(\cdot,V_t)|^2 \phi\, d\|V_t\|\leq \liminf_{l\rightarrow
\infty}\int_{\mathbb R^{n+1}}\frac{|\Phi_{\e_{j_l'}}\ast\delta (\partial\mathcal E_{j_l'}(t))|^2
\phi}{\Phi_{\e_{j_l'}}\ast\|\partial\mathcal E_{j_l'}(t)\|+\e_{j_l'}\Omega^{-1}}\, dx<\infty
\end{equation}
for any $\phi\in\cup_{i\in \mathbb N}\mathcal A_i$.
\end{enumerate}
\end{lemma}
{\it Proof}. Due to \eqref{mecones} and Fatou's Lemma, we have
\begin{equation}
\label{pbi2}
\liminf_{l\rightarrow\infty} \int_{\mathbb R^{n+1}} \frac{|\Phi_{\e_{j_l}}\ast\delta(
\partial\mathcal E_{j_l}(t))|^2\Omega}{\Phi_{\e_{j_l}}\ast\|\partial\mathcal E_{j_l}
(t)\|+\e_{j_l}\Omega^{-1}}\, dx<\infty
\end{equation}
for a.e.~$t\in \mathbb R^+$ 
and for any $T<\infty$, $\sup_{l\in \mathbb N,\,t\in[0,T]}\|\partial\mathcal E_{j_l}(t)\|(\Omega)
<\infty$ due to \eqref{exapp2}.
Suppose we have \eqref{pbi2} and \eqref{mecomake} at $t$. We check that the assumptions of 
Theorem \ref{mainit} are all satisfied for $\{\mathcal E_{j_l}(t)\}_{l=1}^{\infty}$: 
(1) from \eqref{jere}, (2) from above, (3) by \eqref{pbi2} , and (4) from
\eqref{mecomake}. Thus, there exists a further converging 
subsequence of $\{\partial\mathcal E_{j'_l}(t)\}_{l=1}^{\infty}$ 
and a limit $V_t\in {\bf IV}_n(\mathbb R^{n+1})$, where $\{j'_l\}_{l=1}^{\infty}
\subset\{j_l\}_{l=1}^{\infty}$. This convergence is in the sense of varifold, so 
in particular, we have $\lim_{l\rightarrow\infty} \|\partial\mathcal E_{j'_l}(t)\|
=\|V_t\|$. Note that the left-hand side is $\mu_t$ by \eqref{meconeq},
so $\mu_t=\|V_t\|$. This proves (a).
Note that rectifiable (thus integral) 
varifolds are determined by the weight
measure, thus $V_t$ is uniquely determined by $\mu_t$ independent of the 
subsequence $\{j'_l\}_{l=1}^{\infty}$. Let $\{\partial\mathcal E_{j'_l}(t)\}_{l=1}^{\infty}$
be any subsequence with \eqref{pbi3}, then we have already seen that any 
converging further subsequence converges to $V_t$. Since it is unique,
the full sequence $\{\partial \mathcal E_{j'_l}(t)\}_{l=1}^{\infty}$ 
converges to $V_t$. This proves (b). The claim (c) follows from 
Proposition \ref{col}. 
\hfill{$\Box$}
\begin{rem}
Note that we are NOT claiming that $\lim_{l\rightarrow\infty}\partial\mathcal
E_{j_l}(t)=V_t\in {\bf IV}_n(\mathbb R^{n+1})$ for a.e.~$t\in \mathbb R^+$, but only 
the one with uniform bound of \eqref{pbi3}.
\end{rem}
Up to this point, we defined $V_t\in {\bf IV}_{n}(\mathbb R^{n+1})$ for a.e.~$t\in \mathbb R^+$. 
On the complement of such set of time which is $\mathcal L^1$ measure 0,
we still have $\mu_t$. For such $t$, we define an arbitrary varifold with the 
weight measure $\mu_t$. For example, let $T\in {\bf G}(n+1,n)$ be fixed, and
define $V_t(\phi):=\int_{{\bf G}_{n}(\mathbb R^{n+1})} \phi(x,T)\, d\mu_t$
for $\phi\in C_c({\bf G}_{n}(\mathbb R^{n+1}))$. Then we have $\|V_t\|=\mu_t$.
By doing this, we now have a family of varifolds $\{V_t\}_{t\in\mathbb R^+}$
such that $\|V_t\|=\mu_t$ for all $t\in \mathbb R^+$ and $V_t\in {\bf IV}_{n}
(\mathbb R^{n+1})$ for a.e.~$t\in \mathbb R^+$. 
\begin{thm}
\label{brakkeinq}
For all $T>0$, we have
\begin{equation}
\int_0^T\int_{\mathbb R^{n+1}} |h(\cdot,V_t)|^2\Omega\, d\|V_t\|dt<\infty
\label{brng1}
\end{equation}
and for any $\phi\in C_c^{1}(\mathbb R^{n+1}\times \mathbb R^+; \mathbb R^+)$ and
$0\leq t_1<t_2<\infty$, we have
\begin{equation}
\|V_{t}\|(\phi(\cdot,t))\Big|_{t=t_1}^{t_2}
\leq \int_{t_1}^{t_2}\Big(\delta (V_t,\phi(\cdot,t))(h(\cdot,V_t))+\|V_t\|(\frac{\partial\phi}{\partial t}(\cdot,t))\Big)\,
dt.
\label{brng2}
\end{equation}
\end{thm}
{\it Proof}. \eqref{brng1} follows from \eqref{l2gmc}, Fatou's Lemma and
\eqref{mecones}. We prove \eqref{brng2} for time independent $\phi$
first and let $\phi\in C_c^{\infty}(\mathbb R^{n+1};
\mathbb R^+)$ be arbitrary. Since it has a compact support, there exists $c>0$
such that $c\phi(x)<\Omega(x)$ for all $x\in \mathbb R^{n+1}$. 
Due to the linear dependence on $\phi$ of \eqref{brng2}, it suffices 
to prove \eqref{brng2} for $c\phi$ for $C_c^{\infty}$ case, and by suitable
density argument for $C_c^1$ case. Re-writing $c\phi$ as $\phi$, 
we may as well assume that $\phi<\Omega$. Then for all sufficiently large $i\in 
\mathbb N$, we have $\hat\phi:=\phi+i^{-1}\Omega<\Omega$. After fixing $i$,
there exists $m\in \mathbb N$ such that $\hat\phi\in \mathcal A_m$. Fix $0\leq t_1<t_2$ and suppose that 
$l$ is large enough so that $j_l>m$ and $j_l>t_2$. We use \eqref{exapp4} 
with $\hat \phi$. With the notation of \eqref{contime1}, we obtain
\begin{equation}
\label{brng3}
\|\partial\mathcal E_{j_l}(t)\|(\hat\phi)-\|\partial\mathcal E_{j_l}
(t-\Delta t_{j_l})\|(\hat\phi)\leq\Delta t_{j_l}\big( \delta(\partial\mathcal E_{j_l}(t),\hat\phi)(h_{\e_{j_l}}(\cdot,\partial\mathcal E_{j_l}(t))+\e_{j_l}^{\frac18}\big)
\end{equation}
for $t=\Delta t_{j_l},2\Delta t_{j_l},\ldots, j_l 2^{p_{j_l}}\Delta t_{j_l}$. There 
exist $k_1,k_2\in \mathbb N$ such that $t_2\in ((k_2-1)\Delta t_{j_l},k_2
\Delta t_{j_l}]$ and $t_1\in ((k_1-2)\Delta t_{j_l},(k_1-1)\Delta t_{j_l}]$,
where we are assuming that $\Delta t_{j_l}< t_2-t_1$. 
Summing \eqref{brng3} over $t=k_1\Delta t_{j_l},\ldots, k_2\Delta
t_{j_l}$, we obtain
\begin{equation}
\label{brng4}
\|\partial\mathcal E_{j_l}(t)\|(\hat\phi)\Big|_{t=(k_1-1)\Delta t_{j_l}}^{k_2
\Delta t_{j_l}}\leq \sum_{k=k_1}^{k_2}\Delta t_{j_l}\Big(\delta(\partial\mathcal E_{j_l}(k\Delta t_{j_l}), \hat\phi)(h_{\e_{j_l}}(\cdot,\partial\mathcal E_{j_l}(k\Delta t_{j_l})))+\e_{j_l}^{\frac18}\Big).
\end{equation}
Due to the definition of $\hat\phi=\phi+i^{-1}\Omega$, we have
\begin{equation}
\label{brng5}
\|\partial\mathcal E_{j_l}(t)\|(\hat\phi)\Big|_{t=(k_1-1)\Delta t_{j_l}}^{k_2
\Delta t_{j_l}} \geq \|\partial\mathcal E_{j_l}(t_2)\|(\phi)
-\|\partial\mathcal E_{j_l}(t_1)\|(\phi)
-i^{-1} \|\partial\mathcal E_{j_l}(t_1)\|(\Omega).
\end{equation}
As $l\rightarrow\infty$, with \eqref{exapp2}, we obtain
\begin{equation}
\label{brng6}
\limsup_{l\rightarrow\infty} \|\partial\mathcal E_{j_l}(t)\|(\hat\phi)\Big|_{t=(k_1-1)\Delta t_{j_l}}^{k_2
\Delta t_{j_l}} \geq \|V_{t}\|(\phi)\Big|_{t=t_1}^{t_2}
-i^{-1} \|\partial\mathcal E_{0}\|(\Omega) \exp(\frac{\Cr{c_1}^2 t_1}{2}).
\end{equation}
For the right-hand side of \eqref{brng4}, by \eqref{Bdefs} and writing $h_{\e_{j_l}}=
h_{\e_{j_l}}(\cdot,\partial\mathcal E_{j_l}(t))$ and $\partial\mathcal E_{j_l}=
\partial\mathcal E_{j_l}(t)$,
\begin{equation}
\label{brng7}
\delta(\partial\mathcal E_{j_l},\hat\phi)(h_{\e_{j_l}})
=\delta(\partial\mathcal E_{j_l})(\hat\phi h_{\e_{j_l}})
+\int_{{\bf G}_{n}(\mathbb R^{n+1})} S^{\perp}(\nabla\hat\phi)\cdot h_{\e_{j_l}}
\, d(\partial\mathcal E_{j_l}).
\end{equation}
By \eqref{mc0} for all sufficiently large $l$ and
all evaluated at $t=k\Delta t_{j_l}$ and if we write 
\begin{equation}
\label{bshort}
b_{j_l}:=\int_{\mathbb R^{n+1}}\frac{\hat\phi |\Phi_{\e_{j_l}}\ast\delta(\partial\mathcal
E_{j_l})|^2}{\Phi_{\e_{j_l}}\ast\|\partial\mathcal E_{j_l}\|+\e_{j_l}\Omega^{-1}}
\, dx
\end{equation}
for simplicity,
\begin{equation}
\label{brng8}
\begin{split}
&\Big|\delta(\partial\mathcal E_{j_l})(\hat\phi h_{\e_{j_l}})  
+b_{j_l}\Big|\leq  \e_{j_l}^{\frac14}(b_{j_l}+1)
\end{split}
\end{equation}
and by the Cauchy-Schwarz inequality and \eqref{mc01}, we have
\begin{equation}
\begin{split}
\Big|\int_{{\bf G}_{n}(\mathbb R^{n+1})} S^{\perp}(\nabla\hat\phi)\cdot h_{\e_{j_l}}
\, d(\partial\mathcal E_{j_l})\Big| & \leq \Big(\int_{\mathbb R^{n+1}} \hat\phi^{-1}
|\nabla\hat\phi|^2\, d\|\partial\mathcal E_{j_l}\|\Big)^{\frac12}\Big(
\int_{\mathbb R^{n+1}}\hat\phi |h_{\e_{j_l}}|^2\,d\|\partial\mathcal E_{j_l}\|\Big)^{\frac12}  \\ &
\leq c\,\|\partial\mathcal E_{j_l}\|(\Omega)^{\frac12}\Big((1+\e_{j_l}^{\frac14})
b_{j_l}+\e_{j_l}^{\frac14}\Big)^{\frac12},
\end{split}
\label{brng8.5}
\end{equation}
where we estimated as in \eqref{mecon7} and $c$ depends only on 
$\|\phi\|_{C^2}$, $\min_{x\in {\rm spt}\,\phi}\Omega$ and $\Cr{c_1}$
and independent of $i$. Since 
$\sup_{t\in [0,t_2]}\|\partial\mathcal E_{j_l}(t)\|
(\Omega)$ is bounded uniformly, \eqref{brng7}-\eqref{brng8.5} show that
for all sufficiently large $l$, we have
\begin{equation}
\sup_{t\in [t_1,t_2]}\delta(\partial\mathcal E_{j_l}(t),\hat\phi)(h_{\e_{j_l}}(\cdot,\partial\mathcal E_{j_l}(t)))
\leq c
\label{brng8.55}
\end{equation}
where $c$ depends only on $\|\partial\mathcal E_0 \|(\Omega)$,
$t_2$, $\|\phi\|_{C^2}$, $\min_{x\in {\rm spt}\,\phi}\Omega$ and $\Cr{c_1}$. Thus we have
\begin{equation}
\label{brng8.6}
\begin{split}
&\limsup_{l\rightarrow\infty}\sum_{k=k_1}^{k_2}\Delta t_{j_l}\delta(\partial
\mathcal E_{j_l}(k\Delta t_{j_l}),\hat\phi)(h_{\e_{j_l}}(\cdot,\partial\mathcal E_{j_l}
(k\Delta t_{j_l}))) \\
&=\limsup_{l\rightarrow\infty}\int_{t_1}^{t_2}\delta(\partial\mathcal E_{j_l}(t),\hat\phi)(h_{\e_{j_l}}(\cdot,\partial\mathcal E_{j_l}(t)))\, dt \\
& =-\liminf_{l\rightarrow\infty} \int_{t_1}^{t_2}\big(c-\delta(\partial\mathcal E_{j_l}(t),\hat\phi)(h_{\e_{j_l}}(\cdot,\partial\mathcal E_{j_l}(t)))\big)\, dt+c(t_2-t_1) \\
&\leq -\int_{t_1}^{t_2}\liminf_{l\rightarrow\infty}\big(c-\delta(\partial\mathcal E_{j_l}(t),\hat\phi)(h_{\e_{j_l}}(\cdot,\partial\mathcal E_{j_l}(t)))\big)\, dt+c(t_2-t_1) \\
& = \int_{t_1}^{t_2}\limsup_{l\rightarrow\infty}\, \delta(\partial\mathcal E_{j_l}(t),\hat\phi)(h_{\e_{j_l}}(\cdot,\partial\mathcal E_{j_l}(t)))\, dt
\end{split}
\end{equation}
where we used \eqref{brng8.55} and Fatou's Lemma. We estimate 
the integrand of \eqref{brng8.6} from above. Fix $t$. 
Let $\{j_l'\}_{l=1}^{\infty}\subset\{j_l\}_{l=1}^{\infty}$ be a subsequence 
such that the $\limsup$ is achieved, i.e.,
 \begin{equation}
 \label{brngex1}
 \limsup_{l\rightarrow\infty}\, \delta(\partial\mathcal E_{j_l}(t),\hat\phi)(h_{\e_{j_l}}(\cdot,\partial\mathcal E_{j_l}(t)))=\lim_{l\rightarrow\infty}\, \delta(\partial\mathcal E_{j_l'}(t),\hat\phi)(h_{\e_{j_l'}}(\cdot,\partial\mathcal E_{j_l'}(t))).
 \end{equation}
The right-hand side of \eqref{brng7} then have the same property for this
subsequence and
\begin{equation}
\label{brngex2}
\begin{split}
&\lim_{l\rightarrow\infty}\Big( -\delta(\partial \mathcal E_{j_l'})(\hat\phi h_{\e_{j_l'}})
-\int S^{\perp}(\nabla\hat\phi)\cdot h_{\e_{j_l'}}d(\partial\mathcal E_{j_l'})\Big) \\ 
&=\liminf_{l\rightarrow\infty}
\Big(-\delta(\partial \mathcal E_{j_l})(\hat\phi h_{\e_{j_l}})
-\int S^{\perp}(\nabla\hat\phi)\cdot h_{\e_{j_l}}d(\partial\mathcal E_{j_l})\Big).
\end{split}
\end{equation}
Using \eqref{brng8} and \eqref{brng8.5}, the right-hand side of \eqref{brngex2}
may be bounded by $\liminf_{l\rightarrow\infty}2b_{j_l}+c$
from above. The left-hand side of
\eqref{brngex2} is similarly estimated from below by $\limsup_{l\rightarrow
\infty} \frac12 b_{j_l'} -c$. Thus, for any subsequence satisfying \eqref{brngex1}, we have (evaluation at $t$)
\begin{equation}
\label{brngex3}
\limsup_{l\rightarrow\infty} \int_{\mathbb R^{n+1}}\frac{\hat\phi |\Phi_{\e_{j_l'}}\ast\delta(\partial\mathcal
E_{j_l'})|^2}{\Phi_{\e_{j_l'}}\ast\|\partial\mathcal E_{j_l'}\|+\e_{j_l'}\Omega^{-1}}
\, dx
\leq 4\liminf_{l\rightarrow\infty} \int_{\mathbb R^{n+1}}\frac{\hat\phi |\Phi_{\e_{j_l}}\ast\delta(\partial\mathcal
E_{j_l})|^2}{\Phi_{\e_{j_l}}\ast\|\partial\mathcal E_{j_l}\|+\e_{j_l}\Omega^{-1}}
\, dx+c 
\end{equation}
where $c$ is a constant estimated from above in terms of 
$\|\partial\mathcal E_0 \|(\Omega)$,
$t_2$, $\|\phi\|_{C^2}$, $\min_{x\in {\rm spt}\,\phi}\Omega$ and $\Cr{c_1}$.
Define the right-hand side of \eqref{brngex3} as $\tilde M (t)$ in the following.

For any $t$ with $\tilde M(t)<\infty$, by Lemma \ref{pbi1} (b) (note
$\hat\phi\geq i^{-1}\Omega$), the full 
sequence $\{\partial\mathcal E_{j'_l}\}_{l=1}^{\infty}$ converges to 
$V_t\in {\bf IV}_{n}(\mathbb R^{n+1})$ with $\mu_t=\|V_t\|$. 
From $\Omega\leq i\hat\phi$, we also have
\begin{equation}
\label{rng1}
\limsup_{l\rightarrow\infty} \int_{\mathbb R^{n+1}}\frac{\Omega|\Phi_{\e_{j'_l}}\ast\delta
(\partial\mathcal E_{j'_l}(t))|^2}{\Phi_{\e_{j'_l}}\ast\|\partial\mathcal E_{j'_l}(t)\|+
\e_{j'_l}\Omega^{-1}}\, dx\leq  i {\tilde M}(t).
\end{equation}
Set $M:=\|\partial\mathcal E_0\|(\Omega)\exp(\Cr{c_1}^2t_2/2)$ so that
we have
\begin{equation}
\limsup_{l\rightarrow\infty} \sup_{t\in [0,t_2]}
\|\partial\mathcal E_{j_l}(t)\|(\Omega)\leq M.
\label{rng2}
\end{equation}
By \eqref{brngex1}, \eqref{brng7}, \eqref{brng8} and Lemma \ref{pbi1} (c), 
we have
\begin{equation}
\label{brng8.9}
\begin{split}
&\limsup_{l\rightarrow\infty} \, \delta(\partial\mathcal E_{j_l}(t),\hat\phi)(h_{\e_{j_l}}
(\cdot,\partial\mathcal E_{j_l}(t)))=\lim_{l\rightarrow\infty}\, \delta(\partial\mathcal E_{j_l'}(t),\hat\phi)(h_{\e_{j_l'}}(\cdot,\partial\mathcal E_{j_l'}(t)))  \\
& \leq -\int_{\mathbb R^{n+1}}|h(\cdot,V_t)|^2\hat\phi\, d\|V_t\|+\limsup_{l\rightarrow
\infty}\int_{{\bf G}_{n}(\mathbb R^{n+1})} S^{\perp}(\nabla\hat\phi)\cdot h_{\e_{j'_l}}(\cdot,\partial\mathcal E_{j_l'}(t))
\, d(\partial\mathcal E_{j'_l}(t)).
\end{split}
\end{equation}
Let $\epsilon>0$ be arbitrary. Since $V_t\in {\bf IV}_{n}(\mathbb R^{n+1})$, 
there exists $\|V_t\|$ measurable, countably $n$-rectifiable set
$C\subset\mathbb R^{n+1}$ such that 
\begin{equation}
\int_{{\bf G}_{n}(\mathbb R^{n+1})} S^{\perp}(\nabla\phi(x))\, dV_t(x,S)
=\int_{\mathbb R^{n+1}} ({\rm Tan}^{n}( C,x))^{\perp}(\nabla\phi(x))\, 
d\|V_t\|(x)
\label{brng9}
\end{equation}
and $x\longmapsto ({\rm Tan}^{n}(C,x))^{\perp}(\nabla\phi(x))\Omega(x)^{-\frac12}$ is a $\|V_t\|$ measurable function on $\mathbb R^{n+1}$. Hence, 
corresponding to $\epsilon>0$, there exist $g\in C_c^{\infty}(\mathbb R^{n+1};
\mathbb R^{n+1})$ and $m'\in \mathbb N$ such that $g\in \mathcal B_{m'}$ and 
\begin{equation}
\label{brng10}
\int_{\mathbb R^{n+1}}|({\rm Tan}^{n} (C,x))^{\perp}(\nabla\phi(x))-g(x)|^2
\Omega(x)^{-1}\, d\|V_t\|(x)<\epsilon^2.
\end{equation}
Now we compute as (omitting $t$ dependence for simplicity)
\begin{equation}
\label{brng11}
\begin{split}
\int_{{\bf G}_{n}(\mathbb R^{n+1})} & S^{\perp}(\nabla\hat\phi)\cdot h_{\e_{j'_l}}
\, d(\partial\mathcal E_{j'_l})=\int_{{\bf G}_{n}(\mathbb R^{n+1})}
(S^{\perp}(\nabla\hat\phi)-g)\cdot h_{\e_{j'_l}}\, d(\partial\mathcal E_{j'_l}) \\
& +\Big(\int_{\mathbb R^{n+1}} g\cdot h_{\e_{j'_l}}d(\partial\mathcal E_{j'_l})+
\delta(\partial\mathcal E_{j'_l})(g)\Big)-\delta(\partial\mathcal E_{j'_l})(g)
+\delta V_t(g) \\ &+\int_{\mathbb R^{n+1}} h(\cdot,V_t)\cdot
 (g-({\rm Tan}^{n}(C,x))^{\perp}(\nabla\hat\phi))\, d\|V_t\| \\ &
 +\int_{{\bf G}_{n}(\mathbb R^{n+1})} h(\cdot,V_t)\cdot
S^{\perp}(\nabla\hat\phi)\, dV_t(\cdot,S).
\end{split}
\end{equation}
We estimate each term of \eqref{brng11}. We have
\begin{equation}
\label{brng12}
\begin{split}
&\Big|\int_{{\bf G}_{n}(\mathbb R^{n+1})}
(S^{\perp}(\nabla\hat\phi)-g)\cdot h_{\e_{j'_l}}\, d(\partial\mathcal E_{j'_l})
\Big|\leq i^{-1}\int_{\mathbb R^{n+1}} |\nabla\Omega||h_{\e_{j'_l}}|\, d\|\partial
\mathcal E_{j'_l}\| \\ & +
\Big(\int_{{\bf G}_{n}(\mathbb R^{n+1})} |S^{\perp}(\nabla\phi)-g|^2
\Omega^{-1}\, d(\partial\mathcal E_{j'_l})\Big)^{\frac12}
\Big(\int_{\mathbb R^{n+1}}|h_{\e_{j'_l}}|^2 \Omega\, d\|\partial\mathcal E_{j'_l}\|
\Big)^{\frac12} \\
& \leq i^{-1}\Cr{c_1}(\|\partial\mathcal E_{j'_l}\|(\Omega))^{\frac12}
\Big(\int_{\mathbb R^{n+1}}|h_{\e_{j'_l}}|^2\Omega\, d\|\partial\mathcal E_{j'_l}\|
\Big)^{\frac12} \\ & +
\Big(\int_{{\bf G}_{n}(\mathbb R^{n+1})} |S^{\perp}(\nabla\phi)-g|^2
\Omega^{-1}\, d(\partial\mathcal E_{j'_l})\Big)^{\frac12}
\Big(\int_{\mathbb R^{n+1}}|h_{\e_{j'_l}}|^2 \Omega\, d\|\partial\mathcal E_{j'_l}\|
\Big)^{\frac12}. 
\end{split}
\end{equation}
Since $\partial\mathcal E_{j'_l}$ converges to $V_t$ as varifold, 
\begin{equation}
\label{brng13}
\begin{split}
&\lim_{l\rightarrow\infty}
\int_{{\bf G}_{n}(\mathbb R^{n+1})} |S^{\perp}(\nabla\phi)-g|^2
\Omega^{-1}\, d(\partial\mathcal E_{j'_l})=
\int_{{\bf G}_{n}(\mathbb R^{n+1})} |S^{\perp}(\nabla\phi)-g|^2
\Omega^{-1}\, dV_t \\
&=\int_{\mathbb R^{n+1}} |({\rm Tan}^{n} (C,x))^{\perp}(\nabla\phi)-g|^2
\Omega^{-1}\, d\|V_t\|<\epsilon^2
\end{split}
\end{equation}
where we used \eqref{brng10}. 
Using \eqref{mc01} and \eqref{rng1}, \eqref{rng2}, \eqref{brng12} and \eqref{brng13}, we have
\begin{equation}
\label{brng14}
\limsup_{l\rightarrow\infty}
\Big|\int_{{\bf G}_{n}(\mathbb R^{n+1})}
(S^{\perp}(\nabla\hat\phi)-g)\cdot h_{\e_{j'_l}}\, d(\partial\mathcal E_{j'_l})
\Big|\leq \Cr{c_1} M^{\frac12}({\tilde M}(t))^{\frac12} i^{-\frac12}+ (i{\tilde M} (t))^{\frac12}\epsilon.
\end{equation}
By Proposition \ref{cl5} and \eqref{rng1}, we have
\begin{equation}
\label{rng3}
\lim_{l\rightarrow\infty} \Big|\int_{\mathbb R^{n+1}} g\cdot h_{\e_{j'_l}}\,
d(\partial\mathcal E_{j'_l})+\delta(\partial\mathcal E_{j'_l})(g)\Big|=0
\end{equation}
and the varifold convergence shows 
\begin{equation}
\label{rng4}
\lim_{l\rightarrow\infty}\Big|-\delta(\partial\mathcal E_{j'_l})(g)+\delta
V_t(g)\Big|=0.
\end{equation}
For the second last term of \eqref{brng11},
\begin{equation}
\label{rng5}
\begin{split}
&\Big|\int_{\mathbb R^{n+1}} h(\cdot,V_t)\cdot (g-
({\rm Tan}^{n} (C,x))^{\perp}(\nabla\hat\phi))\, d\|V_t\|\Big|
\leq i^{-1}\int_{\mathbb R^{n+1}} |h(\cdot, V_t)| |\nabla\Omega|\, d\|V_t\| \\
& +\int_{\mathbb R^{n+1}} |h(\cdot,V_t)||g-({\rm Tan}^{n} (C,x))^{\perp}(\nabla\phi)|\, d\|V_t\| \\
&\leq i^{-\frac12} \Cr{c_1}M^{\frac12} ({\tilde M}(t))^{\frac12}+(i{\tilde M}(t))^{\frac12}\epsilon
\end{split}
\end{equation}
where we used the Cauchy-Schwarz inequality, 
\eqref{rng1}, \eqref{rng2} (which also hold for the limiting
quantities) and \eqref{brng10}. For the last term of \eqref{brng11}, estimating
as in \eqref{rng5}, 
\begin{equation}
\label{rng6}
\begin{split}
\int_{{\bf G}_{n}(\mathbb R^{n+1})} h(\cdot, V_t)\cdot S^{\perp}(\nabla
\hat\phi)\, dV_t &\leq \int_{{\bf G}_{n}(\mathbb R^{n+1})} h(\cdot, V_t)\cdot S^{\perp}(\nabla
\phi)\, dV_t+i^{-\frac12}\Cr{c_1}M^{\frac12} ({\tilde M}(t))^{\frac12} \\
& = \int_{\mathbb R^{n+1}} h(\cdot, V_t)\cdot \nabla\phi\, d\|V_t\|+
i^{-\frac12}\Cr{c_1}M^{\frac12} ({\tilde M}(t))^{\frac12}
\end{split}
\end{equation}
where we used \eqref{fvf2}. 
Finally, combining \eqref{brng11}, \eqref{brng14}-\eqref{rng6} and 
letting $\epsilon\rightarrow 0$, we obtain
\begin{equation}
\limsup_{l\rightarrow\infty} \int_{{\bf G}_{n}(\mathbb R^{n+1})} 
S^{\perp}(\nabla\hat\phi)\cdot h_{\e_{j'_l}}\, d(\partial\mathcal E_{j'_l})
\leq 3\Cr{c_1}i^{-\frac12}M^{\frac12}({\tilde M}(t))^{\frac12}+ \int_{\mathbb R^{n+1}}
h(\cdot,V_t)\cdot \nabla\phi\, d\|V_t\|.
\label{rng7}
\end{equation}
From \eqref{brng8.9} and \eqref{rng7}, we obtain 
\begin{equation}
\label{rng8}
\limsup_{l\rightarrow\infty}\delta (\partial\mathcal E_{j_l}(t),\hat\phi)
(h_{\e_{j_l}}(\cdot,\partial\mathcal E_{j_l}(t)))\leq\delta(V_t,\phi)(h(\cdot,V_t))
+3\Cr{c_1}i^{-\frac12}(M+{\tilde M}(t)).
\end{equation}
Since $\hat\phi\leq \Omega$, we have by Fatou's Lemma that
\begin{equation}
\int_{t_1}^{t_2}{\tilde M}(t)\, dt\leq 4 \liminf_{l\rightarrow\infty}
\int_{t_1}^{t_2}\int_{\mathbb R^{n+1}} \frac{\Omega|\Phi_{\e_{j_l}}\ast
\delta(\partial\mathcal E_{j_l}(t))|^2}{\Phi_{\e_{j_l}}\ast\|\partial\mathcal E_{j_l}
(t)\|+\e_{j_l}\Omega^{-1}}\, dxdt+c<\infty
\label{rng9}
\end{equation}
by \eqref{mecones}. Thus, by \eqref{brng4}, \eqref{brng6}, \eqref{brng8.6},
\eqref{rng8}, \eqref{rng9} and letting $i\rightarrow\infty$, we obtain 
\eqref{brng2} for time-independent $\phi\in C_c^{\infty}(\mathbb R^{n+1};
\mathbb R^+)$. For time dependent $\phi\in C_c^{\infty}(\mathbb R^{n+1}
\times \mathbb R^+;\mathbb R^+)$, we repeat the same argument. 
We similarly define $\hat\phi$ and use \eqref{exapp4} with $\hat\phi(\cdot,
t)$. Instead of \eqref{brng3}, we obtain a formula with one extra term, namely,
\begin{equation}
\label{rng10}
\begin{split}
\|\partial\mathcal E_{j_l}(s)\|(\hat\phi(\cdot,s))\Big|_{s=t-\Delta t_{j_l}}^t\leq  &\Delta t_{j_l}\big\{\delta(\partial\mathcal E_{j_l}(t),\hat\phi(\cdot,t))(h_{\e_{j_l}}
(\cdot,\partial\mathcal E_{j_l}(t)))+\e_{j_l}^{\frac18}\big\} \\ &
+\|\partial\mathcal E_{j_l}(t-\Delta t_{j_l})\|(\phi(\cdot,t)-\phi(\cdot,t-\Delta t_{j_l})). 
\end{split}
\end{equation}
Note that the last term has $\phi$ instead of $\hat\phi$. 
A similar inequality to \eqref{brng4} will have the summation of the last term 
of \eqref{rng10}. It is not difficult to check using \eqref{meconeq} 
and Lemma \ref{pbi1} (a) that we have
\begin{equation}
\label{rng11}
\begin{split}
&\lim_{l\rightarrow\infty} \sum_{k=k_1}^{k_2}\|\partial\mathcal E_{j_l}
((k-1)\Delta t_{j_l})\|(\phi(\cdot,k\Delta t_{j_l})-\phi(\cdot,(k-1)\Delta t_{j_l})) \\ &
=\lim_{l\rightarrow\infty}\sum_{k=k_1}^{k_2} \|\partial\mathcal E_{j_l}(k\Delta
t_{j_l})\|\big(\frac{\partial\phi}{\partial t}(\cdot,k\Delta t_{j_l})\big) \, \Delta t_{j_l}
=\lim_{l\rightarrow\infty}\int_{t_1}^{t_2}
\|\partial\mathcal E_{j_l}(t)\| \big(\frac{\partial\phi}{\partial t}(\cdot,t)\big)\, dt \\
& =\int_{t_1}^{t_2}\|V_t\|\big(\frac{\partial\phi}{\partial t}(\cdot,t)\big)\, dt,
\end{split}
\end{equation}
where we also used the dominated convergence theorem in the last step.
The rest proceeds by the same argument with error estimates
coming from the time-dependency of $\hat\phi$.  For example, in \eqref{brng8.6},
we need to regard $\hat\phi(\cdot,t)$ as a piecewise constant function with respect to
time variable on $[t_1,t_2]$, namely, in place of $\hat\phi$, we need to have
\begin{equation}
\hat\phi_{j_l}(\cdot, t):=\hat\phi(\cdot, k\Delta
t_{j_l}) \mbox{ if } t\in ((k-1)\Delta t_{j_l},k\Delta t_{j_l}].
\label{rng12}
\end{equation}
For $\delta(\partial\mathcal E_{j_l}(t),\hat\phi_{j_l}(\cdot,t))
(h_{\e_{j_l}}(\cdot,\partial\mathcal E_{j_l}(t)))$ in the last line of \eqref{brng8.6}, if we 
replace $\hat\phi_{j_l}(\cdot,t)$ by $\hat\phi(\cdot,t)$, it only results
in errors of order $\Delta t_{j_l}$ times certain negative power of $\e_{j_l}$ which
remains small and goes to 0 uniformly as $l\rightarrow\infty$. Thus we may 
subsequently proceed just like the time independent case and we have \eqref{brng2} for
$C^{\infty}_c$ case, and by approximation for $C^1_c$ case. 
\hfill{$\Box$}

Now, the proof of Theorem \ref{sthm} is complete: (1) is clear from the construction using $\mathcal E_0=\{E_{0,i}\}_{i=1}^N$,
(2) is by Lemma \ref{pbi1} (a) and (c), (3) and (4) follow from Theorem \ref{brakkeinq}.
We note that the claim of 
Theorem \ref{ktthm} is slightly different from \cite{KT,Tonegawa2} in 
that it is stated for $(x,t)\in \mathbb R^{n+1}\setminus S_t$ here 
instead of ${\rm spt}\,\|V_t\|
\setminus S_t$, allowing a possibility of $O_{(x,t)}\cap {\rm spt}\,\mu$ being empty. But exactly the
same proof of \cite{KT} gives this slightly stronger claim of partial regularity and 
we write the result in this form. 
\section{Proof of Theorem \ref{sthm3}}
\label{secBas}
Let $\mu$ be a measure on $\mathbb R^{n+1}\times \mathbb R^+$ defined as
in Definition~\ref{sptime}.
\begin{lemma}
\label{cmd1}
We have the following properties for $\mu$ and $\{V_t\}_{t\in\mathbb R^+}$. 
\begin{itemize}
\item[(1)] ${\rm spt}\,\|V_t\|\subset \{x\in \mathbb R^{n+1} : (x,t)\in {\rm spt}\,\mu\}$
for all $t>0$.
\item[(2)] ${\rm clos}\,\{(x,t) : x\in {\rm spt}\,
\|V_t\|, V_t\in {\bf IV}_{n}(\mathbb R^{n+1})\}\cap \{(x,t) : t>0\}={\rm spt}\, \mu\cap \{(x,t) : t>0\}$.
\end{itemize}
\end{lemma}
{\it Proof}. Suppose $x\in {\rm spt}\,\|V_t\|$ and $t>0$. Then for any $r>0$, 
there exists some $\phi\in C_c^2(U_{2r}(x);\mathbb R^+)$ with
$\|V_t\|(\phi)>0$. For any $t'\in [0,t)$, by \eqref{brng2} and 
the Cauchy-Schwarz inequality, we have
\begin{equation}
\begin{split}
\|V_t\|(\phi) - \|V_{t'}\|(\phi)&\leq \int_{t'}^t \int_{U_{2r}(x)}
-|h(\cdot,V_s)|^2\phi+\nabla\phi\cdot h(\cdot,V_s)\, d\|V_s\|ds \\ &
\leq \int_{t'}^t \int_{U_{2r}(x)}\frac{|\nabla\phi|^2}{2\phi}\, d\|V_s\|ds
\leq (t-t')\|\phi\|_{C^2} \sup_{s\in [t',t]} \|V_s\|(U_{2r}(x)).
\end{split}
\label{cmd2}
\end{equation}
Choosing $t'$ sufficiently close to $t$, \eqref{cmd2} shows that there exists some $t'<t$ such that $\frac12 \|V_t\|(\phi)\leq \|V_{s}\|(\phi)$ for all $s\in [t',t)$. Thus, $\int_{U_{2r}(x)\times [t',t)}
\phi\, d\mu\geq \frac12 (t-t')\|V_t\|(\phi)>0$. If $(x,t)\notin {\rm spt}\,\mu$, there must 
be some open set $U$ in $\mathbb R^{n+1}\times \mathbb R^+$ with $\mu(U)=0$,
but this is a contradiction to the preceding sentence. Thus we have (1). 

Suppose $(x,t)\in{\rm clos}\,\{(x,t) : x\in {\rm spt}\,
\|V_t\|, V_t\in {\bf IV}_{n}(\mathbb R^{n+1})\} \cap \{(x,t) : t>0\}$. Then there exists a sequence $\{(x_i,t_i)\}_{i=1}^{\infty}$
such that $x_i\in {\rm spt}\,\|V_{t_i}\|$, $t_i>0$ and 
$\lim_{i\rightarrow\infty} (x_i,t_i)=(x,t)$. By (1), $(x_i,t_i)\in {\rm spt}\,\mu$. Since 
${\rm spt}\,\mu$ is a closed set by definition, 
we have $(x,t)\in {\rm spt}\,\mu$, proving $\subset$ of (2). Given $(x,t)\in {\rm spt}\,
\mu$ with $t>0$ and $\epsilon>0$, we have $\mu(B_\epsilon(x)\times (t-\epsilon,
t+\epsilon))>0$. Then, there must be some $t'\in (t-\epsilon,t+\epsilon)$ such that 
$\|V_{t'}\|(B_{\epsilon}(x))>0$ and $V_{t'}\in {\bf IV}_{n}(\mathbb R^{n+1})$. 
If ${\rm spt}\,\|V_{t'}\|\cap B_{\epsilon}(x)=\emptyset$, then 
we would have $\|V_{t'}\|(B_{\epsilon}(x))=0$, a 
contradiction. Thus we have some $x'\in {\rm spt}\,\|V_{t'}\|\cap B_{\epsilon}(x)$
with $V_{t'}\in {\bf IV}_{n}(\mathbb R^{n+1})$ and
$|t'-t|<\epsilon$. Since $\epsilon>0$ is arbitrary, 
this proves $\supset$ of (2). 
\hfill{$\Box$}
\begin{rem}
In (1), it may happen that the left-hand side is strictly smaller than the right-hand side. 
For example, consider a shrinking sphere. At the moment of vanishing, we have $\|V_t
\|=0$ since it is a point and has zero measure, thus 
${\rm spt}\,\|V_t\|=\emptyset$. 
On the other hand, the vanishing point is in ${\rm spt}\,\mu$, and the right-hand
side is not the empty set. We may also encounter a situation where some 
portion of measure vanishes, thus the difference between the left- and right-hand sides 
of (1) may be of positive $\mathcal H^{n}$ 
measure. 
We also point out that, in general, (1) and (2) are not true if $t=0$ is included. We may
have some portion of measure $\|\partial\mathcal E_0\|$ 
vanishing instantly at $t=0$. For example, consider on $\mathbb R^2$ a line segment
with two end points
which is surrounded by one of open partitions. For the first Lipschitz deformation step,
such line segment may be eliminated as we indicated in \ref{ibry}. Thus, even though we have some positive 
measure at $t=0$, ${\rm spt}\,\mu$ may be empty nearby. 
\end{rem}
\bigskip
Let $\eta\in C^{\infty}_c(U_{2};\mathbb R^+)$ be a radially symmetric function 
such that $\eta=1$ on $B_{1}$, $|\nabla\eta|\leq 2$ and $\|\nabla^2 \eta\|\leq 4$.
Then define for $x,y\in\mathbb R^{n+1}$, $s,t\in \mathbb R$ with $s>t$ and $R>0$
\begin{equation}
\label{cmd3}
\begin{split}
&\rho_{(y,s)}(x,t):=\frac{1}{(4\pi(s-t))^{\frac{n}{2}}} \exp\big(-\frac{|x-y|^2}{4(s-t)}\big), \,\,
\hat\rho_{(y,s)}(x,t):=\eta(x-y) \rho_{(y,s)}(x,t), \\
&\hat\rho_{(y,s)}^R(x,t):=\eta\big(\frac{x-y}{R}\big) \rho_{(y,s)}(x,t).
\end{split}
\end{equation}
We often write $\rho_{(y,s)}$ or $\rho$
for $\rho_{(y,s)}(x,t)$ when the meaning is clear from the context and the same for
$\hat\rho_{(y,s)}$ and $\hat\rho_{(y,s)}^R$. The following is a variant of
well-known Huisken's monotonicity formula \cite{Huisken}. 
We include the outline of proof and the reader is advised to 
see \cite[Lemma 6.1]{KT} for more details.
\begin{lemma}
\label{cmd4}
There exists $\Cl[c]{c_5}$ depending only on $n$ with the following property.  
For $0\leq t_1<t_2<s<\infty$, $y\in \mathbb R^{n+1}$ and $R>0$, we have
\begin{equation}
\label{cmd6}
\|V_{t}\|(\hat\rho_{(y,s)}^R(\cdot,t))\Big|_{t=t_1}^{t_2}\leq \Cr{c_5} R^{-2}(t_2-t_1) \sup_{t'\in 
[t_1,t_2]} R^{-n}\|V_{t'}\|(B_{2R}(y)).
\end{equation}
\end{lemma}
{\it Proof}. After change of variables by $\tilde x =(x-y)/R$ and $\tilde t = (t-s)/R^2$, 
we may regard $R=1$ and $(y,s)=(0,0)$. 
A direct computation shows that for any $S\in {\bf G}(n+1,n)$, we have
\begin{equation*}
\frac{\partial\rho}{\partial t}+S\cdot\nabla_x^2 \rho+\frac{|S^{\perp}(\nabla_x\rho)|^2}{
\rho}=0
\end{equation*}
for all $t<0$ and $x\in \mathbb R^{n+1}$. The same computation for $\hat\rho$ has some
extra terms coming from differentiations of $\eta$, and such terms are bounded by 
$c(n)(-t)^{-\frac{n}{2}}\exp(1/4t)$ since ${\rm spt}\,|\nabla\eta|\subset B_2
\setminus U_1$. Thus we have
\begin{equation}
\label{cmd7}
\Big|\frac{\partial\hat\rho}{\partial t}+S\cdot\nabla_x^2 \hat\rho+\frac{|S^{\perp}(\nabla_x\hat\rho)|^2}{
\hat\rho}\Big|\leq \Cr{c_5}\chi_{B_2\setminus U_1}.
\end{equation}
Use $\hat\rho$ in \eqref{brng2} as well as \eqref{cmd7} to find that
\begin{equation}
\label{cmd8}
\|V_{t}\|(\hat\rho_{(0,0)}(\cdot,t))\Big|_{t=t_1}^{t_2}\leq \Cr{c_5}\int_{t_1}^{t_2}\|V_{t'}\|(B_2
\setminus U_1)\, dt'. 
\end{equation}
Then \eqref{cmd8} gives \eqref{cmd6}.
\hfill{$\Box$}
\begin{lemma}
\label{upper}
For any $\lambda >1$, there exists $\Cl[c]{c_f}\in (1,\infty)$ depending only on $n$, $\lambda$, 
$\Omega$
and $\|\partial\mathcal E_0\|(\Omega)$ such that
\begin{equation}
\label{ya1}
\sup_{x\in B_\lambda, r\in (0,1], t\in [\lambda^{-1},\lambda]}r^{-n} \|V_t\|(B_r(x))\leq \Cr{c_f}.
\end{equation}
\end{lemma}
{\it Proof}. We use \eqref{cmd6} with $s=t+r^2$, $t_2=t\in [\lambda^{-1},\lambda]$, $t_1=0$, $R=1$ and $y\in B_{\lambda}$. Then we 
obtain also using $\eta\lfloor_{B_1(y)}=1$ that
\begin{equation}
\label{ya2}
\frac{e^{-\frac14}}{(4\pi r^2)^{\frac{n}{2}}} \|V_t\|(B_r(y))\leq 
\frac{1}{(4\pi t)^{\frac{n}{2}}}  \|V_0\|(B_2(y))+\Cr{c_5} t \sup_{t'\in [0,t]}\|V_{t'}\|
(B_2(y)).
\end{equation}
The quantities on the right-hand side of \eqref{ya2} are all controlled by the stated quantities thus
we obtain \eqref{ya1}. 
\hfill{$\Box$}
\begin{rem}
If $\|\partial\mathcal E_0\|$ satisfies the density ratio upper bound 
\begin{equation}
\sup_{x\in\mathbb R^{n+1}, r\in (0,1]} r^{-n}\|\partial\mathcal E_0\|(B_r(x))<\infty,
\label{ya4}
\end{equation}
then we may obtain 
up to the initial time estimate for \eqref{ya1}.
\end{rem}
The following is essentially Brakke's clearing out lemma \cite[6.3]{Brakke} proved using
Huisken's monotonicity formula. 
\begin{lemma}
\label{cmd9} For any $\lambda>1$, there exist positive constants $\Cl[c]{c_6},\Cl[c]{c_7}\in (0,1)$ 
depending only on $n$, $\lambda$, $\Omega$ and $\|\partial\mathcal E_0\|(\Omega)$
such that the following holds. 
For $(x,t)\in {\rm spt}\,\mu\cap (B_\lambda\times [\lambda^{-1},\lambda])$ and $r\in (0,\frac12]$ with $t-\Cr{c_7}r^2\geq (2\lambda)^{-1}$, 
we have 
\begin{equation}
\label{cmd10}
\|V_{t-\Cr{c_7} r^2}\|(B_r(x))\geq \Cr{c_6} r^{n}.
\end{equation}
\end{lemma}
{\it Proof}. 
By Lemma \ref{cmd1} (2), there exists a sequence $(x_i,t_i)\in {\rm spt}\,\|V_{t_i}\|$
with $\lim_{i\rightarrow\infty}(x_i,t_i)=(x,t)$. We may also have 
$V_{t_i}\in {\bf IV}_{n}(\mathbb R^{n+1})$, thus any neighborhood of $x_i$ contains
some point of integer density of $\|V_{t_i}\|$. Thus we may as well assume that 
$\theta^{n}(\|V_{t_i}\|,x_i)\geq 1$. One uses \eqref{cmd6} with $R=r$, 
$t_1=t-\Cr{c_7}r^2$
($\Cr{c_7}$ to be decided),
$t_2=t_i$, $y=x_i$ and $s=t_i+\epsilon$ to obtain
\begin{equation}
\label{cmd11}
\|V_{s}\|(\hat\rho_{(x_i,t_i+\epsilon)}^r(\cdot,s))\Big|_{s=t-\Cr{c_7}r^2}^{t_i}\leq 
\Cr{c_5}r^{-2}(t_i-t+\Cr{c_7}r^2)\sup_{s\in [t-\Cr{c_7}r^2,t_i]} r^{-n}\|V_s\|(U_{2r}(x_i)).
\end{equation}
 By letting $\epsilon\rightarrow 0+$, $\theta^{n}(\|V_{t_i}\|,x_i)\geq 1$ and \eqref{cmd11} give
\begin{equation}
 \label{cmd12}
 1\leq \|V_{t-\Cr{c_7}r^2}\|(\hat\rho_{(x_i,t_i)}^r(\cdot,t-\Cr{c_7}r^2))
+\Cr{c_5}r^{-2}(t_i-t+\Cr{c_7}r^2)\sup_{s\in [t-\Cr{c_7}r^2,t_i]} r^{-n}\|V_s\|(U_{2r}(x_i)).
\end{equation}
Let $i\rightarrow \infty$ for \eqref{cmd12} to obtain
\begin{equation}
\label{cmd13}
1\leq \|V_{t-\Cr{c_7}r^2}\|(\hat\rho_{(x,t)}^r(\cdot,t-\Cr{c_7}r^2))+\Cr{c_5}\Cr{c_7}
\sup_{s\in [t-\Cr{c_7}r^2,t]}r^{-n}\|V_s\|(U_{2r}(x)). 
\end{equation}
We also have $\|V_{t-\Cr{c_7}r^2}\|(\hat\rho_{(x,t)}^r(\cdot,t-\Cr{c_7}r^2))
\leq (4\pi \Cr{c_7})^{-\frac{n}{2}} r^{-n} \|V_{t-\Cr{c_7}r^2}\|(U_{2r}(x))$. 
Now, given $\lambda$, let $\Cr{c_f}$ be a constant obtained in Lemma \ref{upper}
corresponding to $\lambda$ there equals to $2\lambda$. Suppose we choose
$\Cr{c_7}<(2\lambda)^{-1}$ and $t\geq \lambda^{-1}$ 
so that $t-\Cr{c_7}r^2\geq (2\lambda)^{-1}$. Then by \eqref{ya1} and \eqref{cmd13}, we have
\begin{equation}
\label{cmd14}
1\leq (4\pi \Cr{c_7})^{-\frac{n}{2}} r^{-n} \|V_{t-\Cr{c_7}r^2}\|(U_{2r}(x))+
\Cr{c_5}\Cr{c_7}2^{n} \Cr{c_f}. 
\end{equation}
Choose $\Cr{c_7}$ sufficiently small so that the last term is less than $1/2$.
Then we have a lower bound for $r^{-n}\|V_{t-\Cr{c_7}r^2}\|(B_{2r}(x))$.
By adjusting constants again, we obtain \eqref{cmd10}.
\hfill{$\Box$}
\begin{rem}
If we have \eqref{ya4}, then we may also obtain \eqref{cmd10} up to $t=0$, namely,
we may replace $[\lambda^{-1},\lambda]$ in the statement to $(0,\lambda]$ and 
for $r\in (0,\frac12]$ with $t-\Cr{c_7}r^2\geq 0$. 
\end{rem}
\begin{cor}
\label{ya3}
For any open set $U\subset  B_{\lambda}$ and $t\in (\lambda^{-1},\lambda]$, we have
\begin{equation}
\mathcal H^{n}( \{x \in U : (x,t)\in {\rm spt}\, \mu\}) \leq \limsup_{s\rightarrow t-} 
{\bf B}_{n+1}\Cr{c_6}^{-1}\omega_n \|V_s\|(U).
\label{ya5}
\end{equation}
\end{cor}
{\it Proof}. It is enough to prove the estimate for $K_t:=\{x \in K : (x,t)\in {\rm spt}\, \mu\}$
where $K\subset U$ is compact and arbitrary. For each $x\in K_t$, for all 
sufficiently small $r$, $B_r(x)\subset U$ and by Lemma 
\ref{cmd9}, $\|V_{t-\Cr{c_7}r^2}\|(B_r(x))\geq \Cr{c_6}r^{n}$. Applying the Besicovitch
covering theorem to such family of balls, and recalling the definition of the Hausdorff
measure, we have a disjoint family of balls $\{B_r(x_1),\ldots, B_r(x_J)\}$ such that
($\mathcal H^{n}_{2r}$ is as defined in \cite[Definition 2.1(i)]{Evans-Gariepy})
\begin{equation}
\label{ya6}
{\bf B}_{n+1}^{-1}\mathcal H^{n}_{2r}(K_t)\leq J\omega_{n} r^{n}
\leq \Cr{c_6}^{-1}\omega_n  \sum_{i=1}^J \|V_{t-\Cr{c_7}r^2}\|(B_r(x_i))
\leq \Cr{c_6}^{-1} \omega_n\|V_{t-\Cr{c_7}r^2}\|(U).
\end{equation}
By letting $r\rightarrow 0$ for \eqref{ya6}, we obtain \eqref{ya5}.
\hfill{$\Box$}
\begin{rem} 
Lemma \ref{cmd1} (1) and Corollary \ref{ya3} prove \eqref{thsup1} of Proposition~\ref{sptimeato}.
\end{rem}
\begin{lemma}
\label{ya7}
Let $\{\mathcal E_{j_l}(t)\}_{l=1}^{\infty}$ be a sequence obtained in Proposition
\ref{promecone} and denote the open partitions by $\{E_{j_l,k}(t)\}_{k=1}^N$ for each
$j_l$ and $t\in \mathbb R^+$, i.e., $\mathcal E_{j_l}(t)=\{E_{j_l,k}(t)\}_{k=1}^N$. 
For fixed $k\in \{1,\ldots, N\}$, $0<r<\infty$, $x\in \mathbb R^{n+1}$
and $t>0$ with $t-r^2>0$, suppose 
\begin{equation}
\label{ya8}
\lim_{l\rightarrow\infty}
\mathcal L^{n+1}(B_{2r}(x)\setminus E_{j_l,k}(t))=0
\end{equation}
and 
\begin{equation}
\label{ya10}
\mu(B_{2r}(x)\times [t-r^2,t+r^2])=0.
\end{equation}
Then for all $t'\in (t-r^2,t+r^2]$, we have
\begin{equation}
\label{ya9}
\lim_{l\rightarrow\infty} 
\mathcal L^{n+1}(B_r(x)\setminus E_{j_l,k}(t') )=0.
\end{equation}
\end{lemma}
{\it Proof}. For a contradiction, if \eqref{ya9} were not true for some $t'\in (t-r^2,
t+r^2]$, 
by compactness of $BV$ functions, there exists a subsequence
$\{j'_l\}_{l=1}^\infty$ such that $\chi_{E_{j'_l,k}(t')}$ converges to $\chi_{E_k(t')}$ 
in $L^1(B_{2r}(x))$ and $\mathcal L^{n+1}(B_r(x)\setminus E_k(t'))>0$. By the
lower semicontinuity property, we have
$\|\nabla \chi_{E_k(t')}\|\leq \|V_{t'}\|$. By Lemma \ref{cmd1} (1) and \eqref{ya10},
we have $\|\nabla \chi_{E_k(t')}\|(B_{2r}(x))=0$. Then, $\chi_{E_k(t')}$ is a constant function
on $B_{2r}(x)$ and is identically 1 or 0. Since $\mathcal L^{n+1}(B_r(x)\setminus E_k(t'))>0$, 
$\chi_{E_k(t')}=0$ on $B_{2r}(x)$. Repeating the same argument, we may conclude
that there exist some $k' \in \{1,\ldots, N\}$, $k'\neq k$, and a subsequence 
(denoted again by $\{j'_l\}_{l=1}^{\infty}$) such that $\chi_{E_{j'_l,k'}(t')}$ converges to
$\chi_{E_{k'}(t')}$ and $\mathcal L^{n+1}(B_{2r}(x)\setminus E_{k'}(t'))=0$. Thus, we have
a situation where, at time $t$, $E_{j'_l,k}(t)$ occupies most of $B_{2r}(x)$ while at 
time $t'$, $E_{j'_l,k'}(t')$ occupies most of $B_{2r}(x)$ for all large $l$. 
In particular, for all sufficiently large $l$, we have $\mathcal L^{n+1} (B_{2r}(x)\setminus
E_{j'_l,k}(t))<\omega_{n+1} r^{n+1}/10$ and $\mathcal L^{n+1} (B_{2r}(x)\setminus
E_{j'_l,k'} (t'))<\omega_{n+1} r^{n+1}/10$. The maps
$f_1$ and $f_2$ for the construction of $\{\mathcal E_{j,l}\}$ in Proposition \ref{exapp}
change volume of each open partitions very little at each step 
(note Definition \ref{mrld} (b) for $f_1$, and $f_2$ is diffeomorphism 
which is close to identity, see \eqref{est_F2} and \eqref{est_gradF2}), there exists some $ t_l\in (t,t')$ (or 
$(t',t)$) such that $\frac14 \omega_{n+1} r^{n+1} \leq
\mathcal L^{n+1}(B_r(x)\cap E_{j'_l,k}(t_l))\leq \frac34 \omega_{n+1} r^{n+1}$. By the 
relative isoperimetric inequality, there exists a positive 
constant $c$ depending only on $n$ such that 
\begin{equation}
\|\partial \mathcal E_{j'_l}(t_l)\|(B_r(x))\geq \|\nabla\chi_{E_{j'_l,k}( t_l)}\|(B_r(x))
\geq c r^{n}.
\label{ya12}
\end{equation}
We may assume without loss of generality that $t_l\in 2_{\mathbb Q}$. 
Fix an arbitrary
$\hat t\in 2_{\mathbb Q}\cap (t-r^2, \min\{t,t'\})$. Choose $\phi\in C_c^2(U_{2r}(x);
\mathbb R^+)$ such that $\phi=1$ on $B_r(x)$ and $0\leq \phi\leq 1$ on
$U_{2r}(x)$. Now, we repeat the same argument leading to \eqref{mecon5} with
$t_2=t_l$ and $t_1=\hat t$ to obtain
\begin{equation}
\begin{split}
\liminf_{l\rightarrow\infty}&\Big(\|\partial\mathcal E_{j'_l}(t_l)\|(\phi)-\|\partial\mathcal E_{j'_l}(\hat t)\|
(\phi+i^{-1}\Omega)\Big) \\ & \leq \liminf_{l\rightarrow\infty} \frac12 \int_{\hat t}^{t_l}
\int_{\mathbb R^{n+1}}\frac{|\nabla (\phi+i^{-1}\Omega)|^2}{\phi+i^{-1}\Omega}\, d
\|\partial\mathcal E_{j'_l}(t)\|dt\\
& \leq \liminf_{l\rightarrow\infty} \int_{\hat t}^{t_l}\int_{\mathbb R^{n+1}} \frac{|\nabla\phi|^2}{\phi}+i^{-1}
\Cr{c_1}^2\Omega\, d\|\partial\mathcal E_{j'_l}(t)\|dt\leq i^{-1} \Cr{c_1}^2\int_{\hat t}^{t
+r^2}\|V_t\|(\Omega)\, dt,
\end{split}
\label{ya11}
\end{equation}
where we used the dominated convergence theorem and $\|V_t\|(U_{2r}(x))=0$ 
which follows from \eqref{ya10}. Since $\|\partial\mathcal E_{j'_l}(\hat t)\|(\phi)
\rightarrow \|V_{\hat t}\|(\phi)=0$, \eqref{ya11} proves after letting $i\rightarrow \infty$ that
$\liminf_{l\rightarrow\infty}\|\partial\mathcal E_{j'_l}(t_l)\|(\phi)=0$. But this would be a 
contradiction to \eqref{ya12}. 
\hfill{$\Box$}
\begin{lemma}
\label{ya13}
Let $\{\mathcal E_{j_l}(t)\}_{l=1}^{\infty}$ and $\{E_{j_l,k}(t)\}_{k=1}^N$ be the 
same as Lemma \ref{ya7}. For fixed $k\in \{1,\ldots,N\}$, 
$0<r<\infty$, $x\in \mathbb R^{n+1}$, suppose 
\begin{equation}
\label{ya14}
B_{2r}(x)\subset E_{j_l,k}(0)
\end{equation}
for all $l\in \mathbb N$ and 
\begin{equation}
\label{ya15}
\mu(B_{2r}(x)\times [0,r^2])=0.
\end{equation}
Then, for all $t'\in (0,r^2]$, we have
\begin{equation}
\label{ya16}
\lim_{l\rightarrow\infty} \mathcal L^{n+1} (B_r(x)\setminus E_{j_l,k}(t'))=0.
\end{equation}
\end{lemma}
{\it Proof}. By \eqref{ya14}, we have $\|V_0\|(B_{2r}(x))=0$ and Proposition \ref{cmd1} (1)
and \eqref{ya15} show $\|V_t\|(U_{2r}(x))=0$ for $t\in (0,r^2]$. Then, we may argue 
just like the proof of Lemma \ref{ya7}, where we take $\hat t$ there by $\hat t=0$. 
We omit the proof since it is similar.
\hfill{$\Box$}

The following Lemma \ref{Braclear} is from \cite[3.7, ``Sphere barrier to
external varifolds'']{Brakke}. 
\begin{lemma}
\label{Braclear}
For some $t\in \mathbb R^+$, $x\in \mathbb R^{n+1}$ and $r>0$, 
suppose $\|V_t\|(U_r(x))=0$. Then for $t'\in [t,t+\frac{r^2}{2n}]$, we have
$\|V_{t'}\|(U_{\sqrt{r^2-2n(t' -t)}}(x))=0$.
\end{lemma}
 {\it Proof of Theorem \ref{sthm3}}. 
We may choose a subsequence so that for all $t\in 2_{\mathbb
Q}$, each $\chi_{E_{j_l,k}(t)}$ converges in $L^1_{loc}(\mathbb R^{n+1})$ 
to $\chi_{E_k(t)}$ as $l\rightarrow \infty$.
This is due to the mass bound and $L^1$ compactness of 
BV functions. 
Consider the complement of ${\rm spt}\,\mu \cup ({\rm spt}\,\|V_0\|\times \{0\})$ in $\mathbb R^{n+1}\times \mathbb R^+$ which is open
in $\mathbb R^{n+1}\times\mathbb R^+$, and let $S$ be a connected
component. For any point $(x,t)\in S$, there exists $r>0$ such that 
$B_{2r}(x)\times [t-r^2,t+r^2]\subset S$ if $t>0$, and $B_{2r}(x)
\times [0,r^2]\subset S$ if $t=0$. 
First consider the case $t=0$. Since $B_{2r}(x)$ is in the complement of
${\rm spt}\,\|V_0\|=\Gamma_0$, for some small enough $0<t'\leq r^2$, 
Lemma \ref{Braclear}
shows that ${\rm spt}\,\mu\cap (B_r(x)\times [0,t'])=\emptyset$. Since 
$B_{2r}(x)\subset \mathbb R^n\setminus \Gamma_0$, there exists some $i(x,0)\in 
\{1,\ldots,N\}$ such that $B_{2r}(x)\subset E_{0,i(x,0)}$, thus $B_{2r}(x)\subset 
E_{j_l,i(x,0)}(0)$ for all $l$. Then, by Lemma \ref{ya13}, for some $r'\in (0,r/2)$, we
have $\lim_{l\rightarrow\infty}\mathcal L^{n+1}(B_{r'}(x)\setminus E_{j_l,i(x,0)}(\tilde t))=0$
for all $\tilde t\in (0,(r')^2)$. Similarly, for $t>0$, using Lemma \ref{ya7}, there
exist  $i(x,t)\in \{1,\ldots,N\}$ and $r'\in (0,r/2)$ such that $\lim_{l\rightarrow\infty}
\mathcal L^{n+1}(B_{r'}(x)\setminus  E_{j_l,i(x,t)}(\tilde t))=0$ for all $\tilde t\in (t-(r')^2,
t+(r')^2)$. By the connectedness of $S$, $i(x,t)$ has to be all equal to some $i\in \{1,\ldots,N\}$ on $S$. This also shows that $\chi_{E_{j_l,i}}(t)$ converges to $1$ 
in $L^1$ locally on $\{x : (x,t)\in S\}$ for all $t$. Now, for each $i\in \{1,\ldots,N\}$, 
define $S(i)$ to be the union of all connected component with this property.
Since $E_{0,i}=\{x : (x,0)\in S(i)\}$, each $S(i)$ is 
nonempty. They are open disjoint sets and $\cup_{i=1}^N S(i)=(\mathbb R^{n+1}\times \mathbb R^+)\setminus ({\rm spt}\,\mu \cup ({\rm spt}\,\|V_0\|\times \{0\}))$. Define 
$E_i(t):=\{x : (x,t)\in S(i)\}$. Then it is clear that
$\chi_{E_{j_l,i}(t)}$ locally converges to $\chi_{E_i(t)}$ in $L^1$. Up to this
point, the claims (1)-(5) of Theorem \ref{sthm3} are proved, in particular, (4) follows
from the lower semicontinuity of BV norm. 

To prove (6), let $i=\{1,\ldots, N\}$ and 
$R>0$ be fixed. Without loss of generality, we may assume $x=0$. 
Consider $U_R\cap E_i(t)$ which is open. For $r>0$, set
$A_r:=\{x\in U_{R-r}\cap E_i(t) : {\rm dist}\,(\partial(U_R\cap E_i(t)),x)
<r\}$. Consider a family of closed balls $\{{B_{2r}(x)} : 
x\in A_r\}$. By Vitali's covering theorem, we may choose 
points $x_1,\ldots,x_m\in A_r$ such that $\{B_{2r}(x_j)\}_{j=1}^m$
are mutually disjoint and $A_r\subset\cup_{j=1}^m B_{10r}
(x_j)$. By the definition of $A_r$, there exist $\tilde x_j\in 
U_{r}(x_j)\cap \partial (E_i(t))$ for each $j=1,\ldots,m$. Since $(\partial
(E_i(t))\times\{t\})\subset {\rm spt}\,\mu$, by Lemma \ref{cmd9},
$\|V_{t-\Cr{c_7}r^2}\|(B_r(\tilde x_j))
\geq \Cr{c_6}r^{n}$ for $0<r<r_0$ (with a suitable $\lambda$ chosen). 
Since $B_r(\tilde x_j)\subset 
B_{2r}(x_j)$, $\{B_r(\tilde x_j)\}_{j=1}^m$ are mutually disjoint.
Thus we have 
\begin{equation}
\Cr{c_6}mr^{n}\leq \sum_{j=1}^m \|V_{t-\Cr{c_7} r^2}\|(B_r(
\tilde x_j))=\|V_{t-\Cr{c_7}r^2}\|(\cup_{j=1}^m B_r(
\tilde x_j))\leq \|V_{t-\Cr{c_7}r^2}\|(U_{R+r}).
\end{equation}
On the other hand, 
\begin{equation}
\mathcal L^{n+1}(A_r)\leq m\omega_{n+1} (10r)^{n+1} 
\leq (\Cr{c_6}^{-1} \omega_{n+1} 10^{n+1}  \|V_{t-\Cr{c_7}r^2}\|(U_{R+r}))r.
\end{equation}
For any $x\in (U_{R-r}\cap E_i(t))\setminus A_r$, $U_r(x)\subset E_i(t)$ and 
$\|V_t\|(U_r(x))=0$. Thus by Lemma \ref{Braclear}, 
there exists
$\Cl[c]{c_8}>0$ depending only on $n$ such that 
$B_{r/2}(x)\subset E_i(\tilde t)$ for all $\tilde t \in [t,t+\Cr{c_8} r^2]$.
This means $(U_{R-r}\cap E_i(t))\setminus A_r\subset 
E_i(\tilde t)$ for all $\tilde t\in[t,t+\Cr{c_8} r^2]$. Thus, for such $\tilde t$, 
\begin{equation}
\begin{split}
\mathcal L^{n+1}(U_{R}\cap E_i(t) \setminus E_i(\tilde t))
&\leq \mathcal L^{n+1} ((U_{R}\setminus U_{R-r})\cup A_r)\\
&\leq( (n+1)\omega_{n+1} R^{n}+\Cr{c_6}^{-1} \omega_{n+1} 10^{n+1}  \|V_{t-\Cr{c_7}
r^2}\|(U_{R+r}))r \\
& =: \Cl[c]{c_9}(r) r,
\end{split}
\label{dcon1}
\end{equation}
where $\Cr{c_9}$ is uniformly bounded for small $r$. The estimate \eqref{dcon1} holds for
any $i$ with the same $\Cr{c_9}$. 
$\{E_i(t)\cap U_R\}_{i=1}^N$ is mutually disjoint and 
the union has full $\mathcal L^{n+1}$ measure of $U_R$, and 
so is $\{E_i(\tilde t)\cap U_R\}_{i=1}^N$. Thus, except for a
$\mathcal L^{n+1}$ zero measure set, we have
$E_i(\tilde t )\cap U_R\setminus E_i(t)\subset U_R\cap \cup_{i'\neq
i} E_{i'}(t)\setminus E_{i'}(\tilde t)$. Thus
\begin{equation}
\mathcal L^{n+1} (U_R\cap E_i(\tilde t)\setminus E_i(t))\leq \sum_{i'\neq i}
\mathcal L^{n+1} (U_R\cap E_{i'}(t)\setminus E_{i'}(\tilde t))
\leq (N-1)\Cr{c_9} r.
\label{dcon2}
\end{equation}
\eqref{dcon1} and \eqref{dcon2} prove that
\begin{equation}
\mathcal L^{n+1} (U_R\cap (E_i(t)\triangle E_i(\tilde t )))\leq N\Cr{c_9}r
\end{equation}
for $\tilde t\in [t,t+\Cr{c_8} r^2]$ and $r<r_0$. We may exchange the role of 
$t$ and $\tilde t$ to obtain the similar estimate for $\tilde t<t$. Once this is 
obtained, local $\frac12$-H\"{o}lder continuity for $g$ as defined in (6) follows
for $t>0$ using $(A\triangle B)\triangle (A\triangle C)=B\triangle C$
for any sets $A,B,C$. For $t=0$, 
we cannot estimate as above, but we may still prove 
continuity using Lemma \ref{Braclear}. If we assume an extra property on 
$\mathcal E_0=\{E_{0,i}\}_{i=1}^N$, such as, for each $i=1,\ldots,N$ and $R>0$,
$\mathcal L^{n+1}(\{x\in B_{R-r} \cap E_{0,i} : {\rm dist}\,(x,\partial
E_{0,i})<r)\leq c(R) r$ for all sufficiently small $r$, then we can 
proceed just like above and prove $\frac12$-H\"{o}lder
continuity of $g$ up to $t=0$. 
\hfill{$\Box$}  

\section{Additional comments}
\label{question}

\subsection{Tangent flow}
For Brakke flow $\{V_t\}_{t\in \mathbb R^+}$, at each point $(x,t)$ in space-time, $t>0$, 
there exists a tangent flow (see \cite{Ilmanenp,White0} for the definition and proofs)
which is again a Brakke flow and which tells the local behavior of the flow at that point. 
Just like tangent cones of minimal surfaces, tangent flows have a certain homogeneous
property and one can stratify the singularity depending on the dimensions of 
the homogeneity. In this regard, due to the minimizing step in the construction of
approximate solutions, one may wonder if some extra property of tangent flow 
may be derived.
As far as the approximate solutions are concerned, as indicated in Section \ref{exfe}, 
unstable singularities are likely to break up into more stable ones by Lipschitz deformation. 
There should be some aspects on tangent flow which are affected by the choice of 
$f_1\in {\bf E}(\mathcal E_{j,l}, j)$ in \eqref{ap1} 
as elaborated in Remark \ref{choicef}. 
It is a challenging problem 
to analyze this finer point of the Brakke flow obtained in this paper.

\subsection{A short-time regularity}
\label{astr}
Suppose in addition that $\Gamma_0$ satisfies the following density ratio upper bound
condition. There exist some $\nu\in (0,1)$ and $r_0\in (0,\infty)$
such that $\mathcal H^n(\Gamma_0\cap B_r(x))\leq (2-\nu)\omega_n
r^n$ for all $r\in (0,r_0)$ and $x\in\mathbb R^{n+1}$. Nontrivial examples
with singularities
satisfying such condition are suitably regular 1-dimensional networks 
with finite number of triple junctions, 
since such junctions have density $\frac32$.
Others are  suitably regular 2-dimensional ``soap bubble clusters'' with 
singularities of three surfaces with boundaries 
meeting along a curve, or 6 surfaces with boundaries
meeting at a point and 4 curves. They can have 
densities strictly less than 2. These are interesting classes of examples which 
are also physically relevant.
Under this condition, by using Lemma \ref{cmd4}, one can prove that there exists
$T>0$ such that $\theta^{n}(\|V_t\|,x)=1$ for $\|V_t\|$ almost all $x\in\mathbb R^{n+1}$ and
for almost all $t\in (0,T)$. In other words, there cannot be any points
of integer density greater than or equal to 2.
Thus the solution of the present paper is guaranteed to remain
unit density flow for $t\in (0,T)$. Then Theorem \ref{ktthm} applies and
${\rm spt}\,\mu$ is partially regular as described there for $(0,T)$. 
In the case of $n=1$, this implies further 
that any nontrivial static tangent flow within 
the time interval $(0,T)$ is either a line, or a regular triple junction, both of single-multiplicity. This is precisely
the situation that we may 
apply \cite[Theorem 2.2]{Wick}. The result concludes that
there exists a closed set $S\subset\mathbb R^2\times[0,T)$ of parabolic Hausdorff
dimension at most 1 such that, outside of $S$, ${\rm spt}\,\|V_t\|$ is locally a smooth
curve or a regular triple junction of 120 degree angle moving smoothly by the mean curvature.
We mention that the short-time existence of one-dimensional network flow is 
recently obtained in \cite{Neves}. We allow more general $\Gamma_0$ than \cite{Neves}
but our flow may have singularities of small dimension in general.
Due to the minimizing 
step of the approximate solution, it is likely in the one-dimensional case that any static tangent
flow constructed in this paper is either a line or a regular triple junction
even for later time. This should require a
finer look into the singularities and pose an interesting open question. In any case,
away from space-time region with higher integer multiplicities ($\geq 2$), Brakke flow constructed in this paper is partially
regular as in Theorem \ref{ktthm}. Higher integer 
multiplicities pose outstanding 
regularity questions even for stationary integral varifolds. 

 We also mention 
that there is an initial time regularity property for regular points of $\Gamma_0$ for any $n$
in the following sense.  If $\Gamma_0$ is locally a $C^1$
hypersurface at a point $x$ which is not an interior boundary point of some $E_{0,i}$
(i.e., there exist $i,i'\in \{1,\ldots,N\}$, $i\neq i'$, such that $x\in \partial E_{0,i}
\cap \partial E_{0,i'}$), then 
there exists a space-time neighborhood of $(x,0)$ in 
which the constructed flow is $C^1$ in the parabolic sense up to $t=0$ and 
$C^{\infty}$ for $t>0$. This can be proved by using a $C^{1,\alpha}$ regularity theorem in 
\cite{KT} as demonstrated in \cite[Theorem 2.3(4)]{Takasao} for a phase 
field setting. 

\subsection{Other setting}
If we replace $\mathbb R^{n+1}$ by the flat torus $\mathbb T^{n+1}$, we may simply change
everything by setting quantities periodic on $\mathbb R^{n+1}$ with period 1. We would 
have finite open partitions defined on $\mathbb T^{n+1}$ and all convergence takes place 
accordingly. For general Riemannian manifolds, by adapting definitions and assumptions, similar results should follow with little change.
All the key points of the paper such as the proofs of rectifiability and integrality are local estimates. On
the other hand, if one is interested in the MCF with ``Dirichlet condition'' or ``Neumann
condition'' in a suitable sense, the presence of such boundary condition may pose a nontrivial problem near the 
boundary and further studies are expected. From a geometric point of view in connection
with the Plateau problem, such problem is natural and interesting. As a related 
matter, one aspect that
may puzzle the reader is the finiteness of open partition, i.e., we always fix $N$
of $\mathcal{OP}_{\Omega}^N$
even though we do not see any quantitative statement in the main results concerning
$N$. One may naturally wonder if countably infinite open partition 
$\mathcal{OP}_{\Omega}^{\infty}$ can be allowed. In fact, $N=\infty$ can be dealt with
all the way just before the last step of taking $j_l\rightarrow\infty$. For example, in 
Lemma \ref{ya7}, we want to conclude that a subsequence of 
$\chi_{E_{j_l,k}(t)}$ converges in $L^1_{loc}(\mathbb R^{n+1})$ to some $\chi_{E_k(t)}$
and $\sum_{k=1}^N\chi_{E_k(t)}\equiv 1$ a.e$.$ on $\mathbb R^{n+1}$. 
However, if $N=\infty$, we need to exclude a
possibility that $\sum_{k=1}^{\infty} \chi_{E_k(t)}<1$ on a positive measure set.  This is 
because, even though 
$\sum_{k=1}^{\infty} \chi_{E_{j_l,k}(t)}\equiv 1$ for all $j_l$, 
if there are infinite number of sets, the fear is that
all of them become finer and finer as $j_l$ increases 
and the limit may all vanish. This scenario seems 
unlikely to happen for a.e$.$ $t$, but there has to be some extra argument to 
eliminate such possibility. Since the finite $N$ case is interesting enough, we did not 
pursue $N=\infty$ for the technicality. It is also possible to first find Brakke flow for each $N$ and take a limit $N\rightarrow \infty$. One can argue that there exists a converging
subsequence whose limit is also a Brakke flow as described in the present paper 
and that the limit is nontrivial using the continuity property of the ``grains''. 


\end{document}

%% file: fig_1.eps_tex
\begingroup%
  \makeatletter%
  \providecommand\color[2][]{%
    \errmessage{(Inkscape) Color is used for the text in Inkscape, but the package 'color.sty' is not loaded}%
    \renewcommand\color[2][]{}%
  }%
  \providecommand\transparent[1]{%
    \errmessage{(Inkscape) Transparency is used (non-zero) for the text in Inkscape, but the package 'transparent.sty' is not loaded}%
    \renewcommand\transparent[1]{}%
  }%
  \providecommand\rotatebox[2]{#2}%
  \ifx\svgwidth\undefined%
    \setlength{\unitlength}{774.85102696bp}%
    \ifx\svgscale\undefined%
      \relax%
    \else%
      \setlength{\unitlength}{\unitlength * \real{\svgscale}}%
    \fi%
  \else%
    \setlength{\unitlength}{\svgwidth}%
  \fi%
  \global\let\svgwidth\undefined%
  \global\let\svgscale\undefined%
  \makeatother%
  \begin{picture}(1,0.2274952)%
    \put(0,0){\includegraphics[width=\unitlength]{fig_1.eps}}%
    \put(0.1260709,0.20246168){\color[rgb]{0,0,0}\makebox(0,0)[lt]{\begin{minipage}{0.05506362\unitlength}\raggedright \end{minipage}}}%
    \put(0.1260709,0.20246168){\color[rgb]{0,0,0}\makebox(0,0)[lt]{\begin{minipage}{0.05506362\unitlength}\raggedright \end{minipage}}}%
    \put(0.25091383,0.14803053){\color[rgb]{0,0,0}\makebox(0,0)[lt]{\begin{minipage}{0.06661552\unitlength}\raggedright $E_3$\end{minipage}}}%
    \put(0.1427444,0.08529879){\color[rgb]{0,0,0}\makebox(0,0)[lt]{\begin{minipage}{0.0715393\unitlength}\raggedright $E_4$\end{minipage}}}%
    \put(0.14262203,0.21023204){\color[rgb]{0,0,0}\makebox(0,0)[lt]{\begin{minipage}{0.0613371\unitlength}\raggedright $E_2$\end{minipage}}}%
    \put(0.154371,0.20124292){\color[rgb]{0,0,0}\makebox(0,0)[lt]{\begin{minipage}{0.00662102\unitlength}\raggedright \end{minipage}}}%
    \put(0.13935288,0.03153278){\color[rgb]{0,0,0}\makebox(0,0)[lt]{\begin{minipage}{0.06281331\unitlength}\raggedright (a)\end{minipage}}}%
    \put(0.44001596,0.20246168){\color[rgb]{0,0,0}\makebox(0,0)[lt]{\begin{minipage}{0.05506362\unitlength}\raggedright \end{minipage}}}%
    \put(0.35297615,0.1471332){\color[rgb]{0,0,0}\makebox(0,0)[lt]{\begin{minipage}{0.06160411\unitlength}\raggedright $E_1$\end{minipage}}}%
    \put(0.44001593,0.20246168){\color[rgb]{0,0,0}\makebox(0,0)[lt]{\begin{minipage}{0.05506361\unitlength}\raggedright \end{minipage}}}%
    \put(0.56898867,0.14803053){\color[rgb]{0,0,0}\makebox(0,0)[lt]{\begin{minipage}{0.07273198\unitlength}\raggedright $E_3$\end{minipage}}}%
    \put(0.45875434,0.08529879){\color[rgb]{0,0,0}\makebox(0,0)[lt]{\begin{minipage}{0.06733121\unitlength}\raggedright $E_4$\end{minipage}}}%
    \put(0.458632,0.21023204){\color[rgb]{0,0,0}\makebox(0,0)[lt]{\begin{minipage}{0.06745354\unitlength}\raggedright $E_2$\end{minipage}}}%
    \put(0.45548518,0.03153342){\color[rgb]{0,0,0}\makebox(0,0)[lt]{\begin{minipage}{0.06281335\unitlength}\raggedright (b)\end{minipage}}}%
    \put(0.75382586,0.19997316){\color[rgb]{0,0,0}\makebox(0,0)[lt]{\begin{minipage}{0.05506361\unitlength}\raggedright \end{minipage}}}%
    \put(0.66678603,0.15083942){\color[rgb]{0,0,0}\makebox(0,0)[lt]{\begin{minipage}{0.11947855\unitlength}\raggedright $\tilde E_1$\end{minipage}}}%
    \put(0.75382587,0.19997316){\color[rgb]{0,0,0}\makebox(0,0)[lt]{\begin{minipage}{0.05506362\unitlength}\raggedright \end{minipage}}}%
    \put(0.88279858,0.14967184){\color[rgb]{0,0,0}\makebox(0,0)[lt]{\begin{minipage}{0.12028186\unitlength}\raggedright $\tilde E_3$\end{minipage}}}%
    \put(0.77462915,0.09106992){\color[rgb]{0,0,0}\makebox(0,0)[lt]{\begin{minipage}{0.12520565\unitlength}\raggedright $\tilde E_4$\end{minipage}}}%
    \put(0.77450681,0.20980843){\color[rgb]{0,0,0}\makebox(0,0)[lt]{\begin{minipage}{0.12532799\unitlength}\raggedright $\tilde E_2$\end{minipage}}}%
    \put(0.67950917,0.14272879){\color[rgb]{0,0,0}\makebox(0,0)[lb]{\smash{}}}%
    \put(0.77161599,0.03359763){\color[rgb]{0,0,0}\makebox(0,0)[lt]{\begin{minipage}{0.06281335\unitlength}\raggedright (c)\end{minipage}}}%
    \put(0.03497959,0.1471332){\color[rgb]{0,0,0}\makebox(0,0)[lt]{\begin{minipage}{0.06160411\unitlength}\raggedright $E_1$\end{minipage}}}%
  \end{picture}%
\endgroup%

%% file: fig_2.eps_tex
\begingroup%
  \makeatletter%
  \providecommand\color[2][]{%
    \errmessage{(Inkscape) Color is used for the text in Inkscape, but the package 'color.sty' is not loaded}%
    \renewcommand\color[2][]{}%
  }%
  \providecommand\transparent[1]{%
    \errmessage{(Inkscape) Transparency is used (non-zero) for the text in Inkscape, but the package 'transparent.sty' is not loaded}%
    \renewcommand\transparent[1]{}%
  }%
  \providecommand\rotatebox[2]{#2}%
  \ifx\svgwidth\undefined%
    \setlength{\unitlength}{719.17948146bp}%
    \ifx\svgscale\undefined%
      \relax%
    \else%
      \setlength{\unitlength}{\unitlength * \real{\svgscale}}%
    \fi%
  \else%
    \setlength{\unitlength}{\svgwidth}%
  \fi%
  \global\let\svgwidth\undefined%
  \global\let\svgscale\undefined%
  \makeatother%
  \begin{picture}(1,0.24472335)%
    \put(0,0){\includegraphics[width=\unitlength]{fig_2.eps}}%
    \put(0.14837585,0.03423595){\color[rgb]{0,0,0}\makebox(0,0)[lt]{\begin{minipage}{0.06804814\unitlength}\raggedright (a)\end{minipage}}}%
    \put(0.48416793,0.03423664){\color[rgb]{0,0,0}\makebox(0,0)[lt]{\begin{minipage}{0.06804813\unitlength}\raggedright (b)\end{minipage}}}%
    \put(0.82443021,0.03423592){\color[rgb]{0,0,0}\makebox(0,0)[lt]{\begin{minipage}{0.06804813\unitlength}\raggedright (c)\end{minipage}}}%
    \put(0.12789882,0.2171699){\color[rgb]{0,0,0}\makebox(0,0)[lt]{\begin{minipage}{0.05965258\unitlength}\raggedright \end{minipage}}}%
    \put(0.12789882,0.2171699){\color[rgb]{0,0,0}\makebox(0,0)[lt]{\begin{minipage}{0.05965258\unitlength}\raggedright \end{minipage}}}%
    \put(0.15263615,0.09150984){\color[rgb]{0,0,0}\makebox(0,0)[lt]{\begin{minipage}{0.07653247\unitlength}\raggedright $E_1$\end{minipage}}}%
    \put(0.15250359,0.22550381){\color[rgb]{0,0,0}\makebox(0,0)[lt]{\begin{minipage}{0.07666502\unitlength}\raggedright $E_1$\end{minipage}}}%
    \put(0.46788299,0.21769889){\color[rgb]{0,0,0}\makebox(0,0)[lt]{\begin{minipage}{0.05965262\unitlength}\raggedright \end{minipage}}}%
    \put(0.46788296,0.21769888){\color[rgb]{0,0,0}\makebox(0,0)[lt]{\begin{minipage}{0.05965258\unitlength}\raggedright \end{minipage}}}%
    \put(0.48817074,0.09203882){\color[rgb]{0,0,0}\makebox(0,0)[lt]{\begin{minipage}{0.07362697\unitlength}\raggedright $E_1$\end{minipage}}}%
    \put(0.48803823,0.22603279){\color[rgb]{0,0,0}\makebox(0,0)[lt]{\begin{minipage}{0.0737595\unitlength}\raggedright $E_1$\end{minipage}}}%
    \put(0.80986433,0.21562946){\color[rgb]{0,0,0}\makebox(0,0)[lt]{\begin{minipage}{0.05965258\unitlength}\raggedright \end{minipage}}}%
    \put(0.80986433,0.21562946){\color[rgb]{0,0,0}\makebox(0,0)[lt]{\begin{minipage}{0.05965256\unitlength}\raggedright \end{minipage}}}%
    \put(0.82748861,0.22413305){\color[rgb]{0,0,0}\makebox(0,0)[lt]{\begin{minipage}{0.13657876\unitlength}\raggedright $\tilde E_1$\end{minipage}}}%
    \put(0.13579149,0.21668426){\color[rgb]{0,0,0}\makebox(0,0)[lt]{\begin{minipage}{0.05932609\unitlength}\raggedright \end{minipage}}}%
    \put(0.13579149,0.21668426){\color[rgb]{0,0,0}\makebox(0,0)[lt]{\begin{minipage}{0.05932609\unitlength}\raggedright \end{minipage}}}%
    \put(0.1662823,0.21537117){\color[rgb]{0,0,0}\makebox(0,0)[lt]{\begin{minipage}{0.00713355\unitlength}\raggedright \end{minipage}}}%
    \put(0.47403898,0.21668426){\color[rgb]{0,0,0}\makebox(0,0)[lt]{\begin{minipage}{0.05932609\unitlength}\raggedright \end{minipage}}}%
    \put(0.47403895,0.21668427){\color[rgb]{0,0,0}\makebox(0,0)[lt]{\begin{minipage}{0.05932608\unitlength}\raggedright \end{minipage}}}%
    \put(0.81214084,0.21400311){\color[rgb]{0,0,0}\makebox(0,0)[lt]{\begin{minipage}{0.05932608\unitlength}\raggedright \end{minipage}}}%
    \put(0.81214086,0.21400311){\color[rgb]{0,0,0}\makebox(0,0)[lt]{\begin{minipage}{0.05932609\unitlength}\raggedright \end{minipage}}}%
    \put(0.73207131,0.15232746){\color[rgb]{0,0,0}\makebox(0,0)[lb]{\smash{}}}%
  \end{picture}%
\endgroup%

%% file: fig_5.eps_tex
\begingroup%
  \makeatletter%
  \providecommand\color[2][]{%
    \errmessage{(Inkscape) Color is used for the text in Inkscape, but the package 'color.sty' is not loaded}%
    \renewcommand\color[2][]{}%
  }%
  \providecommand\transparent[1]{%
    \errmessage{(Inkscape) Transparency is used (non-zero) for the text in Inkscape, but the package 'transparent.sty' is not loaded}%
    \renewcommand\transparent[1]{}%
  }%
  \providecommand\rotatebox[2]{#2}%
  \ifx\svgwidth\undefined%
    \setlength{\unitlength}{713.8823412bp}%
    \ifx\svgscale\undefined%
      \relax%
    \else%
      \setlength{\unitlength}{\unitlength * \real{\svgscale}}%
    \fi%
  \else%
    \setlength{\unitlength}{\svgwidth}%
  \fi%
  \global\let\svgwidth\undefined%
  \global\let\svgscale\undefined%
  \makeatother%
  \begin{picture}(1,0.38035597)%
    \put(0,0){\includegraphics[width=\unitlength]{fig_5.eps}}%
    \put(0.75481205,0.36391474){\color[rgb]{0,0,0}\makebox(0,0)[lt]{\begin{minipage}{0.06013206\unitlength}\raggedright \end{minipage}}}%
    \put(0.66808094,0.29367691){\color[rgb]{0,0,0}\makebox(0,0)[lt]{\begin{minipage}{0.00723047\unitlength}\raggedright \end{minipage}}}%
    \put(0.75481205,0.36391474){\color[rgb]{0,0,0}\makebox(0,0)[lt]{\begin{minipage}{0.06013206\unitlength}\raggedright \end{minipage}}}%
    \put(0.61082638,0.17619997){\color[rgb]{0,0,0}\makebox(0,0)[lt]{\begin{minipage}{0.11217941\unitlength}\raggedright $\tilde E_2$\end{minipage}}}%
    \put(0.72384033,0.25443615){\color[rgb]{0,0,0}\makebox(0,0)[lt]{\begin{minipage}{0.10898749\unitlength}\raggedright $\tilde E_1$\end{minipage}}}%
    \put(0.21380384,0.36412227){\color[rgb]{0,0,0}\makebox(0,0)[lt]{\begin{minipage}{0.06013206\unitlength}\raggedright \end{minipage}}}%
    \put(0.12707274,0.29388443){\color[rgb]{0,0,0}\makebox(0,0)[lt]{\begin{minipage}{0.00723047\unitlength}\raggedright \end{minipage}}}%
    \put(0.21380384,0.36412227){\color[rgb]{0,0,0}\makebox(0,0)[lt]{\begin{minipage}{0.06013206\unitlength}\raggedright \end{minipage}}}%
    \put(0.25005698,0.09484557){\color[rgb]{0,0,0}\makebox(0,0)[lt]{\begin{minipage}{0.06767749\unitlength}\raggedright $E_1$\end{minipage}}}%
    \put(0.25005999,0.2504788){\color[rgb]{0,0,0}\makebox(0,0)[lt]{\begin{minipage}{0.06839748\unitlength}\raggedright $E_1$\end{minipage}}}%
    \put(0.24735564,0.03877589){\color[rgb]{0,0,0}\makebox(0,0)[lt]{\begin{minipage}{0.0685951\unitlength}\raggedright (a)\end{minipage}}}%
    \put(0.71944117,0.03877661){\color[rgb]{0,0,0}\makebox(0,0)[lt]{\begin{minipage}{0.0685951\unitlength}\raggedright (b)\end{minipage}}}%
    \put(0.32854317,0.2622281){\color[rgb]{0,0,0}\makebox(0,0)[lb]{\smash{}}}%
    \put(0.83047028,0.17619997){\color[rgb]{0,0,0}\makebox(0,0)[lt]{\begin{minipage}{0.11217941\unitlength}\raggedright $\tilde E_3$\end{minipage}}}%
    \put(0.14042358,0.17358232){\color[rgb]{0,0,0}\makebox(0,0)[lt]{\begin{minipage}{0.06148808\unitlength}\raggedright $E_2$\end{minipage}}}%
    \put(0.36006759,0.17358232){\color[rgb]{0,0,0}\makebox(0,0)[lt]{\begin{minipage}{0.06148808\unitlength}\raggedright $E_3$\end{minipage}}}%
  \end{picture}%
\endgroup%

%% file: fig_3.eps_tex
\begingroup%
  \makeatletter%
  \providecommand\color[2][]{%
    \errmessage{(Inkscape) Color is used for the text in Inkscape, but the package 'color.sty' is not loaded}%
    \renewcommand\color[2][]{}%
  }%
  \providecommand\transparent[1]{%
    \errmessage{(Inkscape) Transparency is used (non-zero) for the text in Inkscape, but the package 'transparent.sty' is not loaded}%
    \renewcommand\transparent[1]{}%
  }%
  \providecommand\rotatebox[2]{#2}%
  \ifx\svgwidth\undefined%
    \setlength{\unitlength}{644.47059076bp}%
    \ifx\svgscale\undefined%
      \relax%
    \else%
      \setlength{\unitlength}{\unitlength * \real{\svgscale}}%
    \fi%
  \else%
    \setlength{\unitlength}{\svgwidth}%
  \fi%
  \global\let\svgwidth\undefined%
  \global\let\svgscale\undefined%
  \makeatother%
  \begin{picture}(1,0.27309239)%
    \put(0,0){\includegraphics[width=\unitlength]{fig_3.eps}}%
    \put(0.22291158,0.03702551){\color[rgb]{0,0,0}\makebox(0,0)[lt]{\begin{minipage}{0.07425316\unitlength}\raggedright (a)\end{minipage}}}%
    \put(0.71763122,0.03702625){\color[rgb]{0,0,0}\makebox(0,0)[lt]{\begin{minipage}{0.07425316\unitlength}\raggedright (b)\end{minipage}}}%
    \put(0.14792409,0.23897742){\color[rgb]{0,0,0}\makebox(0,0)[lt]{\begin{minipage}{0.06509205\unitlength}\raggedright \end{minipage}}}%
    \put(0.09476874,0.17357228){\color[rgb]{0,0,0}\makebox(0,0)[lt]{\begin{minipage}{0.08398262\unitlength}\raggedright $E_1$\end{minipage}}}%
    \put(0.14792411,0.23897741){\color[rgb]{0,0,0}\makebox(0,0)[lt]{\begin{minipage}{0.06509205\unitlength}\raggedright \end{minipage}}}%
    \put(0.35012245,0.17463304){\color[rgb]{0,0,0}\makebox(0,0)[lt]{\begin{minipage}{0.07689471\unitlength}\raggedright $E_1$\end{minipage}}}%
    \put(0.22721804,0.10047632){\color[rgb]{0,0,0}\makebox(0,0)[lt]{\begin{minipage}{0.08063152\unitlength}\raggedright $E_2$\end{minipage}}}%
    \put(0.22707339,0.24816296){\color[rgb]{0,0,0}\makebox(0,0)[lt]{\begin{minipage}{0.08077617\unitlength}\raggedright $E_2$\end{minipage}}}%
    \put(0.74068739,0.23920378){\color[rgb]{0,0,0}\makebox(0,0)[lt]{\begin{minipage}{0.06509205\unitlength}\raggedright \end{minipage}}}%
    \put(0.59319099,0.17379866){\color[rgb]{0,0,0}\makebox(0,0)[lt]{\begin{minipage}{0.14912367\unitlength}\raggedright $\tilde E_1$\end{minipage}}}%
    \put(0.74068739,0.23920378){\color[rgb]{0,0,0}\makebox(0,0)[lt]{\begin{minipage}{0.06509205\unitlength}\raggedright \end{minipage}}}%
    \put(0.84854471,0.17485942){\color[rgb]{0,0,0}\makebox(0,0)[lt]{\begin{minipage}{0.14203575\unitlength}\raggedright $\tilde E_1$\end{minipage}}}%
    \put(0.72053032,0.24838933){\color[rgb]{0,0,0}\makebox(0,0)[lt]{\begin{minipage}{0.14591724\unitlength}\raggedright $\tilde E_2$\end{minipage}}}%
    \put(0.65283583,0.17153381){\color[rgb]{0,0,0}\makebox(0,0)[lb]{\smash{}}}%
  \end{picture}%
\endgroup%

%% file: fig_4.eps_tex
\begingroup%
  \makeatletter%
  \providecommand\color[2][]{%
    \errmessage{(Inkscape) Color is used for the text in Inkscape, but the package 'color.sty' is not loaded}%
    \renewcommand\color[2][]{}%
  }%
  \providecommand\transparent[1]{%
    \errmessage{(Inkscape) Transparency is used (non-zero) for the text in Inkscape, but the package 'transparent.sty' is not loaded}%
    \renewcommand\transparent[1]{}%
  }%
  \providecommand\rotatebox[2]{#2}%
  \ifx\svgwidth\undefined%
    \setlength{\unitlength}{645.49526049bp}%
    \ifx\svgscale\undefined%
      \relax%
    \else%
      \setlength{\unitlength}{\unitlength * \real{\svgscale}}%
    \fi%
  \else%
    \setlength{\unitlength}{\svgwidth}%
  \fi%
  \global\let\svgwidth\undefined%
  \global\let\svgscale\undefined%
  \makeatother%
  \begin{picture}(1,0.42265181)%
    \put(0,0){\includegraphics[width=\unitlength]{fig_4.eps}}%
    \put(0.75634131,0.37272516){\color[rgb]{0,0,0}\makebox(0,0)[lt]{\begin{minipage}{0.06650276\unitlength}\raggedright \end{minipage}}}%
    \put(0.66042146,0.29504597){\color[rgb]{0,0,0}\makebox(0,0)[lt]{\begin{minipage}{0.0079965\unitlength}\raggedright \end{minipage}}}%
    \put(0.75634131,0.37272516){\color[rgb]{0,0,0}\makebox(0,0)[lt]{\begin{minipage}{0.06650276\unitlength}\raggedright \end{minipage}}}%
    \put(0.74898727,0.18430262){\color[rgb]{0,0,0}\makebox(0,0)[lt]{\begin{minipage}{0.13577361\unitlength}\raggedright $\tilde E_2$\end{minipage}}}%
    \put(0.61445224,0.12918152){\color[rgb]{0,0,0}\makebox(0,0)[lt]{\begin{minipage}{0.12406427\unitlength}\raggedright $\tilde E_2$\end{minipage}}}%
    \put(0.69234383,0.34583907){\color[rgb]{0,0,0}\makebox(0,0)[lt]{\begin{minipage}{0.12053418\unitlength}\raggedright $\tilde E_1$\end{minipage}}}%
    \put(0.15801591,0.37295467){\color[rgb]{0,0,0}\makebox(0,0)[lt]{\begin{minipage}{0.06650276\unitlength}\raggedright \end{minipage}}}%
    \put(0.06209607,0.29527548){\color[rgb]{0,0,0}\makebox(0,0)[lt]{\begin{minipage}{0.0079965\unitlength}\raggedright \end{minipage}}}%
    \put(0.15801591,0.37295467){\color[rgb]{0,0,0}\makebox(0,0)[lt]{\begin{minipage}{0.06650276\unitlength}\raggedright \end{minipage}}}%
    \put(0.22909387,0.18421296){\color[rgb]{0,0,0}\makebox(0,0)[lt]{\begin{minipage}{0.07484759\unitlength}\raggedright $E_1$\end{minipage}}}%
    \put(0.19687389,0.34022309){\color[rgb]{0,0,0}\makebox(0,0)[lt]{\begin{minipage}{0.07564386\unitlength}\raggedright $E_1$\end{minipage}}}%
    \put(0.28126129,0.13001142){\color[rgb]{0,0,0}\makebox(0,0)[lt]{\begin{minipage}{0.07553282\unitlength}\raggedright $E_2$\end{minipage}}}%
    \put(0.11776021,0.14983435){\color[rgb]{0,0,0}\makebox(0,0)[lt]{\begin{minipage}{0.06800244\unitlength}\raggedright $E_2$\end{minipage}}}%
    \put(0.3052207,0.1598647){\color[rgb]{0,0,0}\makebox(0,0)[lt]{\begin{minipage}{0.05405213\unitlength}\raggedright $x_0$\end{minipage}}}%
    \put(0.22982437,0.03792659){\color[rgb]{0,0,0}\makebox(0,0)[lt]{\begin{minipage}{0.07586242\unitlength}\raggedright (a)\end{minipage}}}%
    \put(0.71970177,0.03792739){\color[rgb]{0,0,0}\makebox(0,0)[lt]{\begin{minipage}{0.07586242\unitlength}\raggedright (b)\end{minipage}}}%
    \put(0.28491132,0.2602653){\color[rgb]{0,0,0}\makebox(0,0)[lb]{\smash{}}}%
  \end{picture}%
\endgroup%

%% file: deformation_1.eps_tex
\begingroup%
  \makeatletter%
  \providecommand\color[2][]{%
    \errmessage{(Inkscape) Color is used for the text in Inkscape, but the package 'color.sty' is not loaded}%
    \renewcommand\color[2][]{}%
  }%
  \providecommand\transparent[1]{%
    \errmessage{(Inkscape) Transparency is used (non-zero) for the text in Inkscape, but the package 'transparent.sty' is not loaded}%
    \renewcommand\transparent[1]{}%
  }%
  \providecommand\rotatebox[2]{#2}%
  \ifx\svgwidth\undefined%
    \setlength{\unitlength}{1711.75823133bp}%
    \ifx\svgscale\undefined%
      \relax%
    \else%
      \setlength{\unitlength}{\unitlength * \real{\svgscale}}%
    \fi%
  \else%
    \setlength{\unitlength}{\svgwidth}%
  \fi%
  \global\let\svgwidth\undefined%
  \global\let\svgscale\undefined%
  \makeatother%
  \begin{picture}(1,0.42279374)%
    \put(0,0){\includegraphics[width=\unitlength]{deformation_1.eps}}%
    \put(0.29667575,0.21263742){\color[rgb]{0,0,0}\makebox(0,0)[lb]{\smash{}}}%
    \put(0.39250152,0.0406968){\color[rgb]{0,0,0}\makebox(0,0)[lb]{\smash{}}}%
    \put(0.11378885,0.03630299){\color[rgb]{0,0,0}\makebox(0,0)[lb]{\smash{$-r_1$}}}%
    \put(0.34766097,0.03720004){\color[rgb]{0,0,0}\makebox(0,0)[lb]{\smash{$a-\delta$}}}%
    \put(0.69802923,0.03627455){\color[rgb]{0,0,0}\makebox(0,0)[lb]{\smash{$a+\delta$}}}%
    \put(0.54671222,0.0367047){\color[rgb]{0,0,0}\makebox(0,0)[lb]{\smash{$a$}}}%
    \put(0.59796814,0.036669){\color[rgb]{0,0,0}\makebox(0,0)[lb]{\smash{$a^*$}}}%
    \put(0.85197937,0.03756679){\color[rgb]{0,0,0}\makebox(0,0)[lb]{\smash{$r_1$}}}%
    \put(0.11588089,0.40840795){\color[rgb]{0,0,0}\makebox(0,0)[lb]{\smash{$E_1(a)$}}}%
    \put(0.11588089,0.37298667){\color[rgb]{0,0,0}\makebox(0,0)[lb]{\smash{$E_3(a)$}}}%
    \put(0.11588089,0.33653292){\color[rgb]{0,0,0}\makebox(0,0)[lb]{\smash{$E_4(a)$}}}%
    \put(0.06634119,0.29914447){\color[rgb]{0,0,0}\makebox(0,0)[lb]{\smash{$E_2(a)=E_3(a)\cup E_4(a)$}}}%
    \put(0.44873519,0.15526952){\color[rgb]{0,0,0}\makebox(0,0)[lb]{\smash{$\rho _1$}}}%
    \put(0.44873519,0.25528363){\color[rgb]{0,0,0}\makebox(0,0)[lb]{\smash{$\rho _1$}}}%
    \put(0.61986442,0.14031971){\color[rgb]{0,0,0}\makebox(0,0)[lb]{\smash{$\small\xi$}}}%
    \put(0.61986442,0.26463634){\color[rgb]{0,0,0}\makebox(0,0)[lb]{\smash{$\small\xi$}}}%
    \put(0.72240759,0.38203568){\color[rgb]{0,0,0}\makebox(0,0)[lb]{\smash{$\text{slope}= (2\delta)^{-1} \xi$}}}%
    \put(0.80289636,0.32979365){\color[rgb]{0,0,0}\makebox(0,0)[lb]{\smash{$\text{slope}= (2r_1)^{-1} \xi$}}}%
  \end{picture}%
\endgroup%

%% file: grain20160331.bbl
\begin{thebibliography}{99}
\bibitem{Allard}
       W. K. Allard,
       {\em On the first variation of a varifold},
       Ann. of Math. (2) \textbf{95} (1972), 417-491.

\bibitem{Almgren} F. J. Jr. Almgren, {\em Existence and regularity almost everywhere
of solutions to elliptic variational problems with constraints}, Mem. Amer. Math. Soc. {\bf 4}, 
(1976) no.165.

\bibitem{AFP}
L. Ambrosio, N. Fusco, D. Pallara, {\em Functions of bounded variation
and free discontinuity problems}, Oxford Mathematical Monographs. The Clarendon Press, Oxford University Press, New York, 2000.

\bibitem{Andrew}
B. Andrews, P. Bryan, {\em Curvature bound for curve shortening flow via distance comparison and a direct proof of Grayson's theorem}, J. Reine Angew. Math. {\bf 653} (2011), 179-187.

\bibitem{BellettiniT}
G. Bellettini, {\em Lecture notes on mean curvature flow, barriers and singular perturbations}, Lecture Notes. Scuola Normale Superiore di Pisa (New Series), {\bf 12}. Edizioni della Normale, Pisa, 2013.

\bibitem{Bellettini}
G. Bellettini, M. Novaga, {\em Curvature evolution of nonconvex lens-shaped domains},
 J. Reine Angew. Math. {\bf 656} (2011), 17-46.
 
\bibitem{Brakke} 
	 K.~A. Brakke, 
	 {\em The motion of a surface by its mean curvature},
	 Mathematical Notes, vol.~{\bf 20}, Princeton University Press, 1978.
\bibitem{Brakke2}
K.~A. Brakke, {\em The surface evolver}, Experiment. Math. {\bf 1}, (1992), no. 2,
141-165.
\bibitem{Brakke3}
K.~A. Brakke, {\em Grain growth movie}, http://facstaff.susqu.edu/brakke/ 

\bibitem{Bronsard}
L. Bronsard, F. Reitich, 
{\em On three-phase boundary motion and the singular limit of a vector-valued Ginzburg-Landau equation},
Arch. Rational Mech. Anal. {\bf 124} (1993), no. 4, 355-379.

\bibitem{CGG} Y.-G. Chen, Y. Giga, S. Goto, {\em Uniqueness and existence of
viscosity solutions of generalized mean curvature flow equations}, J. 
Differential Geom. {\bf 33} (1991), no. 3, 749-786.

\bibitem{Chen} X. Chen, J.-S. Guo, {\em Motion by curvature of planar curves with end points moving freely on a line}, Math. Ann. {\bf 350} (2011), no. 2, 277-311.
	 
\bibitem{Chen2} X. Chen, J.-S. Guo, {\em Self-similar solutions of a 2-D multiple-phase
curvature flow}, Phys. D {\bf 229} (2007), no.1, 22-34.

\bibitem{Colding} T. H. Colding, W. P. Minicozzi II, {\em Minimal surfaces and mean curvature flow}, Surveys in geometric analysis and relativity, Adv. Lect. Math., {\bf 20}, 73-143, Int. Press, Somerville, MA, 2011.

\bibitem{Ecker}
K. Ecker, {\em Regularity theory for mean curvature flow}, Progress in Nonlinear Differential Equations and their Applications, {\bf 57}. Birkh\"{a}user Boston, Inc., Boston, MA, 2004.

\bibitem{Evans-Gariepy}
L.~C. Evans, R. Gariepy, {\em Measure theory and fine properties of functions,
revised edition}, CRC Press, Boca Raton, FL, 2015.

\bibitem{ES1} L. C. Evans, J. Spruck, {\em Motion of level sets by mean curvature.
I}, J. Differential Geom. {\bf 33} (1991), no. 3, 635-681.

\bibitem{ES2}
L.~C. Evans, J. Spruck, {\em Motion of level sets by mean curvature. IV}, J. Geom. Anal. {\bf 5}
(1995), no. 1, 77-114.

\bibitem{Federer}
Federer, H., \emph{Geometric measure theory}, Springer, 1969.

\bibitem{Freire}
A. Freire,
{\em Mean curvature motion of triple junctions of graphs in two dimensions},
Comm. PDE {\bf 35} (2010), no. 2, 302-327. 

\bibitem{Freire2}
A. Freire, {\em The existence problem for Steiner networks in strictly convex domains}, Arch. Ration. Mech. Anal. {\bf 200} (2011), no. 2, 361-404.

\bibitem{Gage}
M. Gage, R.~S. Hamilton,
{\em The heat equation shrinking convex plane curves},
J. Differential Geom. {\bf 23} (1986), no. 1, 69-96. 

\bibitem{Garcke}
H. Garcke, Y. Kohsaka, D. {\v{S}}ev{\v{c}}ovi{\v{c}}, {\em Nonlinear stability of stationary solutions for curvature flow with triple function}, Hokkaido Math. J. {\bf 38} (2009), no. 4, 721-769.

\bibitem{Giga1}
Y. Giga, {\em Surface evolution equations. A level set approach}, Monographs in Mathematics, {\bf 99}. Birkh\"{a}user Verlag, Basel, 2006.

\bibitem{Grayson}
M. Grayson, {\em The heat equation shrinks embedded plane curves to round points}, 
J. Differential Geom. {\bf 26} (1987), no. 2, 285-314. 


 \bibitem{Huisken}
	 G. Huisken,
	 {\em Asymptotic behavior for singularities of the mean
	 curvature flow}, 
	 J. Differential Geom. \textbf{31} (1990),
	 285-299.
	 
\bibitem{Ikota2005}
R. Ikota, E. Yanagida,
{\em Stability of stationary interfaces of binary-tree type},
Calc. Var. PDE {\bf 22} (2005), no. 4, 375-389.

\bibitem{Ilmanenp} T. Ilmanen, {\em Singularities of mean curvature flow of surfaces}, 
\newline
http://www.math.ethz.ch/$\sim$ilmanen/papers/sing.ps, 1995.

\bibitem{Ilmanen1}
T. Ilmanen, {\em Elliptic regularization and partial regularity for motion by mean curvature},
Mem. Amer. Math. Soc. {\bf 108} (1994), no. 520.

\bibitem{Ilmanen2} T. Ilmanen, {\em Convergence of the Allen-Cahn equation to Brakke's motion by mean curvature}, J. Differential Geom. {\bf 38} (1993), no. 2, 417-461.

\bibitem{Neves}
T. Ilmanen, A. Neves, F. Schulze, {\em On short time existence for the planar network
flow}, arXiv:1407.4756.
	 
 \bibitem{KT} 
	 K. Kasai, Y. Tonegawa,
	 {\em A general regularity theory for weak mean curvature
	 flow}, Calc. Var. PDE. \textbf{50} (2014), 1-68.

\bibitem{Kinderlehrer}
D. Kinderlehrer, C. Liu, {\em Evolution of grain boundaries}, Math. Models Methods Appl. Sci. {\bf 11} (2001), no. 4, 713-729.	 

\bibitem{Lahiri}
A. Lahiri, {\em Regularity of the Brakke flow}, Dissertation, Freie Universit\"{a}t Berlin, 2014.
\bibitem{Magni}
C. Magni, C. Mantegazza, M. Novaga, {\em Motion by curvature of planar networks II}, arXiv:math/0302164, to appear in Annali Sc. Norm. Super. Pisa CL. Sci.

\bibitem{MantegazzaT}
C. Mantegazza, {\em Lecture notes on mean curvature flow}, Progress in Mathematics, {\bf 290}, Birkh\"{a}user/Springer Basel AG, Basel, 2011.
	 
\bibitem{Mantegazza}
C. Mantegazza, M. Novaga, V. M. Tortorelli, {\em Motion by curvature of planar networks},
 Ann. Sc. Norm. Super. Pisa Cl. Sci. (5) {\bf 3} (2004), no. 2, 235-324.
 
\bibitem{Menne} 
U. Menne, {\em Second order rectifiability of integral varifolds of locally
bounded first variation}, J. Geom. Anal. {\bf 23} (2013), no. 2, 709-763.
 
\bibitem{Saez} 
M. S\'{a}ez Trumper, {\em Uniqueness of self-similar solutions to the network flow in
a given topological class}, Comm. PDE \textbf{36} (2011), no. 2, 185-204.

\bibitem{Sch} O. Schn\"{u}rer, A. Azouani, M. Georgi, J. Hell, N. Jangle, A. K\"{o}ller, T. Marxen, S. Ritthaler, 
M. S\'{a}ez, F. Schulze, B. Smith, {\em Evolution of convex lens-shaped networks under 
curve shortening flow}, Trans. Amer. Math. Soc. {\bf 363} (2011), no. 5, 2265-2294.

 \bibitem{Simon}
	 L. Simon,
	 {\em Lectures on geometric measure theory},
	 Proceedings of the Centre for Mathematical Analysis,
	 Australian National University,
	 \textbf{3},
	 Australian National University Centre for Mathematical
	 Analysis, Canberra, 1983.

\bibitem{Solo}
V. A. Solonnikov, {\em Boundary value problems of mathematical physics}, VIII, American Mathematical Society, Providence, R.I., 1975. 	 
	
 \bibitem{Takasao}
	 K. Takasao, Y. Tonegawa,
	 {\em Existence and regularity of mean curvature flow with
	 transport term in higher dimensions}, Math. Ann. \textbf{364} (2016), 857--935.
	 
\bibitem{Tonegawa1}
     Y. Tonegawa, 
     {\em Integrality of varifolds in the singular limit of reaction-diffusion equations},
     Hiroshima Math. J. \textbf{33} (2003), 323--341.
     
\bibitem{Tonegawa2}
      Y. Tonegawa, 
      {\em A second derivative H\"{o}lder estimate for weak mean curvature flow},
      Adv. Cal. Var. \textbf{7} (2014), 91--138.
      	 
\bibitem{Wick} Y. Tonegawa, N. Wickramasekera, {\em The blow up method for Brakke flows: networks near triple junctions}, Arch. Ration. Mech. Anal. (to appear). 

\bibitem{White0} B. White, {\em Stratification of minimal surfaces, mean curvature flows, and harmonic maps}, 
J. Reine Angew. Math. {\bf 488} (1997), 1--35

\end{thebibliography}
